\documentclass[11pt]{amsart}

\usepackage{parskip}
\allowdisplaybreaks
\usepackage{geometry} 
\usepackage{framed} 

\usepackage{tikz}
\tikzset{every picture/.style={line width=1pt}} 

\usepackage{setspace}
\setdisplayskipstretch{2.5}

\usepackage{amsthm} 
\usepackage{amsmath}
\usepackage{amssymb}

\usepackage[symbol]{footmisc}

\usepackage{hyperref}

\usepackage{amsfonts} 
\newcommand\A{\mathcal{A}}
\newcommand\B{\mathcal{B}}
\newcommand\FB{\mathfrak{B}}
\newcommand\C{\mathbb{C}}
\newcommand\CC{\mathcal{C}}

\newcommand\BC{\mathbf{C}}

\newcommand\FF{\mathfrak{F}}
\newcommand\G{\mathcal{G}}
\renewcommand\H{\mathcal{H}}
\newcommand\CK{\mathcal{K}}
\newcommand\FL{\mathfrak{L}}
\newcommand\M{\mathcal{M}}
\newcommand\CN{\mathcal{N}}
\newcommand\N{\mathbb{N}}
\newcommand\R{\mathbb{R}}
\renewcommand\S{\mathcal{S}}
\newcommand\Z{\mathbb{Z}}
\newcommand\CZ{\mathcal{Z}}

\newcommand\w{\omega}

\usepackage{stmaryrd} 

\newcommand\id{\textnormal{id}}
\newcommand\sing{\textnormal{sing}}
\newcommand\reg{\textnormal{reg}}
\newcommand\spt{\textnormal{spt}}

\newcommand\dist{\textnormal{dist}}
\newcommand\graph{\textnormal{graph}}
\newcommand\ext{\mathrm{d}}
\newcommand\del{\partial}
\newcommand\vartan{\textnormal{VarTan}}

\newcommand{\res}{\mathbin{\hspace{0.1em}\vrule height 1.3ex depth 0pt width 0.13ex\vrule height 0.13ex depth 0pt width 1.0ex}} 

\newcommand{\weakly}{\rightharpoonup}
\renewcommand{\div}{\textnormal{div}}

\newtheoremstyle{newtheoremstyle}
{3pt}
{3pt}
{\itshape}
{\parindent}
{\bfseries}
{.}
{0.5em}
{} 
\newtheoremstyle{newtheoremstyledefn}
{3pt}
{3pt}
{}
{\parindent}
{\bfseries}
{.}
{0.5em}
{} 

\theoremstyle{newtheoremstyle}
\newtheorem{theorem}{Theorem}
\newtheorem{lemma}[theorem]{Lemma}
\newtheorem{prop}[theorem]{Proposition}
\newtheorem{corollary}[theorem]{Corollary}

\newtheorem{thmx}{Theorem}

\theoremstyle{newtheoremstyledefn}
\newtheorem{defn}[theorem]{Definition}

\numberwithin{equation}{section} 
\numberwithin{theorem}{section}

\usepackage{fancyhdr}
\begin{document}

\title{The Structure of Stable Codimension One Integral Varifolds near Classical Cones of Density 5/2}

\author{Paul Minter}
\address{Department of Pure Mathematics and Mathematical Statistics, University of Cambridge}
\email{pdtwm2@cam.ac.uk}

\begin{abstract}
	We prove a multi-valued $C^{1,\alpha}$ regularity theorem for the varifolds in the class $\S_2$ (i.e., stable codimension one stationary integral $n$-varifolds admitting no triple junction classical singularities) which are sufficiently close to a stationary integral cone comprised of 5 half-hyperplanes (counted with multiplicity) meeting along a common axis. Such a result is the first of its kind for non-flat cones of higher (i.e. $>1$) multiplicity when branch points are present in the nearby varifolds. For such varifolds, this completes the analysis of the singular set in the region where the density is $<3$, up to a set which is countably $(n-2)$-rectifiable.
	
	Our methods develop the blow-up arguments in \cite{simoncylindrical} and \cite{wickstable}. One key new ingredient of our work is needing to inductively perform successively finer blow-up procedures in order to show that a certain $\epsilon$-regularity property holds at the blow-up level; this is then used to prove a $C^{1,\alpha}$ boundary regularity theory for two-valued $C^{1,\alpha}$ harmonic functions which arise as blow-ups of sequences of such varifolds, the argument for which is carried out in the accompanying work \cite{minter2021}.
\end{abstract}

\maketitle

\tableofcontents

\section{Introduction}\label{sec:intro}

A central problem within geometric analysis concerns understanding the nature of singularities arising in stationary integral varifolds. However, since the first general regularity theorem of W.~Allard (\cite{allard}), surprising little is known. Allard's regularity theorem shows that, whenever a stationary integral varifold (of any dimension and codimension) is sufficiently close, in a varifold sense, to a multiplicity one plane, the varifold is in fact locally expressible as a single-valued $C^{1,\alpha}$ graph over a region of the plane, with estimates on the $C^{1,\alpha}$ norm of the graph; one may then invoke classical quasilinear elliptic regularity theory (\cite{gt}) to infer that the graph is in fact smooth, with estimates on its $C^k$ norm for each positive integer $k$. Results of this nature, where if a stationary integral varifold is ``$\epsilon$-close'' to another, simpler, varifold, then the former can be expressed as a suitable graph of some regularity over the support of the other, we refer to as an $\epsilon$\textit{-regularity theorem}.

As far as $\epsilon$-regularity theorems for stationary integral varifolds go, there are few other results known to hold in the same generality as Allard's regularity theory (e.g. \cite{simonasymptotics}). For us, the other key result is L.~Simon's $\epsilon$-regularity theorem (\cite{simoncylindrical}) for the (multiplicity one) triple junction, i.e. the stationary integral cone comprised of three multiplicity one half-planes with a common boundary (for a given dimension such a varifold is unique up to rotation). Here, the $C^{1,\alpha}$ graph one constructs over the triple junction is comprised of 4 separate functions: for each half-plane in the triple junction, we have a $C^{1,\alpha}$ function defined on a subset of the plane containing it which takes values orthogonal to the plane, and one $C^{1,\alpha}$ function defined on the common axis taking values orthogonal to the axis. 

Simon's $\epsilon$-regularity theorem is in fact a corollary of a more general theory established in \cite{simoncylindrical} which applies to so-called \textit{multiplicity one classes} $\M$ of stationary integral varifolds. The key lemma, \cite[Lemma 2.1]{simoncylindrical}, establishes a dichotomy roughly saying the following: if $\M$ satisfies an additional ``integrability'' hypothesis, then whenever $V\in \M$ is sufficiently close (as varifolds) to a given integral \textit{cylindrical} cone $\BC\in \M$, i.e. $\BC = \BC_0\times\R^m$ for some integer $m$ and $\BC_0$ obeys $\sing(\BC) \subset \{0\}$, then either:
\begin{enumerate}
	\item [(i)] there is a \textit{density gap} in $V$, or
	\item [(ii)] there is some scale $\theta\in (0,1)$ and cone\footnote{In fact $\tilde{\BC}$ also takes the same form as $\BC$, i.e. $\tilde{\BC} = \tilde{\BC}_0\times\R^m$; this fact is crucial for iteration purposes.} $\tilde{\BC}$ close to $\BC$ such that, after a small rotation, the height excess of $V$ relative to $\tilde{\BC}$ in $B_{\theta}(0)$ decays by a fixed factor compared to the height excess of $V$ relative to $\BC$ in $B_1(0)$.
\end{enumerate}
When $\BC$ is a triple junction, topological obstructions and Allard's regularity theorem prevent alternative (i) from occurring; moreover, the new cone $\tilde{\BC}$ found in (ii) is also a triple junction. This enables one to establish that the height excess decays along a geometric sequence of scales (i.e. $1,\theta,\theta^2,\dotsc$) and ultimately establish the $\epsilon$-regularity theorem. Indeed, such an excess decay result is used in the proof of Allard's regularity theorem. One may view such decay of the height excess as the geometric analogue of integral decay required in the Campanato spaces (see \cite{campanato}), and thus this approach is the geometric equivalent of Campanato's regularity theory for functions in Campanato spaces.

In recent years, the ideas developed by L.~Simon have been developed further. Key examples of this include: \cite{beckerkahn}, where density gaps do arise; \cite{polyhedral}, where the (multiplicity one) cone $\BC$ need not be cylindrical; \cite{krummelwick1}, \cite{krummelwick2}, \cite{krummelwick3}, where the ``cone'' $\BC$ is instead the graph (possibly with multiplicity $>1$) of a multi-valued homogeneous harmonic function $\phi$, with varying degrees of homogeneity, including degrees of homogeneity $<1$; \cite{wickstable}, where a higher multiplicity degenerate situation is considered, as well as a situation where the cone $\BC$ is supported on a union of half-hyperplanes meeting along a common axis, where the half-hyperplanes can occur with multiplicity $>1$ (but sheeting still holds away from the axis); \cite{minterwick}, where a degenerate higher multiplicity flat situation in the presence of branch points is analysed, using Almgren's frequency function to establish regularity of the blow-ups. In all these examples one must deal with density gaps, often by showing that they cannot occur under the given assumptions. 

The main result is to establish an $\epsilon$-regularity theorem in a setting where the cone is non-planar and the nearby varifold can have branch points. Let us first set-up some notation, following \cite{minterwick}. Let $\S_2$ denote the class of integral $n$-varifolds $V$ on the open ball $B^{n+1}_2(0)\subset\R^{n+1}$ with $0\in \spt\|V\|$, $\|V\|(B^{n+1}_2(0))<\infty$, and which satisfy the following conditions:
\begin{enumerate}
	\item [$(\S1)$] $V$ is stationary in $B^{n+1}_2(0)$ with respect to the area functional, in the following (usual) sense: for any given vector field $\psi\in C^1_c(B^{n+1}_2(0);\R^{n+1})$, $\epsilon>0$, and $C^2$ map $\phi:(-\epsilon,\epsilon)\times B^{n+1}_2(0)\to B^{n+1}_2(0)$ such that:
	\begin{enumerate}
		\item [(i)] $\phi(t,\cdot):B^{n+1}_2(0)\to B^{n+1}_2(0)$ is a $C^2$ diffeomorphism for each $t\in (-\epsilon,\epsilon)$ with $\phi(0,\cdot)$ equal to the identity map on $B^{n+1}_2(0)$;
		\item [(ii)] $\phi(t,x) = x$ for each $(t,x)\in (-\epsilon,\epsilon)\times\left(B^{n+1}_2(0)\backslash\spt(\psi)\right)$; and
		\item [(iii)] $\left.\del\phi(t,\cdot)/\del t\right|_{t=0} = \psi$;
	\end{enumerate}
	we have that
	$$\left.\frac{\ext}{\ext t}\right|_{t=0}\|\phi(t,\cdot)_\#V\|(B^{n+1}_2(0)) = 0;$$
	equivalently (see \cite[Section 39]{simongmt}),
	$$\int_{B^{n+1}_2(0)\times G_n}\div_S\psi(X)\ \ext V(X,S) = 0$$
	for every vector field $\psi\in C^1_c(B^{n+1}_2(0);\R^{n+1})$, where $G_n$ is the set of $n$-dimensional subspaces in $\R^{n+1}$;
	\item [$(\S2)$] $\reg(V)$ is stable in $B^{n+1}_2(0)$, in the following (usual) sense: for each open ball $B\subset B^{n+1}_2(0)$ with $\sing(V)\cap B = \emptyset$ in the case $2\leq n\leq 6$ or $\H^{n-7+\gamma}(\sing(V)\cap B) = 0$ for every $\gamma>0$ in the case $n\geq 7$, given any vector field $\psi\in C^1_c(B\backslash\sing(V);\R^{n+1})$ with $\psi(X)\perp T_X\reg(V)$ for each $X\in \reg(V)\cap B$,
	$$\left.\frac{\ext^2}{\ext t^2}\right|_{t=0}\|\phi(t,\cdot)_\#V\|(B^{n+1}_2(0))\geq 0$$
	where $\phi(t,\cdot)$, $t\in (-\epsilon,\epsilon)$, are the $C^2$ diffeomorphisms of $B^{n+1}_2(0)$ associated with $\psi$ described in $(\S1)$ above; equivalently (see \cite[Section 9]{simongmt})\footnote{This equivalence requires two-sidedness of $\reg(V)$, which holds in a ball $B$ as above in view of the smallness assumption on the singular set in $B$.} for every such $\Omega$ we have
	$$\int_{\reg(V)\cap B}|A|^2\zeta^2\ \ext\H^n \leq \int_{\reg(V)\cap B}|\nabla\zeta|^2\ \ext \H^n\ \ \ \ \text{for all }\zeta\in C^1_c(\reg(V)\cap B)$$
	where $A$ denotes the second fundamental form of $\reg(V)$, $|A|$ the length of $A$, and $\nabla$ the gradient operator on $\reg(V)$;
	\item [$(\S3)$] $V$ does not contain any triple junction classical singularities.
\end{enumerate}

Note that here by a \textit{triple junction singularity} in $V$ we mean a point $X\in\spt\|V\|$ for which there is a radius $\rho>0$ such that $V\res B^{n+1}_\rho(X)$ is a sum of three (multiplicity one) $C^{1,\alpha}$ submanifolds-with-boundary, which all have a common $C^{1,\alpha}$ boundary, for some $\alpha>0$ (by \cite{krummel2014regularity}, we can in fact assume that the submanifolds are real-analytic and their common boundary is real-analytic also); note that by Simon's triple junction $\epsilon$-regularity theorem (\cite{simoncylindrical}) this is equivalent to $X$ exhibiting a tangent cone which is a sum of 3 (multiplicity one) half-hyperplanes.

To precisely state our main result, let $\BC_0 = \sum^N_{i=1}q_i^{(0)}|H_i^{(0)}|$ be a stationary classical cone in $\R^{n+1}$ with density $\Theta_{\BC_0}(0) = \frac{5}{2}$ and spine $L_{\BC_0} = \{(0,0)\}\times\R^{n-1}$, where $q_i^{(0)}$ are integers $\geq 1$, $H_i^{(0)}$ are distinct half-hyperplanes with $\del H_i^{(0)} = L_{\BC_0}$ for each $i=1,2,\dotsc,N$\footnote{By virtue of the stationarity of $\BC_0$ we must have $q_i^{(0)}\leq 2$ for each $i=1,\dotsc,N$, and the density condition $\Theta_{\BC_0}(0) = \frac{5}{2}$ is equivalent to $\sum^N_{i=1}q_i^{(0)} = 5$. In particular, $N\in \{3,4,5\}$.}; thus $H_i^{(0)} = R_i^{(0)}\times\R^{n-1}$, where $R_i^{(0)} = \{t\mathbf{w}_i^{(0)}:t>0\}$ for distinct unit vectors $\mathbf{w}_1^{(0)},\dotsc,\mathbf{w}_N^{(0)}\in \R^2$.

Let $\sigma_0:= \max\{\mathbf{w}_i^{(0)}\cdot\mathbf{w}_k^{(0)}:i,k=1,2,\dotsc,N,\, i\neq k\}$ and let $N(H_i^{(0)})$ denote the conical neighbourhood of $H_i^{(0)}$ defined by
$$N(H^{(0)}_i):= \left\{(x,y)\in \R^2\times\R^{n-1}: x\cdot\mathbf{w}_i^{(0)}>|x|\sqrt{\frac{1+\sigma_0}{2}}\right\}.$$
Denote by $\tilde{H}_i^{(0)}$ the hyperplane containing $H_i^{(0)}$ and by $(\tilde{H}_i^{(0)})^\perp$ the orthogonal complement of $\tilde{H}^{(0)}_i$ in $\R^{n+1}$.

\begin{thmx}\label{thm:A}
	Let $\BC_0$ be as above. Then there is a constant $\epsilon = \epsilon(\BC_0)$ such that the following holds: if $V\in \S_2$ has $(2+1/8)\w_n \leq \|V\|(B^{n+1}_1(0)) \leq (3-1/8)\w_n$ and
	$$\int_{B^{n+1}_1(0)}\dist^2(X,\spt\|\BC_0\|)\ \ext\|V\|<\epsilon,$$
	then for each $i\in \{1,2,\dotsc,N\}$ there is a function
	$$\gamma_i\in C^{1,\alpha}\left(\overline{L_{\BC_0}\cap B^{n+1}_{1/2}(0)};\tilde{H}_i^{(0)}\cap \{X:\dist(X,L_{\BC_0})<1/16\}\right)$$
	and a function $u_i:\Omega_i\to \A_{q_i^{(0)}}((\tilde{H}_i^{(0)})^\perp)$, where $\Omega_i$ is the connected component of $\tilde{H}^{(0)}_i\cap B_{1/2}(0)^{n+1}\backslash\{x+\gamma_i(x):x\in L_{\BC_0}\cap B^{n+1}_{1/2}(0)\}$ with $(H_i^{(0)}\backslash\{X:\dist(X,L_{\BC_0})<1/16\})\cap B^{n+1}_{1/2}(0)\subset \Omega_i$ such that:
	\begin{enumerate}
		\item [(i)] $u_i\in C^{1,\alpha}(\overline{\Omega}_i;\A_{q_i^{(0)}})$, with $\mathbf{v}(u_i)$ a stationary integral varifold, where $\mathbf{v}(u)$ is the varifold $(\graph(u_i),\theta)$, where if we write $u_i(X) = \sum^{q_i^{(0)}}_{j=1}\llbracket u_i^j(X)\rrbracket$ for $X\in \Omega_i$, then $\graph(u_i) = \{(u_i^j(X),X):X\in \Omega_i, j\in \{1,\dotsc,q^{(0)}_i\}\}$, and the multiplicity function $\theta$ at a point $(u_i^j(X),X)\in \graph(u_i)$ is given by $\theta(u_i^j(X),X) := \#\{k:u_i^k(X) = u_i^j(X)\}$ for each $j=1,\dotsc,q^{(0)}_i$;
		\item [(ii)] for each $i\in \{1,\dotsc,N\}$, $\left.u_i\right|_{\del \Omega_i\cap B^{n+1}_{1/2}(0)} = q_i^{(0)}\llbracket b_i\rrbracket$ for some single-valued $C^{1,\alpha}$ function $b_i:\del\Omega_i\cap B^{n+1}_{1/2}(0)\to (\tilde{H}_i^{(0)})^\perp$, and moreover if $\tilde{b}_i(x):= x+b_i(x)$ for $x\in \del\Omega_i\cap B^{n+1}_{1/2}(0)$, then $\textnormal{image}(\tilde{b}_i) = \textnormal{image}(\tilde{b}_j)$ for all $i,j\in \{1,\dotsc,N\}$;
		\item [(iii)] $V\res B^{n+1}_{1/2}(0) = \sum^N_{i=1}\mathbf{v}(u_i)\res B^{n+1}_{1/2}(0)$;
		\item [(iv)] for each $i\in \{1,\dotsc,N\}$,
		$$\{Z:\Theta_V(Z)\geq 5/2\}\cap B^{n+1}_{1/2}(0) = \{Z:\Theta_V(Z) = 5/2\} \cap B^{n+1}_{1/2}(0) = \tilde{b}_i(\del\Omega_i\cap B^{n+1}_{1/2}(0)).$$
		\end{enumerate}
		Moreover, for each $i\in \{1,\dotsc,N\}$ we have
		$$|u_i|_{1,\alpha;\overline{\Omega}_i}\leq C\left(\int_{B_1^{n+1}(0)}\dist^2(X,\spt\|\BC\|)\ \ext\|V\|(X)\right)^{1/2}.$$
		Here, $C = C(n)\in (0,\infty)$ and $\alpha = \alpha(n) \in (0,1/2)$. In particular, $V$ has a unique tangent cone at every point in $B^{n+1}_1(0)$, and $\{X:\Theta_V(X) = 5/2\}\cap B^{n+1}_{1/2}(0)$ is a connected $C^{1,\alpha}$ $(n-1)$-dimensional submanifold.
\end{thmx}

In \cite{minterwick} the structure of varifolds in $\S_2$ which are near a stationary cone comprised of 4 half-hyperplanes meeting along a common axis is studied; the stationarity condition in fact implies that such a cone must be a sum of two multiplicity one hyperplanes (which could coincide). Thus, the above result solves the next significant case to be studied, namely when the cone is comprised of $5$ half-hyperplanes meeting along a common axis. Indeed, combining Theorem \ref{thm:A} with \cite{minterwick} and \cite{simoncylindrical}, we therefore get for $V\in \S_2$ the following decomposition of the singular set in the region where the density is $<3$:

\begin{thmx}\label{thm:B}
	Let $V\in \S_2$. Then
	$$\spt\|V\|\cap B^{n+1}_1(0) \cap \{\Theta_V<3\} = \Omega\cup\B\cup\mathcal{T}\cup\mathcal{C}\cup K$$
	where:
	\begin{enumerate}
		\item [(i)] $\Omega$ is the set of points $X\in \spt\|V\|\cap B^{n+1}_1(0)\cap \{\Theta<3\}$ such that for some $\delta_X>0$, $\spt\|V\|\cap B_{\delta_X}(X)$ is a smoothly embedded hypersurface;
		\item [(ii)] \textnormal{(\cite{minterwick})} $\B$ is the set of points $X\in \spt\|V\|\backslash\Omega$ such that one tangent cone to $V$ at $X$ is of the form $2|P|$ for some hyperplane $P$; moreover, this is the unique tangent cone to $V$ at $P$, and there is a $\delta_X>0$ such that $V\res B_{\delta_X}(X)$ is given by a $C^{1,1/2}$ two-valued function over a domain a domain in $P$;
		\item [(iii)] \textnormal{(\cite{minterwick})} $\mathcal{T}$ is the set of points $X\in \spt\|V\|$ such that one tangent cone to $V$ at $X$ is of the form $|P_1|+|P_2|$ for a pair of transversely intersecting hyperplanes $P_1,P_2$; moreover, this is the unique tangent cone to $V$ at $X$, and there is a $\delta_X>0$ such that $V\res B_{\delta_X}(X)$ is a union of two transversely intersecting smooth single-valued graphs, one over each hyperplane $P_1,P_2$.
		\item [(iv)] $\mathcal{C}$ is the set of points $X\in\spt\|V\|$ such that one tangent cone to $V$ is the sum of 5 half-hyperplanes meeting along a common axis; moreover, this is the unique tangent cone to $V$ at $X$, and the conclusions of Theorem \ref{thm:A} hold in some ball about $X$;
		\item [(v)] $K = \S_{n-2}$ is the usual $(n-2)$-stratum of the singular set.
	\end{enumerate}
	In particular, when non-empty, $\B$ is countably $(n-2)$-rectifiable (\cite{krummelwick1}; see also \cite{simon2016frequency}), $\mathcal{T}$ and $\mathcal{C}$ are smoothly embedded $(n-1)$-dimensional submanifolds of $B^{n+1}_1(0)$, and $K$ is countably $(n-2)$-rectifiable (\cite{naber2015singular}). 
\end{thmx}

\subsection{Contextual Overview of Higher Multiplicity Singularities}\label{sec:1.1}

In this work we wish to understand the nature of the singular set of stationary integral $n$-varifolds $V$ near certain cones which arise naturally and form the `largest' part of the singular set. In general, one may stratify the singular set (\cite{almgren}; see also \cite{federer1970singular}) into regions based on their tangent cone type, namely, one may always write:
$$\sing(V) = \B\cup \S\cup K$$
where:
\begin{enumerate}
	\item [(a)] $\B$ is the branch set, that is, the singular points $X$ at which at least one tangent cone is a multiplicity $\geq 2$ plane, and there is no (ambient) ball $B$ centred at such that $V\res B$ is a sum of finitely many smoothly embedded minimal submanifolds;
	\item [(b)] $\S$ is the set of singular points in $\sing(V)\backslash\B$ where at least one tangent cone is supported on a union of half-planes meeting along a common $(n-1)$-dimensional axis. A priori we know that such singular points have density $q/2$ for some $q\in \Z_{\geq 3}$, and that $\dim_\H(\S)\leq n-1$ (in fact is countably $(n-1)$-rectifiable by \cite{naber2015singular});
	\item [(c)] $K$ has $\dim_\H(K)\leq n-2$ (again, in fact $K$ is countably $(n-2)$-rectifiable by \cite{naber2015singular}).
\end{enumerate}
It is still an open question whether it is possible to have $\dim_\H(\B) = n$\footnote{An example found by K.~Brakke (\cite[Section 6]{brakke2015motion} demonstrates that it is possible to have $\dim_\H(\B)=n$ when one always for non-zero (generalised) mean curvature, even if it is uniformly arbitrarily small.} It should be noted that in the area-minimising setting simple 1-dimensional comparison arguments show that $\S = \emptyset$, and furthermore for codimension one area-minimisers one also has $\B = \emptyset$.

\begin{figure}[h]
	\centering
	\begin{tikzpicture}
		\draw (0,0) -- (2,0.2) (0,0) -- (1.9,-0.2) (0,0) -- (-1,1) (-0.4,-0.4) -- (-1,-1);
		\draw[dotted] (0,0) -- (-0.4,-0.4);
		\draw[gray] (-2,-0.4) -- (-1.1,0.4) -- (2.8,0.5) -- (2.5,-0.4) -- (-2,-0.4);
		\draw (0,0) -- (-1,1);
	\end{tikzpicture}
	\caption{Example of a twisted cone. The plane containing two of the lines is illustrated in grey.}
	\label{fig:twisted}
\end{figure}
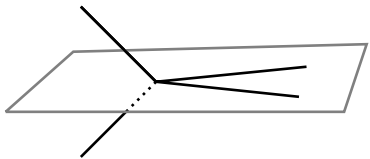

Thus, if one wishes to understand singularities in stationary integral varifolds one is naturally led to study singular points in $\B\cup \S$. The simplest case in this setting is when $X\in \S$ has density $\frac{3}{2}$; then $V$ has a tangent cone at $X$ which is equal to a three multiplicity one half-planes which have a common boundary. L.~Simon (\cite{simoncylindrical}) showed that locally about such $X$, every singularity has density $\frac{3}{2}$ and is in $\S$, and moreover the singular points on a neighbourhood of $X$ form a $C^{1,\alpha}$ $(n-1)$-dimensional submanifold (in fact a smooth submanifold, by \cite{krummel2014regularity}). The next step is to understand $X\in \B\cup \S$ of density $2$; therefore $X$ has a tangent cone which is either a multiplicity two plane or a union of $4$ multiplicity one half-planes with a common boundary (note that a priori both types of tangent cone could occur simultaneously). In codimension one, a stationary union of 4 multiplicity one half (hyper)planes meeting along a common boundary must necessarily be a union of two planes, however this is not true in codimension $>1$ (due to so-called \textit{twisted cones}; see Figure \ref{fig:twisted}. In this general setting, $\epsilon$-regularity theorems for such cones are false as simple examples illustrate (such as scalings of the catenoid; examples with branch points include complex analytic varieties such as $\{(z,w)\in \C^2:z^2=w^3\}$, which are even area-minimising, and stable codimension one examples constructed in \cite{simon2007stable}, \cite{krummel2019existence}).  

As two transverse multiplicity one planes can limit onto a multiplicity two plane, one has to tackle the former case first. By a tremendous effort, an $\epsilon$-regularity theorem has been established by S.~Becker--Kahn and N.~Wickramasekera (\cite{beckerwick}) for stationary integral varifolds close, as varifolds, to a sum of two multiplicity one planes (not necessarily distinct, so even the case of a single multiplicity two plane) in arbitrary dimension and codimension; one only needs to assume a topological assumption, namely that in a sufficiently flat region where the varifold has density $<2$ everywhere, the varifold must decompose as a sum of smoothly embedded sheets. In the codimension one setting this topological assumption is directly implied when one assumes stability of the varifold on the regular part and that the varifold has no triple junction singularities (see \cite{schoen1981regularity} and \cite{wickstable}), and in this case one establishes a strong local structural property for such varifolds (see Theorem \ref{thm:wick1}, Theorem \ref{thm:wick2}), and in particular that the (multiplicity 2) branch set is necessarily countably $(n-2)$-rectifiable (see \cite{simon2016frequency} and \cite{krummelwick1}). Moreover, if one rules out \text{all} classical singularities in a stable codimension one stationary integral varifolds, then branch points do not occur and the singular set is in fact countably $(n-7)$-rectifiable (see \cite{wickstable}); thus, up to a set of codimension 7, every singularity is a limit of classical singularities. These structural results have in fact recently been refined in \cite{minterwick}, where it is shown that one can understand the local structure about a density $Q$ branch point in a stable codimension one stationary integral varifold is there are no classical singularities nearby with density $<Q$.

Thus, after the case of $X\in \B\cup \S$ with density $2$, the next case if to study $X\in \S$ with density $\frac{5}{2}$. Such singular points have a tangent cone which is sum of $5$ half-planes (counted with multiplicity) meeting along a common axis. This is the situation we study in the current work.

\subsection{The Present Work}

As mentioned above, after one has understood the local structure about multiplicity two points in $\B\cup \S$, the next step is understanding the nature of the varifold locally about points in $\S$ of density $\frac{5}{2}$: this is the setting we study here. This is different to previous settings as one needs to deal with both higher multiplicity (and hence branch points) at the same time as non-flatness. Here, we shall study this problem in the context of stable codimension one stationary integral $n$-varifolds which do not contain triple junction singularities; this enables us to invoke the strong structural results of \cite{minterwick} near density 2 points in $\B\cup \S$ which will be crucial for our analysis. In particular, we will be able to use the fact that the  multiplicity two branch set has dimension at most $n-2$ in order to prove that density gaps do not occur; this will be crucial for establishing similar estimates to those seen in \cite{simoncylindrical} (and \cite{wickstable}).

\begin{figure}[h]
		\centering
			\begin{tikzpicture}
				\draw (0,0.8) -- (0,-0.3);
				\draw[ultra thick] (1.1,-0.8) -- (0,-0.3) (-1,-0.8) -- (0,-0.3); 
				\draw[ultra thick] (-3,1) -- (-3,0);
				\draw (-2.2,-1) -- (-3,0) (-3.5,-1) -- (-3,0) (-3.7,-1) -- (-3,0); 
				\draw (-6,-1) -- (-6,-0) (-6.7,1) -- (-6,0) (-6.9,-0.9) -- (-6,0) (-5.3,1) -- (-6,0) (-5.1,0.7) -- (-6,0); 
				\node at (0.75,-0.4) {2};
				\node at (-0.75,-0.4) {2};
				\node at (-2.7,0.5) {2};
				\node at (0,-1.3) {Level 2};
				\node at (-3,-1.3) {Level 1};
				\node at (-6,-1.3) {Level 0};
			\end{tikzpicture}
			\caption{Representative cross-sections $\BC_0$ of the three different \textit{levels} of (tangent) cones we shall be considering. Note that for level 2 cones, the angles between each of the three rays are pre-determined, but not all equal. In each picture the multiplicity two pieces are highlighted.}
			\label{fig:levels}
\end{figure}
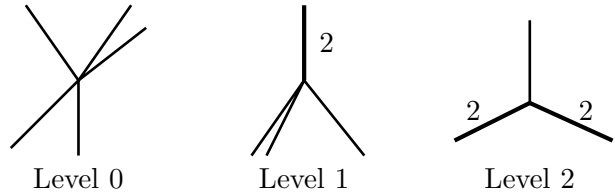

Tangent cones $\BC$ to points in $\S$ of density $\frac{5}{2}$ will, up to a rotation, take the form $\BC = \BC_0\times\R^{n-1}$, where $\BC_0$ is a 1-dimensional stationary integral varifold in $\R^2$ with $\Theta_{\BC_0}(0) = \frac{5}{2}$. Such $\BC_0$ must be supported on at most $5$ multiplicity one rays through the origin; it is easy to check from the stationarity condition that the multiplicity of a given ray in $\BC_0$ is at most $2$. Thus we can divide the different types of $1$-dimensional cross-section $\BC_0$ into three different \textit{classes}, depending on the number of multiplicity two rays, of which there can be $0$, $1$, or $2$. We shall refer to the class where there are $I\in \{0,1,2\}$ multiplicity two rays in $\BC_0$ as the class of \textit{level} $I$ cones; examples of cones in each class are given in Figure \ref{fig:levels}. Note that a level $I$ cone necessarily has $5-I$ rays in the support of its cross-section.

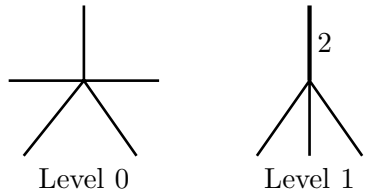
\begin{figure}[b]
	\centering
	\begin{tikzpicture}
				\draw[ultra thick] (0,1) -- (0,0);
				\draw (0,0) -- (0,-1) (0.7,-1) -- (0,0) (-0.7,-1) -- (0,0);
				\node at (0.2, 0.5) {2};
				\node at (0, -1.3) {Level 1}; 
				\draw (-3,0) -- (-3,1) (-3,0) -- (-2,0) (-3,0) -- (-4,0) (-3,0) -- (-2.3,-1) (-3,0) -- (-3.8,-1);
				\node at (-3,-1.3) {Level 0};
			\end{tikzpicture}
		\caption{Examples of decomposable level 0 and level 1 cones.}
		\label{fig:decomposable}
\end{figure}

It is important to note that some level 0 and level 1 cones are \textit{decomposable}, i.e., they are the sum of distinct stationary integral varifolds; in this situation this decomposition always takes the form of a multiplicity one plane and a triple junction. Such examples are illustrated in Figure \ref{fig:decomposable}. A priori this could lead to complications as it allows for density gaps to occur in varifolds arbitrarily close to such cones: just consider rotating one of the planes and ending up with a cone as in Figure \ref{fig:gaps} (which does not have an $(n-1$)-dimensional set of points of density $\geq \Theta_{\BC}(0) = \frac{5}{2}$), or translating the planar part, as shown in Figure \ref{fig:nopoints}. Note that in these examples, there are points of density $\frac{3}{2}$ close to the points of density $\frac{5}{2}$. Our main theorem (Theorem \ref{thm:A}) shows that the presence of points of density $\frac{3}{2}$ is the only obstruction to complete regularity of a stable codimension one integral varifold near a classical cone $\BC$ with $\Theta_{\BC}(0) = \frac{5}{2}$.

\begin{figure}[b]
		\centering
		\begin{minipage}{0.51\textwidth}
		\centering
			\begin{tikzpicture}
				\draw (0,0.5) -- (0.5,0) -- (0.5,-1.5) -- (0,-1) -- (0,0.5);
		\draw (0,-1) -- (-1,-2) -- (-0.5,-2.5) -- (0.5,-1.5);
		\draw (0.5,-1.5) -- (1.5,-2) -- (1.3, -1.7) -- (0.5,-1.3);
		\draw [gray] (0.5,-1.3) -- (0,-1); 
		\draw (-0.75,0.25) -- (1.3,0.45) -- (1.35,-2.3) -- (-0.7,-2.5) -- (-0.75,0.25);
		\draw[dotted] (0.25,-1.2) -- (0.25,0.3) (0.25,-1.2) -- (-0.7,-2.1) (0.25,-1.25) -- (1.3,-1.8);
			\end{tikzpicture}
			\caption{Illustration of density gaps. There is only one density $\frac{5}{2}$ point, yet by rotating the plane the varifold can be arbitrarily close to a cone with a line of density $\frac{5}{2}$.}
			\label{fig:gaps}
		\end{minipage}\hfill
		\begin{minipage}{0.49\textwidth}
		\centering
			\begin{tikzpicture}
				\draw (0,0) -- (0,1.5) (0,0.2) -- (1.5,0.2) (0,0.2) -- (-1.5,0.2) (0,0) -- (1,-1.4) (0,0) -- (-1,-1.4);
			\end{tikzpicture}
			\caption{No points of density $\frac{5}{2}$ but arbitrarily close to a cone which has points of density $\frac{5}{2}$.}
			\label{fig:nopoints}
		\end{minipage}
\end{figure}

The proof of the main theorem employs a blow-up (linearisation) argument. A significant problem that arises in the presence of higher multiplicity (as in level 1 and level 2 cones) is the possibility of branch point singularities nearby in $V$. Such singularities can a priori occur, and using the regularity theory of \cite{minterwick} we are only able to express $V$ over regions in the support, away from the boundary, of such a multiplicity two half-hyperplane in $\BC$ as a $C^{1,\alpha}$ two-valued graph (in fact $\alpha=1/2$: see \cite{simon2016frequency}). Thus, when we perform the blow-up procedure we shall end up with a $C^{1,\alpha}$ two-valued harmonic function defined on the support of each multiplicity two half-hyperplane of $\BC$. The difference between the single-valued harmonic setting (which occurs in \cite{simoncylindrical} and \cite{wickstable}) is that in the multi-valued harmonic setting we do not have a reflection principle or $C^{1,\alpha}$ boundary regularity theory which we can apply to establish regularity up-to-the-boundary for these two-valued $C^{1,\alpha}$ harmonic functions. Establishing that such results do in fact hold in our setting will take up the majority of our work. Indeed, to achieve this we shall adapt the arguments seen in \cite{wickstable}: we first prove a suitable regularity claim for homogeneous degree one blow-ups in order to perform a reflection argument to classify the homogeneous degree one blow-ups, and subsequently use this to prove the regularity claim for general blow-ups (using an argument based on the \textit{reverse Hardt--Simon inequality}).

Many of our estimates will be integral estimates, and to pass from integral decay estimates to regularity statements one typically uses Campanato regularity theory (\cite{campanato}); thus one needs to establish a variant of the Campanato regularity theory for multi-valued functions in order to conclude the regularity statements for the blow-ups; this is done in the accompanying work \cite{minter2021}. The main difficulty in the approach outlined above is establishing the regularity up-to-the-boundary for homogeneous degree one blow-ups, which is again achieved by an argument based on the (reverse) Hardt--Simon inequality. To do this, one needs a suitable $\epsilon$-regularity property for the blow-ups similar to that seen in \cite[$(\B7)$, Section 4]{wickstable}. Loosely speaking, in \cite{wickstable} this property can be thought of as saying that no classical singularities occur at the blow-up (i.e. linear) level if classical singularities do not occur at the varifold level. In our setting, this $\epsilon$-regularity property takes a different form: it tells us that if a blow-up, relative to a level $I\in \{1,2\}$ cone, has a graph which is sufficiently close to a union of $>5-I$ half-hyperplanes meeting along a common axis, then in fact the blow-up is $C^{1,\alpha}$ up-to-the-boundary. Put another way, when a blow-up off a level $I$ cone is close to a level $<I$ cone, it is $C^{1,\alpha}$ up-to-the-boundary.

Establishing this $\epsilon$-regularity property for blow-ups is what will in fact take the majority of our work. To prove it we are led naturally to performing a fine blow-up procedure, as is performed in \cite{wickstable} in the flat-setting. In order to say something about the varifolds in this procedure, one must already have established the \textit{varifold} $\epsilon$-regularity theorem (i.e. Theorem \ref{thm:A}) for cones of a lower level. Thus we have an inductive procedure: first prove Theorem \ref{thm:A} for level 0 cones (where the blow-up functions are comprised of only single-valued harmonic functions for which we have a boundary regularity theory), and use this to prove the $\epsilon$-regularity property for blow-ups relative to a level 1 cone, which we can then use to prove Theorem \ref{thm:A} for level 1 cones, and so on. The fine blow-up procedure will be crucial for another reason, namely in establishing a \textit{fine} $\epsilon$-regularity theorem at the varifold level, which will then be used in proving the Theorem \ref{thm:A}. This extra technicality arises because it is possible for a sequence of cones of level $I$ to converge to a cone of level $>I$, and thus one needs to deal with this case at the same time.

The final technicality which arises in this setting is when trying to prove the $\epsilon$-regularity property for blow-ups relative to a level 2 cone. When one performs a fine blow-up procedure, it is possible that in the fine blow-up we still have a two-valued harmonic function: geometrically this is because there are two multiplicity 2 half-hyperplanes in the level 2 cone, and it is possible to have a sequence of level 1 cones converging to it; as such, only one multiplicity two piece ``splits'' in the fine blow-up procedure. As such, one needs to establish a boundary regularity theory for the two-valued harmonic function which arises in the \textit{fine} blow-up class. In order to follow the same procedure as above for the other, simpler, cases, we would need to establish a $\epsilon$-regularity theorem for the fine blow-up class, which requires performing an even finer blow-up procedure. This is what we refer to as an \textit{ultra fine blow-up}. Once such a procedure has been performed, we are left with a class of functions where all two-valued harmonic functions have ``split'' into single-valued harmonic functions, the boundary regularity of which is classical.

We note that, up to the presence of density gaps, our ideas extend readily to any classical cones which are comprised of half-hyperplanes of multiplicity at most $2$. As such, it seems likely the corresponding results to those in \cite{polyhedral} (i.e. allowing for certain $4-way$ and $5$-way junctions) will be true here. When the half-hyperplanes can have multiplicities $>2$, it seems that our arguments have the potential to be extended inductively if one has available suitable regularity theorems near suitable higher multiplicity planes and dimension bounds on the branch sets of corresponding densities. Indeed, in the case where the cone has density $q+\frac{1}{2}$ (as opposed to $2+\frac{1}{2}$ as it is here) armed with such results one could work in a multiplicity $q$ class (see Section \ref{sec:MC}) and inductively extend the finer blow-up procedures. To this end, this work is split into two papers: \cite{minter2021} studies the blow-up (i.e. linear) regularity theory in a more general setting, whilst this paper focuses on developing the non-linear regularity theory, namely proving the properties of the blow-up classes from the varifold level, bringing in the results of \cite{minter2021} in order to prove Theorem \ref{thm:A}.

As a final point of note, we shall heavily rely on the ideas and techniques seen in \cite{simoncylindrical} and \cite{wickstable}. Where possible, we shall avoid unnecessary repetition and refer the reader to these works for the full details of an argument if they are extremely similar, and instead just detail how the proof differs in this setting, allowing us to focus on the key new ideas.

\begin{center}
	\Small ACKNOWLEDGEMENTS
\end{center}

The author would like to thank his PhD supervisor, Neshan Wickramasekera, for numerous discussions surrounding both this work and broader works within the regularity theory of stationary integral varifolds which were invaluable to the current paper. This work was supported by the UK Engineering and Physical Sciences Research Council (EPSRC) grant EP/L016516/1 for the University of Cambridge Centre for Doctoral Training, the Cambridge Centre for Analysis.

\section{Notation and Preliminaries}

\subsection{Basic Notation}

We work in $\R^{n+1}$ throughout. Often we will work with coordinates relative to the $(n-1)$-dimensional spine of a cone, in which case we write $X = (x,y)\in \R^2\times\R^{n-1}$, where the $\R^2$ factor will be coordinates for the cross-section of the cone and the $\R^{n-1}$ factor will be coordinates along its spine; with this notation we write $r(X):= |x|$ and $R(X) := |X|$.

For $x\in \R^{m}$ and $\rho>0$ we write $B^m_{\rho}(x):= \{y\in \R^{m}: |y-x|<\rho\}$ for the open ball of radius $\rho$ centred at $x$. When $m=n+1$ we will often just write $B_\rho(x)$ for $B^{n+1}_\rho(x)$. For simplicity we often write $B^m_\rho$ for $B^m_\rho(0)$. When $A\subset\R^{m}$, we define $\dist(x,A):= \inf_{y\in A}|y-x|$.

We define the \textit{homothety} at $x\in \R^{n+1}$ by scale $\rho>0$ to be the map $\eta_{x,\rho}:\R^{n+1}\to \R^{n+1}$ given by $\eta_{x,\rho}(y):= \rho^{-1}(y-x)$. We also define $\tau_{x}:= \eta_{x,1}$ to be the translation by $x$. For $s\geq 0$ we write $\H^s$ for the $s$-dimensional Hausdorff measure, and $\dim_\H(A)\equiv \dim(A)$ for the Hausdorff dimension of a subset $A\subset\R^{n+1}$. For $A,B\subset\R^{n+1}$ the \textit{Hausdorff distance}, $d_\H$, between $A$ and $B$ is defined by
$$d_\H(A,B):= \max\left\{\sup_{x\in A}\dist(x,B),\ \sup_{x\in B}\dist(x,A)\right\}.$$
For us a \textit{hyperplane} $P$ will be any $n$-dimensional affine subspace of $\R^{n+1}$ and we write $\pi_P$ for the orthogonal projection $\R^{n+1}\to P$. 

We write $G(n+1,n)$ for the \textit{Grassmannian} of $n$-dimensional subspaces of $\R^{n+1}$. An $n$\textit{-varifold} $V$ on an open subset $U\subset\R^{n+1}$ is a Radon measure on $U\times G(n+1,n)$; we write $\|V\|$ for the \textit{weight} or \textit{mass measure} of $V$, which is the Radon measure on $U$ defined by
$$\|V\|(A):= V(A\times G(n+1,n))\ \ \ \ \text{for }A\subset U.$$
We define the \textit{support} of a varifold $V$ by $\spt\|V\|$. We equip the set of $n$-varifolds on $U$ with the \textit{varifold topology}, which is simply the usual topology on Radon measures. It is standard that any countably $n$-rectifiable set $M$ defines an $n$-varifold on $U$, denoted $|M|$, via
$$|M|(A):= \H^n(\{x:(x,T_xM)\in A\})\ \ \ \ \text{for }A\subset U\times G(n+1,n)$$
where $T_x M$ is the approximate tangent space of $M$ at $x$ (which is defined $\H^n$-a.e. on $M$). We say that an $n$-varifold $V$ is an \textit{integral} $n$\textit{-varifold} if we can write $V = \sum^\infty_{j=1}c_j|M_j|$ for some $(c_j)_j\subset\N$ and $n$-rectifiable sets $M_j$. For $f:U\to \tilde{U}$ a $C^1$ function with $f|_{\spt\|V\|\cap U}$ proper, we write $f_\# V$ for the \textit{image varifold}, or \textit{pushforward}, of $V$ under $f$. We write $T_x V$ for the \textit{approximate tangent plane} of $V$ at $x$, which we know exists $\H^n$-a.e. in $\spt\|V\|$. We define the \textit{regular part} of $V$, denoted $\reg(V)$, to be the set of points $x\in \spt\|V\|$ such that $\exists\rho>0$ for which $\spt\|V\|\cap B_\rho^{n+1}(x)$ is a smooth embedded submanifold in $B_\rho^{n+1}(0)$. We then write $\sing(V)$ for the \textit{interior singular set} of $V$, i.e., $\sing(V):= \left(\spt\|V\|\backslash\reg(V)\right)\cap B^{n+1}_2(0)$.

We shall exclusively work with integral $n$-varifolds $V$ on $B^{n+1}_2(0)$; minor modifications of our arguments can be made to extend our results to more general settings, although we do not present these here to avoid the additional technical complications which arise.

\subsection{Some Varifold Preliminaries}\label{sec:basics}

Recall the definition of a \textit{stationary} integral varifold as defined in $(\S1)$ of the Introduction. By a suitable choice of $\psi$ in the first variation formula (see \cite{simongmt}) we see that if $V$ is a stationary integral varifold in $B^{n+1}_2(0)$, then for any $z\in B^{n+1}_2(0)$ and $0<\sigma<\rho<2-|z|$ one has the \textit{monotonicity formula}
\begin{equation}\label{eqn:monotonicity}
	\frac{\|V\|(B_\rho(z))}{\rho^n} - \frac{\|V\|(B_\sigma(z))}{\sigma^n} = \int_{B_\rho(z)\backslash B_\sigma(z)}\frac{|X^\perp|^2}{|X|^{n+2}}\ \ext\|V\|(X)
\end{equation}
where $X^\perp := X - \pi_{T_XV}(X)$ is the projection of $X$ onto the orthogonal complement $T^\perp_X V$. The monotonicity formula implies that for each $z\in B^{n+1}_2(0)$ the function $\rho\mapsto\frac{\|V\|(B_\rho(z))}{\rho^n}$ is a monotonically non-decreasing for $\rho\in (0,2-|x_0|)$. In particular, the \textit{density} of $V$ at $z$, i.e.
$$\Theta_V(z):= \lim_{\rho\downarrow 0}\frac{\|V\|(B_\rho(z))}{\w_n\rho^n}$$
is well-defined everywhere in $\spt\|V\|$ (here, $\w_n = \H^n(B_1^n(0))$ is the volume of the $n$-dimensional unit ball in $\R^n$). It also follows that $\Theta_V(z)$ is an upper-semicontinuous function of both $z$ and $V$ (with respect to the Euclidean topology on $\R^{n+1}$ and the varifold topology, respectively). Taking $\sigma\downarrow 0$ in the monotonicity formula (\ref{eqn:monotonicity}) we find
\begin{equation}\tag{2.3'}\label{eqn:monotonicity2}
	\frac{\|V\|(B_\rho(z))}{\w_n\rho^n} - \Theta_V(z) = \frac{1}{\w_n}\int_{B_\rho(z)}\frac{|X^\perp|^2}{|X|^{n+2}}\ \ext\|V\|(X).
\end{equation}
We shall refer to the integral on the right-hand side of (\ref{eqn:monotonicity2}) as the \textit{mass drop}. Finding a suitable bound for the mass drop for $V$ sufficiently close to certain cones will be crucial to our $L^2$ estimates later on.

For $X\in \spt\|V\|$, we write $\vartan_X(V)$ for the set of all \textit{tangent cones} to $V$ at $X$, i.e. the set of all varifold limits $\BC \equiv \lim_{j\to\infty}(\eta_{X,\rho_j})_\# V$ for some $\rho_j\downarrow 0$. We know from standard compactness theorems for varifolds (see \cite{simongmt}) that each $\BC\in \vartan_X(V)$ is a stationary integral varifold in $\R^{n+1}$, which by the monotonicity formula is a \textit{cone}, i.e. $(\eta_{0,\rho})_\#\BC = \BC$ for every $\rho>0$. Due to this homogeneity property of a tangent cone $\BC$, the set of points in $\spt\|\BC\|$ under which $\BC$ is translation invariant, i.e.
$$S(\BC) := \{x\in \spt\|\BC\|: (\tau_x)_\#\BC = \BC\}$$
is a subspace of $\R^{n+1}$, called the \textit{spine} of $\BC$. Moreover, from the upper semi-continuity of the density it follows that $S(\BC) = \{x\in \spt\|\BC\|:\Theta_{\BC}(x) = \Theta_{\BC}(0)\}$. We can therefore always find a rotation $q$ of $\R^{n+1}$ such that\footnote{By this product notation, $\BC_0\times\R^k$, we mean the varifold whose support is $\spt\|\BC_0\|\times\R^k$ and with density function $\Theta_{\BC_0\times\R^k}(x,y) = \Theta_{\BC_0}(x)$.}
$$q_\# \BC = \BC_0\times \R^{\dim(S(\BC))}$$
where $\BC_0$ is a stationary integral cone in $\R^{n+1-\dim(S(\BC))}$.

\begin{defn}
	Let $V$ be a stationary integral $n$-varifold in $\B^{n+1}_2(0)$. We say that $X\in \sing(V)$ is a \textit{branch point} if at least one tangent cone to $V$ at $X$ is supported on a hyperplane, yet there is no neighbourhood of $X$ on which $\spt\|V\|$ is a union of finitely many embedded submanifolds. We write $\B$ for the set of branch point singularities of $V$, and $\B_q:= \B\cap \{\Theta_V= q\}$ for the branch points of density $q\in\{2,3,\dotsc\}$.
\end{defn}
Note that Allard's regularity theorem (\cite{allard}) tells us that branch points have density $\geq 2$.

For $j\in \{0,1,\dotsc,n-1\}$ we define the $j^{\text{th}}$\textit{-stratum} of $\sing(V)$, denoted $\S_j$, by
$$\S_j:= \{x\in \sing(V): \dim(S(\BC))\leq j\ \ \forall \BC\in \vartan_X(V)\}.$$
Almgren's stratification theorem (\cite{almgren}) tells us that $\dim(\S_j)\leq j$ for each such $j$ (in fact $\S_j$ is countably $j$-rectifiable for each $j$ by Naber--Valtorta (\cite{naber2015singular})). Therefore, we can write $\sing(V)$ as a disjoint union
$$\sing(V):= \B\cup \tilde{\B}\cup \bigcup^{n-1}_{j=0}(\S_j\backslash \S_{j-1})$$
where $\S_{-1} = \emptyset$ and $\tilde{\B}$ denotes those singular points which are not branch points yet one tangent cone is supported on a hyperplane; necessarily by standard quasilinear elliptic PDE theory (\cite{gt}), we have $\dim(\tilde{\B})\leq n-2$. Every point in $\S_j\backslash \S_{j-1}$ has the property that every tangent cone has spine dimension at most $j$, \textit{and} that there is at least one tangent cone with spine dimension equal to $j$. Thus from Almgren's stratification theorem we see that $\B$ is the only part of the singular set which could have dimension $>n-1$; indeed branch point singularities are the primary difficulty in understanding the singular set.

For $\BC\in \vartan_X(V)$ it is clear that whenever $\sing(\BC)\neq\emptyset$, i.e. $\BC$ is not supported on a hyperplane, we have $S(\BC)\subset\sing(\BC)$. We shall say that a stationary integral cone $\BC$ is a \textit{cylindrical cone} if $\sing(\BC) = S(\BC)$, and as such we can write (up to a rotation) $\BC = \BC_0\times\R^{\dim(S(\BC))}$, where $\BC_0\subset\R^{n+1-\dim(S(\BC))}$ is a stationary integral cone with $\sing(\BC_0) = \{0\}$, i.e. $\BC_0$ has an isolated singularity. Moreover, we shall say that $\BC$ is a \textit{classical tangent cone} if $\dim(S(\BC)) = n-1$. It follows easily that a classical tangent cone is necessarily cylindrical with $\BC_0$ being a finite collection of rays through the origin, which in particular means that $\Theta_{\BC}(0) = q/2$ for some $q\in \Z_{\geq 3}$, and so $\BC$ is comprised of finitely many half-hyperplanes with some integer multiplicities meeting along a common axis. Thus, $\S_{n-1}\backslash\S_{n-2}$ is the set of non-branch point singularities which have at least one classical tangent cone arising as a tangent cone.

The only subsets of $\sing(V)$ which can have dimension $\geq n-1$ are $\B$ and $\S_{n-1}\backslash \S_{n-2}$, which necessarily have density taking values in $\{\frac{3}{2},2,\frac{5}{2},\dotsc\}$. L.~Simon's $\epsilon$-regularity theorem already provides us with the appropriate understanding when $X\in \S_{n-1}\backslash\S_{n-2}$ has density $\frac{3}{2}$. In the case of $V\in \S_2$, \cite{minterwick} understands the case when $\Theta_V(X) = 2$ and $X\in \B\cup (\S_{n-1}\backslash\S_{n-2})$, and thus the next case to understand is when $X\in \S_{n-1}\backslash\S_{n-2}$ has $\Theta_V(X) = \frac{5}{2}$; this is Theorem \ref{thm:A} and the aim of our work here.

\subsection{Two-Valued Functions}

Two-valued functions will be used to model the behaviour of $V\in \S_2$ near multiplicity two branch points. We recall the key definitions and properties here. More information on multi-valued functions can be found in \cite{almgren}, \cite{de2010almgren}, whilst more specifics for two-valued functions -- which we will make use of -- can be found in \cite{simon2016frequency}, \cite{krummelwick1}, and \cite{krummelwick3}.

We write $\A_2(\R^m)$ for the space of unordered pairs $x = \{x_1,x_2\}$, where $x_1,x_2\in \R^m$ ($x_1,x_2$ are not necessarily distinct). We make $\A_2(\R^m)$ into a metric space by endowing it with the metric
$$\G(x,y):= \min\left\{\sqrt{|x_1-y_1|^2 +|x_2-y_2|^2}, \sqrt{|x_1-y_2|^2+|x_2-y_1|^2}\right\}.$$
For each $x = \{x_1,x_2\}\in \A_2(\R^m)$ we set $|x|:= \G(x,\{0,0\}) \equiv \sqrt{|x_1|^2+|x_2|^2}$. It is important to note that since there is no well-defined notion of \lq\lq addition'' for unordered pairs $\A_2(\R^m)$ is not a vector space in any natural way, however we shall sometimes abuse notation and write $f+g$ when $f$ is single-valued and $g = \{g_1,g_2\}$ is two-valued to mean the two-valued function $\{g_1+f,g_2+f\}$.

\begin{defn}
	For $U\subset\R^n$ open, a \textit{two-valued function} $u$ is a map $u:U\to \A_2(\R^m)$. We write $u(X):= \{u_1(X),u_2(X)\}$ for each $X\in U$, for some $u_1(X),u_2(X)\in \R^m$.
\end{defn}

Every two-valued function is determined uniquely by two functions: the (single-valued) \textit{average part} $u_a:U\to \R^m$ and the (two-valued) \textit{symmetric part} $u_s:U\to \A_2(\R^m)$, defined by
$$u_a(X):=\frac{u_1(X)+u_2(X)}{2}\ \ \ \ \text{and}\ \ \ \ u_s(X):= \left\{\pm\frac{u_1(X)-u_2(X)}{2}\right\}$$
where $u(X) = \{u_1(X),u_2(X)\}$. In general we say that a two-valued function $v$ is \textit{symmetric} if $v_a\equiv 0$, in which case we can write $v(X) \equiv v_s(X) = \{\pm \phi(X)\}$ for all $X\in U$, for some single-valued $\phi$.

Since we have a metric on $\A_2(\R^m)$ we can define notions of continuity and differentiability for two-valued functions, and thus we can define (metric) spaces of two-valued functions such as $C^0(U;\A_2(\R^m))$, $C^{1,\alpha}(U;\A_2(\R^m))$, $L^p(U;\A_2(\R^m))$, and so on; we omit the details and refer the reader to e.g. \cite{simon2016frequency}.

For $u\in C^1(U;\A_2(\R^m))$ we define sets
$$\CZ_u:= \{X\in U: u_1(X) = u_2(X)\}$$
and
$$\CK_u:= \{X\in U: u_1(X) = u_2(X)\ \text{and}\  Du_1(X) = Du_2(X)\}.$$

\begin{defn}
	The \textit{branch set} $\B_u$ for a two-valued function $u$ is the set of $Y\in U$ for which there is no $\rho\in (0,\dist(Y,\del U))$ such that on $B_\rho(Y)$ we can write $u(X) = \{u_1(X),u_2(X)\}$ for some (single-valued) $C^1$ functions $u_1,u_2:B_\rho(Y)\to \R^m$.
\end{defn}

Clearly we have $\B_u\subset\CK_u\subset\CZ_u$ and moreover
$$\CZ_u = \CZ_{u_s}\equiv \{X\in U: u_s(X) = \{0,0\}\}$$
$$\CK_u = \CK_{u_s}\equiv \{X\in U: u_s(X) = \{0,0\}\ \ \text{and}\ \ Du_s(X) = \{0,0\}\}.$$

\subsection{Two-Valued Harmonic Functions}\label{sec:two-valued_harmonic}

Two-valued harmonic functions play the same role for $V\in\S_2$ near multiplicity two planes as single-valued harmonic functions do for arbitrary stationary integral varifolds near multiplicity one planes, namely they provide the appropriate linear theory in order to understand blow-ups.

\begin{defn}
	Let $\alpha\in (0,1]$ and $U\subset\R^n$ be open. Then we say $u\in C^{1,\alpha}(U;\A_2(\R^m))$ is \textit{locally harmonic} in $U\backslash \B_u$, or is \textit{two-valued harmonic} in $U$, if for every $B_\rho(X_0)\subset U\backslash \B_u$, there is a (unique) pair of single-valued harmonic functions $u_1,u_2: B_\rho(X_0)\to \R^m$ such that $u(X) = \{u_1(X),u_2(X)\}$ for all $X\in B_{\rho}(X_0)$.
\end{defn}

It is possible to show that whenever $v$ is a \textit{symmetric} two-valued $C^{1,\alpha}$ harmonic function, then (see \cite{simon2016frequency}):
\begin{enumerate}
	\item [(i)] either $v\equiv \{0,0\}$ on $U$ or $\dim(\CK_v)\leq n-2$; moreover either $\B_v =\emptyset$ or $\dim(\B_v) = n-2$ with $\H^{n-2}(\B_v) >0$ (moreover $\B_v$ is countably $(n-2)$-rectifiable, from \cite{krummelwick1});
	\item [(ii)] in fact $v\in C^{1,1/2}(U;\A_2(\R^m))\cap W^{2,2}_{\text{loc}}(U;\A_2(\R^m))$ with the estimate
	\begin{equation}\label{eqn:2harmonic_bounds}
		\sup_{B_{\rho/2}(Y)}|v| + \rho\sup_{B_{\rho/2}(Y)}|Dv| + \rho^{3/2}[Dv]_{1/2;B_{\rho/2}(Y)} \leq C\rho^{-n/2}\|v\|_{L^2(B_\rho(Y))}
	\end{equation}
	for every $B_\rho(Y)$ with $\overline{B_\rho(Y)}\subset U$, where $C = C(n,m)$.
\end{enumerate}

In particular, for any two-valued $C^{1,\alpha}$ harmonic function $u$ we see that $\dim(\B_u)\leq n-2$, and thus the average part $u_a$ is always a single-valued harmonic function on \textit{all} of $U$.

One crucial difference between single-valued and two-valued harmonic functions for this work is that there is currently no known general boundary regularity theory for two-valued harmonic functions, unlike in the single-valued case where we have classical boundary regularity results from elliptic PDE theory (see e.g. \cite{gt}, \cite{morrey1966multiple}). Even a reflection principle is unclear (it should be stressed that the reason the usual reflection principle for harmonic functions is so powerful is because it is not necessary to assume any control on the derivatives at the boundary). Instead, we will have to establish the boundary regularity in a different way, appealing to Campanato-style results and integral estimates, which we establish by classifying the homogeneous degree one two-valued harmonic functions defined on a half-plane and using Hardt--Simon inequality arguments (see e.g. \cite[Section 4]{wickstable} and \cite{minter2021}).

\subsection{Two-Valued Stationary Graphs}\label{sec:two-valued_stationary_graphs}

Two-valued stationary graphs will provide the graphical representation for $V\in \S_2$ near a multiplicity two plane. Let $\alpha\in (0,1]$ and let $U\subset\R^n$ be open. For a two-valued function $u\in C^{1,\alpha}(U;\A_2(\R^m))$ the \textit{graph} of $u$ is defined by
$$\graph(u):= \{(X,Y)\in U\times\R^m: Y = u_1(X)\ \ \text{or}\ \ Y= u_2(X)\}.$$
We can associate to $\graph(u)$ an $n$-dimensional varifold $V_u := (\graph(u),\theta_u)$ where the multiplicity function $\theta_u:\graph(u)\to \Z_{\geq 1}$ is defined by
$$\theta_u(X,Y) := \begin{cases}
2 & \text{if }u_1(X) = u_2(X)\\
1 & \text{if }u_1(X) \neq u_2(X)
\end{cases}$$
for $Y\in \{u_1(X),u_2(X)\}$.

Note that $V_u$ is determined by the function $x\mapsto (x,u(x))\in \R^n\times\A_2(\R^m)$, and so if we define the two-valued Jacobian function by $J(x):= \{J_1(x),J_2(x)\}$ where
$$J_i(x):= \det\left[(\delta_{pq}+D_p u_i(x)\cdot D_q u_i(x))_{p,q}\right]^{1/2}$$
for $i\in \{1,2\}$, then from the area formula (see \cite{de2013multiple}) we have
\begin{equation}\label{eqn:area_formula}
	\int_{A\times \R^m}g(x)\ \ext\|V_u\|(x) = \int_A\sum^2_{i=1}g(x,u_i(x))J_i(x)\ \ext x
\end{equation}
for any measurable $A\subset U$ and bounded compactly supported Borel function $g:U\times\R^m\to \R$. Moreover we have $1\leq J_i \leq 1+C|Du|^2$, where $C = C(n,m)$.

\begin{defn}
	We say that $u\in C^{1,\alpha}(U;\A_2(\R^m))$ is a \textit{stationary} two-valued graph in $U\times\R^m$ if $V_u:= (\graph(u),\theta_u)$ is a stationary varifold in $U\times\R^m$.
\end{defn}

In this setting we can define analogous sets to the two-valued harmonic setting, namely $\CZ_{\graph(u)}$ and $\CK_{\graph(u)}$, on the graph level instead of the domain level via:
$$\CZ_{\graph(u)} := \theta^{-1}(\{2\})$$
$$\CK_{\graph(u)} := \{Z\in \CZ_{\graph(u)}: \exists\ \text{a multiplicity two tangent plane at } Z\}$$
and then in this setting we define $\sing(\graph(u))$ to be the set of points $Z\in \graph(u)$ such that there is no $\rho>0$ such that $\graph(u)\cap B_\rho^{n+m}(Z)$ is a finite union of smoothly embedded submanifolds.

Again from \cite{simon2016frequency} we know that $\exists\epsilon = \epsilon(n,m)$ such that if $\|u\|_{C^{1,\alpha}}<\epsilon(n,m)$ then either $\CK_{\graph(u)} = \graph(u)$ or $\dim(\CK_{\graph(u)})\leq n-2$. Moreover $\sing(\graph(u)) = \emptyset$ or $\dim(\sing(\graph(u)))=n-2$ with $\H^{n-2}(\sing(\graph(u)))>0$. We also have in fact that $u_a\in C^{1,1}(U;\R^m)$ and $u_s\in C^{1,1/2}(U;\A_2(\R^m))\cap W^{2,2}_{\text{loc}}(U;\A_2(\R^m))$ with the estimates:
\begin{equation}\label{eqn:two-valued_estimate_1}
	\sup_{B_{\rho/2}(Y)}|u_a| + \rho\sup_{B_{\rho/2}(Y)} |Du_a| + \rho^2\sup_{B_{\rho/2}(Y)}|D^2u| \leq C\rho^{-n/2}\|u_a - u_a(Y)\|_{L^2(B_\rho(Y))}
\end{equation}
\begin{equation}\label{eqn:two-valued_estimate_2}
	\sup_{B_{\rho/2}(Y)}|u_s| + \rho\sup_{B_{\rho/2}(Y)}|Du_s| + \rho^{3/2}[Du_s]_{1/2;B_{\rho/2}(Y)} \leq C\rho^{-n/2}\|u_s\|_{L^2(B_\rho(Y))}
\end{equation}
for each $B_\rho(Y)$ with $\overline{B_\rho(Y)} \subset U$, where $C = C(n,m)$. Moreover if $U = B_1(0)$ then we have for every $X\in B_{1/2}(0)$ with $d(X):= \dist(X,\CK_u) \leq \frac{1}{4}$,
\begin{equation}
	|u_s(X)| + d(X)|Du_s(X)| + d(X)^2|D^2u_s(X)| \leq Cd(X)^{3/2}\|u_s\|_{L^2(B_1(0))}.
\end{equation}
Moreover from \cite{krummelwick3}, we know that the branch set of $V_u$ is always countably $(n-2)$-rectifiable.

\subsection{Some Regularity Results for Stable Codimension One Varifolds}\label{sec:regularity}

From now on we shall be focused on the class $\S_2$, i.e. stationary integral varifolds in $B^{n+1}_2(0)$ which have stable regular part (in the sense of $(\S2)$) and contain no triple junction singularities. We first recall the two key results for the class $\S_2$ from \cite{minterwick}. 

\begin{theorem}[\cite{minterwick},  Theorem C, Branched Case]\label{thm:wick1}
	Fix $\delta\in (0,1)$. Then $\exists\epsilon = \epsilon(n,\delta)\in (0,1)$ such that, whenever $V\in \S_2$ obeys:
	\begin{enumerate}
		\item [(i)] $2-\delta\leq \w_n^{-1}\|V\|(\R\times B_1^n(0)) < 3-\delta$;
		\item [(ii)] $\hat{E}_V<\epsilon$, where $\hat{E}_V$ is the (one-sided) height excess of $V$ relative to the hyperplane $\R^n\times\{0\}$, i.e.
		$$\hat{E}^2_V := \int_{B_1^n(0)\times\R} |x^{n+1}|^2\ \ext\|V\|;$$
	\end{enumerate}
	then we have $V\res (B^n_{1/2}(0)\times \R) = \graph(u)$, i.e., locally $V$ is expressible as the graph of a (stationary) two-valued $C^{1,1/2}$ function $u:B^n_{1/2}(0)\times\R\to \A_2(\R)$ satisfying $\|u\|_{C^{1,1/2}(B_{1/2}(0))}\leq C\hat{E}_V$, where $C = C(n)$.
\end{theorem}

\begin{theorem}[\cite{minterwick}, Theorem C, Transverse Case]\label{thm:wick2}
	Fix $\delta\in (0,1)$ and a cone of the form $\BC = |P_1| + |P_2|$ where $P_1,P_2$ are distinct hyperplanes in $\R^{n+1}$. Then $\exists\epsilon = \epsilon(n,\delta,\BC)$ such that whenever $V\in \S_2$ obeys:
	\begin{enumerate}
		\item [(i)] $\w_n^{-1}\|V\|(B_1^{n+1}(0)) < 3-\delta$;
		\item [(ii)] $\hat{Q}_{V,\BC}<\epsilon$, where $\hat{Q}_{V,\BC}$ is the two-sided height excess of $V$ relative to $\BC$, i.e.
		$$\hat{Q}_{V,\BC}^2 := \int_{B_1^{n+1}(0)} \dist^2(X,\spt\|\BC\|)\ \ext\|V\| + \int_{B^{n+1}_{1/2}(0)\backslash \{r_\BC(X)>1/16\}} \dist^2(X,\spt\|V\|)\ \ext\|\BC\|$$
	\end{enumerate}
	where $r_\BC(X):= \dist(X,S(\BC))$, then we have $V\res B^{n+1}_{1/2}(0) = |\graph(u_1)| + |\graph(u_2)|$, where for $j=1,2$, $u_j:P_j\cap B^{n+1}_{1/2}(0)\to P_j^\perp$ is a $C^2$ function satisfying $\|u_j\|_{C^2(P_j\cap B^{n+1}_{1/2}(0))} \leq C\hat{Q}_{V,\BC}$, where $C = C(n)$.
\end{theorem}

Theorem \ref{thm:wick1} and Theorem \ref{thm:wick2} are sharp with respect to their respective hypotheses. In Theorem \ref{thm:wick1} we use a one-sided height excess along with a lower bound on the mass; this is because the support of a plane is indecomposable, and thus as long as we have the mass lower and upper bound, smallness of the one-sided excess implies closeness as varifolds. In Theorem \ref{thm:wick2} we need to work with a two-sided height excess since the $\BC$ in question is decomposable, and so we need to prevent $V$ looking like just one of the two planes (as it could be a multiplicity one or two version of the single plane). This is also why we do not need to assume any mass lower bound in Theorem \ref{thm:wick2}, since the upper bound is enough to know $V$ is multiplicity one away from the spine of $\BC$ once we know it is close to all of $\BC$.

From Theorem \ref{thm:wick1} and \cite{krummelwick3} we then have the following important corollary:

\begin{corollary}\label{cor:multiplicity_2_branch_points}
	Let $V\in\S_2$. Then $\B_2$, the set of branch points of density $2$ in $V$, is countably $(n-2)$-rectifiable; in particular $\dim_\H(\B_2)\leq n-2$.
\end{corollary}

\subsection{A Unique Continuation Property}

Later, we will need to construct a two-valued stationary graph relative to a plane by patching together two-valued stationary graphs over smaller open regions. To do this, we will use the following unique continuation principle for $C^{1,\alpha}$ stationary two-valued graphs:

\begin{lemma}\label{lemma:unique_continuation}
	Let $U\subset\R^n$ be open, and suppose $u_1,u_2\in C^{1,\alpha}(U;\A_2(\R))$ are both stationary two-valued graphs. Then if there is an open subset $V\subset U$ for which $\left. u_1\right|_V \equiv \left. u_2\right|_V$, then we have $u_1\equiv u_2$.
\end{lemma}

\begin{proof}
	Firstly, it is well-known, albeit hard to find in the literature (see e.g. \cite[Lemma 2.9]{hiesmayr2020bernstein}) that the varifold $\mathbf{v}(u)$ associated to a $C^{1,\alpha}$ stationary two-valued graph $u:U\to \A_2(\R)$ is stable. We know from \cite{simon2016frequency} that moreover $\dim_\H(\B_{\mathbf{v}(u_i)})\leq n-2$ for $i=1,2$, and thus Almgren's stratification of the singular set gives
	$$\sing(\mathbf{v}(u_i)) = S_i\cup B_i$$
	where for each $i=1,2$, $B_i$ is a relatively closed set with $\dim_\H(B_i)\leq n-2$ and each $X\in \S_i$ has the property that locally about $X$, $\mathbf{v}(u_i)$ is a union of two smoothly embedded transverse hypersurfaces (note that triple junctions do not occur in two-valued stationary graphs).
	
	Now suppose there is a point $X\in U\backslash \pi(B_1\cup B_2)$ where $u_1(X)\neq u_2(X)$; here, $\pi:U\times\R\to U\times\{0\}\cong U$ is the orthogonal projection. It suffices to find a contradiction to this, as then we have $\left. u_1\right|_{U\backslash\pi(B_1\cup B_2)} \equiv \left. u_2\right|_{U\backslash \pi(B_1\cup B_2)}$, which then implies that $u_1\equiv u_2$ by continuity and the fact that $\dim_\H(\pi(B_1\cup B_2))\leq n-2$. We may also assume that $U$ is bounded by restricting to an appropriate ball containing $X$ and an open subset of $V$. In particular, this dimension bound on $B_1\cup B_2$ implies that $U\backslash \pi(B_1\cup B_2)$ is connected, and so it is path connected (as it is open), and thus choosing $X_0\in V\backslash\pi(B_1\cup B_2)$ we may find a path $\gamma:[0,1]\to U\backslash \pi(B_1\cup B_2)$ with $\gamma(0) = X_0$ and $\gamma(1) = X$.
	
	Set $A:= \{t\in [0,1]: u_1 = u_2 \text{ on a neighbourhood of }\gamma(t)\}$. Then we know that, by our assumption on $V$, that there is some $\epsilon>0$ for which $[0,\epsilon)\subset A$. Let $t_0:=\sup\{t:[0,t)\subset A\}$; clearly $[0,t_0)\subset A$. We claim that if $t_0<1$ we have a contradiction. Indeed, we know from continuity that $u_1(t_0) = u_2(t_0)$ and, by construction, we know that there is a $\rho>0$ with $B_\rho(\gamma(t_0))\subset U\backslash\pi(B_1\cup B_2)$ and moreover for which on $B_\rho(\gamma(t_0))$ we have $u_i = \llbracket u^{(1)}_i\llbracket + \llbracket u^{(2)}_i\rrbracket$ for $i=1,2$, where $u_i^{(j)}$ is a smooth solution of the minimal surface equation. However, by definition of $A$ and $t_0$, we know that $u_1^{(i)}$ agrees with $u_2^{(j_i)}$ on an open subset of $B_\rho(\gamma(t_0))$ for some $j_i\in \{1,2\}$, and so by unique continuation of solutions to the minimal surface equation, they must agree on all of $B_\rho(\gamma(t_0))$, i.e. there exists $t>t_0$ for which $t\in A$, contradicting the definition of $t_0$. Thus we must have $t_0=1$, and so by continuity we must have $u_1(\gamma(1)) = u_2(\gamma(1))$, i.e. $u_1(X) = u_2(X)$, providing the necessary contradiction and proving the result.
\end{proof}

Making the obvious modifications in the above proof and using instead the results from Section \ref{sec:two-valued_harmonic} and the unique continuation principle for single-valued harmonic functions, we can similarly prove a unique continuation principle for two-valued $C^{1,\alpha}$ harmonic functions (see also \cite{minter2021}):

\begin{lemma}\label{lemma:unique_continuation_harmonic}
	Let $U\subset\R^n$ be open, and suppose that $u_1,u_2\in C^{1,\alpha}(U;\A_2(\R))$ are both two-valued harmonic functions. Then if there is an open subset $V\subset U$ on which $\left.u_1\right|_V \equiv \left.u_2\right|_V$, then we have $u_1\equiv u_2$.
\end{lemma}

\subsection{Classes of Varifolds}

We now describe the set up for the proof of Theorem \ref{thm:A}. Note that we know $\S_2$ is a closed class, i.e. any limit point of this class also belongs to $\S_2$.

From the regularity theories already described, for $V\in \S_2$ we already understand the singular set well in the region $\{\Theta_V<\frac{5}{2}\}$. In the region $\{\frac{5}{2}\leq \Theta_V<3\}$, the only points we do not yet understand which could create an $(n-1)$-dimensional singular set are those $X\in \sing(V)$ for which $\exists \BC\in \vartan_X(V)$ a classical tangent cone with $\Theta_{\BC}(0) = \frac{5}{2}$. Our aim in proving Theorem \ref{thm:A} is essentially to understand the behaviour of $V\in \M$ near such a $\BC$, thus completing the analysis of the top-dimensional part of the singular set in the region $\{\Theta_V<3\}$.

So take such a cone $\BC$, and rotate so that we can without loss of generality write $\BC = \BC_0\times\R^{n-1}$; thus $\BC_0$ is a stationary $1$-dimensional cone in $\R^2$ which has $\Theta_{\BC_0}(0) = \frac{5}{2}$, and thus $\BC_0$ is comprised of $5$ rays from the origin (counted with multiplicity), which could coincide. 

\textbf{Remark:} For $\BC_0\subset\R^2$ a 1-dimensional stationary integral cone with $\Theta_{\BC_0}(0) = \frac{5}{2}$, if we write $\{n_1,\dotsc,n_k\}$ for the unit vectors in the (outward) directions of the rays of $\spt\|\BC_0\|$ and $\{\theta_1,\dotsc,\theta_k\}\subset\Z_{\geq 1}$ for the multiplicity of each ray respectively, then the stationary condition requires that $\sum^k_{i=1}\theta_i n_i = 0$, whilst the density condition implies that $\sum^k_{i=1}\theta_i = 5$. It follows immediately from these two facts that $\theta_i\in \{1,2\}$ for each $i$, and thus $k\in\{3,4,5\}$. A simple calculation shows that when $k=3$, up to an orthogonal rotation of $\R^2$, $\spt\|\BC_0\|$ is completely determined, whilst when $k\in \{4,5\}$ there is a $(k-3)$-parameter family of possible cones, which is a closed family except for the possibility of two or more rays coinciding, giving rise to a cone with a fewer number of distinct rays; in particular, the collection of all cones supported on at most $p$ rays, for each $p\in \{3,4,5\}$, is a closed class.

The above remark implies that $\BC_0$ must be comprised of either: (i) five multiplicity 1 rays, or (ii) three multiplicity 1 rays and one multiplicity 2 ray, or (iii) one multiplicity 1 ray and two multiplicity 2 rays. Each case will need its own consideration in the proof of Theorem \ref{thm:A}, and indeed we will need to know that Theorem \ref{thm:A} is true for $\BC_0$ supported on $>p$ rays when proving the result when $\BC_0$ is only supported on $p$ rays.

\begin{defn}
	Fix $I\in \{0,1,2\}$. We say that a cone $\BC$ is \textit{level }$I$ if it is a classical cone with $\Theta_\BC(0) = \frac{5}{2}$ and $\spt\|\BC\|$ is comprised of $5-I$ distinct half-hyperplanes (equivalently, $\BC$ contains $I$ half-hyperplanes of multiplicity $2$). We write $\FL_I$ for the set of cones of level $I$.
\end{defn}

\textbf{Remark:} We make no stationarity assumption on the cones in $\FL_I$.

Set $\FL:= \FL_0\cup \FL_1\cup \FL_2$; the set $\FL$ comprises of all cones which we are interested in for the proof of Theorem \ref{thm:A}. We write $\FL_S\subset\FL$ for the set of cones in $\FL$ which are also stationary as varifolds.

There are two height excesses which we shall need for the proof of Theorem \ref{thm:A}, namely\footnote{A mass lower bound is of course natural for our setting, as we do not want to consider the situation where $V$ is close to a subcone to $\BC$, such as a multiplicity one hyperplane or triple junction. As long as the total mass of $V$ is in some $(2\w_n+\delta,3\w_n-\delta)$, we will be fine.\label{footnote-1}}:
\begin{defn}
	For $V\in \S_2$ and $\BC\in \FL$, the \textit{one-sided height excess} of $V$ relative to $\BC$ is
	$$E_{V,\BC}^2 := \int_{B_1}\dist^2(X,\spt\|\BC\|)\ \ext\|V\|$$
	and the \textit{two-sided height excess} of $V$ relative to $\BC$ is
	$$Q_{V,\BC}^2:=\int_{B_1}\dist^2(X,\spt\|\BC\|)\ \ext\|V\| + \int_{B_{1/2}\backslash\{r_{\BC}<1/16\}}\dist^2(X,\spt\|V\|)\ \ext\|\BC\|$$
	where $r_\BC(X) := \dist(X,S(\BC))$ (which is just $|x|$ if $X = (x,y)\in \R^2\times\R^{n-1}$ and $S(\BC) = \{0\}^2\times\R^{n-1}$).
\end{defn}
As mentioned in the above footnote (\footref{footnote-1}), the height excess can only tell us how close the support of $V$ is to $\spt\|\BC\|$, and so to ensure a varifold is close in the varifold topology to a given cone, we need an assumption on the mass. This leads us to define a class of nearby varifolds for each $\BC\in \FL$:

\begin{defn}
	For $\BC\in \FL$ and $\epsilon>0$, define $\CN_\epsilon(\BC)$ to be the class of $V\in \S_2$ which have $\frac{5}{2}\w_n - 1\equiv\|\BC\|(B_1(0)) - 1\leq \|V\|(B_1(0)) \leq \|\BC\|(B_1(0)) + 1 \equiv \frac{5}{2}\w_n + 1$ and $E_{V,\BC}<\epsilon$.
\end{defn}

For $\epsilon>0$, we also define the class $\FL_\epsilon(\BC)$ of nearby cones to a given cone $\BC\in \FL$ in the following manner: $\BC^\prime\in \FL_\epsilon(\BC)$ if $S(\BC^\prime) = S(\BC)$ and, after performing a rotation so that $\BC = \BC_0\times\R^{n-1}$ and $\BC^\prime = \BC_0^\prime\times\R^{n-1}$, if we write $\BC_0 = \sum^5_{i=1}|\ell_i|$ for some rays $\ell_i$ through $0\in \R^2$, then $\BC_0^\prime = \sum^5_{i=1}(q_i)_\#|\ell_i|$ for some $q_i\in SO(2)$ with $|q_i-\id|<\epsilon$, for $\id:\R^2\to \R^2$ the identity map.

\textbf{Remark:} If $\BC\in \FL_I$, then there is $\epsilon = \epsilon(\BC)>0$ such that $\FL_\epsilon(\BC)$ only contains cones of level at most $I$. Moreover, if $\BC\in \FL\backslash\FL_S$ is not stationary, then there is $\epsilon = \epsilon(\BC)>0$ such that $\CN_{\epsilon}(\BC) = \emptyset$.

We now prove that varifolds in $\CN_\epsilon(\BC)$ are close, as varifolds, to $\BC$ for $\epsilon>0$ sufficiently small, in the following sense:

\begin{lemma}\label{lemma:convergence}
	Fix $\BC\in \FL_S$. If $V_i\in \CN_{\epsilon_i}(\BC)$ with $\epsilon_i\downarrow 0$, we have $V_i\weakly \BC$ (i.e. as varifolds in $B_1$).
\end{lemma}

\begin{proof}
	It suffices to show that every subsequence of $(V_i)_i$ has a further subsequence which converges to $\BC$ as varifolds; so let us suppose, without relabelling the sequence, we have already passed to some subsequence and that this is $(V_i)_i$. Then, from the compactness properties of the class $\S_2$, we may pass to another subsequence to ensure that $V_i\weakly V$ for some $V\in \S_2$, which moreover has $\frac{5}{2}\w_n-1\leq \|V\|(B_1(0))\leq \frac{5}{2}\w_n + 1$. By definition of $\CN_{\epsilon_i}(\BC)$ and the fact that varifold convergence implies local convergence of the supports with respect to the Hausdorff distance, we see that $\spt\|V\|\cap B_1\subset \spt\|\BC\|$. In particular, as $\BC$ comprises of half-hyperplanes and $\spt\|V\|$ has no boundary in $B^{n+1}_2(0)$, we see that necessarily $\spt\|V\| \cap B_1$ consists of some subcollection of half-hyperplanes in $\spt\|\BC\|\cap B_1$. Thus, as $V$ is integral, each half-hyperplane in $V$ arises with a fixed integer multiplicity. From the mass bound on $V$ we know that $\Theta_V(0)<3$; the form of $V$ (i.e. supported on half-hyperplanes) then implies $\Theta_V(0) \leq \frac{5}{2}$. It is then straightforward to check case by case that for $\BC$ of this specific form that we must have, as $\|V\|(B_1)>2\w_n$, that $V = \BC$; this then completes the proof.
\end{proof}

For a given classical cone $\BC$, we introduce a notation for functions defined over $\spt\|\BC\|$ which ``respect the multiplicity'' of $\BC$. These are the functions which we will use to approximate nearby varifolds as graphs over $\spt\|\BC\|$.

\begin{defn}\label{defn:graphs}
	Fix $\alpha\in (0,1)$, $I\in \{0,1,2\}$, and $\BC\in \FL_I$; write $\BC = \sum^{5-2I}_{i=1}\llbracket H_i\rrbracket + 2\sum^I_{j=1}\llbracket G_j\rrbracket$, where $(H_i)_{i=1}^{5-2I}$, $(G_i)_{i=1}^I$ are the half-hyperplanes $\spt\|\BC\|$ is comprised of.  We say $u\in C^{1,\alpha}(\BC)$, written $u:\spt\|\BC\|\to \A_{\BC}(\spt\|\BC\|^\perp)$, if the following holds:
	\begin{enumerate}
		\item [(a)] For each $i=1,\dotsc,5-2I$, $\left. u\right|_{H_i} =u_i$, where $u_i\in C^{1,\alpha}(\overline{H}_i,H_i^\perp)$, and for each $j=1,\dotsc,I$, $\left. u\right|_{G_j} = v_j$, where $v_j\in C^{1,\alpha}(\overline{G}_j,\A_2(G_j^\perp))$;
		\item [(b)] There is a $C^{1,\alpha}$ function $w:S(\BC)\to \R^2$ such that for each $i,j$ we have $\left. u_i\right|_{S(\BC)} = w^{\perp_{H_i}}$ and $\left. v_j\right|_{S(\BC)} = w^{\perp_{G_j}}$, where $\perp_{H}$ denotes the orthogonal projection onto the normal direction to $H$.
	\end{enumerate}
	We also write $u\in C^{1,\alpha}(\BC\res U)$, for $U\subset\R^{n+1}$ obeying $\dist(U,S(\BC))>0$, to mean that $u = \left. v\right|_{\spt\|\BC\|\cap U}$ for some $v\in C^{1,\alpha}(\BC)$, and we write $u\in C^{1,\alpha}(\BC\res\{r_{\BC}>0\})$ to mean a function $u$ obeying (a) above, except we only require $u_i\in C^{1,\alpha}(H_i,H_i^\perp)$ and $v_j\in C^{1,\alpha}(G_i,\A_2(G^\perp_j))$.
\end{defn}

\textbf{Remark:} Condition (b) tells us that each two-valued function $v_j$ necessarily has boundary values on $\del H_j$ which are determined by a single-valued function.

Throughout this work, we will almost always assume that we have rotated the system so that $S(\BC) = \{0\}^2\times\R^{n-1}$ when $\BC\in \FL$.

\subsection{Multiplicity Two Classes}\label{sec:MC}

The aim of this section is to show that, for $\epsilon = \epsilon(\BC)>0$ sufficiently small, $\CN_\epsilon(\BC)$ is contained in a so-called \textit{multiplicity two class}; this is a natural extension of the notion of a multiplicity one class originally introduced in \cite{simoncylindrical}. Working in such a class provides us with more powerful estimates than those available otherwise, as we know that if the support of a varifold is close to a hyperplane, then the multiplicity of the hyperplane must be $2$ or $1$, allowing us to apply Theorem \ref{thm:wick1} or Allard's regularity theorem; as such, we will be in places able to argue in the style of \cite[Sections 2 and 3]{simoncylindrical} as opposed to those in \cite[Section 10]{wickstable} (which is significantly more involved).

\begin{defn}
	We say a class $\M_2$ is a \textit{multiplicity two class} if:
	\begin{enumerate}
		\item [(i)] Elements of $\M_2$ are pairs $(V,U_V)$, where $U_V\subset\R^{n+1}$ is open, $V$ is a stationary integral $n$-varifold in $U_V$, with stable regular part (in the sense of $(\S2)$) in $U_V$ and with no triple junction singularities;
		\item [(ii)] $\M_2$ is closed under rotations and suitable homotheties, i.e. if $(V,U_V)\in \M_2$ then for any orthogonal rotation $q$ of $\R^{n+1}$, $X\in U_V$, and $\rho\in (0,\dist(X,\del U_V))$ we have $((q\circ\eta_{X,\rho})_\# V, (q\circ \eta_{X,\rho})(U_V))\in \M_2$;
		\item [(iii)] If $(V_j,U_{j})\subset \M_2$ and $U\subset\R^{n+1}$ is open such that $U\subset U_j$ for all sufficiently large $j$ and $\sup_{j\geq 1}\|V_j\|(K)<\infty$ for each compact $K\subset U$, then there is a subsequence $(V_{j'})_{j'}$ and $(V,U_V)\in \M_2$ such that $U\subset U_{V}$, $V_{j'}\res U\weakly V\res U$, and moreover $\left.\Theta_V\right|_{\reg(V)\cap U}\leq 2$.
	\end{enumerate}
\end{defn}
\textbf{Remark:} When $U_V$ is contextually clear, we shall write for simplicity $V\in \M_2$ instead of $(V,U_V)\in \M_2$. Moreover, from (iii) we see that $\left.\Theta_V\right|_{\reg(V)} \leq 2$.

It follows from Allard's regularity theorem and Theorem \ref{thm:wick1} that we have the following $\epsilon$-regularity theorem for multiplicity two classes:

\begin{theorem}\label{thm:multiplicity_two_class}
	Let $\Lambda>0$ and let $\M_2$ be a multiplicity two class. Then there exists a constant $\beta = \beta(\M_2,\Lambda)>0$ such that the following is true: if $(V,U_V)\in \M_2$, $\rho>0$, $B_\rho(X_0)\subset U_V$, $\|V\|(B_\rho(X_0))\leq \Lambda$, $\spt\|V\|\cap B_{3\rho/4}(X_0)\neq\emptyset$, and $\rho^{-n-2}\int_{B_\rho(X_0)}\dist^2(X,P)\ \ext\|V\|(X)<\beta$ for some $n$-dimensional hyperplane $P\subset\R^{n+1}$, then either:
	\begin{enumerate}
		\item [(i)] There is a $C^2$ map $u:P\cap B_{3\rho/4}(X_0)\to P^\perp$ such that $V\res B_{\rho/2}(X_0) = |\graph(u)|\res B_{\rho/2}(X_0)$, $\graph(u)\subset \spt\|V\|$;
		\item [(ii)] There is a $C^{1,1/2}$ map $u:P\cap B_{3\rho/4}(X_0)\to \A_2(P^\perp)$ such that $V\res B_{\rho/2}(X_0) = \mathbf{v}(u)\res B_{\rho/2}(X_0)$, $\graph(u)\subset\spt\|V\|$;
	\end{enumerate}
	moreover, in either case we have (for some $C = C(n)$):
	$$\rho^{-2}\sup|u|^2 + \sup|Du|^2 \leq C\rho^{-n-2}\int_{B_\rho(X_0)}\dist^2(X,P)\ \ext\|V\|.$$
\end{theorem}

\begin{proof}
	If this were false, then we could find sequences $(V_k)_k\subset\M_2$, $(\rho_k)_k$, $(X_k)_k$ with $B_{\rho_k}(X_k)\subset U_{V_k}$, $\|V_k\|(B_{\rho_k}(X_k))\leq \Lambda$, $\spt\|V_k\|\cap B_{3\rho_k/4}(X_k)\neq \emptyset$, and
	$$\rho_k^{-n-2}\int_{B_{\rho_k}(X_k)}\dist^2(X,P_k)\ \ext\|V_k\|\to 0$$
	for some sequence of hyperplanes $(P_k)_k$, such that neither conclusion holds for every $k$. Now, for each $k$ we can find a rotation $q_k$ of $\R^{n+1}$ with $q_k(P_k) = \{0\} \times \R^n$, and then by definition of a multiplicity two class we know that $\tilde{V}_k:= (q_k\circ\eta_{X_k,\rho_k})_\# V\in \M_2$ for each $k$. Moreover, by construction we have $\hat{E}_{\tilde{V}_k}\to 0$, where $\hat{E}^2_{\tilde{V}_k} := \int_{B_1}|x^1|^2\ \ext\|\tilde{V}_k\|$. Also, as $\|\tilde{V}\|(K)\leq \Lambda$ for each compact $K\subset B_1$, we can pass to a subsequence to ensure that $\tilde{V}_k\weakly V\in \M_2$. By construction we necessarily have $\spt\|V\|\cap B_1(0)\subset \{0\}\times\R^n$, and thus as $V\in \M_2$, we have $V\res B_1 = \theta|\{0\}\times B^n_1(0)|$, for some constant $\theta\in \{1,2\}$. But if $\theta=1$ we contradict Allard's regularity theorem, and if $\theta = 2$ we contradict Theorem \ref{thm:wick1}.
\end{proof}

An important observation to containing $V\in \CN_\epsilon(\BC)$ in a multiplicity two class if that we can control the density of points in $V$ close to the spine of $\BC$. The following simple lemma enables us to do this.

\begin{lemma}\label{lemma:mass_bounds}
	Fix $\delta_0>0$. Then there exists $\epsilon_0 = \epsilon_0(n,\delta_0)$ such that if $\delta\geq \delta_0$ and $V$ is a stationary $n$-varifold in $B_1(0)$ which satisfies $\w_n^{-1}\|V\|(B_1(0))\leq \frac{5}{2}+\delta$, then for any $X\in B_{\epsilon_0}(0)$ and any $\rho\in (0,1-|X|)$ we have
	$$\frac{\|V\|(B_\rho(X))}{\w_n\rho^n}\leq \frac{5}{2} + 2\delta.$$
\end{lemma}

\begin{proof}
	From the monotonicity formula it follows that:
	$$\frac{\|V\|(B_\rho(X))}{\w_n\rho^n} \leq \frac{\|V\|(B_{1-|X|}(X))}{\w_n(1-|X|)^n} \leq \frac{\|V\|(B_1(0))}{\w_n}\cdot\frac{1}{(1-|X|)^n} \leq \left(\frac{5}{2}+\delta\right)\cdot \frac{1}{(1-\epsilon_0)^n}$$
	and since $y\mapsto \frac{\frac{5}{2}+2y}{\frac{5}{2}+y}$ is increasing for $y>0$, it suffices to take $\epsilon_0$ obeying $(1-\epsilon_0)^{-n}\leq \frac{\frac{5}{2}+2\delta_0}{\frac{5}{2}+\delta_0}$.
\end{proof}

Thus it is crucial that our cone has half-integer density as opposed to integer density. We can now prove:

\begin{theorem}\label{thm:M2C}
	Let $\BC\in \FL_S$. Then there exists $\epsilon_1 = \epsilon_1(\BC)\in (0,1)$ and a multiplicity two class $\M_2 = \M_2(\BC)$ such that for all $\epsilon\leq \epsilon_1$, $\CN_{\epsilon}(\BC)\subset \M_2$, in the sense that there is a fixed $U\supset B^{n+1}_{3/4}(0)$ for which $(V,U)\in \M_2$ for each $V\in \CN_\epsilon(\BC)$.
\end{theorem}

\begin{proof}
	We follow a similar argument to that seen in \cite[Corollary 3]{simoncylindrical}. Firstly, note the trivial inclusion $\CN_\epsilon(\BC)\subset\CN_{\epsilon_1}(\BC)$ for $\epsilon\leq \epsilon_1$, and so it suffices to prove the containment for $\CN_{\epsilon_1}(\BC)$. Next note from Lemma \ref{lemma:convergence} that if $\epsilon_1 = \epsilon_1(\BC)$ is sufficiently small we have for all $V\in \CN_{\epsilon_1}(\BC)$ that $\w_n^{-1}\|V\|(B_1(0))\leq \frac{5}{2}+1/8$; thus applying Lemma \ref{lemma:mass_bounds} with $\delta_0 = 1/8$ we see that we can find some $\rho_0 = \rho_0(n)$ such that for any $X\in B_{\rho_0}(0)$ and $\rho\in (0,1-|X|)$ we have $(\w_n\rho^n)^{-1}\|V\|(B_\rho(X))\leq \frac{5}{2}+1/4$.
	
	By translating parallel to $S(\BC)$ (which we can without loss of generality assume is $\{0\}^2\times\R^{n-1}$) we can also arrange that the same argument holds at any point $Y\in \{0\}^2\times B^{n-1}_{7/8}(0)$; thus we can choose $\epsilon_1 = \epsilon_1(\BC)$ sufficiently small such that for all $X\in B^2_{\rho_0}(0)\times B^{n-1}_{7/8}(0)$ and $\rho\in (0,1-|X|)$ we have $(\w_n\rho^n)^{-1}\|V\|(B_\rho(X)) \leq \frac{5}{2}+1/4$. Moreover, from Lemma \ref{lemma:convergence}, Allard's regularity theorem, and Theorem \ref{thm:wick1}, we can find $\epsilon_* = \epsilon_*(\BC) < \epsilon_1$ such that if $\epsilon\leq\epsilon_*$, then on $\{|x|>\rho_0/4\}\cap B^{n-1}_{7/8}(0)$ we can express $V\in \CN_{\epsilon}(\BC)$ as a sum of single-valued and two-valued $C^{1,1/2}$ (stationary) graphs defined on $\spt\|\BC\|\cap\{|x|>\rho_0/4\}\cap B^{n+1}_{7/8}(0)$.
	
	Now define $\M_2$ to be all pairs $(V,U)$, where the varifolds $V$ are either of the form $(q\circ\eta_{Y,\rho})_\# \tilde{V}$ or a varifold limit of varifolds of this form, where $Y\in B^{n+1}_{3/4}(0)$, $\rho>0$, $q$ a rotation of $\R^{n+1}$, and $\tilde{V}\in \CN_{\epsilon_*}(\BC)$, and $U = B^{2}_{13/16}(0)\times B^{n-1}_{13/16}(0)$. Then the above conditions imply that $\M_2$ is a multiplicity two class which contains $\CN_{\epsilon_1}(\BC)$ in the sense stated in the theorem, which completes the proof.
\end{proof}

\textbf{Remark:} The above argument does not actually depend on $\BC\in \FL_S$: no sequence of cones in $\FL_S$ can converge to a cone which has a half-hyperplane of multiplicity $\geq 3$, as this would contradict the stationarity condition. As such, in the proof of Theorem \ref{thm:M2C} one could find a multiplicity two class $\M_2$, dependent only on the dimension (and the number of half-hyperplanes in $\BC$ counted with multiplicity, or equivalently $\Theta_{\BC}(0) = 5/2$) and an $\epsilon = \epsilon(n)>0$ such that $\CN_{\epsilon}(\BC)\subset\M_2$, in the above sense, for \textit{any} $\BC\in \FL_S$. As such, the constant $\beta$ from Theorem \ref{thm:multiplicity_two_class} when applied later to the classes $\CN_{\epsilon}(\BC)$  can be chosen such that $\beta = \beta(n,\Theta_{\BC}(0)) = \beta(n)$, i.e. independent of $\BC$ and only dependent on the dimension. We shall make use of this later to ensure our constants do not depend on $\BC$ explicitly.

For the rest of this work, we shall always assume that $\epsilon = \epsilon(n)$ is sufficiently small so that $\CN_{\epsilon}(\BC)$ is contained in a multiplicity two class; this is ensured by Lemma \ref{thm:M2C}. By rescaling, we can also assume that if $V\in \CN_{\epsilon}(\BC)$ then $(V,B^{n+1}_1(0)) \in \M_2$.

\subsection{Density Gaps}

The last important property we record here for later is that, for each $V\in \CN_\epsilon(\BC)$ and $Z\in S(\BC)\cap B_1$, points of sufficiently large density in $V$ accumulate at $Z$. Geometrically, this ensures that $V$ is not losing symmetries that $\BC$ has, and analytically we will be able to use this property to show that certain $L^2$ estimates hold on balls centred at $S(\BC)$, which ultimately we used to establish the $C^{1,\alpha}$ boundary regularity for blow-ups.

\begin{defn}
	Let $\delta>0$ and $\BC\in \FL$. We say that $V\in \S_2$ has no $\delta$\textit{-density gaps} with respect to $\BC$ if for each $y\in S(\BC)\cap B_1(0)$ we have $\{\Theta_V\geq \Theta_{\BC}(0) = \frac{5}{2}\}\cap B_\delta(y)\neq\emptyset$.
\end{defn}

\begin{lemma}\label{lemma:gaps}
	Fix $\delta>0$ and $\BC\in \FL$. Then there exists $\epsilon = \epsilon(\BC,\delta) \in (0,1)$ such that each $V\in \CN_{\epsilon}(\BC)$ has has no $\delta$-density gaps with respect to $\BC$.
\end{lemma}

\begin{proof}
	Without loss of generality rotate so that $S(\BC) = \{0\}^2\times\R^{n-1}$. We will in fact prove more: we will show that there is $\epsilon_* = \epsilon_*(\BC)$ sufficiently small such that if $V\in \CN_{\epsilon}(\BC)$ then $\H^{n-1}$-every (in fact every, by upper semi-continuity of the density) two-dimensional slice $\R^2\times\{y\}$ contains a point of density $\geq \frac{5}{2}$. This of course proves the result, as for any $\delta>0$ one may apply Lemma \ref{lemma:convergence} to find $\epsilon = \epsilon(\BC,\delta)\in (0,\epsilon_*)$ such that if $V\in \CN_{\epsilon}(\BC)$ then on $\{|x|>\delta/4\}\cap B^{n+1}_1(0)$ we have that $V$ is a sum of single-valued $C^2$ and two-valued $C^{1,1/2}$ graphs (by Allard's regularity theorem and Theorem \ref{thm:wick1}), and thus has density $<5/2$ in this region; thus any point of density $\geq 5/2$ must lie in $B_\delta(S(\BC))\cap B^{n+1}_1(0)$.
	
	We first claim that, for $\epsilon = \epsilon(\BC)\in (0,1)$ sufficiently small, for every $Y = (0,y)\in \{0\}^2\times B^{n-1}_1(0)$ and $V\in \CN_\epsilon(\BC)$ we have
	\begin{equation}\tag{$\star$}
		(\sing(V)\backslash\CC_2)\cap (\R^2\times\{y\}) \cap B^{n+1}_1(0) \neq\emptyset
	\end{equation}
	where by $\CC_2$ we mean the set of density 2 immersed classical singular points in $V$. Indeed, to see this, first choose $\epsilon = \epsilon(\BC)$ such that if $V\in \CN_{\epsilon}(\BC)$ then on $\{(x,y)\in B^{n+1}_1(0): |x|>100\}$ we can express $V$ as a sum of two-valued $C^{1,1/2}$ stationary graphs and $C^2$ single-valued graphs defined on appropriate subsets of the half-hyperplanes in $\spt\|\BC\|$; such $\epsilon$ exists by virtue of Lemma \ref{lemma:convergence}, Allard's regularity theorem, and Theorem \ref{thm:wick1}. Now, if $(\star)$ fails with this $\epsilon$, then we can find some $(0,y)\in \{0\}^2\times B^{n-1}_1(0)$ and $V\in \CN_{\epsilon}(\BC)$ for which
	$$\sing(V)\cap (\R^2\times\{y\})\cap B^{n+1}_1(0) \subset \CC_2.$$
	However, we know that $\CC_2\subset\sing(V)$ is open, and thus by a simple compactness argument we see that we must be able to find $\rho>0$ such that
	$$\sing(V)\cap (\R^2\times B_\rho^{n-1}(y))\cap B^{n+1}_{1/2}(0)\subset\CC_2.$$
	This means in particular that $M:=\spt\|V\|\cap (\R^2\times\{y\})$ is a smoothly immersed 1-dimensional submanifold in $(\R^2\times\{y\})\cap B^{n+1}_{1/2}(0)\cong B^2_{1/2}(0)$; by choice of $\epsilon$ and the fact that there are no (multiplicity 2) branch points in $V$ in $(\R^2\times B_\rho(y))\cap B^{n+1}_{1/2}(0)$, we see that $M$ has five (counted with multiplicity) connected components in $\del B^{2}_{1/2}(0)$. However, such an immersed 1-dimensional submanifold must have an even number of such boundary components, providing the contradiction and establishing $(\star)$.
	
	From $(\star)$ we can now prove the result. Choosing $\epsilon_* = \epsilon_*(\BC)$ so that $(\star)$ holds with $\epsilon_*$. But note that, $\H^{n-1}(\sing(V)\cap\{\Theta_V<5/2\}\backslash\CC_2) = 0$; indeed, by Almgren's stratification and Theorem \ref{thm:wick2}, and Corollary \ref{cor:multiplicity_2_branch_points}, we have
	$$\H^{n-1}(\sing(V)\cap\{\Theta_V<5/2\}\backslash\CC_2) \leq \H^{n-1}(\B_2) + \H^{n-1}(\S_{n-2}) = 0.$$
	Hence $(\star)$ implies that $\H^{n-1}(\{y\in B^{n-1}_1(0): \Theta_V<5/2\text{ on }\R^2\times\{y\}\}) = 0$, which shows that on $\H^{n-1}$-a.e. slice $\R^2\times\{y\}$ there is a point $X$ with $\Theta_V(X)\geq \frac{5}{2}$. Thus we are done.
\end{proof}

\textbf{Remark:} Just like for Theorem \ref{thm:M2C}, one sees that the above argument does not depend on the base cone $\BC$, but just that no sequence of cones in $\FL_S$ can limit onto a cone with multiplicity $>2$ on a hyperplane. Hence we see that in fact the constant $\epsilon$ in Lemma \ref{lemma:gaps} can be chosen to only depend on $n$ and $\delta$.

Here and elsewhere in the paper, we suppress the dependence on the dimension $n$ of any constant which also depends on a cone $\BC\in \FL$; this is simply because $n$ can be recovered from $\BC$ from the dimension of its spine. Thus we could write $C = C(\BC)$ or $C = C(\BC_0,n)$ if $\BC = \BC_0\times\R^{n-1}$, and we shall opt for the former. Thus when a constant if written to depend on $n$, this is to stress that it \textit{does not} depend on the form of base cone $\BC$.

\subsection{Outline of the Proof of Theorem \ref{thm:A}}

Theorem \ref{thm:A} will be established by proving a suitable excess decay statement. In the simplest case, as will be true for level 0 cones, this just says that there is a scale $\theta = \theta(n)$ such that if $\epsilon = \epsilon(\BC)$ is sufficiently small, then whenever $V\in \CN_{\epsilon}(\BC)$ one can find another cone $\tilde{\BC}$ of the same form as $\BC$ (i.e. the same level) for which the (one-sided) excess of $V$ relative to $\tilde{\BC}$ at scale $\theta$ has decayed by a factor of $\frac{1}{2}$ relative to the (one-sided) excess of $V$ relative to $\BC$ at scale 1, i.e.
$$\theta^{-n-2}\int_{B_\theta}\dist^2(X,\spt\|\tilde{\BC}\|)\ \ext\|V\| \leq \frac{1}{2}E_{V,\BC}^2.$$
By iterating this excess decay statement one will be able to deduce Theorem \ref{thm:A} in a relatively standard manner. When $\BC$ is not level 0 however there is a slight technicality regarding this excess decay statement, namely that the new cone $\tilde{\BC}$ need not be of the same form as $\BC$: it could be of a lower level. When the excess at scale $\theta$ of the new cone $\tilde{\BC}$ is comparable (i.e. up to a dimensional constant) to that of some other cone $\BC^\prime$ which \textit{is} the same level as $\BC$, this is still not a problem: by decreasing the scale $\theta$ we can still get decay with respect to a cone of the same form as the original base cone $\BC$. The difficulty really appears when this excess is not comparable: this is the situation when the varifold $V$ is actually much closer to a cone of lower level as $\BC$. To deal with this case, we shall need a variant of Theorem \ref{thm:A} under the assumption that $V$ is \textit{significantly closer} to the lower level cone than any cone of the same level as $\BC$: this we refer to as a \textit{fine }$\epsilon$\text{-regularity theorem}. This will be established during the respective fine blow-up procedure, which will be needed to establish regularity of the (coarse) blow-up class, and excess decay statement, anyway. Thus, in general our excess decay statement will be a dichotomy: either we get excess decay with respect to a cone of the same level, or the excess decays with respect to a cone of lower level and moreover the varifold is significantly closer to this lower level cone than any cone of the original level -- in which the fine $\epsilon$-regularity theorem will say that the varifold actually already has the structure provided by Theorem \ref{thm:A}. It will be by iterating this excess decay dichotomy that Theorem \ref{thm:A} will be established.

The excess decay statement will be established by a blow-up argument, which is most conveniently phrased as a contradiction argument. So fix a base cone $\BC^{(0)}\in \FL_S$ and consider sequences $V_k\in \CN_{\epsilon_k}(\BC^{(0)})$, $\BC_k\in \FL_{\epsilon_k}(\BC^{(0)})$, where $\epsilon_k\to 0$; rotate everything so that $\BC^{(0)} = \BC^{(0)}_0\times \R^{n-1}$. Let us look at each case individually.

\textbf{Case 1: $\BC^{(0)}\in \FL_0$ is level 0.} This will be the simplest case as it will only involve a coarse blow-up process. We can without loss of generality assume that all the $\BC_k$ are level 0 and essentially follow the ideas in \cite{simoncylindrical}, since we have ruled out density gaps (Lemma \ref{lemma:gaps}) and our varifolds lie within a multiplicity two class (Theorem \ref{thm:M2C}). By fixing a suitable sequence $\tau_k\downarrow 0$, on $\{|x|>\tau_k\}$ we can (essentially) write $V_k$ as a single-valued graph, $u_k$, over $\spt\|\BC_k\|\cap\{|x|>\tau_k\}$; thus $u_k$ is comprised of $5$ single-valued $C^2$ functions, one on some subset of each half-hyperplane in $\spt\|\BC_k\|$, which each solve the minimal surface equation over their respect domains of definition and which have disjoint graphs. We can control the $L^2$ norm of each function over $\{|x|>\tau_k\}$ by the (one-sided) excess, $E_k:= E_{V_k,\BC_k}$. By using the stationarity of each $V_k$ along with the fact that density gaps have been ruled out, we will be able to establish $L^2$ estimates on $V_k$ and $u_k$ analogous to those seen in \cite[Theorem 3.1]{simoncylindrical}. We will be able to do this, avoiding the complications necessary the argument when the base cone $\BC^{(0)}$ has integer density (seen in \cite[Section 10]{wickstable}) because we are able to work in a multiplicity two class (thus $V_k$ can never be arbitrarily close in a ball near the spine $S(\BC^{(0)})$ to a hyperplane of multiplicity $\geq 3$). In particular these $L^2$ estimates enable us to prove that no excess concentrates along $S(\BC^{(0)})$, meaning that when we consider the \textit{blow-up sequence}, $v_k:= u_k/E_k$, we will be able to extract a limit $v$ in $L^2(\BC\res B_1)$ (as opposed to getting just a limit in $L^2_{\text{loc}}(\BC\res(B_1\backslash S(\BC^{(0)})))$, i.e. global convergence in $L^2$ as opposed to locally away from the spine); this limit will be comprised of 5 functions, one over each half-hyperplanes in $\BC^{(0)}$, which are harmonic in the interior (i.e. away from the boundary of the half-hyperplane, that is, the spine of $\BC^{(0)}$). Initially however we have no control on the derivatives of $u$ up to the spine $S(\BC^{(0)})$ (i.e. no control on the derivatives of a given component of $v$ up to the boundary of the half-hyperplane on which it is defined). However, we can use the strong $L^2$ convergence to pass the $L^2$ estimates established for the $V_k$ to the blow-up level; these will give that $v$ is in fact $C^{0,\alpha}$ regular up-to-the-boundary for some $\alpha\in (0,1)$, and moreover that the boundary values of $v$ are $C^{2,\alpha}$ regular. Thus, we can apply classical boundary regularity theory for harmonic functions (e.g. \cite{morrey1966multiple}, \cite{gt}) to deduce that in fact $u$ is $C^{2,\alpha}$ up to the boundary on each half-hyperplane in $\spt\|\BC^{(0)}\|$. Thus the derivatives of $v$ at the boundary define a new cone, $\tilde{\BC}$ (it turns out that the derivatives of $v$ parallel to $S(\BC^{(0)})$ are the same, i.e. independent of the half-hyperplane the component of $v$ is defined on, and so $\tilde{\BC}$ is still $5$ half-hyperplanes meeting along a common axis). Passing this cone $\tilde{\BC}$ back to the varifold level, $V_k$, by rescaling each half-hyperplane by $E_k$, then shows that in fact the excess decay statement does hold for some suitable cone, providing the desired contradicting to establish the excess decay lemma (and moreover it will be in the ``simple'' form where the new cone is of the same form as $\BC^{(0)}$, i.e. level 0, so the simpler argument outlined above will prove Theorem \ref{thm:A} when $\BC^{(0)}$ is level 0).

\textbf{Case 2: $\BC^{(0)}\in \FL_1$ is level 1.} The key difference in this setting to the level 0 case is that we no longer know what level the $\BC_k$ are: they could be level 0 or level 1. Of course, we can pass to a subsequence to assume without loss of generality that either $\BC_k$ is level 0 for all $k$ or level $1$ for all $k$. Let us first focus on the case where all the $\BC_k$ are level 1, i.e. the same level as $\BC^{(0)}$; we will see how to deal with the case when all the $\BC_k$ are level 0 through this case.

When all the $\BC_k$ are level 1, we can follow a similar argument as in the level 0 setting. Indeed, we are still in a multiplicity two class and have ruled out the possibility of density gaps. Thus we may fix a suitable sequence $\tau_k\downarrow 0$ and write $V_k$ as a graph over $\spt\|\BC_k\|\cap \{|x|>\tau_k\}$; the difference is that now we are forced to apply Theorem \ref{thm:wick1} over the multiplicity two half-hyperplane in $\BC_k$, and thus over one half-hyperplane in $\spt\|\BC_k\|$, the function $u_k$ is a two-valued $C^{1,1/2}$ function. This difference does not significantly impact the proofs of the key $L^2$ integral estimates from \cite[Theorem 3.1]{simoncylindrical}, and thus we may still preform the same blow-up procedure, constructing a limit $v = \lim_k v_k \equiv \lim_k E_k^{-1}u_k$ where the convergence is strong in $L^2$ over all of $\spt\|\BC^{(0)}\|\cap B^{n+1}_{3/4}$. The limit $v$ is therefore a smooth harmonic function over multiplicity one half-hyperplanes in $\BC^{(0)}$ and a $C^{1,1/2}$ two-valued harmonic function over the multiplicity 2 half-hyperplane in $\BC^{(0)}$. Just as in the level 0 case, we will still be able to show that, as a single-valued or two-valued function, each component of $v$ is $C^{0,\alpha}$ up-to-the-boundary of each half-hyperplane. The next difference comes from the fact that now we are only able to show that, over a given half-hyperplane in $\spt\|\BC^{(0)}$, the boundary values of the \textit{average} of the corresponding function are $C^{2,\alpha}$; thus as the average is always harmonic (if it is single-valued harmonic then the average is just itself, if it is two-valued harmonic then it is the average of the two values, which we know is harmonic) we hence get that the average is always a $C^{2,\alpha}$ function up-to-the-boundary. This deals with the boundary regularity of $v$ on each half-hyperplane, except the one which is multiplicity two in $\BC^{(0)}$ and thus for which $v$ is represented by a two-valued function; in this case we have the regularity of the average part, but we only know the symmetric part is $C^{0,\alpha}$ up-to-the-boundary. We will however be able to show that the boundary values for any two-valued function actually agree, i.e. they are given by $\{f,f\}$, for some function $f$. Thus, the symmetric part always has zero boundary values; this is crucial since it geometrically means that when a two-valued half-hyperplane splits into two multiplicity one half-hyperplane, they must maintain the same axis, meaning that examples such as those shown in Figure \ref{fig:gaps} do not arise (hence we are using the fact that there are no density gaps in a crucial way here). 

Since there is no known general boundary regularity theory for two-valued $C^{1,1/2}$ harmonic functions, we will need to establish this in the current setting. Our method is to classify the homogeneous degree one blow-ups using methods based on the Hardt--Simon inequality and Campanato regularity theory (similar to those seen in \cite[Section 4]{simoncylindrical} and \cite[Section 4]{wickstable}). The one ingredient we are missing to carry out this classification is a property which plays the role of \cite[$(\B7)$]{wickstable}; it should be noted that \cite[$(\B7)$]{wickstable} as stated does not hold in our setting. We will be able to establish a similar property however, roughly saying the following: whenever a blow-up has a graph which is sufficiently close (in $L^2$) to a union of $5$ (distinct) multiplicity one half-hyperplanes meeting along a common axis, i.e. a level 0 cone, then in fact the blow-up must be $C^{1,\alpha}$ up-to-the-boundary. This is a type of $\epsilon$-regularity property for the blow-up class. To establish it, we shall study the corresponding \textit{fine} blow-up process, and the proof will requiring knowing the validity of Theorem \ref{thm:A} for level 0 cones which we have already discussed. Proving the corresponding boundary regularity statement for the functions in the fine blow-up class will be possible because the two-valued function ``splits'' into two single-valued functions, and thus the functions in the fine blow-up class will be comprised of 5 single-valued harmonic functions, for which we have a boundary regularity theory. It is in this way that we establish the regularity of the fine blow-up class, hence prove the $\epsilon$-regularity property for the original (coarse) blow-up class, and hence prove the boundary regularity for the (coarse) blow-up class, giving rise to an excess decay statement.

However, there is a issue: since the symmetric part of the two-valued function in the (coarse) blow-up need not vanish, it can have non-zero derivative at the origin and hence the new cone for which we get the excess decay need not be level $1$: this is a problem for iteration as the whole analysis above was performed under the assumption that the $\BC_k$ was level $1$. This is when we need to use another result which comes from the fine blow-up procedure: the \textit{fine $\epsilon$-regularity theorem}. This roughly says that there is a fixed dimensional constant $\beta$ such that if the (two-sided) excess relative to a level 0 cone is significantly smaller than the (two-sided) excess relative to \textit{any} level 1 cone, then in fact we already have a regularity conclusion for $V$ similar to that of Theorem \ref{thm:A}; intuitively this corresponds to the case where $V$ consists of $5$ separate multiplicity one sheets, but two of them happen to be very close. So, if the excess decays with respect to a level 0 cone $\BC^\prime$, one may ask: is the excess relative to $\BC^\prime$, at some fixed scale $\theta = \theta(n)$, significantly smaller than the excess relative to any level 1 cone at scale $\theta$? If so, one does not need to iterate the excess decay inequality further, as we have the desired regularity conclusion on some smaller ball (say, $B_{\theta/2}$) from the fine $\epsilon$-regularity theorem. If this does not hold however, one may find a level $1$ cone, $\tilde{\BC}$, for which the excess of $V$ relative to $\BC^\prime$ at scale $\theta$ is bounded below by $\beta$ times the excess of $V$ relative to $\tilde{\BC}$ at scale $\theta$; hence one may use that the excess decays relative to $\BC^\prime$ to see that the excess does in fact decay with respect to a level 1 cone, namely $\tilde{\BC}$ (some constants may change, but as they are only dependent on the dimension this is fine). In this case we therefore have a suitable excess decay with respect to cones of the same level, which we can iterate a further step and repeat the process. This provides a suitable ``excess decay dichotomy'', either when we iterate the excess decay infinitely many times in the usual fashion, or we stop at some scale and apply the fine $\epsilon$-regularity theorem. This is then enough to conclude Theorem \ref{thm:A} for level 1 cones. 

\textbf{Case 3: $\BC^{(0)}\in \FL_2$ is level 2.} Broadly speaking, we follow the same ideas as in setting where $\BC^{(0)}\in \FL_1$ is level 1, but with a significant extra technicality resulting from the fact that now the $\BC_k$ can be any level. As before, let us first focus on the case where all the $\BC_k$ are level 2, i.e. the same level as $\BC^{(0)}$, and see how in dealing with this case we will also develop the necessary tools to deal with the other cases.

We follow the same general blow-up argument as before, except now we will have that over two of the half-hyperplanes in $\spt\|\BC^{(0)}\|$ we have that the function $v$ is given be two-valued $C^{1,1/2}$ harmonic functions in the interior for which we need to establish the boundary regularity theory for. For this we wish to establish the corresponding $\epsilon$-regularity property for the blow-up class again, which this time will take the form: whenever a blow-up has a graph which is sufficiently close (in $L^2$) to a cone of level $<2$, then in fact the blow-up must be $C^{1,\alpha}$ up-to-the-boundary. One can then study the corresponding fine blow-up class; however, when we are taking a fine blow-up relative to a sequence of \textit{level 1} cones, only one of the two-valued functions will ``split'' into two single-valued functions, meaning that over some half-hyperplane the functions in the fine blow-up class are still represented by a two-valued $C^{1,\alpha}$ function; hence we can no longer apply standard elliptic boundary regularity arguments to deduce the boundary regularity of functions in the fine blow-up class. Thus, in order to prove this we take the same approach we have used previously for deducing boundary regularity for the (coarse) blow-up classes: we will use an argument based on the (reverse) Hardt--Simon inequality to classify the homogeneous degree one elements of the fine blow-up class. However, for this to work we will also need an $\epsilon$-regularity property for the fine blow-up class: to prove this we will need to perform an even finer blow-up process, which we call an \textit{ultra fine blow-up}. This is carried out in a similar way to the fine blow-up, and functions in the ultra fine blow-up class will consist of 5 single-valued harmonic functions; hence we can apply standard elliptic boundary regularity theory to deduce the boundary regularity of functions in the ultra fine blow-up class, which in turn allows us to deduce an \textit{ultra fine} $\epsilon$-regularity theorem for our varifolds (this uses Theorem \ref{thm:A} for level 0 cones), which in turn allows us to prove the $\epsilon$-regularity property for the fine blow-up class, which in turn allows us to deduce the boundary regularity of the fine blow-up class, which in turn allows us to prove a fine $\epsilon$-regularity theorem for varifolds converging to a level 2 cone which are significantly closer to a sequence of level 1 cones (this uses Theorem \ref{thm:A} for level 1 cones), which allows us to prove the $\epsilon$-regularity property for the (coarse) blow-up class, \textit{in the case where the graph of the blow-up is close to a level 1 cone}. When the graph of a function in the coarse blow-up class is close to a level 0 cone, passing to the varifold level, one needs to ask, similarly to what we saw in setting where $\BC^{(0)}\in \FL_1$ was level 1, whether there is a level 1 cone which has excess comparable to the level 0 cone; if so, one may reduce to the setting where the fine blow-up is taken relative to a sequence of level 1 cones as above. If this is not the case, then we are actually in a setting where we can take a fine blow-up relative to a sequence of level 0 cones; in which case elements of the fine blow-up class are in fact made from $5$ single-valued harmonic functions, for which the boundary regularity is simple. Thus Theorem \ref{thm:A} for level 0 cones can be used to prove the $\epsilon$-regularity property for the coarse blow-up class. Combining all of the above then proves the boundary regularity of the coarse blow-up class.

Given all the above analysis, one can then prove a suitable excess decay dichotomy, similar to the case where the base cone $\BC^{(0)}$ was level 1. If all the $\BC_k$ are level 2, the above is enough to deduce that one can find some other cone $\BC^\prime$ for which the excess decays; however we do not know if this cone is level 2, and so we do not know how to iterate this. If however eventually always the $\BC_k$ are level $<2$, one asks: is the excess relative to $\BC_k$ significantly smaller (again, by a fixed dimensional constant) than \textit{every} level $2$ cone? If not, then one can replace the sequence $\BC_k$ by a suitable sequence of level $2$ cones and deduce that the excess decays again. If however this is true, we are in the realm of the fine blow-up process; hence one may apply the suitable fine $\epsilon$-regularity theorem established above (the proof of which requires the full ultra fine blow-up process) to deduce the desired regularity already holds. One may then iterate such a statement, i.e. we either stop at some finite scale and have the desired statement of Theorem \ref{thm:A} already, or we get that the excess decays along a geometric sequence of scales relative to a sequence of level 2 cones; then the result can be concluded in the usual fashion. This will prove Theorem \ref{thm:A} for level 2 cones, and hence combining with all the above will prove Theorem \ref{thm:A} in full.

There is one small additional technicality regarding proving the boundary regularity of the fine blow-up class when $\BC^{(0)}\in \FL_2$ is level 2 which we mention now. The construction of the fine blow-up class depends on a choice of parameter $M>1$, and it turns out that the fine blow-up class for a fixed parameter $M$ is not closed under simple operations, such as domain rescalings. This means that our general arguments for boundary regularity don't quite hold in this setting; however, it turns out that the closure of the fine blow-up class under these operations is contained within another fine blow-up class for a \textit{fixed} parameter $M^\prime = M^\prime(M,n)$; this turns out to be sufficient for our purposes.

\section{Proper Blow-Up Classes}\label{sec:coarse_regularity}

In this section we set up the general language we shall use for our blow-up classes, and state their properties; this is so that we can simply refer back to the properties later on when we prove them in various settings, and to give the reader a point of reference. We will also state the main regularity results, which are established from the presented properties and a suitable multi-valued Campanato theory, all of which is carried out in the accompanying work \cite{minter2021}.

Fix a base cone $\BC^{(0)}\in \FL_S\cap \FL_I$, where $I\in \{0,1,2\}$; let us rotate so that $\BC^{(0)} = \BC^{(0)}_0 \times\R^{n-1}$. Functions in the various blow-up classes will be defined over the half-hyperplanes in $\spt\|\BC^{(0)}\|$, and thus it is convenient to rotate each half-hyperplane and view all the components of such a function as a function defined on a single half-hyperplane. Indeed, for each ray $\ell_i$ in $\spt\|\BC^{(0)}_0\|$, find a rotation $q_i$ of $\R^2$ which maps $\ell_i$ to $\{(x^1,x^2)\in \R^2: x^1=0, x^2>0\}$; then define a rotation $Q_i$ of $\R^{n+1}$ by $Q_i(x,y) = (q_i(x),y)$; then $Q_i$ is a rotation which fixes $S(\BC^{(0)})$ and rotates the half-hyperplane $H_i = \ell_i\times \R^{n-1}$ in $\spt\|\BC^{(0)}\|$ to $H:=\{(x^1,\dotsc,x^{n+1})\in \R^{n+1}: x^1=0, x^2>0\}$. Moreover, any function $v_i:H_i\to H_i^\perp$ can be rotated to a function $\tilde{v}_i:H\to H^\perp$ by $\tilde{v}_i(x):= Q_i v(Q_i^{-1}x)$.

The following will be our general definition of a (coarse) blow-up class:

\begin{defn}\label{defn:proper-blow-up-class}
	We say a collection of functions $\FB(\BC)$ is a \textit{proper} (\text{coarse}) \textit{blow-up class} over $\BC\in \FL_S\cap \FL_I$ if it obeys the following properties:
	\begin{enumerate}
	\item [$(\FB1)$] Each element $v\in \FB(\BC)$ takes the form $v = (v^1,\dotsc,v^{5-I})$, where $v^i\in L^2(B_1(0)\cap H;\A_{q_i}(H^\perp))\cap W^{1,2}_{\text{loc}}(B_1(0)\cap H; \A_{q_i}(H^\perp))$, where $q_1,\dotsc,q_{5-2I} = 1$ and $q_{5-2I+1},\dotsc,q_{5-I} = 2$;
	\item [$(\FB2)$] (Interior regularity). If $v\in \FB(\BC)$, then $v^i$ is a $q_i$-valued harmonic function for each $i=1,\dotsc,5-I$, which is smooth if $q_i=1$ and $C^{1,1/2}$ if $q_i=2$;
	\item [$(\FB3)$] (Boundary estimates). If $v\in \FB(\BC)$ and $z\in B_1(0)\cap \del H$, then for each $\rho\in (0,\frac{3}{8}(1-|z|)]$ we have
	$$\int_{B_{\rho/2}(z)\cap H}\sum^{5-I}_{i=1}\frac{|v^i(x)-\kappa^i(z)|^2}{|x-z|^{n+3/2}}\ \ext x \leq C\rho^{-n-3/2}\int_{B_\rho(z)\cap H}\sum^{5-I}_{i=1}|v^i(x)-\kappa(z)|^2\ \ext x$$
	where $\kappa:B_1(0)\cap \del H\to \R^2$ is a smooth single-valued function which obeys
	$$\sup_{B_{5/16}(0)\cap \del H}|\kappa|^2 \leq C\int_{B_{1/2}(0)\cap H}|v|^2$$
	and $\kappa^i$ denotes the projection of $\kappa$ onto the normal direction to $H_i$;
	\item [$(\FB4)$] (Hardt--Simon inequality) For $v\in \FB(\BC)$, $z\in B_1(0)\cap \del H$, and $\rho\in (0,\frac{3}{8}(1-|z|)]$, we have:
	$$\int_{B_{\rho/2}(z)\cap H}\sum^{5-I}_{i=1}R_z^{2-n}\left(\frac{\del}{\del R_z}\left(\frac{v^i-v^i_a(z)}{R_z}\right)\right)^2 \leq C\rho^{-n-2}\int_{B_\rho(z)\cap H}\sum^{5-I}_{i=1}|v^i-\ell_{v^i,z}|^2$$
	where $R_z(x):= |x-z|$ and $\ell_{v^i,z}(x):= v^i_a(x) + (x-z)\cdot Dv^i_a(z)$ is the first-order linear approximation to the average part of $v^i$ at $z$;
	\item [$(\FB5)$] (Closure properties). If $v\in \FB(\BC)$, then:
	\begin{enumerate}
		\item [$(\FB5\text{I})$] For each $z\in B_1(0)\cap \del H$ and $\sigma\in (0,\frac{3}{8}(1-|z|)]$, if $v\not\equiv 0$ in $B_\sigma(z)\cap H$ then $v_{z,\sigma}(\cdot):= \|v(z+\sigma(\cdot))\|^{-1}_{L^2(B_1(0)\cap H)}v(z+\sigma(\cdot))\in \FB(\BC)$;
		\item [$(\FB5\text{II})$] $\|v-\ell_v\|^{-1}_{L^2(B_1(0)\cap H}(v-\ell_v)\in \FB(\BC)$ whenever $v-\ell_v\not\equiv 0$ in $B_1(0)\cap H$, where $v-\ell_v = (v^1-\ell_{v^1},\dotsc,v^{5-I}-\ell_{v^{5-I}})$ and $\ell_{v^i} = \ell_{v^i,0}$ (from $(\FB4)$);
	\end{enumerate}
	\item [$(\FB6)$] (Compactness property). If $(v_m)_m\subset\FB(\BC)$, then there is a subsequence $(m')\subset (m)$ and a function $v\in \FB(\BC)$ such that $v_{m'}\to v$ strongly in $L^2_{\text{loc}}(B_1(0)\cap \overline{H})$ and weakly in $W^{1,2}_{\text{loc}}(B_1(0)\cap H$);
	\item [$(\FB7)$] ($\epsilon$-regularity property). There exist constants $\alpha = \alpha(n)$ and $\epsilon = \epsilon(\BC)$ such that whenever $v\in \FB(\BC)$ has $v^i_a(0) = 0$, $Dv^i_a(0) = 0$ for each $i=1,\dotsc,5-I$, and $\|v\|_{L^2(B_1(0)\cap H)} = 1$, then the following is true: if $v_* = (v_*^1,\dotsc,v_*^{5-I})$ is such that for each $i$, $v^i_*:H\to \A_{q_i}(H^\perp)$, with $q_1,\dotsc,q_{5-2I} = 1$, $q_{5-2I+1},\dotsc,q_{5-I} = 2$, $\graph(v^i_*)$ is a union of $q_i$ half-hyperplanes with boundaries meeting along $\del H$, and has $(v^i_*)_a \equiv 0$ for each $i$ but $v^i_*\not\equiv 0$ for at least one $i>5-2I$, then if 
	$$\int_{B_1(0)\cap H}\G(v,v_*)^2<\epsilon$$
	then we have $\left. v\right|_{B_{1/2}(0)\cap H}\in C^{1,\alpha}(B_{1/2}(0)\cap \overline{H})$.
	\end{enumerate}
	Here, $C = C(n)$ is simply a dimension constant to be chosen.
\end{defn}

Note that we are only ever subtracting a single-valued function from a (possibly) two-valued function in all the above. Here $\G(v,v_*)^2 = \sum^{5-I}_{i=1}\G(v^i,v_*^i)^2$.

\textbf{Remark:} If $v = (v^1,\dotsc,v^{5-I})\in \FB(\BC)$, it follows from $(\FB2)$ that $v^i_a$ is a smooth harmonic function for each $i$, and moreover from $(\FB3)$ it follows that the boundary values of $v_a^i$ equal $\kappa^i$, which is a smooth function, As such, by standard elliptic boundary regularity theory we always have $v^i_a\in C^\infty(B_1(0)\cap\overline{H};H^\perp)$, and hence it makes sense to talk above $\left. v^i_a\right|_{B_1(0)\cap \del H}$ and $\left. Dv^i_a\right|_{B_1(0)}$, as in $(\FB4)$ and $(\FB5\text{II})$. The smoothness of $\kappa$ will in fact come from the integral estimates in $(\FB3)$ coupled with further identities arising from the stationarity condition of the varifolds. In fact, $(\FB3)$ allows us to deduce that each $v\in \FB(\BC)$ is actually $C^{0,\beta}(B_1\cap\overline{H})$, for some $\beta = \beta(n)$ (see \cite{minter2021}). The fact that $\beta$ does not depend on $\BC$, and only the dimension, is due to the fact that the constant $C$ only depends on the dimension $n$ (it is for this reason as well that the final regularity constant $\alpha$ in Theorem \ref{thm:A} is independent of $\BC$).

\textbf{Note:} In the case $I=0$, $(\FB7)$ is automatically satisfied.

The main boundary regularity result for $\FB(\BC)$ is the following:

\begin{theorem}[\cite{minter2021}, Theorem 3.1]\label{thm:coarse_reg}
	Let $\BC\in \FL_S\cap \FL$, where $I\in \{0,1,2\}$. Then there exists $\gamma = \gamma(n)\in (0,1/2)$ such that if $v\in \FB(\BC)$, then $v\in C^{1,\gamma}(B_{1/8}(0)\cap \overline{H})$. Moreover, $\left. v_s\right|_{\del H} \equiv 0$, the branch set of $v$ (including any boundary branch points) is countably $(n-2)$-rectifiable, and we have the estimate:
	$$\rho^{-n-2}\int_{B_\rho(z)\cap H}\G(v,\ell_z)^2 \leq C\rho^{2\gamma}\int_{B_{1/2}(0)\cap H}|v|^2$$
	for every $\rho\in (0,1/8]$ and $z\in \overline{H}\cap B_{1/8}$; here, $\ell_z = (\ell_z^1,\dotsc,\ell_z^{5-I})$ is $\ell^i_z(x):= v^i_a(z) + (x-z)\cdot Dv_a^i(z) + (x-z)\cdot \llbracket \pm Dv^i_s(z)\rrbracket$, and $C = C(n)\in (0,\infty)$. Furthermore, $(\FB3)$ implies that $\left.v^i\right|_{B_1\cap \del H} = \kappa^i$ for each $i$.
\end{theorem}

The refer the reader to \cite{minter2021} for a proof of this result. We remark that whilst the boundary regularity conclusion is stated on $\overline{H}\cap B_{1/8}$, this can of course be improved to $\overline{H}\cap B_1$ using $(\FB5\text{I})$.

Proper coarse blow-up classes will be constructed when we blow-up a sequence of varifolds $V_k$ converging to $\BC^{(0)}$ relative to a sequence of cones $\BC_k$ which are the \textit{same level} as $\BC^{(0)}$. However, we will also have to perform other blow-up procedures, known as \textit{fine} and \textit{ultra fine} blow-ups, when the sequence of cones $\BC_k$ are of \textit{strictly lower} level than $\BC^{(0)}$. This will only be necessary under certain closeness assumptions, namely when the $V_k$ are significantly closer to a lower level cone than \textit{any} cone of the same level as $\BC^{(0)}$. This requires the introduction of a parameter, $M>1$, in the construction of fine (and ultra fine) blow-up classes. However, this parameter $M$ is not well-behaved under the closure properties, $(\FB5)$, and there is no guarantee that the functions detailed in $(\FB5)$ will lie in a blow-up class constructed with the same parameter (we will be able to establish all other properties $(\FB1) - (\FB7)$). However, we will see that the constant $M$ can only increase by fixed dimensional constant; this observation is enough to prove the corresponding boundary regularity theorem as in Theorem \ref{thm:coarse_reg} for fine blow-up classes. This is result is also established in \cite{minter2021}. Thus, we remark:

\begin{theorem}\label{thm:fine-reg}
	The conclusions of Theorem \ref{thm:coarse_reg} also hold for any fine blow-up class, $\FB^F_{p,q;M}(\BC)$, where $M>1$, with $\gamma = \gamma(n,M)\in (0,1/2)$ and $C = C(n,M)\in (0,\infty)$.
\end{theorem}

See Section \ref{sec:fine_construction} for an explanation of the notation $\FB_{p,q;M}^F(\BC)$ used here.

\section{The Coarse Blow-Up Class}\label{sec:coarse_construction}

In this section we shall construct the coarse blow-up class for our setting and show that it obeys properties $(\FB1) - (\FB6)$ from Section \ref{sec:coarse_regularity}. As usual, we fix throughout a level $I\in \{0,1,2\}$ and a base cone $\BC^{(0)}\in \FL_S\cap \FL_I$, and rotate to assume without loss of generality that $\BC^{(0)} = \BC^{(0)}_0\times\R^{n-1}$. Recall that we write $X = (x,y)\in \R^2\times\R^{n-1}$ for coordinates, and $r = |x|$, $R= |X|$. As noted before, the \textit{coarse} blow-up class will be constructed by suitable scaling limits of graphs approximating a sequence of varifolds $V_k\in \CN_{\epsilon_k}(\BC^{(0)})$, where $\epsilon_k\downarrow 0$, relative to a sequence of cones $\BC_k\in \FL_{\epsilon_k}(\BC^{(0)}) \cap \FL_I$, i.e. the cones $\BC_k$ have the \textit{same} level as $\BC^{(0)}$.\footnote{Note that we cannot perform small rotations of $V_k$ and $\BC_k$ to assume that $\BC_k\equiv \BC^{(0)}$ as even though $\BC_k$ and $\BC^{(0)}$ are the same level, the angles between their half-hyperplanes need not agree. This is different to the situation where we are blowing up relative to a hyperplane and the $\BC_k$ are all hyperplanes, as is the case during the (coarse) blow-up procedures in \cite{wickstable} and \cite{minterwick}. Note however we can always perform a small rotation to ensure that $S(\BC_k) = S(\BC^{(0)})$ for all $k$; indeed, this small rotation is already taken into consideration in the definition of $\FL_\epsilon(\BC^{(0)})$.}

\subsection{Approximate Graphical Representation and Initial Estimates}

First we need to construct functions defined on $\spt\|\BC_k\|$ which represent the varifold sequence $V_k$ on a large set. This will be possible away from a fixed $\tau$-neighbourhood of the spine $S(\BC_k)\equiv S(\BC^{(0)})$ for $k$ sufficiently large using Allard's regularity theorem and Theorem \ref{thm:wick1}. However close to the spine, it is less clear whether an approximate representation is possible, even if the excess on a small ball relative to some hyperplane is small, since it is a priori possible that different sheets of the $V_k$ (from the different half-hyperplanes) come close, and so the multiplicity of the close hyperplane could be $>2$. It is possible to deal with this problem (and indeed in a different situation this problem is overcome in \cite[Section 10]{wickstable}), however we do not need to worry about this: $\CN_{\epsilon_k}(\BC^{(0)})$ is contained in a multiplicity 2 class for all $k$ sufficiently large (Theorem \ref{thm:M2C}) and so in this situation the multiplicity of the nearby plane will always be at most $2$. This is one significantly simplification which is possible in the case where the base cone has half-integer density as opposed to full-integer density. As such, one is able to follow arguments similar to that seen in \cite{simoncylindrical} to prove the following:

\begin{lemma}\label{lemma:coarse_graphical_rep}
	Let $\BC^{(0)}\in \FL_S\cap \FL_I$ be as above, and fix $\tau\in (0,1/40)$. Then, there exists $\epsilon_0 = \epsilon_0(\BC^{(0)},\tau)$ such that if $\BC\in \FL_{\epsilon_0}(\BC^{(0)})\cap \FL_I$ and $V\in \CN_{\epsilon_0}(\BC^{(0)})$, then there is an open subset $U\subset\spt\|\BC\|\cap B_1$ with the following properties:
	\begin{enumerate}
		\item [(i)] $U_\tau:= \{(x,y)\in \spt\|\BC\|\cap B_{3/4}: |x|>\tau\}\subset U$;
		\item [(ii)] There exists a function $u$ with domain $U$ such that $\left. u\right|_{U_\tau}\in C^{1,1/2}(\BC\res U_\tau)$ and moreover, for each point $x\in U$ there is a $\rho>0$ such that $\left. u\right|_{U\cap B_\rho(x)}$ is given by either a $C^2$ single-valued function or a $C^{1,1/2}$ two-valued function, valued in the normal direction of the half-hyperplane of $\spt\|\BC\|$ which contains $U\cap B_\rho(x)$.
	\end{enumerate}
	Moreover, the function $u$ obeys:
	\begin{enumerate}
		\item [(a)] $V\res (B_{3/4}\cap \{|x|>\tau\}) = \mathbf{v}(u)\res (B_{3/4}\cap \{|x|>\tau\})$;
		\item [(b)] $\sup_U r^{-1}|u| + \sup_U |Du|\leq \beta$, where $\beta = \beta(n)$ is the constant from Theorem \ref{thm:multiplicity_two_class};
		\item [(c)]
		$$\int_{B_{3/4}\backslash\graph(u)}r^2\ \ext\|V\| + \int_{U\cap B_{3/4}}r^2|Du|^2 \leq CE_{V,\BC}^2$$
		where $C = C(n)$ is independent of $\BC^{(0)}$, $\BC$, $V$, and $\tau$.
	\end{enumerate}
\end{lemma}

\textbf{Remark:} We lose some of the additional structure when compared with the corresponding lemma in \cite{simoncylindrical}, and have the more complicated condition (ii), as we do not have control near the axis of whether the graph will be single-valued or two-valued.

\begin{proof}
	Firstly, choose $\epsilon^\prime = \epsilon^\prime(n)$ so that Theorem \ref{thm:M2C} holds, and let $\epsilon\in (0,\epsilon^\prime)$; hence we may assume $\CN_{\epsilon}(\BC)\subset\M_2$ for some multiplicity two class $\M_2$, in the sense of Theorem \ref{thm:M2C} (and $\M_2$ only depends on $n$). In particular, by Theorem \ref{thm:multiplicity_two_class} we deduce the existence of a constant $\beta = \beta(n)$ such that whenever $V\in \CN_{\epsilon}(\BC^{(0)})$ has $\rho^{-n-2}\int_{B_\rho(x)}\dist^2(X,P)\ \ext\|V\|(X)<\beta^2$ for any ball $B_\rho(x)\subset B_1$ and hyperplane $P$, there is either a single-valued $C^2$ function or two-valued $C^{1,1/2}$ function $u$ with domain $P\cap B_{\rho/2}(x)$ which represents $V$ and moreover obeys
	$$\rho^{-1}\sup|u| + \sup|Du|\leq \beta.$$
	So now fix $V\in \CN_{\epsilon}(\BC^{(0)})$. For each $\rho\in (0,1]$ and $\zeta\in \R^{n-1}$, define a toroidal-region $T_{\rho}(\zeta)$ centered at $(0,\zeta)\in\R^{n+1}$ by:
	$$T_{\rho}(\zeta):= \{(x,y)\in \R^{n+1}:(|x|-\rho)^2 + |y-\zeta|^2 < (\rho/8)^2\}.$$
	Now let $U$ denote the union of all $T_{|\xi|}(\zeta)\cap H$ obeying the following conditions: $(\xi,\zeta)\in \spt\|\BC\|\cap B_{1/2}$, $H$ is a half-hyperplane in $\spt\|\BC^{(0)}\|$, and there is a function $u_{|\xi|,\zeta;H}$ defined on $B_{|\xi|/16}(T_{|\xi|}(\zeta))\cap H$ which is either $C^2$ single-valued or $C^{1,1/2}$ two-valued which obeys
	$$V\res (T_{|\xi|}(\zeta)\cap \tilde{B}_{|\xi|,\zeta;H}) = \mathbf{v}(u_{|\xi|,\zeta;H})\res (T_{|\xi|}(\zeta)\cap \tilde{B}_{|\xi|,\zeta;H})$$
	and
	$$|\xi|^{-1}\sup|u_{|\xi|,\zeta;H}| + \sup|Du_{|\xi|,\zeta;H}|\leq \beta/2$$
	where by $\tilde{B}_{|\xi|,\zeta;H}$ we mean the open ball centred on $H$ whose intersection with $H$ is precisely equal to $T|_{|\xi|}(\zeta)\cap H$\footnote{This is similar to \cite[Section 16]{wickstable}.} We can then define a function $u$ on all of $U$ by:
	$$\left. u\right|_{T_{|\xi|}(\zeta)\cap H}:= \left.u_{|\xi|,\zeta;H}\right|_{T_{|\xi|}(\zeta)\cap H}.$$
	By unique continuation of single-valued $C^2$ and two-valued $C^{1,1/2}$ stationary graphs (Lemma \ref{lemma:unique_continuation}) it follows that for $\epsilon = \epsilon(\BC^{(0)},\tau)$ sufficiently small, properties (i), (ii), (a), and (b) from the lemma statement hold. So all that remains to be checked is (c).
	
	Note that if $(\xi,\zeta)\in \spt\|\BC\|\cap B_{3/4}\cap \del U$, with $|\xi|>0$, then we must necessarily have
	\begin{equation}\label{E:coarse1}
		\int_{B_{3|\xi|/16}(T_{|\xi|}(\zeta))}\dist^2(X,\spt\|\BC\|)\ \ext\|V\|\geq (3|\xi|/16)^{n+2}\beta^2
	\end{equation}
	since otherwise we could apply Theorem \ref{thm:multiplicity_two_class} to extend the definition of $u$ to a neighbourhood of $(\xi,\zeta)$, which then contradicts the definition of $U$. Moreover, since $|\xi|<\tau < 1/40$ and $|\zeta|<3/4$, we have
	\begin{equation}\label{E:coarse-extra}
	\begin{split}
	\int_{U\cap B_{10|\xi|}(0,\zeta)}r^2\ \ext\H^n & \leq (10|\xi|)^2 \cdot \H^n(U\cap B_{10|\xi|}(0,\zeta))\\
	& \leq (10|\xi|)^{n+2}\cdot \frac{\|\BC\|(B_{1/4}(0,\zeta))}{(1/4)^n} \leq (10|\xi|)^{n+2}\cdot 4^n \cdot \frac{5}{2}\w_n
	\end{split}
	\end{equation}
	i.e.
	\begin{equation}\label{E:coarse2}
		\int_{U\cap B_{10|\xi|}(0,\zeta)}r^2\ \ext\H^n \leq C|\xi|^{n+2}
	\end{equation}
	where $C = C(n)$. So using the fact that on $U\cap B_{10|\xi|}(0,\zeta)$ we have $|Du|\leq \beta$, combining (\ref{E:coarse1}) and (\ref{E:coarse2}) we arrive at
	\begin{equation}\label{E:coarse3}
	\int_{U\cap B_{10|\xi|}(0,\zeta)}r^2|Du|^2\ \ext\H^n \leq C\beta^2|\zeta|^{n+2} \leq \tilde{C}\int_{B_{3|\xi|/16}(T_{|\xi|}(\zeta))}\dist^2(X,\spt\|\BC\|)\ \ext\|V\|
	\end{equation}
	where $\tilde{C} = \tilde{C}(n)$; we therefore know that this holds whenever $(\xi,\zeta)\in \spt\|\BC\|\cap B_{3/4}\cap \del U$ and $|\xi|>0$. Then, as we have the trivial cover
	$$\{X = (x,y)\in U\cap B_{3/4}:\dist(X,B_{3/4}\cap\del U) < |x|/2\}\subset \bigcup_{(x,y)\in \spt\|\BC\|\cap B_{3/4}\cap \del U}B_{2|x|}(0,y)$$
	and since $B_{2|\xi_1|}(0,\zeta_1)\cap B_{2|\xi_2|}(0,\zeta_2) = \emptyset$ implies that $B_{|\xi_1|/4}(T_{|\xi_1|}(\zeta_1))\cap B_{|\xi_2|/4}(T_{|\xi_2|}(\zeta_2)) = \emptyset$, by the Vitali covering lemma we may extract a countably collection of balls $(B_{2|x_j|}(0,y_j))_j$ where $(x_j,y_j)\in \spt\|\BC\|\cap B_{3/4}\cap \del U$ such that these balls are pairwise disjoint and, if $A:= \{X = (x,y)\in U\cap B_{3/4}:\dist(X,B_{3/4}\cap \del U)<|x|/2\}$, that
	$$A \subset \cup_j B_{10|x_j|}(0,y_j)$$
	which then implies by (\ref{E:coarse3}) that
	\begin{align*}
		\int_{A}r^2|Du|^2 & \leq \tilde{C}\sum_j\int_{B_{3|x_j|/16}(T_{|x_j|}(y_j))}\dist^2(X,\spt\|\BC\|)\ \ext\|V\|\\
		& \leq \tilde{C}\int_{B_1}\dist^2(X,\spt\|\BC\|)\ \ext\|V\|.
	\end{align*}
	But also, if we set $B:= \{X = (x,y)\in U\cap B_{3/4}: \dist(X,B_{3/4}\cap \del U)\geq |x|/2\}$, then for any $X\in B$ we can apply either (i) standard $L^2$ estimates for single-valued solutions to the minimal surface equation, or (ii) the $L^2$ estimates for $C^{1,1/2}$ two-valued stationary graphs in Section \ref{sec:two-valued_stationary_graphs} (namely (\ref{eqn:two-valued_estimate_1}) and (\ref{eqn:two-valued_estimate_2})) to deduce that
	\begin{equation}\label{E:coarse4}
	\int_{\spt\|\BC\|\cap B_{\rho/2}(X)} r^2|Du|^2 \leq C\int_{B_\rho(X)}|u|^2\ \ \ \ \text{for all }\rho\in (0,|x|/8),
	\end{equation}
	where $C^\prime = C^\prime(n)$. Thus, by considering the cover
	$$B\subset \bigcup_{(x,y)\in B}B_{|x|/32}(x,y)$$
	and applying (a suitable simple adaptation of) the Besicovitch covering lemma, one may find finitely many subcollections $\{\Gamma_1,\dotsc,\Gamma_N\}$ of points in $B$, where $N = N(n)$, such that $B\subset \cup_{i=1}^N\cup_{(x,y)\in \Gamma_i}B_{|x|/32}(x,y)$ and if $(x_1,y_1), (x_2,y_2)\in \Gamma_i$ then $B_{|x_1|/16}(x_1,y_1)\cap B_{|x_2|/16}(x_2,y_2) = \emptyset$. Combining this with (\ref{E:coarse4}) when $\rho = |x|/16$, we arrives at
	$$\int_{B}r^2|Du|^2 \leq NC^\prime \int_{U\cap B_1}|u|^2.$$
	Using properties of the Jacobian (see \ref{eqn:area_formula}) for the two-valued case) we know
	$$\int_{U\cap B_1}|u|^2 \leq C_1\int_{B_1}\dist^2(X,\spt\|\BC\|)\ \ext\|V\|$$
	for some $C_1 = C_1(n)$, and hence combining all the above we see that (as $A\cup B = U\cap B_{3/4}$)
	$$\int_{U\cap B_{3/4}}r^2|Du|^2 \leq CE_{V,\BC}^2$$
	for some $C = C(n)$; this completes one half of (c). For the other half, note that if $(\xi,\zeta)\in \spt\|\BC\|\cap B_{3/4}\cap\del U$, then in a similar way to (\ref{E:coarse-extra}), except now using the monotonicity formula for $V$ and the fact that $\|V\|(B_1)\leq \left(\frac{5}{2}\w_n+1\right)$, we have
	\begin{equation*}
		\begin{split}
			\int_{B_{10|\xi|}(0,\zeta)}r^2\ \ext\|V\| & \leq (10|\xi|^2)\|V\|(B_{10|\xi|}(0,\zeta))\\
			& \leq (10|\xi|^{n+2})\cdot\frac{\|V\|(B_{1/4}(0,\zeta))}{(1/4)^n}\leq (10|\zeta|)^{n+2}\cdot 4^n\cdot\left(\frac{5}{2}\w_n + 1\right)
		\end{split}
	\end{equation*}
	i.e.
	$$\int_{B_{10|\xi|}(0,\zeta)}r^2\ \ext\|V\| \leq C_2|\xi|^{n+2}$$
	where $C_2 = C_2(n)$. Thus combining this with (\ref{E:coarse1}) we get
	$$\int_{B_{10|\xi|}(0,\zeta)}r^2\ \ext\|V\| \leq C\beta^{-2}\int_{B_{3|\xi|/16}(T_{|\xi|}(\zeta))}\dist^2(X,\spt\|\BC\|)\ \ext\|V\|.$$
	By construction, $\spt\|V\|\cap B_{3/4}\backslash \graph(u)\subset \spt\|V\|\cap(\bigcup B_{2|\xi|}(0,\zeta))$, where the union is taken over $(\xi,\zeta)\in \spt\|\BC\|\cap B_{3/4}\cap \del U$. So applying the Vitali covering lemma again, in the same way as above we deduce that
	$$\int_{B_{3/4}\backslash\graph(u)}r^2\ \ext\|V\| \leq CE_{V,\BC}^2$$
	which completes the proof.
\end{proof}

Next we establish key $L^2$ estimates for the graphical representation $u$ from Lemma \ref{lemma:coarse_graphical_rep}.

\begin{lemma}[Coarse $L^2$ Estimates]\label{lemma:L2_coarse}
	Let $\BC^{(0)}\in \FL_S\cap \FL_I$ be as above, and fix $\tau\in (0,1/40)$. Then there exists $\epsilon_0 = \epsilon_0(\BC^{(0)},\tau)\in (0,1)$ such that the following is true: if $\epsilon\in (0,\epsilon_0)$, $\BC\in \FL_{\epsilon}(\BC^{(0)})\cap \FL_I$, $V\in \CN_{\epsilon}(\BC)$, and $U,u$ are as in Lemma \ref{lemma:coarse_graphical_rep}, then for every $Z = (\xi,\eta)\in B_{3/4}$ with $\Theta_V(Z)\geq \Theta_{\BC^{(0)}}(0) = \frac{5}{2}$ we have:
	\begin{enumerate}
		\item [(i)] $\dist(Z,S(\BC)) \leq CE_{V,\BC}$\textnormal{;}
		\item [(ii)] Writing $(e_j)_{j=1}^{n+1}$ for the standard basis vectors on $\R^{n+1}$, 
		$$\int_{B_{3/4}}\sum^{n+1}_{j=3}|e^{\perp_{T_XV}}_j|^2\ \ext\|V\|(X) \leq CE_{V,\BC}^2;$$
		\item [(iii)] $$\int_{B_{3/4}}\frac{\dist^2(X,\spt\|\BC\|)}{|X-Z|^{n-1/2}}\ \ext\|V\|(X) \leq CE_{V,\BC}^2;$$
		\item [(iv)] Writing $\xi^\perp(X)$ for the projection of $(\xi,0)$ onto $T_X^\perp\BC$ (which is just $\xi\cdot \mathbf{n}_X$, for $\mathbf{n}_X$ the unit normal in $\R^2$ to the ray in the cross-section $\BC_0$ in $\BC$ whose corresponding half-hyperplane in $\BC$ contains $X$),
		$$\int_{U_\tau}\frac{|u(X)-\xi^\perp(X)|^2}{|X + u(X)-Z|^{n+3/2}} \leq CE_{V,\BC}^2;$$
		\item [(v)] $$\int_{U\cap B_{3/4}}R^{2-n}\left|\frac{\del(u/R)}{\del R}\right|^2 \leq CE_{V,\BC}^2.$$
	\end{enumerate}
	Here, $C = C(n)$.
\end{lemma}

\textbf{Remark:} Let us briefly discuss the significance of each inequality in Lemma \ref{lemma:L2_coarse}. (i) tells us that ``good'' singular points, i.e. those of sufficiently high density, are not just $\tau$-close to $S(\BC^{(0)})$, but are significantly closer; this will be used at various points when we need to combine estimates at different points of density $\geq \frac{5}{2}$. (ii) is a bound on the tilt excess, \textit{in directions parallel to }$S(\BC^{(0)})$, in terms of the height excess. As such, it will control in $L^2$ the derivatives parallel to the spine of the coarse blow-ups; this will be needed to prove the regularity of the boundary values of the coarse blow-up. (iii), along with the absence of density gaps (Lemma \ref{lemma:gaps}) will show that the height excess cannot accumulate along the spine, giving strong $L^2$ convergence of the blow-up sequence globally on $\spt\|\BC^{(0)}\|\cap B_1$ instead of just locally away from $S(\BC^{(0)})$. (iv) is almost a type of $L^2$ bound on the full derivative $Du$, except we need to subtract a small power in the denominator, i.e. we have $n+3/2$ as opposed to $n+2$. Thus, the inequality is not strong enough to achieve global control on the $W^{1,2}$ norm of the blow-up sequence, and so we can only control the $W^{1,2}$ norm locally away from $S(\BC^{(0)})$. Nonetheless, (iv) we still show that the boundary values of the blow-up are always determined by a \textit{single} function, and that the blow-ups are always $C^{0,\alpha}$ up-to-the-boundary. Finally, (iv) is the Hardt--Simon inequality: it will be key for studying the boundary regularity of the blow-ups. We shall not need the Hardt--Simon inequality in the level 0 setting, as there we will be able to instead appeal to classical boundary regularity theory for harmonic functions.

\begin{proof}
	We only outline the proof, pointing out how the corresponding argument in \cite[Section 3]{simoncylindrical} can be modified to this setting. As usual, we are always working with $\epsilon \in (0,\epsilon_*)$, where $\epsilon_* = \epsilon_*(n)$ is sufficiently small so that $\CN_{\epsilon_*}(\BC)$ is contained with a multiplicity two class (as in Theorem \ref{thm:M2C}).
	
	\textbf{Step 1: $Z=0$.} Let us first consider the case where $Z = 0$ has $\Theta_V(0) \geq \Theta_{\BC^{(0)}}(0) = \frac{5}{2}$. Then the monotonicity formula gives (in both distributional sense and for a.e. $\rho\in (0,1)$):
	\begin{align*}
		n\rho^{n-1}\int_{B_\rho}\frac{|X^\perp|^2}{|X|^{n+2}}\ \ext\|V\|(X) & = \frac{\ext}{\ext \rho}\left[\rho^n\int_{B_\rho}\frac{|X^\perp|^2}{|X|^{n+2}}\ \ext\|V\|(X)\right] - \rho^n\frac{\ext}{\ext\rho}\int_{B_\rho}\frac{|X^\perp|^2}{|X|^{n+2}}\ \ext\|V\|(X)\\
		& \leq \frac{\ext}{\ext\rho}\left(\|V\|(B_\rho)-\Theta_V(0)\cdot\w_n\rho^n\right)\\
		& \leq \frac{\ext}{\ext\rho}\left(\|V\|(B_\rho) - \|\BC^{(0)}\|(B_\rho)\right)
	\end{align*}
	where in the first inequality we have used that the last term on the first line is positive, and in the second inequality we have used that $\frac{\ext}{\ext\rho}(\Theta_V(0) \w_n\rho^n) = n\rho^{n-1}\w_n\Theta_V(0) \geq n\rho^{n-1}\w_n\Theta_{\BC^{(0)}} = \frac{\ext}{\ext\rho}(\Theta_{\BC^{(0)}}(0)\cdot\w_n\rho^n)$, and $\Theta_{\BC^{(0)}}(0)\cdot\w_n\rho^n = \|\BC^{(0)}\|(B_\rho)$. Now choosing $\psi:\R\to [0,1]$ a decreasing $C^1$ function with $\left.\psi\right|_{(-\infty,7/8)}\equiv 1$ and $\left.\psi\right|_{(15/16,\infty)}\equiv 0$, we can multiply the above by $\psi^2(\rho)$ and, noting that $\psi^2(\rho)\frac{\ext}{\ext\rho}\|V\|(B_\rho) = \frac{\ext}{\ext \rho}\int_{B_\rho}\psi^2(|X|)\ \ext\|V\|$ for a.e.\;$\rho\in (0,1)$ (and similarly for the $\psi^2(\rho)\frac{\ext}{\ext\rho}\|\BC^{(0)}\|(B_\rho)$ term), if integrate over $\rho\in (0,1)$ we get:
	\begin{equation*}\
	n\int^1_0\psi^2(\rho)\rho^{n-1}\left(\int_{B_\rho}\frac{|X^\perp|^2}{|X|^{n+2}}\ \ext\|V\|(X)\right)\ \ext\rho \leq \int_{B_1}\psi^2(|X|)\ \ext\|V\|(X) - \int_{B_1}\psi^2(|X|)\ \ext\|\BC\|(X)
	\end{equation*}
	and hence, as $\psi|_{[3/4,7/8]} \equiv 1$, we get
	\begin{equation}\label{E:L2-1}
		n\left(\frac{3}{4}\right)^{n-1}\cdot\frac{1}{8}\int_{B_{3/4}}\frac{|X^\perp|^2}{|X|^{n+2}}\ \ext\|V\|(X) \leq \int_{B_1}\psi^2(R)\ \ext\|V\|(X) - \int_{B_1}\psi^2(R)\ \ext\|\BC\|(X)
	\end{equation}
	where recall that we write $R = |X|$ and $r=|x|$. By the same application of the first variation formula as in \cite[Lemma 3.4, (3)]{simoncylindrical}, we also have
	\begin{equation}\label{E:L2-2}
	\begin{split}
	\int_{B_1}&\left(1+\frac{1}{2}\sum^{n+1}_{j=3}|e^\perp_{j}|^2\right)\psi^2(R)\ \ext\|V\|\\
	&\hspace{3em} \leq C\int_{B_1}|(x,0)^\perp|^2\left(\psi^2(R)+\psi'(R)^2\right)\ \ext\|V\| - 2\int_{B_1}r^2R^{-1}\psi(R)\psi'(R)\ \ext\|V\|.
	\end{split}
	\end{equation}
	Now, if $(x,y) \in \graph(u)$, then $(x,y) = (x',y) + u(x',y)$ for some $x'\in \spt\|\BC\|$, and we have $(x,0)^\perp = u(x',y) + (P_{(x,y)} - Q_{(x',y)})(x,0)$, where $P_{(x,y)}$ and $Q_{(x',y)}$ denote the orthogonal projections onto $T_{(x,y)}^\perp V$ and $T_{(x',y)}^\perp\BC$ respectively; note that is true for $\H^n$-a.e. such point, as $V$ has a (unique) tangent plane at every branch point, and thus everywhere except a set which is $\H^{n-1}$-null; if $(x,y)$ is a point where $V$ is locally expressed as a two-valued function, if $(x,y)$ is a branch point it does not matter which choice of value of $u$ we take, and otherwise if $(x,y)\in \reg(V)$, we simply take the value of $u$ which locally describes $V$ about $(x,y)$. But also note that $|P_{(x,y)}-Q_{(x',y)}| \leq C|Du(x',y)|$, where $C = C(n)$, and hence using this in (\ref{E:L2-2}) we get (noting $\psi(R) \equiv 0$ for $R\geq 15/16$):
	\begin{equation}\label{E:L2-3}
	\begin{split}
	\int_{B_1}\left(1+\frac{1}{2}\sum^{n+1}_{j=3}|e^\perp_{j}|^2\right)\psi^2(R)\ \ext\|V\| \leq C\int_{U\cap B_{15/16}}&|u|^2 + r^2|Du|^2 + C\int_{B_{15/16}\backslash\graph(u)}r^2\ \ext\|V\|\\
	& - 2\int_{\graph(u)\cap B_1}r^2R^{-1}\psi(R)\psi'(R)\ \ext\|V\|
	\end{split}
	\end{equation}
	where as the last integral only takes place over $R\in (7/8,15/16)$, this region will be graphical and so certainly we change the domain of integration to $B_1\cap \graph(u)$ as opposed to $B_1$ to stress this; of course, over any region where $u$ is two-valued, by $|u|^2$ we mean $|u_1|^2 + |u_2|^2$, etc.
	
	By the same simple 1-dimensional integration argument as in \cite[Lemma 3.4, (6)]{simoncylindrical}, as $\spt\|\BC^{(0)}\|$ is comprised of half-hyperplanes we readily see that
	$$\int_{B_1}\psi^2(R)\ \ext\|\BC\| = -2\int_{B_1}r^2R^{-1}\psi(R)\psi'(R)\ \ext\|\BC\|$$
	(and indeed the right hand side is $\geq -2\int_{U\cap B_1}r^2R^{-1}\psi(R)\psi'(R)\ \ext\|\BC\|$), and similarly for $\graph(u)\cap B_1$ we have
	$$\int_{\graph(u)\cap B_1}r^2R^{-1}\psi(R)\psi'(R)\ \ext\|V\| = \int_{U\cap B_1}r_u^2R^{-1}_u\psi(R_u)\psi'(R_u)\cdot \sqrt{g}$$
	where $r_u^2 = |x|^2 + |u(x,y)|^2 \equiv r^2 + |u(x,y)|^2$, $R_u^2 = |x|^2 + |u(x,y)|^2 + |y|^2 \equiv R^2 + |u(x,y)|^2$, and $g$ is the volume element of $\graph(u)$; note that we know $1\leq \sqrt{g} \leq 1 + C|Du|^2$. Of course, over regions where $u$ is two-valued, we need to understand these terms as a sum of two terms taking the same form, one for each value of $u$.  Thus combining this with (\ref{E:L2-3}) we get
	\begin{align*}
	\frac{1}{2}\int_{B_{3/4}}&\sum^{n+1}_{j=3}|e^\perp_j|^2\ \ext\|V\| + \int_{B_1}\psi^2(R)\ \ext\|V\| - \int_{B_1}\psi^2(R)\ \ext\|\BC\|\\
	& \leq C\int_{U\cap B_{15/16}}|u|^2 + r^2|Du|^2 + C\int_{B_{15/16}\backslash\graph(u)}r^2\ \ext\|V\|\\
	& \hspace{2em} - 2\left(\int_{U\cap B_1}r_u^2R_u^{-1}\psi(R_u)\psi'(R_u)\sqrt{g} - \int_{U\cap B_1}r^2R^{-1}\psi(R)\psi'(R)\ \ext\|\BC\|\right)\\
	& \leq \tilde{C}\int_{U\cap B_{15/16}}|u|^2 + r^2|Du|^2 + C\int_{B_{15/16}\backslash \graph(u)}r^2\ \ext\|V\|
	\end{align*}
	where $\tilde{C} = \tilde{C}(n)$. Applying Lemma \ref{lemma:coarse_graphical_rep} (the same argument goes through working on $B_{15/16}$ as opposed to $B_{3/4}$, up to changing (dimensional) constants) and combining the above with (\ref{E:L2-1}) we arrive at
	\begin{equation}\label{E:L2-4}
		\int_{B_{3/4}}\frac{|X^\perp|^2}{|X|^{n+2}}\ \ext\|V\| + \int_{B_{3/4}}\sum_{j=3}^{n+1}|e^\perp_j|^2\ \ext\|V\| \leq CE_{V,\BC}^2
	\end{equation}
	where $C = C(n)$; in particular this establishes (ii) when $Z=0$.
	
	By exactly the same argument based on the first variation formula as \cite[Lemma 3.4]{simoncylindrical}, one may derive
	\begin{equation}\label{E:L2-extra2}
	\int_{B_{3/4}}\frac{\dist^2(X,\spt\|\BC\|)}{|X|^{n+(2-\alpha)}}\ \ext\|V\|(X) \leq C\int_{B_1}\zeta^2\frac{|X^\perp|^2}{|X|^{n+(2-\alpha)}} + \frac{\dist^2(X,\spt\|\BC\|)}{|X|^{n-\alpha}}|\nabla\zeta|^2\ \ext\|V\|(X)
	\end{equation}
	where $C = C(n,\alpha)$ and $\zeta\in C^\infty(\R^{n+1})$ obeys $\zeta|_{B_{7/8}(0)} \equiv 1$ and $\zeta|_{\R^{n+1}\backslash B_1}\equiv 0$, with $|\nabla\zeta|\leq 16$. Clearly taking $\alpha = 1/2$ in the above and using (\ref{E:L2-4}) (again, we can change the domains of integration to $B_{7/8}$ by re-running the argument on the larger ball), we see that we arrive at (iii) when $Z=0$.
	
	To complete Step 1 of the proof, notice that for $X = (x,y)\in (\B\cup \reg(V))\cap \graph(u)$ (which is $\|V\|$-a.e. point in $\spt\|V\|$), in exactly the same way as in the argument leading up to \cite[Lemma 3.4, (11)]{simoncylindrical} we have
	$$(x,y)^\perp = -R^2\left(\frac{\del}{\del R}\left(\frac{u(x',y)}{R}\right)\right)^\perp$$
	where $(x',y)\in U$ is such that $x'$ is the nearest point projection of $x$ onto the cross-section $\BC_0$, and at a branch point the choice of value of $u$ does not matter and away from the branch set we mean the value of $u$ which locally expresses $V$ about this point. Therefore, by reducing $\beta = \beta(n)$ if necessary to ensure that $\|P_{(x,y)}-Q_{(x',y)}\|\leq 1/2$, we get
	$$|(x,y)^\perp| \geq \frac{1}{2}R^2\left|\frac{\del}{\del R}\left(\frac{u(x',y)}{R}\right)\right|$$
	and thus combining this with (\ref{E:L2-4}) we arrive at (v).
	
	\textbf{Step 2: (i) and translating $Z$.}
	We first claim that there is $\epsilon^\prime = \epsilon^\prime(n)$ and $\vartheta = \vartheta(n)\in (0,1)$ such that if $\BC^{(0)}\in \FL_S\cap \FL_I$, $\BC\in \FL_{\epsilon^\prime}(\BC^{(0)})\cap \FL_I$, $V\in \CN_{\epsilon^\prime}(\BC^{(0)})$, $Z = (\xi,\zeta)\in \sing(V)\in B_{3/4}$ with $\Theta_V(Z)\geq \frac{5}{2}$, then any $X = (x,y)$ obeying $|x|\geq \vartheta^{-1}(|\xi|+\dist(X,\spt\|\BC\|))$ has
	\begin{equation}\label{E:L2-extra}
	\dist(X,\spt\|(\tau_Z)_\#\BC\|) = |(x,y) - (x',y) - \xi^\perp| + R
	\end{equation}
	where $x'$ is the nearest point projection of $x$ onto $\BC_0$ (in particular $|(x,y)-(x',y)| = \dist(X,\spt\|\BC\|)$), and $\xi^\perp$ is the projection of $(\xi,0)$ onto $T^\perp_{(x',0)}\BC$, and $|R|\leq C|x|^{-1}|\xi|^2$, where $C = C(n)$ . Indeed, if this were not true then taking $\epsilon^\prime = 1/k$ and $\vartheta = 1/k$, we could find sequences $\BC^{(0)}_k\in \FL_S\cap \FL_I$, $\BC_k\in \FL_{1/k}(\BC^{(0)}_k)\cap \FL_I$, $V\in \CN_{1/k}(\BC^{(0)}_k)$, $Z_k = (\xi_k,\zeta_k)\in \sing(V_k)\cap B_{3/4}$ with $\Theta_{V_k}(Z_k) \geq \frac{5}{2}$, and $X_k = (x_k,y_k)$ obeying $|x_k|\geq k(|\xi_k| + \dist(X_k,\spt\|\BC_k\|))$, yet the conclusion fails. In particular, this tells us that, after passing to a subsequence, that $\BC^{(0)}_k \weakly \BC^{(0)}\cap \FL_S$ and $x_k/|x_k| \to a\in \spt\|\BC^{(0)}_0\|\cap S^1$. But $|\xi_k|/|x_k| \to 0$, and so as all the cones are of the same level, it follows that the equality must be true for all $k$ sufficiently large (just by a simple geometric argument and calculation, easiest seen by rescaling by $|x_k|$; the error term at this scale is $C(|x_k|^{-1}|\xi|)^2$ by Taylor's theorem, and allows for degeneration of a level $I$ cone to a cone of a higher level), which provides the contradiction. Note that (\ref{E:L2-extra}) readily implies
	\begin{equation}\label{E:L2-extra3}
	|\xi^\perp| \leq \dist(X,\spt\|\BC\|) + \dist(X,\spt\|(\tau_Z)_\#\BC\|) + |R|
	\end{equation}
	We now claim that there is a constant $\delta = \delta(n)>0$ such that the following holds: for each $\rho\in (0,1/4)$, there is a constant $\epsilon_0 = \epsilon_0(n,\rho)$ such that if $\BC^{(0)}\in \FL_S\cap \FL_I$, $V\in \CN_{\epsilon_0}(\BC^{(0)})$, $\BC\in \FL_{\epsilon_0}(\BC^{(0)})$, $a\in \R^2$, and $Z = (\xi,\zeta)\in \sing(V)\cap B_{1/2}(0)$ is not a density $2$ branch point or density $2$ classical singularity, then
	\begin{equation}\label{E:L2-5}
		\|V\|(\{X\in B_{\rho}(Z)\cap \{|x|>\rho/10\}: |a^\perp|\geq \delta|a|\text{ and }|x|\geq\vartheta^{-1}(|\xi|+\dist(X,\spt\|\BC\|))\}) \geq \delta\rho^n 
	\end{equation}
	where here $a^\perp$ at $X$ is $a^{\perp_{T_x\BC_0}}$, the orthogonal projection of $a$ onto $T^\perp_x\BC_0$, and $\vartheta = \vartheta(n)$ is as in (\ref{E:L2-extra}). Indeed, if this does not hold then for each $\delta>0$, there is $\rho>0$ such that with $\epsilon_j = 1/j$, there exists $\BC^{(0)}_j\in \FL_S\cap \FL_I$, $V_j\in \CN_{\epsilon_j}(\BC^{(0)}_j)$, $\BC_j\in \FL_{\epsilon_j}(\BC^{(0)})$, $a_j\in S^1$, and $Z_j\in \sing(V_j)\cap B_{1/2}$ which is not a density 2 branch point or density $2$ classical singularity such that
	$$\|V_j\|(\{X\in B_\rho(Z_j)\cap\{|x|>\rho/10\}: |a^\perp_j|\geq\delta\text{ and }|x|\geq \vartheta^{-1}(|\xi_j| + \dist(X,\spt\|\BC_j\|))\}) < \delta\rho^n.$$
	After passing to a subsequence, we may assume that $\BC^{(0)}_j\weakly \BC^{(0)}\in \FL_S\cap \FL_I$, $\BC^j\weakly \BC^{(0)}$, $V_j\weakly \BC^{(0)}$, $Z_j\to Z$ for some $Z\in \{0\}\times\overline{B}^{n-1}_{1/2}(0)$ (by Lemma \ref{lemma:coarse_graphical_rep}), and $a_j\to a\in S^1$ such that
	\begin{equation}\label{E:L2-6}
		\|\BC^{(0)}\|(\{X\in B_\rho(Z)\cap\{|x|>\rho/10\}: |a^\perp|\geq \delta\text{ and }|x|\geq\vartheta^{-1}\dist(X,\spt\|\BC^{(0)}\|)\}) < 2\delta\rho^n.
	\end{equation}
	Thus we have shown that if (\ref{E:L2-5}) is false, then for every $\delta>0$ there is a $\BC^{(0)}\in \FL_S$, $\rho>0$, $Z\in \{0\}\times\overline{B}^{n-1}_{1/2}(0)$ and $a\in S^1$ such that (\ref{E:L2-6}) holds. By translating by $Z$ and rescaling by $\rho$, we may without loss of generality assume that $Z=0$ and $\rho=1$ (indeed, $\BC^{(0)}$ is translation invariant by $Z$ and using how the quantities scale). Thus, taking $\delta_j = 1/j$, we can find a new sequence $\BC^{(0)}_j\in \FL_S\cap \FL_I$ and $a_j\in S^1$ such that
	$$\|\BC^{(0)}_j\|(\{X\in B_1\cap\{|x|>1/10\}:|a^\perp_j|\geq 1/j\text{ and }|x|\geq\vartheta^{-1}\dist(X,\spt\|\BC^{(0)}_j\|)\}) < 2/j$$
	and so again passing to a subsequence to ensure that $\BC^{(0)}_j\weakly \BC^{(0)}\in \FL_S\cap \FL_I$, we get
	$$\|\BC^{(0)}\|(\{X\in B_1\cap\{|x|>1/10\}: |a^\perp|>0 \text{ and }|x|\geq\vartheta^{-1}\dist(X,\spt\|\BC^{(0)}\|)\}) = 0$$
	i.e. $a^\perp = 0$ for $\H^{n}$-a.e. $X\in B_1\cap\{|x|>1/10\}$ which obeys $|x|\geq\vartheta^{-1}\dist(X,\spt\|\BC^{(0)}\|)$, which is obviously false as $a\in S^1$ is fixed and $\BC^{(0)}$ has cross-section whose unit vectors span $\R^2$.
	
	Now for $\rho\in (0,1/4)$ (to be chosen only depending on $n$), let $\epsilon_0 = \epsilon_0(n,\rho)$ be as in (\ref{E:L2-5}); we know that if $\epsilon_0$ is sufficiently small then Lemma \ref{lemma:coarse_graphical_rep} will apply on $B_{3/4}\cap\{|x|>\rho/10\}$, and so up to dimensional constants we may pass between $\ext\|V\|$ and $\ext\|\BC^{(0)}\|$ on this region. Then, for $Z = (\xi,\zeta)\in \sing(V)$ such that $\Theta_V(Z)\geq \Theta_{\BC^{(0)}}(0) = \frac{5}{2}$, take $a = \xi$ in (\ref{E:L2-5}) to obtain, for some set $S\subset B_\rho(Z)\cap \{X:|x|\geq\vartheta^{-1}(|\xi|+\dist(X,\spt\|\BC\|))\}$ with $\|V\|(S)\geq \delta\rho^n$,
	\begin{equation}\label{E:L2-7}
	\delta^2|\xi|^2\cdot \delta\rho^n \leq \int_S \delta^2|\xi|^2\ \ext\|V\| \leq \int_S |\xi^\perp|^2\ \ext\|V\|.
	\end{equation}
	Using (\ref{E:L2-extra3}) in (\ref{E:L2-7}) we get
	\begin{equation}\label{E:L2-8}
	\begin{split}
	\delta^3\rho^n|\xi|^2 & \leq 4\int_{B_\rho(Z)}\dist^2(X,\spt\|(\tau_Z)_\#\BC\|)\ \ext\|V\|(X) + 4\int_{B_\rho(Z)}\dist^2(X,\spt\|\BC\|)\ \ext\|V\|\\
	& \hspace{5em} + 4C\int_{B_\rho(Z)\cap\{|x|>\rho/10\}\cap \{X:|x|\geq\vartheta^{-1}(|\xi|+\dist(X,\spt\|\BC\|))\}}|x|^{-2}|\xi|^4\ \ext\|V\|  
	\end{split}
	\end{equation}
	We now need to deal with these terms individually; let us start with the first. Note that by the triangle inequality we trivially have
	\begin{equation}\label{E:L2-9}
	|\dist(X,\spt\|(\tau_Z)_\#\BC\|) - \dist(X,\spt\|\BC\|)| \leq |\xi|
	\end{equation}
	as $(\tau_Z)_\#\BC = (\tau_{(\xi,0)})_\#\BC$, and thus as we know for each $\epsilon>0$ there is a $\delta = \delta(\epsilon)$ with $\delta(\epsilon)\to 0$ as $\epsilon\to 0$ such that if $V\in \CN_{\epsilon}(\BC^{(0)})$ then $|\xi|<\delta(\epsilon)$, it follows that
	$$4^{-n-2}\int_{B_{1/4}(Z)}\dist^2(X,\spt\|(\tau_Z)_\#\BC\|)\ \ext\|V\| \leq C\int_{B_1}\dist^2(X,\spt\|\BC\|)\ \ext\|V\| + C|\xi|^2 \leq C(\epsilon+\delta(\epsilon))$$
	and so we can, for $\epsilon = \epsilon(n)$ sufficiently small, apply the results of Step 1 (namely (\ref{E:L2-extra2})) to $(\eta_{Z,1/4})_\#V$ in place of $V$ to obtain
	\begin{equation*}
	\begin{split}
	\rho^{-n-3/2}\int_{B_\rho(Z)}\dist^2(X,\spt\|(\tau_Z)_\#\BC\|)\ \ext\|V\| & \leq C\int_{B_1}\dist^2(X,\spt\|(\tau_Z)_\#\BC\|)\ \ext\|V\|\\
	& \leq C\int_{B_1}\dist^2(X,\spt\|\BC\|)\ \ext\|V\| + C|\xi|^2.
	\end{split}
	\end{equation*}
	The second term in (\ref{E:L2-8}) we leave as it is. For the third term directly compute, passing to $\BC^{(0)}$,
	\begin{align*}
	\int_{B_\rho(Z)\cap\{|x|>\rho/10\}\cap \{X:|\xi|\leq \vartheta |x|\}}|x|^{-2}|\xi|^4\ \ext\|V\| \leq C|\xi|^4\rho^{n-2}\cdot (\rho^{-1}+|\xi|^{-1}).
	\end{align*}
	Hence combining everything in (\ref{E:L2-8}), we get 
	$$\rho^n|\xi|^2 \leq C\int_{B_1}\dist^2(X,\spt\|\BC\|)\ \ext\|V\| + C\rho^n|\xi|^2\left(\rho^{3/2} + |\xi|^2\rho^{-3}+|\xi|\rho^{-2}\right)$$
	where $C = C(n)$. Hence, choosing $\rho = \rho(n)$ so that $C\rho^{3/2} < 1/4$, and then choosing $\tau = \tau(n)$ such that $C(\tau^2\rho^{-3}+\tau\rho^{-2}) < 1/4$, we get that if we choose $\epsilon < \epsilon_*$, where $\epsilon_* = \epsilon_*(n)$ is such that if $V\in \CN_{\epsilon_*}(\BC^{(0)})$, then $|\xi|<\tau$ for each $Z = (\xi,\zeta)\in\sing(V)\cap B_{3/4}$ with $\Theta_{V}(Z)\geq 5/2$, then we get for such $\epsilon = \epsilon(\BC^{(0)},n)$,
	$$|\xi|^2 \leq C\int_{B_1}\dist^2(X,\spt\|\BC\|)\ \ext\|V\|.$$
	as desired. Combining this with (\ref{E:L2-9}) then also gives
	\begin{equation}\label{E:L2-10}
	\int_{B_1}\dist^2(X,\spt\|(\tau_Z)_\#\BC\|)\ \ext\|V\|(X)\leq C\int_{B_1}\dist^2(X,\spt\|\BC\|)\ \ext\|V\|(X) \equiv CE_{V,\BC}^2.
	\end{equation}
	\textbf{Step 3: Conclude.} Combining (\ref{E:L2-10}) with Step 1 (namely (iii), applied with $(\eta_{Z,1/4})_\#V$ in place of $V$) and using (\ref{E:L2-9}) and (i) (whose truth is established in Step 2) we get
	$$\int_{B_{1/4}(Z)}\frac{\dist^2(X,\spt\|\BC\|)}{|X-Z|^{n-1/2}}\ \ext\|V\| \leq CE_{V,\BC}^2$$
	which readily establishes (iii). Now we just need to establish (iv). Note that by (\ref{E:L2-extra}) that, for any $P\in \graph(u)\cap (\reg(V)\cup \B)$, where if we write $P = (x,y)+u(x,y)$ (for some choice of value $u$ if $u$ is two-valued about $(x,y)$) then $(x,y)\in U_\tau$, then
	$$\dist(P,\spt\|(\tau_Z)_\#\BC\|) = |u(x,y) -\xi^\perp| + R$$
	where now we have $|x|>\tau$, and so $|R|\leq C\tau^{-1}|\xi|^2$. Using the same application of (iii) to $(\eta_{Z,1/4})_\#V$, we have
	$$\int_{U_\tau\cap B_{1/4}(Z)}\frac{|u(x,y)-\xi^\perp|^2}{|X + u(X) -Z|^{n+3/2}} \leq CE_{V,\BC}^2$$
	in the same way as \cite[Theorem 3.1]{simoncylindrical} (up to using (\ref{eqn:area_formula}) to change the domain of integration over regions where $V$ is represented by a two-valued function); this uses the fact that we have established already that $|\xi|^2\leq C\hat{E}_{V,\BC}^2$, and so if we choose $\epsilon_0<\tau^2$, then $|\xi|<C\tau^2$. This completes the proof.
\end{proof}

Combining Lemma \ref{lemma:L2_coarse} (iii) -- (iv) with the Lemma \ref{lemma:gaps}, we can establish the following two inequalities which will be of critical importance:

\begin{corollary}\label{corollary:non-concentration}
	Let $\tau,\delta\in (0,1/10)$. Then there exists $\epsilon_1 = \epsilon_1(\BC^{(0)},\tau)\in (0,1)$ such that the following is true: if $\epsilon\leq \min\{\epsilon_1,\delta\}$, $\BC\in \FL_{\epsilon}(\BC^{(0)})\cap \FL_I$, $V\in \CN_{\epsilon}(\BC^{(0)})$, and $U,u$ are as in Lemma \ref{lemma:coarse_graphical_rep}, then:
	\begin{enumerate}
		\item [(i)] For each $\rho\in (0,1/4)$, if we allow $\epsilon_1$ to depend on $\rho$ also, we get that for each $Z = (\xi,\zeta)\in \spt\|V\|\cap B_{3/8}$ with $\Theta_V(Z)\geq \Theta_{\BC^{(0)}}(0) = \frac{5}{2}$, we have
		$$\int_{B_{\rho/2}(Z)}\frac{\dist^2(X,\spt\|(\tau_Z)_\#\BC\|)}{|X-Z|^{n+3/2}}\ \ext\|V\| \leq C\rho^{-n-3/2}\int_{B_\rho(Z)}\dist^2(X,\spt\|(\tau_Z)_\#\BC\|)\ \ext\|V\|;$$
		\item [(ii)] Writing $r_\delta := \max\{r,\delta\}$,
		$$\int_{B_{1/2}}\frac{\dist^2(X,\spt\|\BC\|)}{r_\delta^{1/2}}\ \ext\|V\| \leq CE_{V,\BC}^2;$$
	\end{enumerate}
	here, $C = C(n)$.
\end{corollary}

\begin{proof}
	Given Lemma \ref{lemma:L2_coarse}, the proof of (ii) is identical to that of \cite[Corollary 3.2]{simoncylindrical}; as such we omit the details and refer the reader to \cite{simoncylindrical}. To see (i), we wish to apply Lemma \ref{lemma:L2_coarse} with $(\eta_{Z,\rho})_\#V$ in place of $V$ (and $0$ in place of $Z$); indeed, we have
	$$\int_{B_1}\dist^2(X,\spt\|\BC\|)\ \ext\|(\eta_{Z,\rho})_\#V\| \leq 2\rho^{-n-2}\int_{B_\rho(Z)}\dist^2(X,\spt\|\BC\|)\ \ext\|V\| + C\rho^{-2}|\xi|^2$$
	from which is follows (as we can bound the other term in $Q_{V,\BC}$ similarly) from Lemma \ref{lemma:L2_coarse}(i) that if $\epsilon$ is sufficiently small (depending on $\rho$) we can apply Lemma \ref{lemma:L2_coarse} to deduce (i) holds.
\end{proof}

\textbf{Remark:} Note that Corollary \ref{corollary:non-concentration}(ii) gives us that, for any $\delta\in (0,1)$, there exists $\epsilon_1 = \epsilon_1(\BC^{(0)},\tau)$ such that if $\epsilon\leq \min\{\epsilon_1,\delta\}$, then
$$\int_{\{X:|x|<\delta\}\cap B_{1/2}}\dist^2(X,\spt\|\BC\|)\ \ext\|V\| \leq C\delta^{1/2}\int_{B_1}\dist^2(X,\spt\|\BC\|)\ \ext\|V\|$$
where $C = C(n)$ is in particular independent of $\delta$; this will be used to show that the height excess does not concentrate along the spine when we take coarse blow-ups, and hence will give strong convergence in $L^2$ on all of $\spt\|\BC^{(0)}\|\cap B_{1/2}$.

\subsection{Construction of the Coarse Blow-Up Class}\label{sec:coarse-construction}

Fix $I\in \{0,1,2\}$ and $\BC^{(0)}\in \FL_S\cap \FL_I$. Let $(\epsilon_k)_k$ and $(\delta_k)_k$ be arbitrary sequences obeying $0<\epsilon_k <\delta_k\to 0$. Let $V_k\in \CN_{\epsilon_k}(\BC^{(0)})$ and $\BC_k\in \FL_{\epsilon_k}(\BC^{(0)})$ be arbitrary sequences. Now let $(\tau_k)_k$ be a sequence with $\tau_k\downarrow 0$ sufficiently slower so that $\tau_k^{-1}E_{V_k,\BC_k}\to 0$ and the conclusions of Lemma \ref{lemma:coarse_graphical_rep} hold with $\epsilon_k$, $\tau_k$, $\BC_k$, and $V_k$ in place of $\epsilon$, $\tau$, $\BC$, and $V$, respectively.

Then Lemma \ref{lemma:coarse_graphical_rep} gives the existence of functions $u_k\in C^{1,1/2}(\BC_k\res U_k)$, where $U_k:= U_{\tau_k}\cap B_{3/4}$, with, the estimates in Lemma \ref{lemma:L2_coarse} and Corollary \ref{corollary:non-concentration} holding for all $k$ sufficiently large. Now, one may work relative to the base cone $\BC^{(0)}$ to provide a common domain of definition, i.e. we may find a function $\psi_k\in C^2(\spt\|\BC^{(0)}\|\cap B_1\backslash\{|x|<\tau_k/2\};\spt\|\BC^{(0)}\|^\perp)$ such that $\spt\|\BC_j\|\cap B_1\backslash\{|x|<\tau_k\}\subset\graph(\psi_j)$ and $|\psi_k|_{C^2}\leq C\epsilon_k\to 0$; note that $\psi_k$ defines a corresponding region to $U_k$ in $\spt\|\BC_k\|$ in $\spt\|\BC^{(0)}\|$, which for simplicity we shall also denote by $U_k$. So let us now define the \textit{blow-up sequence} by
$$v_k:= E_{V_k,\BC_k}^{-1}u_k(x+\psi_k(x)).$$
For notational simplicity, we shall write $E_{k}:= E_{V_k,\BC_k}$; we also extend $v_k$ by $0$ to all of $\spt\|\BC^{(0)}\|\cap B_1$. It then follows, by Lemma \ref{lemma:coarse_graphical_rep} and standard elliptic estimates (over half-hyperplanes in $\spt\|\BC^{(0)}\|$ where $v_k$ is represented by a single-valued function) and from the estimates in Section \ref{sec:two-valued_stationary_graphs} and Theorem \ref{thm:wick1} (over half-hyperplanes in $\spt\|\BC^{(0)}\|$ where $v_k$ is represented by a two-valued function) that for each compact subset $K\subset B_1\backslash\{|x|=0\}$, for all $k$ sufficiently large (depending on $K$):
$$|v_k|_{C^{1,1/2}(\BC^{(0)}\res K} \leq C$$
for some $C = C(n,K)$ independent of $k$; hence $v_j$ converges in $C^1$ on each compact subset of $\spt\|\BC^{(0)}\|\cap B_1\backslash\{|x|=0\}$, to a limit function $v\in C^{1,1/2}(\BC^{(0)}\res B_1)$ which obeys $\Delta_{\BC^{(0)}}v = 0$, by which we mean over each half-hyperplane in $\spt\|\BC^{(0)}\|$, $v$ is either a single-valued harmonic function or two-valued $C^{1,1/2}$ harmonic function.

\begin{defn}
	We call $v$ constructed in the above manner a \textit{coarse blow-up} of $(V_k)_k$ relative to $(\BC_k)_k$; we write $\FB(\BC^{(0)})$ for the collection of all such coarse blow-ups.
\end{defn}

Let us remark some basic properties of coarse blow-ups from Lemma \ref{lemma:L2_coarse} and Corollary \ref{corollary:non-concentration}. Firstly, from Corollary \ref{corollary:non-concentration}(ii) we have for any fixed $\delta\in (0,1)$, for all $k$ sufficiently large:
$$\int_{B_{1/2}}\frac{\dist^2(X,\spt\|\BC_k\|)}{r^{1/2}_{\delta}}\ \ext\|V\| \leq CE_{V_k,\BC_k}^2$$
which in particular gives
$$\int_{B_{1/2}\cap U_k\cap \{|x|<\delta\}}|u_k|^2 \leq C\delta^{1/2}E_{V_k,\BC_k}^2$$
i.e.
$$\int_{B_{1/2}\cap U_k\cap \{|x|<\delta\}}|v_k|^2 \leq C\delta^{1/2}$$
which in fact tells us that we have $v_k\to v$ strongly in $L^2(\BC\res B_{1/2})$ (indeed, for each $0<\tau<\delta$ we have $v_k\to v$ in $C^2$ on $B_{1/2}\cap \{\tau < |x| <\delta\}$, and thus $\int_{B_{1/2}\cap\{\tau<|x|<\delta\}}|v|^2 \leq C\delta^{1/2}$ for each $\tau>0$; let $\tau\downarrow 0$ to get $\int_{B_{1/2}\cap \{|x|<\delta\}}|v|^2$ for each $\delta>0$, from which the claim follows). Next, consider any $Y\in S(\BC^{(0)})\cap B_{1/2}$. By Lemma \ref{lemma:gaps}, we know that we may choose a sequence $Z_k \equiv (\xi_k,\eta_k)\in \spt\|V_k\|\cap B_{3/4}(0)$ with $\Theta_{V_k}(Z_k)\geq \Theta_{\BC^{(0)}}(0) = \frac{5}{2}$ with $Z_k\to Y$. Thus, for each $\rho\in (0,1/8]$, from Corollary \ref{corollary:non-concentration}(i) (in the same way as for Lemma \ref{lemma:L2_coarse}(iv)), we have for all $k$ sufficiently large (depending on $\rho$):
$$\int_{U_k\cap B_{\rho/2}(Z_k)}\frac{|u_k-\xi_k^\perp|^2}{|X + u_k(X)-Z|^{n+3/2}}\ \ext\H^n \leq C\rho^{-n-3/2}\int_{U_k\cap B_{\rho}(Z_k)}|u_k - \xi_k^\perp|^2.$$
Note that by Lemma \ref{lemma:L2_coarse}(i) that $|\xi_k| \leq CE_{k}$, and hence $E_{k}^{-1}(\xi_k,0)$ converges to some limit: let us call this limit $\kappa(Y)$ (we shall see momentarily that this is only dependent on $Y$, and not the choice of sequence $(Y_k)_k$). Then, dividing the above inequality by $E_k^2$ and passing to the limit, using the fact that $v_k\to v$ strongly in $L^2(\BC\res B_{1/2})$, we get
$$\int_{B_{\rho/2}(Y)}\frac{|v-\kappa^\perp(Y)|^2}{|X-Z|^{n+3/2}}\ \ext(\H^n\res\spt\|\BC^{(0)}\|) \leq C\rho^{-n-3/2}\int_{B_\rho(Y)}|v-\kappa^\perp(Y)|^2\ \ext(\H^n\res\spt\|\BC^{(0)}\|).$$
Note that this shows that $\kappa(Y)$ is independent of the sequence $(Z_k)_k$, as this integral on the left hand side needs to be finite, which uniquely determines $\kappa(Y)$ (using the fact that the unit normals to $\spt\|\BC^{(0)}_0\|$ span $\R^2$). Also, note that applying Lemma \ref{lemma:L2_coarse}(i) to $(\eta_{0,1/2})_\#V_k$ instead of $V_k$ shows that $\sup_{B_{5/16}\cap S(\BC^{(0)})}|\kappa|^2 \leq C\int_{B_{1/2}}|v|^2$.

Note that in terms of the properties of a proper coarse blow-up class as in Section \ref{sec:coarse_regularity}, we have now established that the class of functions $\FB(\BC)$ satisfies (after rotating the functions onto the fixed half-hyperplane $H$) $(\FB1)$, $(\FB2)$, $(\FB3)$ (except for the smoothness of $\kappa$). Note that, in a similar way to how we established $(\FB3)$ above (namely through Corollary \ref{corollary:non-concentration}, which itself came from an application of Lemma \ref{lemma:L2_coarse} with $(\eta_{Z,\rho})_\# V$ in place of $V$), we can show that, for each $\rho\in (0,1/8)$ and $k$ sufficiently large (depending on $\rho$), if we apply Lemma \ref{lemma:L2_coarse}(v) to $(\eta_{Z_k,\rho})_\#V_k$ and pass to the limit, we get
$$\int_{B_{\rho/2}(Y)}R_Y^{2-n}\left(\frac{\del}{\del R_Y}\left(\frac{v}{R_Y}\right)\right)^2\ \ext(\H^n\res \spt\|\BC^{(0)}\|) \leq C\rho^{-n-2}\int_{B_\rho(Y)}|v|^2\ \ext(\H^n\res\spt\|\BC^{(0)}\|).$$
This is almost $(\FB4)$: this will follow from the above by applying the closure properties in $(\FB5)$, noting that if $\ell$ is a homogeneous degree one function in $|x|$, then $\frac{\del}{\del |x|}(\ell/|x|)= 0$. Thus, we are left with establishing properties $(\FB5)$, $(\FB6)$ $(\FB7)$, and showing that the function $\kappa$ is smooth.

\subsection{Further Properties of the Coarse Blow-Up Class}\label{sec:coarse_properties}

In this section we shall prove that the coarse blow-up class, $\FB(\BC^{(0)})$, satisfies properties $(\FB5)$, $(\FB6)$, and that the function $\kappa$ in $(\FB3)$ is smooth; the remainder of the paper will then be devoted to proving that $(\FB7)$ is also satisfied, and using this with various other intermediary regularity results to prove Theorem \ref{thm:A}.

Fix $I\in \{0,1,2\}$, $\BC^{(0)}\in \FL_S\cap \FL_I$, and let $v\in \FB(\BC^{(0)})$; we shall always write $(V_k)_k$, $(\BC_k)_k$, $\epsilon_k$, $\delta_k$, and $\tau_k$ for sequences as in Section \ref{sec:coarse-construction} which give rise to $v$.

Let us start by showing $(\FB5\text{I})$. Let $z\in S(\BC^{(0)})\cap B_1$ and $\sigma\in (0,\frac{3}{8}(1-|z|)]$. Then set $\tilde{V}_k:= (\eta_{z,\sigma})_\#V_k$; note that we still have $\tilde{V}_k\weakly \BC^{(0)}$. Then, since $(\eta_{z,\sigma})_\#\BC_k = \BC_k$ for all $k$, it is straightforward to check that the coarse blow-up of $\tilde{V}_k$ relative to $\BC_k$ is $\tilde{v}(\cdot) = \|v(z+\sigma(\cdot))\|^{-1}_{L^2(B_1)}v(z+\sigma(\cdot))$; thus $(\FB5\text{I})$ holds.

Now let us prove $(\FB6)$. Let $(v_\ell)_\ell\subset\FB(\BC^{(0)})$ and for each $\ell$ let $(V_k^\ell)_k\subset\S_2$ and $(\BC_k^\ell)_k\subset \FL_I$ be the sequences whose coarse blow-up is $v_\ell$. Now inductively choose integers $k_\ell$ such that:
\begin{enumerate}
	\item [(a)] $k_1<k_2<\cdots$;
	\item [(b)] $E_{V^\ell_{k_\ell},\BC^\ell_{k_\ell}} < \min\{\ell^{-1},\epsilon(\BC^{(0)},\ell^{-1})\}$, where $\epsilon_0$ is the constant from Lemma \ref{lemma:coarse_graphical_rep};
	\item [(c)] $\|E^{-1}_{V_{k_\ell}^\ell,\BC^\ell_{k_\ell}}u_{\ell,k_\ell} - v_\ell\|_{L^2(B_1)} < \ell^{-1}$, where $u_{\ell,k_\ell}$ is the function from Lemma \ref{lemma:coarse_graphical_rep} with $\tau = \ell^{-1}$; such a function exists by the construction of the coarse blow-up $v_\ell$.
\end{enumerate}

Now from Lemma \ref{lemma:coarse_graphical_rep}(c), we know that for each compact subset of $B_1\cap \spt\|\BC^{(0)}\|\backslash S(\BC^{(0)})$, for all sufficiently large $\ell$ we have $\int_K |v_\ell|^2 + |Dv_\ell|^2 \leq C(n,K)$, and so if $v$ is the coarse blow-up of an appropriate subsequence $(V_{k_{\ell^\prime}}^{\ell^\prime})_{\ell^\prime}$ of $(V_{k_\ell}^\ell)_\ell$, it is then straight forward to see that, after perhaps passing to another subsequence, that on each compact subset $K$ we have $v_{\ell^\prime}\to v$ strongly in $L^2(K)$ and weakly in $W^{1,2}(K)$. Finally, as usual, from the non-concentration estimate of the $L^2$ norm about the spine, we see that this $L^2$ convergence also holds on each compact subset of $\spt\|\BC^{(0)}\|$; this proves $(\FB6)$.

Let us now prove that $\kappa$ is smooth. We already know (see Section \ref{sec:coarse-construction}) that for each $Y\in B_{1/2}\cap S(\BC^{(0)})$ we have for all $\rho\in (0,1/8)$,
\begin{equation}\label{E:kappa-0}
\int_{B_{\rho/2}(Y)}\frac{|v-\kappa^\perp(Y)|^2}{|X-Y|^{n+3/2}} \leq C\rho^{-n-3/2}\int_{B_\rho(Y)}|v-\kappa^\perp(Y)|^2
\end{equation}
where $\kappa$ is a single-valued function. Note that this inequality by itself is enough to prove that the average $v_a$ is $C^{0,\alpha}$ up-to-the-boundary by standard Campanato estimates (see \cite{minter2021}) with boundary values exactly determined by $\kappa$ (as the integral on the left hand side is finite); moreover, once we have that $\kappa$ is smooth, the above inequality will imply that $\kappa^{\perp_{H_i}} = \left. (v_i)_a\right|_{B_{1/2}\cap \del H_i}$ for $H_i$ a half-hyperplane in $\spt\|\BC^{(0)}\|$ and $v_i \equiv \left. v\right|_{H_i}$, and so standard boundary regularity theory will imply that $v_a$ is smooth up-to-the-boundary on each half-hyperplane; this means we can make sense of $v_a(0)$ and $Dv_a(0)$, which we make use of in $(\FB5\text{II})$.

We will follow the argument seen in \cite[Lemma 1]{simoncylindrical} (a variant of which can be found in \cite[Lemma 12.2]{wickstable}) to show that $\kappa$ is smooth. Let $\zeta = \zeta(r,y)$ be a smooth function such that $\zeta(r,y)\equiv 0$ on $\R^{n+1}\backslash B_{3/8}$ with $\frac{\del\zeta}{\del r} = 0$ on a neighbourhood of $\{|x|=0\}$; more precisely, such that there is $\tau_*>0$ such that $D_i\zeta = 0$ for $i=1,2$ on $|x|<2\tau_*$.

For $i=1,2$ and $p=1,\dotsc,n-1$, apply the first variation formula with test function $\Phi:= e_i\zeta_p$ and $V = V_k$ to get
$$\int_{B_{1/2}}\nabla^{V_k}x^i\cdot \nabla^{V_k}\zeta_{y_p}\ \ext\|V_k\| = \int_{B_{1/2}}e_i\cdot\nabla^{V_k}\zeta_{y_p}\ \ext\|V\|  = 0$$
where for notational simplicity we have written $\zeta_{y_p}\equiv \frac{\del\zeta}{\del y^p}$; the second integral here vanishes because $\zeta$ has compact support and because $D_{e_i}\zeta\equiv 0$ on a neighbourhood of $0$. Set $U_k:= \spt\|\BC_k\|\cap B_1\backslash\{|x|<\tau_k\}$ and note that on $B_{3/4}\cap \spt\|V_k\|\backslash\{|x|<\tau_k\}$ is expressible as a function $u^k$ as in Lemma \ref{lemma:coarse_graphical_rep}. Let us write $G_k:= \graph(\left. u^k\right|_{U_k\cap B_{1/2}})$. Now clearly we have, for any $\tau>0$, for all $k$ sufficiently large, if we write $(g^{ij}_k)_{ij}$ for the matrix representing the orthogonal projection $R^{n+1}\to T_X V^k$ (which is defined $\H^n$-a.e. point $X$ in $\spt\|V_k\|$),
\begin{align*}
	\int_{B_{1/2}\backslash G_k}\left|\nabla^{V_k}x^i\cdot\nabla^{V_k}\zeta_{y_p}\right|\ \ext\|V_k\| & = \int_{B_{1/2}\backslash G_k}\left|\sum^{n+1}_{j=3}(\delta_{ij}-g^{ij})D_{y^{j-2}}\zeta_{y_p}\right|\ \ext\|V_k\|\\
	& \leq \int_{B_{1/2}\backslash G_k}\left(\sum^{n+1}_{j=3}|e_j^\perp|^2\right)^{1/2} |D\zeta_{y_p}|\ \ext\|V_k\|\\
	& \leq \sup_{B_{1/2}}|D\zeta_{y_p}|\cdot\sqrt{\|V_k\|(B_{1/2}\backslash G_k)}\cdot\left(\int_{B_{1/2}}\sum^{n+1}_{j=3}|e_j^\perp|^2\ \ext\|V_k\|\right)^{1/2}\\
	& \leq C\sup_{B_{1/2}}|D^2\zeta|\cdot\tau^{1/2}\cdot E_{V_k,\BC_k};
\end{align*}
here, in the first equality we have used the fact that, if $k$ is sufficiently large (depending on $\tau$), we have $B_{1/2}\backslash G_k\subset B_{1/2}\cap\{|x|<\tau\}$ and thus $D_{x_i}\zeta = 0$ on $B_{1/2}\backslash G_k$ for $i=1,2$ (if $\tau<\tau_*$), so these terms do not appear in the sum in the integrand; in the first inequality we have simply used the Cauchy--Schwarz inequality for vectors in $\R^{n-1}$; in the second inequality we have used the Cauchy--Schwarz inequality on $L^2$; for the final inequality we have used Lemma \ref{lemma:L2_coarse}(ii) to bound the integral and then we have used that, for fixed $\tau>0$, since $V_k\weakly \BC^{(0)}$ by Lemma \ref{lemma:convergence} we have $\|V_k\|(B_{1/2}\cap\{|x|<\tau\}) \to \|\BC^{(0)}\|(B_{1/2}\cap\{|x|<\tau\}) = C\tau$ (as the cross-section is one-dimensional and consists of rays of length $<\tau$ on this set) for some $C = C(n)$, and $B_{1/2}\backslash G_k\subset B_{1/2}\cap\{|x|<\tau\}$ for all $k$ sufficiently large. Therefore we conclude, for each $\tau\in (0,1/2)$ and all $k$ sufficiently large (depending on $\tau$):
$$\int_{B_{1/2}\backslash G_k}\left|\nabla^{V_k}x^i\cdot\nabla^{V_k}\zeta_{y_p}\right|\ \ext\|V_k\| \leq C\sup|D^2\zeta|E_k\cdot\tau^{1/2}.$$

Now write $(\w_k^j)_{j=1}^{5-I}$ for the unit vectors in the direction of the rays making up the cross-section of $\BC_k$. Set $U_k(\tau):= \spt\|\BC_k\|\cap B_1\backslash\{|x|<\tau\}$ and write $U^j_k(\tau)$ for the intersection of $U_k(\tau)$ with the half-hyperplane in $\BC_k$ whose ray in the cross-section is in the direction $\w_k^{j}$. Also, write $G_k(\tau):= \graph(\left.u_k\right|_{U_k(\tau)\cap B_{1/2}})$ and $G^j_k(\tau):= \graph(\left. u_k\right|_{U_k^j(\tau)\cap B_{1/2}})$ (this could be determined by a single-valued or two-valued function).

We are left with estimating
$$\int_{B_{1/2}\cap G_k}\nabla^{V_k}x^i\cdot\nabla^{V_k}\zeta_{y_p}\ \ext\|V_k\| \equiv \sum^{5-I}_{j=1}\int_{G^j_k(\tau)}\nabla^{V_k}x^i\cdot\nabla^{V_k}\zeta_{y_p}\ \ext\|V_k\|.$$
Rotate so that $\w_k^1 = e_1$ (we will see that the expression we find is invariant under rotations, and so we can rotate back at the end). Let us begin with the case $i=1$, i.e. when our deformation is in the direction parallel to $\w_k^1$. Then we have, for $\H^n$-a.e. point in $G^1_k(\tau)$,
$$\nabla^{V_k}x^1\cdot\nabla^{V_k}\zeta_{y_p} = h^{11}_k\frac{\del\zeta_{y_p}}{\del x^1} + \sum^{n+1}_{j=3}h^{1,j}_{k}\cdot\frac{\del\zeta_{y_p}}{\del y^{j-2}}$$
where $(h^{ij}_k)_{ij}$ is the inverse of the matrix Jacobian matrix for the graph of $u^1_k$, i.e. the component of $u_k$ over the half-hyperplane determined by $\w_k^1$; note that if $u^1_k$ is two-valued, then the above expression is taken to be for the corresponding component of $u^1_k$ if the point is a non-branch point and the (unique) value at a branch point. Thus, using the area formula as in (\ref{eqn:area_formula}) we see that
\begin{align*}
	\int_{G_k^1(\tau)}\nabla^{V_k}&x^1\cdot\nabla^{V_k}\zeta_{y_p}\ \ext\|V_k\|\\
	& = \int^1_\tau\int_{\R^{n-1}}\sum_\ell\left((h^\ell)_k^{11}\frac{\del}{\del x^1}\zeta_{y_p}\left(\sqrt{|x^1|^2 + |(u_k^1)^\ell|^2},\, y\right)\right)\\
	& \hspace{10em} + \left.\sum^{n+1}_{j=3}(h^\ell)^{1j}_k\frac{\del}{\del y^{j-2}}\zeta_{y_p}\left(\sqrt{|x^1|^2+|(u_k^1)^\ell|^2},\, y\right)\right)\sqrt{h^\ell_k}\ \ext y\ext x^1
\end{align*}
where here $h^\ell_k$ is the determinant of the corresponding $(h^{ij}_k)_{ij}$; the sum in $\ell$ is over the values of $u^1_k$, so if $u^1_k$ is single-valued there is just one term and if it is two-valued there are two. We may then estimate, just as in \cite{simoncylindrical}, using in place of usual quasilinear elliptic estimates the estimates in Section \ref{sec:two-valued_stationary_graphs} whenever $u^1_k$ is two-valued, to get for all $k$ sufficiently large (depending on $\tau$)
$$\left|\int_{G_k^1(\tau)}\nabla^{V_k}x^1\cdot\nabla^{V_k}\zeta_{y_p}\ \ext\|V_k\|\right| \leq C\left(\sup|D\zeta| + \sup|D^2\zeta|\right)\cdot\int_{U^1_k(\tau/2)}|u^1_k|^2.$$

Now for $i=2$, i.e. deformations perpendicular to $\w_k^1$ in the cross-section, we see in a similar fashion that (see \cite[Lemma 1, (26)]{simoncylindrical})
$$\int_{G_k^1(\tau)}\nabla^{V_k}x^2\cdot\nabla^{V_k}\zeta_{y_p}\ \ext\|V_k\| = \int_{U_k^1(\tau)} \sum_\ell \nabla (u^1_k)^\ell\cdot\nabla\zeta_{y_p} + \beta_k$$
where $\beta_k$ is a term with $E_k^{-1}\beta_k \to 0$ as $k\to\infty$. All of these expressions are invariant under rotations, and so they hold without the assumption that $\w_k^1 = e_1$, and so these hold for each $G^j_k(\tau)$. So summing these results over $i=1,2$ and over $j=1,\dotsc,5-I$ we get
$$\sum^{5-I}_{j=1}\int_{U^j_k(\tau)}\theta_j\nabla (u^j_k)_{a}\cdot\nabla\zeta_{y_p} = R_k + S E_{k}$$
where $|S|\leq C\tau^{1/2}$ and $E_{k}^{-1}R_k\to 0$ as $k\to \infty$; here, we have used that $\sum_\ell \nabla(u^1_k)^\ell = \theta_k (u^j_k)_a$, for $\theta_j = \Theta_{\BC_k}(\w_k^j)$ is the multiplicity of the respective half-hyperplane in $\BC_k$ (which is independent of $k$), i.e. $\theta_j$ is $2$ whenever $u^j_k$ is two-valued, and is $1$ otherwise. We stress here that $(u^j_k)_a$ is the \textit{average} part of $u^j_k$, and so the index $a$ does not represent a derivative.

If we divide this expression by $E_k$, for fixed $\tau\in (0,\tau_*)$ we may let $k\to\infty$ to obtain
\begin{equation}\label{E:kappa-1}
\int_{\spt\|\BC^{(0)}\|\cap \{|x|>\tau\}}\Theta_{\BC^{(0)}}(z)\nabla v_a(z)\cdot\nabla\zeta_{y_p}(z) = S
\end{equation}
where $|S|\leq C\tau^{1/2}$; so letting $\tau\downarrow 0$, we get
$$\int_{\spt\|\BC^{(0)}\|}\Theta_{\BC^{(0)}}(z)\nabla v_a(z)\cdot \nabla\zeta_{y_p}(z) = 0;$$
note here we have used the fact that, by Lemma \ref{lemma:L2_coarse}(ii), we have that $\int_{B_{1/2}}|D_yv^k|^2 \leq C$ for all $k$ sufficiently large, where $C = C(n)$, and so $D_{y_q}v\in L^2(\BC^{(0)}\res B_{1/2})$ for all $q\in \{1,2,\dotsc,n-1\}$ (i.e. all the directions parallel to the spine) and that since $D_{x_i} \zeta \equiv 0$ for $i=1,2$, when taking the dot-product in (\ref{E:kappa-1}), on the region $\{|x|<\tau_*\}$ the only derivatives of $v_a$ which occur are those parallel to the spine, and on the region $\{|x|\geq\tau_*\}$ we have $C^{1,1/2}$ convergence, and thus we may pass to the limit. Thus, we have
$$\int_{B_{3/4}}\nabla v_{a}\cdot \nabla\zeta_{y_p}\ \ext\|\BC^{(0)}\| = 0$$
for any $\zeta\in C^\infty_c(B_{1/4})$ with $D_{x_i}\zeta \equiv 0$ on a neighbourhood of $\{r=0\}$, and for any $p\in \{1,\dotsc,n-1\}$. We stress that we know $v_a$ is \text{always} a smooth single-valued harmonic function, away from $S(\BC^{(0)})$. Now, if we perform a rotation on each half-hyperplane in $\BC^{(0)}$, rotating them to the fixed hyperplane $H = \{(r,y)\in \R^n: r>0\}$, changing the domain of integration of the integral to this fixed hyperplane, and then integrate by parts, using the fact that $\BC^{(0)}\in \FL_S$ is stationary and so the sum of the unit normals over half-hyperplanes in $\spt\|\BC^{(0)}\|$, counted with multiplicity, vanishes, the above expression can be written as
$$\int_H \tilde{v}\cdot\nabla\zeta_{y_p} = 0$$
where $\tilde{v}(r,y):= \sum^{5-I}_{j=1}\Theta_{\BC^{(0)}}(\w_j)\cdot v_a(r\w_j,y)\w_j^\perp$, for $(\w_j)^{5-I}_{j=1}$ the unit vectors in the direction of the rays in $\BC^{(0)}_0$; we stress here that this function is vector-valued, as we have only changed the domain of integration. Since $\zeta$ as above is arbitrary, we can now follow the argument leading to \cite[Lemma 1, (28)]{simoncylindrical} to see that, if $\hat{v}$ denotes the even reflection of $\tilde{v}$ over $\del H$ to a function on all of $\R^n$, then $D_{y^p}\hat{v}$ is a smooth harmonic function on all of $\R^n$, with the desired estimates on its values and derivatives at $0$ (which follow from standard estimates for harmonic functions). In particular, the function $\left.\hat{v}\right|_{\del H}$, which are the boundary values of the original function $v$, is smooth on all of $B_{1/2}\cap\del H$ with the same estimates; indeed, we have for $Y\in B_{1/2}\cap \del H$ that $\hat{v}(Y) = \sum^{5-I}_{j=1}\theta_j(\kappa\cdot \w_j)^\perp \w_j^\perp$ (note that, by the remark after (\ref{E:kappa-0}), we know that the boundary values of $v_a$ agree with $\kappa$ as $v_a$ is continuous up-to-the-boundary), which is a smooth function and has supremum bounded by $C\int_{B_{1/2}}|v|^2$. Hence, we see that, on $\del H$, the function $F:= \sum_{j=1}^{5-I}\theta_j(\kappa\cdot\w_j^\perp)\w_j^\perp$, where $\theta_j\in \{1,2\}$, is smooth. However, since the normal directions of the cross-section, i.e. $\w_j^{\perp}$, span $\R^2$, this sum $F$ being smooth readily implied that $\kappa$ itself must be a smooth function; hence we have shown that $\kappa$ is smooth, proving $(\FB3)$ in its totality.

Finally, now that $Dv_a(0)$ makes sense, let us prove $(\FB5\text{II})$. Note that from $(\FB3)$ we also know that it is a single point, $\kappa(0)\in \R^2$, which determines all the values of $v_a(0)$ along each half-hyperplane via the projections of $(\kappa(0),0)$ to each half-hyperplane. So, note that for any $y\in \R^{n+1}$ with $v-y\not\equiv 0$ in $B_1$ and each $\sigma\in (0,1)$ ( the introduction of the parameter $\sigma\in (0,1)$ is a necessary step to ensure that we can still apply our graphical representation over suitable balls), the function $\|v(\sigma(\cdot))-y\|_{L^2(B_1)}^{-1}(v(\sigma(\cdot))-y)$ is the coarse blow-up of the sequence $(\eta_{\sigma E_k y,\, \sigma})_\#V_k$ relative to $\BC_k$, and so belongs to $\FB(\BC^{(0)})$; note that by $v-y$ in this setting we mean the function which over a half-hyperplane $H_i$ in $\BC^{(0)}$ is given by $v^i - y^{\perp_{H_i}}$, for $y^{\perp_{H_i}}$ the projection onto the orthogonal complement of $H_i$. Then by applying $(\FB6)$, with $\sigma\uparrow 1$, we get that $\|v-y\|_{L^2(B_1)}^{-1}(v-y)\in \FB(\BC^{(0)})$. If we apply this with $y = (\kappa(0),0)$, we see that this does subtract the relevant value of $v^i_a(0)$ from $v^i$.

In order to remove a linear function, we need to perform two rotations: the first will rotate the sequence $V_k$ to remove the derivatives of the function which are parallel to the spine (note that from $(\FB3)$ we once again know that for each component of $v$ these derivatives are determined by the same function, appropriately projected), and then second we rotate the individual half-hyperplanes in the $\BC_k$ to remove the derivatives in the direction of the corresponding ray of each half-hyperplane. (Note this should be compared with \cite[Section 4.2, ($\FB5\text{II}$)]{beckerkahn}.)

So consider a function of the form $\kappa^{\perp_{T_{(\cdot)}\BC^{(0)}}} + \psi$, where $\kappa\in S(\BC^{(0)})^\perp$ and $\psi:\spt\|\BC^{(0)}\|\cap\{|x|>0\}\to \spt\|\BC^{(0)}\|^\perp$ a linear function on each half-hyperplane, i.e. $\psi$ is of the form $\psi(X) = \sum^{n-1}_{j=1}y^jc_j^{\perp_{T_X\BC^{(0)}}} + |x|\phi(x/|x|)$ for some collection of vectors $c_1,\dotsc,c_{n-1}\in S(\BC^{(0)})^\perp$ and function $\phi:\{\w_1,\dotsc,\w_5\}\to \spt\|\BC^{(0)}\|^\perp$, where $\w_i$ is the unit vector in the direction of the ray determining the $i^{\text{th}}$ half-hyperplane in $\BC^{(0)}$; this is the type of function we want to subtract off from $v$. By the same argument as in Section \cite[Section 2]{simoncylindrical}, there is a sequence of cones of the form $\mathbf{D}_k:= R^k_\#\tilde{\mathbf{D}}_k$, where $\tilde{\mathbf{D}}_k\in \FL_I$ and $R^k$ is a rotation of $\R^{n+1}$ such that $|R^k-\id|\to 0$. Now let $d_k$ be the function which represents $\tilde{\mathbf{D}}_k$ as a graph over $\BC^{(0)}$, and then let $\tilde{\BC}_k\in \FL_I$ be the cone determined by the function $c_k + E_k d_k$, where here $c_k$ is the function which represents $\BC_k$ over $\BC^{(0)}$. Then, considering $\tilde{V}_k:= (\tau_{E_k \kappa} \circ (R^k)^{-1})_\#V_k$, we see that the coarse blow-up of $\tilde{V}_k$ relative to $\tilde{\BC}_k$ is exactly $\|v-\kappa^{\perp_{T_{(\cdot)}\BC^{(0)}}} - \psi\|_{L^2}^{-1}(v- \kappa^{\perp_{T_{(\cdot)}\BC^{(0)}}} -\psi)$, as desired. Taking this with the $\kappa$ and $\psi$ determined by $v_a(0)$, we get $(\FB5\text{II})$.

\section{Level 0: Proof of Main Theorem}\label{sec:L0excess}

Here we now prove Theorem \ref{thm:A} when the base cone is level 0, i.e. $\BC^{(0)}\in \FL_S\cap \FL_0$. We are able to do this now because our coarse blow-ups consist only of single-valued harmonic functions, and thus we do not require any additional properties other than $(\FB3)$ to deduce the desired boundary regularity.

Our the technical lemma toward proving Theorem \ref{thm:A} is the following excess decay statement.

\begin{lemma}\label{lemma:ed-level-0}
	Suppose $\BC^{(0)} = \BC^{(0)}_0\times\R^{n-1}\in \FL_S\cap \FL_0$, and fix $\theta \in (0,1)$. Then there exists $\epsilon = \epsilon(\BC^{(0)},\theta)\in (0,1)$ such that the following is true: if $V\in \CN_{\epsilon}(\BC^{(0)})$ and $\BC\in \FL_\epsilon(\BC^{(0)})$, then there exists $\tilde{\BC}$ and an orthogonal rotation $\Gamma$ of $\R^{n+1}$ such that:
	\begin{enumerate}
		\item [(i)] $|\Gamma-\id|\leq CE_{V,\BC}$;
		\item [(ii)] $\dist(\spt\|\tilde{\BC}\|\cap B_1,\spt\|\BC\|\cap B_1) \leq CE_{V,\BC}$;
		\item [(iii)] 
		$$\theta^{-n-2}\int_{B_\theta} \dist^2(X,\spt\|\Gamma_\#\tilde{\BC}\|)\ \ext\|V\|(X) \leq C\theta E_{V,\BC}^2;$$
	\end{enumerate}
	here, $C = C(n)$.
\end{lemma}

\begin{proof}
	We will prove this by contradiction; so suppose that the lemma does not hold (for $C = C(n)$ to be chosen): therefore, we may find sequences $\epsilon_k\to 0$, $V_k\in \CN_{\epsilon_k}(\BC^{(0)})$, $\BC_k\in \FL_{\epsilon_k}(\BC^{(0)})$, such that there lemma does not hold true for this choice of $\theta$ and $\BC^{(0)}$; it suffices to show for infinitely many $k$, the lemma does hold.
	
	For $i=1,2,\dotsc,n-1$, set $Y_i:= \frac{\theta}{2}e_{2+i}\in S(\BC^{(0)})$. Lemma \ref{lemma:gaps} tells us that for each $k\geq 1$ and $i\in \{1,\dotsc,n-1\}$ we may find sequences $Z_{i,k}\in \spt\|V_k\|\cap B_1$ such that $\Theta_{V_k}(Z_{i,k})\geq \frac{5}{2}$ and $Z_{i,k}\to Y_i$ as $k\to\infty$. Since $\{Y_1,\dotsc,Y_{n-1}\}$ span an $(n-1)$-dimensional subspace, it must be the case that, for all $k$ sufficiently large, the $\{Z_{1,k},\dotsc,Z_{n-1,k}\}$ also span an $(n-1)$-dimensional subspace of $\R^{n+1}$; call this subspace $\Sigma_k$. We may then choose for each $k$ sufficiently large a rotation $\Gamma_k$ of $\R^{n+1}$ such that $\Gamma_k(\Sigma_k) = S(\BC^{(0)})$ and such that $|\Gamma_k - \id|$ is minimal across all rotations $\Gamma$ which obey $\Gamma(\Sigma)= S(\BC^{(0)})$. From Lemma \ref{lemma:L2_coarse}(i), we know that $\dist(Z_{i,k},S(\BC^{(0)}))\leq CE_k$ for each $i$ and all $k$ sufficiently large, where $C = C(n)$; here we have written $E_k\equiv E_{V_k,\BC_k}$ in the usual way. Thus, this shows that for all $k$ sufficiently large,
	$$|\Gamma_k-\id|\leq CE_k$$
	where $C = C(n)$. Now, setting $\tilde{V}_k:= (\Gamma_k)_\#V_k$, note that by the triangle inequality,
	\begin{align*}
		E_{\tilde{V}_k,\BC_k}^2 & := \int_{B_1}\dist^2(X,\spt\|\BC_k\|)\ \ext\|(\Gamma_k)_\#V_k\|\\
		& \leq 2\int_{B_1}\dist^2(X,\spt\|(\Gamma_k)_\#\BC_k\|)\ \ext\|(\Gamma_k)_\#V_k\| + \left(5\w_n+2\right)\dist^2(\spt\|(\Gamma_k)_\#\BC_k\|\cap B_1,\spt\|\BC_k\|\cap B_1)\\
		& \leq 2\int_{B_1}\dist^2(X,\spt\|\BC_k\|)\ \ext\|V_k\| + CE^2_{k}\\
		& = CE_k^2
	\end{align*}
	for some $C = C(n)$. Thus, if $k$ is sufficiently large we can apply Lemma \ref{lemma:coarse_graphical_rep} and the analysis of Section \ref{sec:coarse_construction} to $\tilde{V}_k$. So, let $\tilde{v}\in \FB(\BC^{(0)})$ be the coarse blow-up of (a subsequence of) $\tilde{V}_k$ relative to $\BC_k$. By construction we have $\tilde{v}(Y_i) = 0$ for $I=1,\dotsc,n-1$ and also $\tilde{v}(0) = 0$ (which we can also arrange by construction, by initially translating each $V_k$ by some $Z_k$ with $\Theta_{V_k}(Z_k)\geq\frac{5}{2}$; by Lemma \ref{lemma:L2_coarse}(i) we can also estimate the height excess of this translated $V_k$ in terms of the original $V_k$). Write $\tilde{v}^1,\dotsc,\tilde{v}^5$ for the $5$ (single-valued) $C^2$ functions which determine $v$. The above tells us, as each $v^i$ is $C^2$ (in fact smooth) up-to-the-boundary by Theorem \ref{thm:coarse_reg}, that for each $\ell=1,\dotsc,5$ and $i=1,\dotsc,n-1$, we may find $S_{\ell,i}\in B_{\theta/2}\cap S(\BC^{(0)})$ such that $\frac{\del \tilde{v}^\ell}{\del y_i}(S_{\ell,i}) = 0$. Given this, the $C^{0,1/2}$ regularity of $D_y \tilde{v}$ provided by Theorem \ref{thm:coarse_reg} (in fact we can get improved estimates in this case as $\tilde{v}$ is smooth) gives that
	$$|D_y \tilde{v}^\ell(0)|^2 \leq C\theta\int_{B_{1/2}}|\tilde{v}|^2.$$
	We can now use $Dv(0)$ to define a new sequence of cones, $\tilde{\BC}_k$, for which the excess improves, as follows. Let $\{H_1,\dotsc,H_5\}$ denote the half-hyperplanes in $\spt\|\BC^{(0)}\|$, with unit vectors $\w_1,\dotsc,\w_5$ in $\R^2$ parallel respective ray in the cross-section of $\BC^{(0)}$ determined by $H_1,\dotsc,H_5$, respectively. Thus, a function over $H_i$ is expressible as $(r,y)\mapsto v^i(r\w_i,y)\w_i^\perp$. Now consider for each $i=1,\dotsc,5$, the half-hyperplane $P_i:H_i\to H_i^\perp$ determined by the graph of $P_i(r,y):= rD_r\tilde{v}^i(0) + y\cdot D_y\tilde{v}^i(0)$, and the half-hyperplane with axis that of $S(\BC^{(0)})$, namely $c_i:H_i\to H_i^\perp$ given by $c_i(r,y):= rD_r\tilde{v}^i(0)$. The above estimates clearly give
	$$|P_i(r,y)-c_i(r,y)|^2 \leq C\theta\cdot|y|^2\int_{B_{1/2}}|\tilde{v}|^2$$
	and moreover from the regularity of $\tilde{v}$ provided by Theorem \ref{thm:coarse_reg} we know
	$$\theta^{-n-2}\int_{H_i\cap B_\theta}|v^i-P_i|^2 \leq C\theta\int_{B_{1/2}}|\tilde{v}|^2.$$
	These two inequalities clearly give that
	$$\theta^{-n-2}\int_{B_\theta}|v-c|^2 \leq C\theta\int_{B_{1/2}}|\tilde{v}|^2$$
	where $c:\spt\|\BC^{(0)}\|\to \spt\|\BC^{(0)}\|^\perp$ is the function determined by $H_i$ by $c_i$; notice that by construction all the $c_i$ have the same boundary values, and so $\graph(c)$ determines another level $0$ cone. Indeed, $c$ is exactly what we use to determine the new sequence of cones; suppose $\BC_k$ has half-hyperplanes $H^k_1,\dotsc,H^k_5$, with $H^k_i\to H_i$; then, $H^k_i$ is determined as a graph over $H_i$ of some function, say $g^k_i$. We then set $\tilde{g}^k_i:= g^k_i + E_kc_i^k$; the graph of $\tilde{g}^k$ determines a new sequence of level 0 cones, which we denote by $\tilde{\BC}_k$. It is clear that if we write $\tilde{u}_k$, $u_k$ for the graphs provided by Lemma \ref{lemma:coarse_graphical_rep} of $\tilde{V}_k$ over $\BC_k$ and $\tilde{\BC}_k$ respectively, then
	$$u_k(X+\tilde{g}^k(X)) = \tilde{u}_k(X+\tilde{g}^k(X)) -  cE_{k} + \beta_k$$
	where $E_k^{-1}\beta_k\to 0$ as $k\to\infty$. Thus we see
	$$\frac{1}{E_k^2}\cdot\theta^{-n-2}\int_{B_\theta}\dist^2(X,\spt\|\tilde{\BC}_k\|)\ \ext\|\tilde{V}_k\| \to \theta^{-n-2}\int_{B_\theta}|v-c|^2 \leq C\theta$$
	which implies that the conclusions of the lemma hold for all $k$ sufficiently large, providing the contradiction.
\end{proof}

Now we use Lemma \ref{lemma:ed-level-0} to prove Theorem \ref{thm:A} when $\BC^{(0)}\in \FL_S\cap \FL_0$ is level 0.

\begin{theorem}\label{thm:mainL0}
	Theorem \ref{thm:A} is true whenever $\BC^{(0)}\in \FL_S\cap \FL_0$.
\end{theorem}

\begin{proof}
	Fix $\BC^{(0)}\in \FL_S\cap \FL_0$; without loss of generality we can rotate so that $\BC^{(0)} = \BC^{(0)}_0\times\R^{n-1}$. Choose $\theta = \theta(n)$ so that $C\theta< 1/4$, where $C = C(n)$ is the constant from Lemma \ref{lemma:ed-level-0}. Thus we know that there is a $\epsilon_0 = \epsilon_0(\BC^{(0)})$ for which Lemma \ref{lemma:ed-level-0} holds.
	
	So let $\epsilon \in (0,\epsilon_0)$ to be a constant eventually chosen to only depend on $\BC^{(0)}$. Suppose that $V\in \CN_{\epsilon}(\BC^{(0)})$. We claim that if $\epsilon = \epsilon(\BC^{(0)})$ is small enough we can iterate Lemma \ref{lemma:ed-level-0} to prove a sequence of rotations $(\Gamma^j)_j$ and cones $\BC^i\cap \FL_0$ such that for each $j$,
	\begin{enumerate}
		\item [(i)] $|\Gamma^j-\Gamma^{j-1}|^2\leq C_*4^{-j}E^2_{V,\BC^{(0)}}$;
		\item [(ii)] $\dist^2(\spt\|\BC^j\|\cap B_1,\spt\|\BC^{j-1}\|\cap B_1) \leq C_*4^{-j}E^2_{V,\BC^{(0)}}$;
		\item [(iii)]
		$$(\theta^j)^{-n-2}\int_{B_{\theta^j}}\dist^2(X,\spt\|(\Gamma^j)_\#\BC^j\|)\ \ext\|V\|(X) \leq 4^{-j}E^2_{V,\BC^{(0)}}$$
	\end{enumerate}
	for suitable $C_* = C_*(n)$; here, $\Gamma^0 = \id$ and $\BC^0 = \BC^{(0)}$. We shall prove this by induction on $j$. The $j=1$ case follows immediately from Lemma \ref{lemma:ed-level-0}, applying it with $V$ and $\BC^{(0)}$ in place of $\BC$. Now suppose that (i) -- (iii) hold for $j=1,\dotsc,k$; we shall construct $\Gamma^{k+1}$ and $\BC^{k+1}$ to prove (i) -- (iii) hold for $j=k+1$. To do this, we shall show that we can apply Lemma \ref{lemma:ed-level-0} with $\BC^k$ and $V^k := (\eta_{0,\theta^k}\circ(\Gamma^k)^{-1})_\# V$ in place of $\BC$ and $V$, respectively. In order to do this, we need to choose $\epsilon$ is sufficiently small, independent of $k$, to ensure that $V^{k+1}\in \CN_{\epsilon_0}(\BC^{(0)})$ and $\BC^{k+1}\in \CN_{\epsilon_0}(\BC^{(0)})$ whenever (i) -- (iii) hold for $j=1,\dotsc,k$.
	
	Note first that from (ii), we have from the triangle inequality and the form of $\BC^{(0)}$,
	$$\dist(\spt\|\BC^k\|\cap B_1,\spt\|\BC^{(0)}\|\cap B_1) \leq C_*E_{V,\BC^{(0)}}\sum^k_{i=1}2^{-i} \leq C_*E_{V,\BC^{(0)}}$$
	and so if $C_*\epsilon<\epsilon_0$, then $E_{V,\BC^{(0)}}<\epsilon$ will ensure that $\BC^k\in \FL_{\epsilon_0}(\BC^{(0)})$ for all $k$.
	
	To show $V^{k}\in \CN_{\epsilon_0}(\BC^{(0)})$, from (i) for $j=1,\dotsc,k$ it is clear again by the triangle inequality that
	$$|\Gamma_k - \id| \leq C_*Q_{V,\BC^{(0)}}\sum^k_{i=1}2^{-i} \leq C_*Q_{V,\BC^{(0)}}.$$
	Moreover, (iii) gives
	$$E_{V^k,\BC^k}^2 \leq 4^{-k}E_{V,\BC^{(0)}}^2$$
	which then gives by the triangle inequality,
	\begin{align*}
	\int_{B_1}\dist^2(X,\spt\|\BC^{(0)}\|)\ \ext\|V^k\| & \leq 2\int_{B_1}\dist^2(X,\spt\|\BC^k\|)\ \ext\|V^k\| + 6\w_n\dist^2(\spt\|\BC^k\|\cap B_1,\spt\|\BC^{(0)}\|\cap B_1)\\
	& \leq 2E_{V_k,\BC_k}^2 + 6\w_n C_*^2E_{V,\BC^{(0)}}^2\\
	& \leq (4^{1-k} + 6\w_n C_*^2)E_{V,\BC^{(0)}}^2
	\end{align*}
	Thus, if $\epsilon = \epsilon(\BC^{(0)})$ is suitably such that $(4+6\w_n C_*^2)\epsilon<\epsilon_0$, we have $V^k\in \CN_{\epsilon_0}(\BC^{(0)})$, and thus we may apply Lemma \ref{lemma:ed-level-0} to $V^k$ and $\BC^k$ to produce a rotation $\Gamma$ and a cone $\BC^{k+1}$ such that:
	\begin{enumerate}
		\item [(1)] $|\Gamma-\id|^2\leq CE^2_{V^k,\BC^k}$;
		\item [(2)] $\dist^2(\spt\|\BC^{k+1}\|\cap B_1, \spt\|\BC^k\|\cap B_1)\leq CE^2_{V^k,\BC^k}$;
		\item [(3)] $$\theta^{-n-2}\int_{B_\theta}\dist^2(X,\spt\|\Gamma_\#\BC^{k+1}\|)\ \ext\|V^k\| \leq C\theta E_{V^k,\BC^k}^2;$$
	\end{enumerate}
	once again, $C = C(n)$ is just the constant from Lemma \ref{lemma:ed-level-0}. Now set $\Gamma_{k+1} := \Gamma\circ\Gamma_k$; since $E^2_{V^k,\BC^k} \leq 4^{-k}E^2_{V,\BC^{(0)}}$ and $C\theta < 1/4$, we see that as long as $C_*>4C$, we get that (i) -- (iii) hold again, with the \textit{same} $C_*$; thus we have reset the constant at each stage, and so our choice of $\epsilon$ does allow us to inductively prove (i) -- (iii) hold for all $j$. Fix $\epsilon_1 = \epsilon_1(\BC^{(0)})$ such that all the above holds when $V\in \CN_{\epsilon_1}(\BC^{(0)})$.
	
	Now for any $Z\in \sing(V)\cap B_{1/2}$ with $\Theta_V(Z)\geq \frac{5}{2}$, if $V_Z := (\eta_{Z,1/2})_\# V$ note that, using Lemma \ref{lemma:L2_coarse}(i) and (\ref{E:L2-9}), that
	$$(1/2)^{-n-2}\int_{B_{1/2}}\dist^2(X,\spt\|(\tau_Z)_\#\BC\|)\ \ext\|V\| \leq 2^{n+2}(1+6C\w_n)E_{V,\BC}^2$$
	where $C = C(n)$ is from Lemma \ref{lemma:L2_coarse}; thus if $\epsilon$ is sufficiently small, we see that $V_Z\in \CN_{\epsilon_1}(\BC^{(0)})$ holds for each such $Z$, and so all the above arguments hold for $V_Z$. Thus for such $Z$ we can deduce that there is a sequence of rotations $\Gamma_Z^j$ of $\R^{n+1}$ and a sequence of cones $\BC^j_Z\in \FL_0$ satisfying (i) -- (iii) above, and thus we can find limits $\Gamma^j_Z\to \Gamma_Z$, $\BC^j_Z\weakly \BC_Z\in \FL_0$, such that
	\begin{enumerate}
		\item [(I)] $|\Gamma_Z-\id|\leq CE_{V_Z,\BC^{(0)}}$;
		\item [(II)] $\dist(\spt\|\BC_Z\|\cap B_1, \spt\|\BC^{(0)}\|\cap B_1) \leq CE_{V_Z,\BC^{(0)}}$ (i.e. $\BC_Z\in \FL_{CE_{V_Z,\BC^{(0)}}}(\BC^{(0)})$);
		\item [(III)] there exists $\alpha = \alpha(n) \in (0,1/2)$ such that for each $\rho\in (0,\theta)$,
		$$\rho^{-n-2}\int_{B_\rho}\dist^2(X,\spt\|(\Gamma_Z)_\#\BC_Z\|)\ \ext\|V_Z\| \leq C\rho^{2\alpha}E_{V_Z,\BC^{(0)}}^2.$$
	\end{enumerate}
	Indeed, clearly (i), (ii) prove that the sequences $\Gamma_Z^j$ and $\BC_Z^j$ are Cauchy sequences and thus converge, and (III) follows by a standard interpolation argument based on (iii); we can even specify $\alpha := \frac{1}{2}\log_{\theta^{-1}}(2)$. In particular, (II) shows that for $\epsilon$ sufficiently small we must have that $\BC_Z$ is level 0, and (III) shows that $(\Gamma_Z)_\#\BC_Z$ is the unique tangent cone to $V$ at $Z$, and thus we see that every singular point of density $\geq \frac{5}{2}$ has a unique tangent cone which is in $\FL_0\cap \FL_S$, and so in particular must have density exactly $\frac{5}{2}$. Thus for $\epsilon_2 = \epsilon_2(\BC^{(0)})$ sufficiently small, $\{\Theta_V > 5/2\}\cap B_{1/2} = \emptyset$.
	
	We now claim that for each $y\in B^{n-1}_{1/2}(0)$ there is at most one singular point of density $\frac{5}{2}$ in the slice $(\R^2\times\{y\})\cap B_{1/2}$; in fact we will be able to show that if such a slice has one point of density $\frac{5}{2}$, all other points must have density at most $2$. We argue this by contradiction. Suppose that we have $Z_1,Z_2\in \sing(V)\cap B_{1/2}\cap (\R^2\times\{y\})$ with $\Theta_V(Z_1) = \frac{5}{2}$ and $\Theta_V(Z_2)>2$; set $\sigma = |Z_1-Z_2|$. From Lemma \ref{lemma:coarse_graphical_rep} we may assume that $|Z_1|, |Z_2| < \theta/4$; thus $\sigma < \theta/2$. Hence we may apply (III) above at $Z_1$ to get
	$$(2\sigma)^{-n-2}\int_{B_{2\sigma}}\dist^2(X,\spt\|(\Gamma_{Z_1})_\#\BC_{Z_1}\|)\ \ext\|V_{Z_1}\| \leq C(2\sigma)^{2\alpha}E_{V_{Z_1},\BC^{(0)}}^2.$$
	Using this with (I) and (II), we then get
	$$E_{(\eta_{0,2\sigma})_\# V_{Z_1},\BC^{(0)}} \leq CE_{V_{Z_1},\BC^{(0)}}$$
	where $C = C(n)$ is independent of $\sigma$. Thus, if $\epsilon$ is sufficiently small, only depending on $\BC^{(0)}$, this implies that we may apply Lemma \ref{lemma:coarse_graphical_rep} to $(\eta_{0,2\sigma})_\# V_{Z_1}$, expressing it as a sum of single-valued and two-valued stationary graphs over the region $\{|x|>1/4\}$; but this is a contradiction, as by assumption we know it has a singularity (determined by $Z_2$) of density $>2$ on $\{|x|=1/2\}\cap B_{1}$. Thus we see that whenever we have $\sing(V)\cap \{\Theta_V\geq \frac{5}{2}\}\cap (\R^2\times\{y\})\neq\emptyset$, then in fact there is a unique $Z$ (depending on $y$) such that
	\begin{equation}\label{E:thm-A-level-0-1}
	\sing(V)\cap\{\Theta_V>2\}\cap (\R^2\times\{y\})\cap B_{1/2} = \{Z\}.
	\end{equation}
	But we know from (the proof of) Lemma \ref{lemma:gaps} that in fact for $\epsilon = \epsilon(n)$ sufficiently small, there is at least one point of density $\geq\frac{5}{2}$ in each slice $(\R^2\times\{y\})\cap B_{1/2}$; thus we can find, for $\epsilon = \epsilon(\BC^{(0)})$, a function $w:\{0\}^2\times\R^{n-1}\to \R^2$ such that
	$$\{Z\in B_{1/2}:\Theta_V(Z)\geq 5/2\} = \graph(w).$$
	To finish the proof, we need to show that $w$ has the desired regularity and find the functions describing $V$ away from the set of points of density $\frac{5}{2}$. Indeed, if $\Theta_V(Z)\geq 5/2$, if we define $\tilde{V}_Z := (\eta_{0,\rho}\circ\Gamma^{-1}_Z)_\#V_Z$, then (III) applied at $Z$ gives
	$$E_{\tilde{V}_Z,\BC_Z}^2 \leq C\rho^{2\alpha}E_{V_Z,\BC^{(0)}}^2$$
	and hence, we may apply Lemma \ref{lemma:L2_coarse} to $\tilde{V}_{Z}$ when $\epsilon$ is sufficiently small, to see that (in particular, from Lemma \ref{lemma:L2_coarse}(i)) for any $\tilde{Y}\in \sing(\tilde{V}_Z)\cap B_{1/2}$ with $\Theta_{\tilde{V}_Z}(\tilde{Y})\geq \frac{5}{2}$,
	$$\dist^2(\tilde{Y},S(\BC_Z)) \leq CE_{\tilde{V}_Z,\BC_Z}^2.$$
	Unpacking this, as $S(\BC_Z) = S(\BC^{(0)})$, it is equivalent to
	$$\rho^{-1-\alpha}\dist(Y,\Gamma_Z(S(\BC^{(0)}))) \leq CE_{V_Z,\BC^{(0)}}$$
	for every $Y\in \sing(V_Z)\cap B_{\rho/2}$ with $\Theta_{V_Z}(Y)\geq\frac{5}{2}$; indeed, taking $Y\to Z$ (with $\Theta_V(Y)\geq 5/2$) we can take $\rho\downarrow 0$ in the above, to see that this implies that $w$ is differentiable, and the tangent plane at $Z = w(z)$ in $\graph(w)$ is $\Gamma_Z(S(\BC^{(0)}))$. So currently we have: $\sup|w| \leq CE_{V,\BC^{(0)}}$ (from Lemma \ref{lemma:L2_coarse}(i)) and that $w$ is differentiable everywhere, with moreover that $\sup|Dw| \leq CE_{V,\BC^{(0)}}$ (from (I)). So all that remains to be shown for $w$ is that it is $C^{1,\alpha}$, with the desired bound on the H\"older semi-norm of $Dw$; from the form of the tangent plane above, this amounts to showing a $C^{0,\alpha}$ bound on $Z\mapsto \Gamma_Z$. So fix $Z_1,Z_2\in \sing(V)\cap \{\Theta_V\geq \frac{5}{2}\}\cap B_{1/2}$ and set $\sigma_*:= |Z_1-Z_2|$. Then, if $\sigma_*<\theta/4$, it follows from (III) that if $V_* := (\eta_{0,4\sigma_*}\circ\Gamma_{Z_1}^{-1})_\#V_{Z_1}$, then
	$$E^2_{V_*,\BC_{Z_1}} \leq C\sigma_*^{2\alpha}E_{V,\BC^{(0)}}^2.$$
	We can then repeat the previous iteration scheme which established (I), (II), and (III), with $V_*$, $\BC_{Z_1}$, and $Z_*:= (4\sigma_*)^{-1}\Gamma_{Z_1}^{-1}(Z_2-Z_1)$ in place of $V$, $\BC_{Z_1}$, and $Z$ (again, for $\epsilon$ sufficiently small independent of $Z_1,Z_2$); hence we find some rotation $\Gamma^*_{Z_*}$ and a cone $\BC^*_{Z_*}$ such that $|\Gamma^*_{Z_*}-\id|\leq CE_{V_*,\BC_{Z_1}}$ and that (from the equivalent of (III)) that $(\Gamma^*_{Z_*})_\#\BC^*_{Z_*}$ is the unique tangent cone to $V_*$ at $Z_*$. However, we know that $(\Gamma_{Z_2})_\#\BC_{Z_2}$ is the unique tangent cone to $V$ at $Z_2$, and so unravelling the transformations which gave rise to $V_*$, equating the unique tangent cones we see that one needs $(\Gamma^*_{Z_*})_\#\BC^*_{Z_*} = (\Gamma_{Z_1}^{-1})_\#((\Gamma_{Z_2})_\#\BC_{Z_2})$, which implies from the form of $\BC^{(0)}$ and (II) that $\Gamma^*_{Z_*} = \Gamma_{Z_1}^{-1}\circ\Gamma_{Z_2}$. Thus, from the property (I), we have
	$$|\Gamma_{Z_2}-\Gamma_{Z_1}| = |\Gamma_{Z_1}\circ\Gamma^*_{Z_*}| = |\Gamma^*_{Z_*}-\id| \leq CE_{V_*,\BC_{Z_1}} \leq C\sigma_*^\alpha E_{V,\BC^{(0)}} \equiv [CE_{V,\BC^{(0)}}]\cdot|Z_1-Z_2|^\alpha$$
	which is the desired $C^{0,\alpha}$ bound when $|Z_1-Z_2|<\theta/4$. However, if $|Z_1-Z_2|\geq \theta/4$, we simply need to iterate the above inequality (using the triangle inequality) at most $N = N(1/\theta) = N(n)$ times to recover the desired inequality, and hence we see that $w\in C^{1,\alpha}$ with the desired bounds; in particular, we have now seen that $B_{1/2}^{n+1}\cap \{\Theta_V>2\}\equiv B^{n+1}_{1/2}\cap \{\Theta_V = 5/2\}$ forms an embedded $C^{1,\alpha}$ submanifold of $B^{n+1}_{1/2}(0)$, each point of which has a tangent cone $\BC_Z\in \FL_0\cap \FL_S$ which obeys $\dist(\spt\|\BC_Z\|\cap B_1,\spt\|\BC^{(0)}\|\cap B_1) \leq CE_{V,\BC^{(0)}}$.
	
	Now all that is left is to prove the existence of the remaining functions defined over the hyperplanes determined by the half-hyperplanes in $\BC^{(0)}$; note that the functions $\gamma_i$ are determined by projecting $\graph(w)$ into the respective hyperplanes, and so from the above it follows that each $\gamma_i$ is a $C^{1,\alpha}$ function; indeed, if we rotate so that one half-hyperplane $H_i$ in $\BC^{(0)}$ is $H_i = \{(0,x^2,y): x^2>0, y\in \R^{n-1}\}$, and $\pi$ denotes the orthogonal projection of $\R^{n+1}$ onto the hyperplane $\{x^1=0\}$ determined by $H_i$, then the function $\gamma_i$ is given by $\gamma_i(y):= (0,w_2(y),y)$, which is still $C^{1,\alpha}$ with the same bounds. If $\epsilon$ is sufficiently small we know that we have $\graph(\gamma_i)\subset \{x^1 = 0, |x^2|<1/16\}$, and so this does split $B_{1/8}\cap H_i$ into two connected components; let $\Omega$ be the component containing $B_{1/8}\cap \{x^2>1/16\}$. But applying the estimates from (III) for points within $\theta/4$ from $\{\Theta_V=5/2\}$, we see that we can locally represent $V$ over the domain by single-valued or two-valued stationary graphs (in fact here they will always be single-valued as all half-hyperplanes in $\BC^{(0)}$ are multiplicity 1); for those points $Y\in \spt\|V\|$ with $\dist(Y,\{\Theta_V = 5/2\})>\theta/4$, the same holds for $\epsilon$ sufficiently small (as $V$ is close to a given hyperplane on $B_{\theta/8}(Y)$, as $(\theta/8)^{-n-2}\int_{B_{\theta/8}(Y)}\dist^2(X,\spt\|\BC^{(0)}\|)^2\ \ext\|V\| \leq (\theta/8)^{-n-2}E_{V,\BC^{(0)}}^2)$. As such, using the unique continuation of single-valued stationary graphs, we construct the desired functions $u_i$ over each half-hyperplane, which moreover have the boundary values determined by $\left. u_i\right|_{\del\Omega \cap B_{1/8}(0)}(0,w_2(y),y) = w_1(y)$; using standard boundary regularity for quasilinear elliptic equations (e.g. \cite{morrey1966multiple}) we therefore deduce that $u_i$ is $C^{1,\alpha}(\overline{\Omega})$, as desired. This completes the proof of Theorem \ref{thm:A} when $\BC^{(0)}\in \FL_S\cap \FL_0$.
\end{proof}

\textbf{Note:} For the last step, even if the $u_i$ are two-valued $C^{1,\alpha}$ stationary graphs, the regularity up-to-the-boundary can still be established by appealing to the Campanato regularity theory for multi-valued functions (\cite{minter2021}) as we still have integral decay estimates at the boundary, provided by those of $w$.

\textbf{Remark:} In fact, in this setting $V$ takes the form of a $C^{1,\alpha}$ classical singularity, so one may apply \cite{krummel2014regularity} to deduce that in fact $V$ is smooth (in fact, real-analytic as in $\R^{n+1}$ with the usual metric), with the points of density $\frac{5}{2}$ being a smooth (real-analytic) submanifold.

\section{The Fine Blow-Up Class}\label{sec:fine_construction}

In this section we provide the construction of the appropriate \textit{fine} blow-up class; such a procedure was originally introduced in \cite{wickstable} to study classical singularities arising in the corresponding coarse blow-up class for that setting. The construction is performed when we have a varifold (or sequence of) which is close to a level $I\in \{1,2\}$ cone, but in fact $V$ is significantly closer to a cone of level $<I$; one would like to say in such a situation that one can deduce some regularity conclusion on $V$ from that of level $<I$ cones rather than that of the level $I$ cone; this degenerate situation will be of crucial importance in proving Theorem \ref{thm:A} when the base cone is level $1$ or $2$, as one may have exactly this situation where cones of a lower level can converge to one of a higher level.

We will need to work with the two-sided excess in this section, as we need to ensure that $V$ is close to all half-hyperplanes in the cone of lower level. We also introduce the following notation: for $V\in \S_2$, we write
$$Q^*_V:=\inf_{\BC\in \FL_1\cup \FL_2}Q_{V,\BC}$$
i.e. $Q_V^*$ is the optimal excess relative to cones of level $1$ or $2$; equivalently, the infimum is taken over all $\BC\in \FL$ which have support consisting of at most $4$ half-hyperplanes. We will only need $Q^*_V$ when $\BC^{(0)}\in \FL_2$ is level 2 (that is, when multiple degenerations can occur simultaneously). 

Let us fix some notation now for our cones. Fix $\BC^{(0)}\in \FL_S\cap \FL_I$, where $I\in \{1,2\}$; so $\BC^{(0)}$ has $I$ multiplicity two half-hyperplanes. For suitable $\epsilon = \epsilon(\BC^{(0)})$, if $\BC\in \FL_{\epsilon}(\BC^{(0)})$, near a multiplicity $q\in \{1,2\}$ half-hyperplane in $\BC^{(0)}$, $\BC$ must have $q$ half-hyperplanes nearby (counted with multiplicity). In particular, if $\BC\in \FL_{I'}$, then $I'\leq I$ and $r:= I-I'$ of the multiplicity two half-hyperplanes in $\BC^{(0)}$ have split into distinct (multiplicity one) half-hyperplanes in $\BC$ (we will see later that $r>0$ in our setting).

Now fix $\BC^c\in \FL_{\epsilon}(\BC^{(0)})\cap \FL_I$ and $\BC\in \FL_{\epsilon}(\BC^{(0)})\cap \FL_{I'}$; $\BC^c$ will represent a nearby cone of the same level as $\BC^{(0)}$ which before we took the coarse blow-up with respect to, and $\BC$ will be another nearby cone which $V$ will have much smaller excess relative to when compared to $\BC^c$. Thus, $\BC^{(0)}$ has $p\equiv 5-2I$ multiplicity one half-hyperplanes, which will be close to $p$ half-hyperplanes in $\BC$, $r$ multiplicity two half-hyperplanes which are close to $2r$ multiplicity one half-hyperplanes in $\BC$ (these are the multiplicity two half-hyperplanes which ``split'' in $\BC$), and $q-r$ multiplicity two half-hyperplanes in $\BC^{(0)}$, which are close to $q-r$ multiplicity two half-hyperplane sin $\BC$. We introduce the following notation for this:
\begin{enumerate}
	\item [(i)] we denote by $H_1^{(0)},\dotsc,H^{(0)}_p$ the multiplicity one half-hyperplanes in $\BC^{(0)}$ and by $H_1,\dotsc,H_p$ the multiplicity one half-hyperplanes $H_1,\dotsc,H_p$ in $\BC$;
	\item [(ii)] we write $\tilde{G}^{(0)}_1,\dotsc,\tilde{G}^{(0)}_r$ for the $r$ multiplicity two half-hyperplanes in $\BC^{(0)}$ such that $\tilde{G}^{(0)}_i$ splits into multiplicity one half-hyperplanes, $\tilde{H}^1_i, \tilde{H}^2_i,$, in $\BC$;
	\item [(iii)] we write $G^{(0)}_1,\dotsc,G^{(0)}_{q-r}$ for the multiplicity two half-hyperplanes in $\BC^{(0)}$ which do not split in $\BC$, and write $G_1,\dotsc,G_{q-r}$ for the corresponding (multiplicity two) half-hyperplanes in $\BC$.
\end{enumerate}
Thus we we have
$$\BC^{(0)} = \sum^p_{i=1}|H^{(0)}_i| + 2\sum^r_{i=1}|\tilde{G}^{(0)}_i| + 2\sum^{q-r}_{i=1}|G^{(0)}_i|$$
and
$$\BC = \sum^p_{i=1}|H_i| + \sum^{r}_{i=1}|\tilde{H}^1_i| + |\tilde{H}^i_2| + 2\sum^{q-r}_{i=1}|G_i|$$
where $H_i$ is achieved from $H^{(0)}_i$ by a small rotation, and similarly $H^1_i$ and $H_i^2$ are achieved by two distinct small rotations from $|\tilde{G}^{(0)}_i|$, and similarly $G_i$ from $G^{(0)}_i$; our notation is chosen so that ``$H$'' always represents a multiplicity one half-hyperplane, ``$G$'' a multiplicity two half-hyperplane, with a ``$\sim$'' representing a half-hyperplane which splits or arises from a split. We also denote the corresponding half-hyperplanes in $\BC^c$ by $H_i^c$, $\tilde{G}_i^c$, and $G_i^c$, i.e.
$$\BC^c = \sum^p_{i=1}|H_i^c| + 2\sum^r_{i=1}|\tilde{G}_i^c| + 2\sum^{q-r}_{i=1}|G_i^c|.$$
Next, choose unit vectors $(\w_i)_{i=1}^p$, $(\tilde{\vartheta})_{i=1}^r$, and $(\vartheta)_{i=1}^{q-r}\subset\R^2$ such that
$$H_i^{(0)} = \{(r\w_i,y): r>0\},\ \ \ \ \tilde{G}^{(0)}_i = \{(r\tilde{\vartheta}_i,y): r>0\},\ \ \ \ G_i^{(0)} = \{(r\vartheta_i,y): r>0\}$$
and similarly choose unit vectors $(\w_i^c)_{i}$, $(\tilde{\vartheta}^c_i)_i$, and $(\vartheta^c_i)_i$ determining $H_i^c$, $\tilde{G}_i^c$, and $G_i^c$, respectively.

Now for $\tau>0$ and $\epsilon = \epsilon(\BC^{(0)},\tau)>0$ sufficiently small, we can find linear functions defined over the half-hyperplanes in $\BC^c$ whose graphs coincide with the half-hyperplanes in $\BC$ on the region $\{|x|>r\}$. That is, we can find single-valued linear functions $(h_i)^p_{i=1}$, $(\tilde{g}_i^j)_{i=1,\dotsc,r;\; j=1,2}$, and $(g_i)_{i=1}^{q-r}$, with
$$h_i:H_i^c\to (H_i^c)^\perp,\ \ \ \ \tilde{g}^j_i: \tilde{G}^c_i\to (\tilde{G}^c_i)^\perp,\ \ \ \ g_i:G^c_i\to (G^c_i)^\perp$$
such that $\graph(h_i|_{\{|x|>\tau\}}) = H_i\cap \{|x|>\tau\}$, $\graph(\tilde{g}^j_i|_{|x|>\tau}) = \tilde{H}^j_i\cap \{|x|>\tau\}$, and $\graph(g_i|_{\{|x|>\tau\}}) = G_i\cap \{|x|>\tau\}$. We then explicitly write
$$h_i(r\w_i,y) = \lambda_i r(\w_i^c)^{\perp_{H_i^c}},\ \ \ \ \tilde{g}^j_i(r\tilde{\vartheta}_i,y) = \tilde{\lambda}^j_i r(\tilde{\vartheta}^c_i)^{\perp_{\tilde{G}_i^c}},\ \ \ \ g_i(r\vartheta_i,y) = \mu_ir(\vartheta^c_i)^{\perp_{G_i^c}}$$
where for $\w\in S^1\subset\R^2$ we write $\w^\perp$ to be $(\w,0)^\perp$; here $\lambda_i,\tilde{\lambda}^j_i,\mu_i\in \R$, and the unit normals are chosen in an anti-clockwise manner, i.e. after rotating the unit vector $\w$ to $(1,0)\in \R^2$, the unit normal is $(0,1)$.

In this section we will be working under various sets of hypotheses. The first are the following:
\vspace{0.5em}
\begin{leftbar}
	\textbf{Hypothesis (H):} For appropriately small $\epsilon, \gamma\in (0,1)$ to be determined depending only on $\BC^{(0)}$, we have
	\begin{enumerate}
		\item [(H1)] $\BC^{(0)}\in \FL_S\cap \FL_I$ for some $I\in \{1,2\}$ and $\BC^c\in \FL_{\epsilon}(\BC^{(0)})\cap \FL_I$;
		\item [(H2)] $V\in \CN_{\epsilon}(\BC^{(0)})$;
		\item [(H3)] $\BC\in \FL_{\epsilon}(\BC^{(0)})$;
		\item [(H4)] $Q_{V,\BC}^2 < \gamma E_{V,\BC^c}^2$.
	\end{enumerate}
\end{leftbar}

\textbf{Remark 1:} Note that there exists $\epsilon = \epsilon(\BC^{(0)})$ such that if Hypothesis (H) hold with any $\gamma\in (0,1)$, then
$$\max_{i,j}\{|\lambda_i|,|\tilde{\lambda}^j_i|,|\mu_i|\} \leq c_1E_{V,\BC^c}$$
where $c_1 = c_1(n)\in (0,\infty)$. Indeed, by Lemma \ref{lemma:coarse_graphical_rep} we know that for $\epsilon = \epsilon(\BC^{(0)})$ sufficiently small, $V$ can be represented by a sum of single-valued and two-valued functions in the region $\{|x|>1/8\}\cap B_{7/8}$ over the half-hyperplanes in $\BC^c$, for which we then get, for example,
\begin{align*}
	c(n)\cdot\lambda_i^2 & \leq \int_{B_{3/4}\cap\{|x|>1/4\}}\dist^2(X,\spt\|\BC\|)\ \ext\|\BC^c\res H^c_i\|\\
	& \leq 4\int_{B_{3/4}\cap \{|x|>1/4\}}\dist^2(X,\spt\|\BC\|)\ \ext\|V\| + 4\int_{B_{3/4}\cap \{|x|>1/4\}}\dist^2(X,\spt\|\BC^c\|)\ \ext\|V\|\\
	& \leq 4E_{V,\BC}^2 + 4E_{V,\BC^c}^2 \leq 4(1+\gamma)E_{V,\BC^c}^2 \leq 8E^2_{V,\BC^c}.
\end{align*}
Now, for $V,\BC^{(0)},\BC^c$, and $\BC$ satisfying Hypothesis (H), we also assume the following for suitable values of $M = M(n)>1$:
\vspace{0.5em}
\begin{leftbar}
	\textbf{Hypothesis ($\star$):} We have
	$$E^2_{V,\BC^c} < M\inf_{\tilde{\BC}\in \FL_I}E^2_{V,\tilde{\BC}}.$$
\end{leftbar}
Hypothesis ($\star$) therefore tells us that $\BC^c$ is close to the best approximating level $I$ cone to $V$.

\textbf{Remark 2:} If Hypothesis (H) and Hypothesis ($\star$) hold for sufficiently small $\epsilon = \epsilon(\BC^{(0)}) \in (0,1)$ and $\gamma = \gamma(\BC^{(0)})\in (0,1)$, then, we have $\BC\in \FL_{I'}$ for some $I'<I$, and moreover there is an $i\in \{1,\dotsc,r\}$ such that $c_2E_{V,\BC^c} \leq |\tilde{\lambda}^1_i - \tilde{\lambda}^2_i|$; here $c_2 = c_2(n)$. In particular, we have
$$c_2E_{V,\BC^c} \leq \max_{i,j}|\tilde{\lambda}^j_i|.$$
Indeed, if $\gamma<1/2M$ then from (H4) and Hypothesis ($\star$) we have $E_{V,\BC} < \frac{1}{2}\inf_{\tilde{\BC}\in \FL_I}E_{V,\tilde{\BC}}$, which implies that $\BC\not\in \FL_I$ and thus $\BC\in \FL_{I'}$ for some $I'<I$. To see the inequalities, create a new level $I$ cone, $\tilde{\BC}\in \FL_I$, from $\BC$ by replacing for each $i\in \{1,\dotsc,r\}$ the two multiplicity one half-hyperplanes $\tilde{H}^1_i$, $\tilde{H}^2_i$ by a single multiplicity two half-hyperplane given by their average, i.e. taking 
$$\tilde{G}_i := \{(r\theta_i,y):r>0\}$$
where $\theta_i$ is the unit vector in $\R^2$ determined by the angle bisector of the angles determined by $\tilde{\vartheta}_i^1$ and $\tilde{\vartheta}_i^2$, then we set
$$\tilde{\BC}:= \sum^p_{i=1}|H_i| + 2\sum^r_{i=1}|\tilde{G}_i| + 2\sum^{q-r}_{i=1}|G_i|;$$
note by construction $\tilde{\BC}\in \FL_I$ and as a graph over $\BC^c$ on the region $\{|x|>1/4\}$, $\tilde{G}_i$ is determined by $\tilde{g}_i(r\tilde{\vartheta}_i,y) = \frac{1}{2}\left(\tilde{\lambda}^1_i + \tilde{\lambda}^2_i\right)r(\tilde{\vartheta}_i^c)^\perp$. Since $\tilde{\BC}\in \FL_I$, by Hypothesis ($\star$) we have $E_{V,\BC^c} < ME_{V,\tilde{\BC}}$, and since triangle inequality gives
$$\dist^2(X,\spt\|\tilde{\BC}\|) \leq 2\dist^2(X,\spt\|\BC\|) + 2\dist^2(\spt\|\BC\|\cap B_1,\spt\|\tilde{\BC}\|\cap B_1)$$
for $X\in B_1$, we see that by integrating this over $X\in \spt\|V\|$ that, for $\epsilon = \epsilon(\BC^{(0)})$ sufficiently small, $E_{V,\tilde{\BC}}^2 \leq 2E_{V,\BC}^2 + 2(6\w_n)\sum^r_{i=1}|\tilde{\lambda}^1_i-\tilde{\lambda}^2_i|^2$, and so we see that
$$M^{-1}E_{V,\BC^c}^2 < 2\gamma E_{V,\BC^c}^2 + 12\w_n\sum^r_{i=1}|\tilde{\lambda}^1_i-\tilde{\lambda}^2_i|^2$$
which shows that, for any $M = M(n)$, provided $2\gamma < (2M)^{-1}$, we get that $(24M\w_n)^{-1}E_{V,\BC^c}^2 \leq \sum^r_{i=1}|\tilde{\lambda}^1_i-\tilde{\lambda}^2_i|^2$. The claim follows from this.

In particular, we see that for suitably chosen $\epsilon$, $\gamma$, and $M$, if Hypothesis (H) and Hypothesis ($\star$) hold, then
\begin{itemize}
	\item If $\BC^{(0)}$ is level 1, then $\BC$ is level 0;
	\item If $\BC^{(0)}$ is level 2, then $\BC$ is level 1 or level 0.
\end{itemize}

Throughout our arguments we will have to take different values for the constant $M$ in Hypothesis $(\star)$. The reason for this is that we cannot guarantee that Hypothesis $(\star)$ holds, with the same $M$, when we perform rescalings and translations of $V$. However, we will see that $M$ will only ever change by a fixed constant factor depending only on $n$. An upper bound on this constant factor we shall name $M_0 = M_0(n)$, and is given by:
$$M_0:= \max\left\{\frac{3}{2}, \frac{2^{4n+20}\w_n^2c_1^2}{\bar{C}_1},\frac{2^{3n+20}\w_n}{\bar{C}_1}\right\}$$
where $c_1 = c_1(n)$ is the constant from Remark 1, and $\bar{C}_1 = \bar{C}_1(n) := \int_{B_{1/2}^n\backslash\{x^2>1/16\}}|x^2|^2\ \ext\H^n(x^2,y)$.

Finally, for $V$, $\BC^{(0)}$, $\BC^c$, and $\BC$ as in Hypothesis (H), for small $\beta\in (0,1/2)$ to be determined depending only on $\BC^{(0)}$, we will also need to consider the following:
\vspace{0.5em}
\begin{leftbar}
	\textbf{Hypothesis ($\dagger$):} Either:
	\begin{enumerate}
		\item [(i)] $\BC^{(0)}\in \FL_I$ and $\BC\in \FL_{I-1}$, where $I\in \{1,2\}$;
		\item [(ii)] $\BC^{(0)}\in \FL_2$,  $\BC\in \FL_0$, and moreover they obey $Q_{V,\BC}^2 < \beta(Q^*_V)^2$.
	\end{enumerate}
\end{leftbar}
\textbf{Remark 3:} If $V,\BC^{(0)},\BC^c,\BC$ are as in Hypothesis (H) and satisfy Hypothesis ($\dagger$)(ii), then there is a constant $c_3 = c_3(n)$ such that for all $\epsilon,\gamma,\beta$ sufficiently small (depending on $\BC^{(0)}$):
$$|\tilde{\lambda}^1_i-\tilde{\lambda}^2_i|\geq 2c_3Q_V^*$$
for all $i=1,\dotsc,r$. Indeed, this follows in the same way as in Remark 2, except instead of replacing \textit{all} splitting multiplicity one half-hyperplane pairs by single multiplicity two half-hyperplanes, we only replace a single pair at a time, a run the same argument to this new cone formed by just collapsing a given pair to a single half-hyperplane of multiplicity two.

\subsection{The Fine Graphical Representation and Initial Estimates}\label{sec:fine_estimates}

The following lemma regarding multiplicity two classes is the first crucial observation for the construction of the fine blow-up class.

\begin{lemma}\label{lemma:splitting}
	Let $\M_2$ be a multiplicity two class and let $\Lambda>0$. Then there exists constants $\epsilon = \epsilon(\M_2,\Lambda)>0$ and $\gamma = \gamma(\M_2,\Lambda)>0$ such that the following is true: if $(V,U_V)\in \M_2$, $\rho>0$, $B_\rho(X_0)\subset U_V$, $\|V\|(B_\rho(X_0))\leq \Lambda$, $\spt\|V\|\cap B_{3\rho/4}(X_0)\neq\emptyset$, and $\rho^{-n-2}\int_{B_{\rho}(X_0)}\dist^2(X,P_1+P_2)\ \ext\|V\|(X)<\epsilon$ for some pair of disjoint hyperplanes $P_1,P_2$, and moreover if
	$$\rho^{-n-2}\int_{B_{\rho}(X_0)}\dist^2(X,P_1+P_2)\ \ext\|V\| < \gamma\inf_{P}\rho^{-n-2}\int_{B_{\rho}(X_0)}\dist^2(X,P)\ \ext\|V\|(X)$$
	where the infimum is taken over all affine hyperplanes $P$, then there are $C^2$ functions $u_i:P_i\cap B_{3\rho/4}(X_0)\to P_i^\perp$ such that $V_j\res B_{5\rho/8}^{n+1}(X_0) = V_1 + V_2$, where $V_i = |\graph(u_i)\cap B_{5\rho/8}^{n+1}(X_0)|$; moreover, $\spt\|V_1\|\cap \spt\|V_2\| = \emptyset$, $\|u_i\|^2_{C^2(P_i\cap B_{3\rho/4}(X_0))}\leq C\rho^{-n-2}\int_{B_\rho(X_0)}\dist^2(X,P_1+P_2)\ \ext\|V\|$, and
	$$\rho^{-n-2}\int_{B_{5\rho/8}^{n+1}(X_0)}\dist^2(X,P_1+P_2)\ \ext\|V\| = \sum^2_{i=1}\rho^{-n-2}\int_{B_{5\rho/8}^{n+1}(X_0)}\dist^2(X,P_i)\ \ext\|V_i\|(X).$$
\end{lemma}

\begin{proof}
	We argue by contradiction. If the result is not true, then we can find sequences $\epsilon_k\downarrow 0$, $\gamma_k\downarrow 0$, $(V_k,U_k)\in \M_2$, $\rho_k>0$, $B_{\rho_k}(X_k)\subset U_k$, with $\|V_k\|(B_{\rho_k}(X_k))\leq\Lambda$, $\spt\|V_k\|\cap B_{3\rho_k/4}(X_k)\neq\emptyset$, $\rho_k^{-n-2}\int_{B_{\rho_k}(X_k)}\dist^2(X,P_k^1+P_k^2)\ \ext\|V_k\|(X)<\epsilon_k$ for some pair of hyperplanes $P^1_k$ and $P_k^2$, and
	\begin{equation}\label{E:splitting-1}
	\left(\inf_P\int_{B_{\rho_k}(X_k)}\dist^2(X,P)\ \ext\|V_k\|(X)\right)^{-1}\int_{B_{\rho_k}(X_k)}\dist^2(X,P_k^1+P_k^2)\ \ext\|V_k\| < \gamma_k
	\end{equation}
	yet the conclusions do not hold. First let us translate and rescale, i.e. consider instead $\tilde{V}_k:= (\eta_{X_k,\rho_k})_\#V_k$, so that we may assume that $X_k = 1$ and $\rho_k=1$ for all $k$. Then, if we have $\inf_P\int_{B_1(0)}\dist^2(X,P)\ \ext\|V_k\|(X)\not\to 0$, then the result holds for all sufficiently large $k$ by Theorem \ref{thm:wick2}, so we may assume that $\inf_P\int_{B_1(0)}\dist^2(X,P)\ \ext\|V_k\|(X)\to 0$. So choose a hyperplane $P_k$ such that
	$$\int_{B_1(0)}\dist^2(X,P_k)\ \ext\|V_k\|(X) < \frac{3}{2}\inf_P\int_{B_1(0)}\dist^2(X,P)\ \ext\|V_k\|(X).$$
	By performing a rotation $\Gamma_k$, we may without loss of generality assume that $P_k = \{0\}\times\R^n$ for all $k$; by passing to a subsequence we may also assume that $\Gamma_k\to\id$. Note that $V_k\weakly V$, where $\spt\|V\| = \{0\}\times B^n_1(0)$. In particular, as $\M_2$ is a multiplicity two class, we must have $V = \theta|\{0\}\times B^n_1(0)|$, where $\theta\in \{1,2\}$. We cannot however have $\theta =1$, as then by Allard's regularity theorem $V_k$ would be expressible as a single smooth graph, and this would contradict (\ref{E:splitting-1}) for all $k$ sufficiently large. So $V_k\weakly 2|\{0\}\times B^n_1(0)|$. In particular, we may apply Theorem \ref{thm:wick1} and apply a blow-up procedure (relative to a fixed hyperplane now) to see that the generated blow-up $v = (v^1,v^2)$ must have, by (\ref{E:splitting-1}), that $v^1$ and $v^2$ have graphs given by disjoint affine hyperplanes; in particular, using the local $C^{1,1/2}$ convergence to the blow-up, we see that there are no points of density $2$ in $V_k\res B_{7/8}^{n+1}(0)$ for all $k$ sufficiently large, and so the two-valued graphical representation provided by Theorem \ref{thm:wick1} is in fact given by two single-valued functions over $\{0\}\times B^{n}_{7/8}(0)$. The conclusions now follow.
\end{proof}

In particular, Lemma \ref{lemma:splitting} tells us that for each $\tau>0$, there is an $\epsilon_0 = \epsilon_0(\BC^{(0)},\tau)>0$ and $\gamma_0 = \gamma_0(\BC^{(0)},\tau)>0$ such that if $V,\BC^{(0)},\BC^c,\BC$ obey Hypothesis (H) and Hypothesis ($\dagger$) with these $\epsilon_0$ in place of $\epsilon$ and $\gamma_0$ in place of $\gamma$, then if $\BC^{(0)}\in \FL_I$ and $\BC\in \FL_{I-1}$ (here, $I\geq 1$ by Remark 2) we are able to express $V$ as a graph over $\BC$ in the region $\{|x|>\tau\}$ such that over the two multiplicity one half-hyperplanes which have split in $\BC$ from a multiplicity two half-hyperplane in $\BC^{(0)}$, $V$ is represented by two single-valued functions, as opposed to a two-valued function.

Our goal will now be to find suitable $\epsilon,\gamma,\beta$, depending only on $\BC^{(0)}$ and $\tau$, such that under Hypothesis (H), Hypothesis ($\star$), and Hypothesis ($\dagger$) we can not only express $V$ as a graph of a function relative to $\BC$, but also that the function obeys the same integral estimates as we saw in Section \ref{sec:coarse_construction}, except now with an upper bound in terms of the excess $E_{V,\BC}$.

\begin{theorem}\label{thm:fine_representation}
	Let $\tau\in (0,1/40)$ and $\BC^{(0)}\in \FL_S\cap \FL_I$, with $I\in \{1,2\}$. Then there exist constants $\epsilon_1 = \epsilon_1(\BC^{(0)},\tau)\in (0,1)$, $\gamma_1 = \gamma_1(\BC^{(0)},\tau)\in (0,1)$, and $\beta_1 = \beta_1(\BC^{(0)},\tau)\in (0,1)$ such that the following is true: let $V$, $\BC^{(0)}$, $\BC^c$, and $\BC$ satisfy Hypothesis (H), Hypothesis $(\star)$, and Hypothesis $(\dagger)$ with $\epsilon_1,\gamma_1,\beta_1$, and $\frac{3}{2}M_0^4$ in place of $\epsilon,\gamma,\beta$, and $M$, respectively, and suppose $\Theta_V(0) \geq \Theta_{\BC^{(0)}}(0) = \frac{5}{2}$. Then we have:
	\begin{align*}
		\ \ \ \ \textnormal{(a)}\ \ & V\res B_{3/4}\cap \{|x|>\tau\} = \mathbf{v}(u)\res\{|x|>\tau\}\text{, where }u\in C^{1,1/2}(\BC\res B_{3/4}^{n+1}(0)\cap \{|x|>\tau\});\\
		& \text{equivalently, we can find }p+2r\text{ single-valued functions }u_1,\dotsc,u_p\text{, }\tilde{u}^1_1, \tilde{u}^2_1,\dotsc,\tilde{u}^1_r,\tilde{u}^2_r\text{, and }\\
		& q-r \text{ two-valued functions }v_1,\dotsc,v_{q-r}\text{, each with their graph being stationary and pairwise}\\
		&\text{disjoint, such that}\\
		&\ \ \ V\res B_{3/4}^{n+1}(0)\cap \{|x|>\tau\} = \sum^p_{i=1}|\graph(h_i + u_i)| + \sum_{i,j}|\graph(\tilde{g}^j_i + \tilde{u}^j_i)| + \sum^{q-r}_{i=1}\mathbf{v}(g_i+v_i)\\
		& \text{where }u_i\in C^2(H_i^c\cap B_{3/4}\cap \{|x|>\tau\}; (H^c_i)^\perp)\text{, }\tilde{u}^j_i\in C^2(\tilde{G}^c_i\cap B_{3/4}\cap \{|x|>\tau\}; (\tilde{G}^c_{i})^\perp)\text{, and}\\
		& v_i\in C^{1,1/2}(G^c_i\cap B_{3/4}\cap \{|x|>\tau\};\A_2((G^c_i)^\perp));\\
		\textnormal{(b)}\ \ & \int_{B_{5/8}^{n+1}(0)}\frac{|X^\perp|^2}{|X|^{n+2}}\ \ext\|V\|\leq CE_{V,\BC}^2;\\
		\textnormal{(c)}\ \ & \int_{B_{5/8}^{n+1}(0)}\sum^{n+1}_{j=3}|e_j^\perp|^2\ \ext\|V\|\leq CE_{V,\BC}^2;\\
		\textnormal{(d)}\ \ & \int_{B^{n+1}_{5/8}(0)}\frac{\dist^2(X,\spt\|\BC\|)}{|X|^{n+3/2}}\ \ext\|V\| \leq CE_{V,\BC}^2;\\
		\textnormal{(e)}\ \ & \int_{B_{5/8}^{n+1}(0)}R^{2-n}\left(\frac{\del(u/R)}{\del R}\right)^2\ \ext\H^n\leq CE_{V,\BC}^2;
	\end{align*}
	here, $C = C(n)$.
\end{theorem}

\textbf{Remark:} Theorem \ref{thm:fine_representation}(e) will give rise to the Hardt--Simon inequality for functions in the fine blow-up class. Such an estimate will only be needed in the setting where $\BC^{(0)}\in \FL_2$ is level 2 and $\BC\in \FL_1$ is level 1, as then one of the functions in the fine blow-up class is two-valued $C^{1,1/2}$ harmonic in the interior, and so to establish its boundary regularity we will need to invoke the methods mentioned in Section \ref{sec:coarse_regularity}.

\begin{proof}
	Let us first prove (a). We already know from Lemma \ref{lemma:splitting} that (a) will hold, with appropriately chosen $\epsilon_1,\gamma_1$, when Hypothesis ($\dagger$)(i) holds, i.e. when $\BC$ is one less level than $\BC^{(0)}$. So let us now assume that $\BC^{(0)}\in \FL_2$ and $\BC\in \FL_0$, so that Hypothesis ($\dagger$) becomes $Q^2_{V,\BC}<\beta (Q^*_V)^2$; we shall argue (a) by contradiction, with our proof using the fact that we already know (a) holds when $\BC\in \FL_1$ and $\BC^{(0)}\in \FL_2$. 
	
	If (a) does not hold in this setting, we may then find sequences $\epsilon_k,\gamma_k,\beta_k\downarrow 0$, $(V_k)_k\subset\S_2$, $(\BC_k)_k\subset\FL_0$, and $(\BC_k^c)_k\subset\FL_2$ such that (a) does not hold, yet Hypothesis (H), Hypothesis ($\star$), and Hypothesis ($\dagger$) hold with $\epsilon_k,\gamma_k,\beta_k, \frac{3}{2}M_0^4$, in place of $\epsilon,\gamma,\beta, M$. Now choose, for each $k$, $\tilde{\BC}_k\in \FL_1\cup \FL_2$ such that
	$$Q_{V_k,\tilde{\BC}_k}^2 < \frac{3}{2} (Q^*_{V_k})^2;$$
	In particular, we have $Q_{V_k,\BC_k}^2 < \beta_k (Q_{V_k}^*)^2 \leq \beta_k Q^2_{V_k,\tilde{\BC}_k}$, and from Remark 3, we have for all $k$ sufficiently large, for $i=1,2$,
	$$|\tilde{\lambda}^1_{i;k} - \tilde{\lambda}^2_{i;k}| \geq 2c_3Q^*_{V_k} > \frac{4c_3}{2}Q_{V_k,\tilde{\BC}_k}$$
	where we recall that $(\tilde{\lambda}^1_{i;k},\tilde{\lambda}^2_{i;k})_{i=1,2}$ are the gradients of the multiplicity one rays in the cross-section of $\BC_k$ relative to the multiplicity two ray in $\BC_k^c$ from which they have split. Thus, we see that if $\tilde{\BC}_k\in \FL_2$ for infinitely many $k$, then by Hypothesis ($\star$) and Lemma \ref{lemma:splitting}\footnote{Note that in fact one can apply essentially the same argument as seen in Lemma \ref{lemma:splitting} to this non-flat setting, replacing the blow-up argument there with a blow-up argument based on the coarse blow-up constructed in Section \ref{sec:coarse_construction}.} we have that we can express $V_k\res B_{3/4}\cap \{|x|>\tau\}$ as a (single-valued) $C^2$ graph over $\BC_k\res B_{3/4}\cap \{|x|>\tau\}$ for all $k$ sufficiently large, providing the contradiction. Thus, we may assume (after passing to a subsequence) that $\tilde{\BC}_k\in \FL_1$ for all $k$.
	
	Now, let $\epsilon_1 = \epsilon_1(\BC^{(0)},\tau/2)$ and $\gamma_1 = \gamma_1(\BC^{(0)},\tau/2)$ be the constants from the theorem in the setting where Hypothesis $(\dagger$)(i) holds. If for infinitely many $k$ we have $Q_{V_k,\tilde{\BC}_k}^2\geq \gamma_1 E^2_{V_k,\BC_k^c}$, by the same argument as above (but now relative to $\BC_k^c\in \FL_2$ as opposed to $\tilde{\BC}_k$) we achieve the same contradiction. Thus we may assume that $Q^2_{V_k,\tilde{\BC}_k} < \gamma_1 E^2_{V_k,\BC^k_c}$ for all but finitely many $k$. But then for all $k$ sufficiently large, we have that the sequences $(V_k)_k$, $(\BC_k^c)_k$, $(\tilde{\BC}_k)_k$ obey Hypothesis (H), Hypothesis $(\star)$, and Hypothesis ($\dagger$) with the constants $\epsilon_1,\gamma_1$, and hence we may apply Theorem \ref{thm:fine_representation}(a) to these sequences, giving the existence of a function $u_k$ representing $V_k$ in $B_{7/8}^{n+1}(0)\cap \{|x|<\tau/2\}$ over $\tilde{\BC}_k$; moreover, the function $u_k$ obeys $|u_k|^2_{C^2(B_{7/8}\cap \{|x|>\tau/2\})}\leq CE_{V_k,\tilde{\BC}_k}^2$ (see Lemma \ref{lemma:splitting}). If one then defines $v_k:= E_{V_k,\tilde{\BC}_k}^{-1}u_k$ (where again, just as we did in the setting of coarse blow-ups, we reparameterise $u_k$ over a fixed domain in $\BC^{(0)}$), we see that $v_k$ converges to some $v\in C^{1,1/2}(\BC\res B_{7/8}\cap\{|x|>\tau/2\})$, where the convergence is locally in $C^{1,1/2}$ on $\spt\|\BC^{(0)}\|\cap B_{7/8}\cap \{|x|>\tau/2\}$. Moreover, by construction we have $Q_{V_k,\tilde{\BC}_k}\leq CE_{V_k,\tilde{\BC}_k}$ for some $C = C(n)$, and
	$$E_{V_k,\tilde{\BC}_k}^{-1}Q_{V_k,\BC_k} \to 0\ \ \ \ \text{and}\ \ \ \ E_{V_k,\tilde{\BC}_k}^{-1}|\hat{\lambda}^1_k-\hat{\lambda}_k^2|\to \hat{c}>0$$
	where $(\hat{\lambda}^1_k,\hat{\lambda}^2_k)$ are the gradients of the multiplicity one half-hyperplanes relative to the multiplicity two half-hyperplane in $\tilde{\BC}_k$ from which they split (this follows in the same way as in Remark 2). Thus, we see that in fact $v$ must be given by $5$ multiplicity one half-hyperplanes in the region $B_{7/8}\cap \{|x|>\tau/2\}$. In particular, this gives that for all $k$ sufficiently large, $V_k\res B_{13/16}\cap \{|x|>3\tau/4\}$ has no points of density $2$; thus the conclusion follows now in the same way as in Lemma \ref{lemma:splitting}.
		
	Now let us turn out attention to the conclusions (b) -- (e). This will follow in the same way as in Lemma \ref{lemma:L2_coarse}, provided that we can extend our graphical representation function $u$ from (a) to a larger subset $U\subset\spt\|\BC^{(0)}\|$ such that $u$ and $G:= \graph(u)$ obey the estimates from Lemma \ref{lemma:coarse_graphical_rep}(c), i.e.
	\begin{equation}\label{E:fine-rep-0}
	\int_{B_{3/4}\backslash G}r^2\ \ext\|V\| + \int_{U\cap B_{3/4}}r^2|Du|^2 \leq CE_{V,\BC}^2
	\end{equation}
	where here the excess term is for $\BC$, not $\BC^c$. We will be able to do this here in a simpler manner, instead of the more complicated argument seen in \cite[Section 10]{wickstable}, thanks to Lemma \ref{lemma:splitting}, which requires the knowledge that we are in a multiplicity two class; our argument here is similar to that seen in \cite[Lemma 6.20]{wickramasekera2004rigidity}.
	
	So let us again introduce the toroidal regions $T_\rho(\zeta):= \{(x,y)\in \R^{n+1}:(|x| - \rho)^2 + |y-\zeta|^2<(\rho/8)^2\}$ and $\tilde{T}_\rho(\zeta):= \{(x,y)\in \R^{n+1}: (|x|-\rho)^2 + |y-\zeta|^2 < (\rho/8)^2/2\}$, where here $(0,\zeta)\in \{0\}^2\times\R^{n-1}$. Take any $\rho\in (0,1/2)$, and for $\delta = \delta(n)$ and $\gamma = \gamma(n)$ sufficiently small (to be chosen) consider the four alternatives:
	\begin{enumerate}
		\item [(i)] $\rho^{-n-2}\int_{T_\rho(\zeta)}\dist^2(X,\spt\|\BC^{c}\|)\ \ext\|V\|<\gamma$ and $$\int_{\tilde{T}_\rho(\zeta)}\dist^2(X,\spt\|\BC\|)\ \ext\|V\| > \delta\int_{T_\rho(\zeta)}\dist^2(X,\spt\|\BC^c\|)\ \ext\|V\|;$$
		\item [(ii)] $\rho^{-n-2}\int_{T_\rho(\zeta)}\dist^2(X,\spt\|\BC^{c}\|)\ \ext\|V\|\geq \gamma$ and $$\int_{\tilde{T}_\rho(\zeta)}\dist^2(X,\spt\|\BC\|)\ \ext\|V\| > \delta\int_{T_\rho(\zeta)}\dist^2(X,\spt\|\BC^c\|)\ \ext\|V\|;$$
		\item [(iii)] $\rho^{-n-2}\int_{T_\rho(\zeta)}\dist^2(X,\spt\|\BC\|)\ \ext\|V\| <\gamma$ and $$\int_{\tilde{T}_\rho(\zeta)}\dist^2(X,\spt\|\BC\|)\ \ext\|V\| \leq \delta\int_{T_\rho(\zeta)}\dist^2(X,\spt\|\BC^c\|)\ \ext\|V\|;$$
		\item [(iv)] $\rho^{-n-2}\int_{T_\rho(\zeta)}\dist^2(X,\spt\|\BC\|)\ \ext\|V\| \geq \gamma$ and $$\int_{\tilde{T}_\rho(\zeta)}\dist^2(X,\spt\|\BC\|)\ \ext\|V\| \leq \delta\int_{T_\rho(\zeta)}\dist^2(X,\spt\|\BC^c\|)\ \ext\|V\|.$$
	\end{enumerate}
	i.e. (i) and (ii) are two alternatives when Hypothesis (H4) fails, and (iii) and (iv) are two alternatives when it holds. Clearly in the case of (i), if $\gamma = \gamma(n)$ is sufficiently small, one may follow the argument in Lemma \ref{lemma:coarse_graphical_rep}, applying Theorem \ref{thm:multiplicity_two_class} and using the bounds provided in (i) to establish (\ref{E:fine-rep-0}) on such regions. In alternative (iii), one may apply Lemma \ref{lemma:splitting} to deduce the same result on these such regions. When either (ii) or (iv) holds, we trivially get $\int_{\tilde{T}_\rho(\zeta)}r^2\ \ext\|V\|\leq C\int_{T_\rho(\zeta)}\dist^2(X,\spt\|\BC^c\|)\ \ext\|V\|$, and so we take $U$ such that $U\cap \tilde{T}_\rho(\zeta)= \emptyset$. Thus, if we define $u$ over the regions determined by (i) and (iii), then as before in Lemma \ref{lemma:coarse_graphical_rep}, we can prove (\ref{E:fine-rep-0}), and thus the result.
\end{proof}

Before stating the next corollary, recall the following basic inequality regarding cone translates, which we saw in (\ref{E:L2-9}): if $Z = (\xi,\zeta)\in \R^2\times\R^{n-1}$ and $S(\BC) = \{0\}^2\times\R^{n-1}$, then
\begin{equation}\label{eqn:cone_translates}
\left|\dist(X,\spt\|(\tau_Z)_\#\BC\|) - \dist(X,\spt\|\BC\|)\right| \leq |\xi|.
\end{equation}

\begin{corollary}\label{cor:fine_estimates}
	Let $\BC^{(0)}\in \FL_S\cap \FL_I$, where $I\in \{1,2\}$. Then there exists $\epsilon_0 = \epsilon_0(\BC^{(0)})$, $\gamma_0 = \gamma_0(\BC^{(0)})$, and $\beta_0 = \beta_0(\BC^{(0)})$ such that the following holds: if $V$, $\BC^{(0)}$, $\BC^c$, and $\BC$ satisfy Hypothesis (H), Hypothesis $(\star)$, and Hypothesis $(\dagger)$ with $\epsilon_0,\gamma_0,\beta_0$, and $\frac{3}{2}M_0^3$ in place of $\epsilon,\gamma,\beta$, and $M$ respectively, then for each $Z = (\xi,\zeta)\in \spt\|V\|\cap (\R^2\times B^{n-1}_{3/8}(0))$ with $\Theta_V(Z)\geq \Theta_{\BC^{(0)}}(0) = \frac{5}{2}$, we have the following:
	\begin{enumerate}
		\item [(a)] $|\xi| \leq CE_{V,\BC}$;
		\item [(b)] For any $\rho\in (0,1)$, if we allow $\epsilon_0,\gamma_0,\beta_0$ to depend on $\rho$ also, we have
		\begin{align*}
		\int_{B^{n+1}_{5\rho/8}(Z)}&\frac{\dist^2(X,\spt\|(\tau_Z)_\#\BC\|)}{|X-Z|^{n+3/2}}\ \ext\|V\|(X) \leq C\rho^{-n-3/2}\int_{B^{n+1}_\rho(Z)}\dist^2(X,\spt\|(\tau_Z)_\#\BC\|)\ \ext\|V\|(X);
		\end{align*}
	\end{enumerate}
	here, $C = C(n)$ is in particular independent of $\rho$.
\end{corollary}

\textbf{Remark:} Unlike in \cite[Corollary 10.2]{wickstable}, we do not need to bound a specific linear combination of $\xi^1,\xi^2$ as the normal directions to the rays in $\BC^{(0)}_0$ will always span $\R^2$, and so we can determine a function on $\R^2$ by the projections onto these rays.

\begin{proof}
	We first argue that for any $\delta\in (0,1)$, there exists $\epsilon_0(\BC^{(0)},\delta)$ and $\gamma_0 = \gamma_0(\BC^{(0)},\delta)$ such that if Hypothesis (H) holds for $V,\BC^{(0)},\BC,\BC^c$ with $\epsilon_0,\gamma_0$ in place of $\epsilon$, $\gamma$, respectively, then
	\begin{equation}\label{eqn:delta_bound}
		|\xi| \leq \delta E_{V,\BC^c}.
	\end{equation}
	Indeed, if this does not hold, then we can find $\delta>0$ and sequences $\epsilon_k,\gamma_k\downarrow 0$, $(V_k)_k$, $(\BC_k)_k$, $(\BC_k^c)_k$, satisfying Hypothesis (H) with $\epsilon_k,\gamma_k$ in place of $\epsilon,\gamma$ respectively, yet there is some $Z_k = (\xi_k,\zeta_k)\in \spt\|V_k\|\cap (\R^2\times B^{n-1}_{3/8})$ with $\Theta_{V_k}(Z_k)\geq \frac{5}{2}$ and $|\xi_k|\geq \delta E_{V_k,\BC^c_k}$. Now, let $v$ be the coarse blow-up of $(V_k)_k$ relative to $(\BC_k^c)_k$, as described in Section \ref{sec:coarse_construction}; thus we have functions $u_k$ defined on $\spt\|\BC^{(0)}\|\cap \{|x|>\tau_k\} \equiv U_{\tau_k}$, for some suitable sequence $\tau_k\downarrow 0$, and from Section \ref{sec:coarse-construction} we know $v_k:= E_{V_k,\BC_k^c}^{-1}u_k \to v\in L^2(\BC^{(0)}\res B_1)$ with the convergence being strong in $L^2(B_1)$. By Lemma \ref{lemma:L2_coarse} we have $|\xi_k|\leq CE_{V_k,\BC_k^c}$ and
	$$\int_{U_{\tau_k}}\frac{|u_k-\xi_k^\perp|^2}{|X-Z_k|^{n+3/2}}\ \ext\H^n \leq CE_{V_k,\BC^c_k}^2$$
	where here $C = C(n)$. Thus, we may pass to a subsequence to ensure that $\zeta_k\to \zeta\in \overline{B}^{n-1}_{3/8}(0)$, $E_{V_k,\BC_k^c}^{-1}\xi_k\to\xi$, where by assumption we know $|\xi|\in [\delta,C]$. Thus, $Z_k\to Z = (0,0,\zeta)$, and so dividing by $E^2_{V_k,\BC_k^c}$ in the above inequality and taking $k\to\infty$ we see
	\begin{equation}\label{E:fine-est-1}
	\int_{\spt\|\BC^{(0)}\|\cap B_1}\frac{|v-\xi^\perp|^2}{|X-Z|^{n+3/2}}\leq C < \infty.
	\end{equation}
	Moreover, as by Hypothesis (H) we know $E_{V_k,\BC_k^c}^{-1}Q_{V_k,\BC_k} \to 0$, this tells us that $v$ must be a linear function over each half-hyperplane, and moreover that $v$ must vanish along the axis $\{0\}^2\times\R^{n-1}$; in particular, $v(Z) = 0$. But then finiteness of the integral in (\ref{E:fine-est-1}) implies that $0 = v(Z) = \xi^\perp$, i.e. the projection of the fixed vector $(\xi,0)$, which obeys $|\xi| \in [\delta,C]$, onto the normal direction of each half-hyperplane in $\BC^{(0)}$ vanishes. But this is impossible unless $\xi = 0$, since we know that the unit vectors in the directions of the rays in the cross-section $\BC^{(0)}_0$ span $\R^2$; thus we have a contradiction to $|\xi|\geq \delta>0$. Thus (\ref{eqn:delta_bound}) holds.
	
	We remark now that when Hypothesis ($\dagger$) is more than a geometric condition, i.e. when $\BC^{(0)}\in \FL_2$ and $\BC\in \FL_0$, we will require a stronger inequality than (\ref{eqn:delta_bound}), namely that $|\xi|\leq \delta Q_V^*$ for any $\delta>0$ (with $\epsilon_0,\gamma_0$ chosen depending on $\delta$). We will be able to prove this in the same way as above once we have established the corollary of the simpler case when $\BC\in \FL_1$ and $\BC^{(0)}\in \FL_2$. So let us first focus on the case where Hypothesis $(\dagger)$(i) holds, i.e. $\BC\in \FL_{I-1}$.
	
	Let $\epsilon_0,\gamma_0,\beta_0$ be the constants given in Theorem \ref{thm:fine_representation} taken with, say, $\tau = 1/16$. Let $\rho\in (0,1/4]$. To prove Corollary \ref{cor:fine_estimates}, we will apply Theorem \ref{thm:fine_representation} with $\tau = 1/16$ and with $(\eta_{Z,\rho})_\#V$ in place of $V$ (with the same $\BC,\BC^c$) for any $Z = (\xi,\zeta)\in \spt\|V\|\cap (\R^2\times B^{n-1}_{3/8})$ with $\Theta_V(Z)\geq \frac{5}{2}$. Thus, we need to choose $\epsilon_0,\gamma_0,\beta_0$ independently of $Z$. Firstly note that the fact $(\eta_{Z,\rho})_\#V\in \CN_{\epsilon_1}(\BC^{(0)})$, where $\epsilon_1 = \epsilon_1(\BC^{(0)})$ is as in Theorem \ref{thm:fine_representation}, when $V\in \CN_{\epsilon_0}(\BC^{(0)})$ and $\epsilon_0 = \epsilon_0(\epsilon_1,\BC^{(0)})$ is sufficiently small, follows from taking $\BC^c = \BC^{(0)}$ in the coarse estimate $|\xi|\leq CE_{V,\BC^c}$ from Lemma \ref{lemma:L2_coarse} (which is strengthened in (\ref{eqn:delta_bound}) above).
	
	Write $\tilde{V}:= (\eta_{Z,\rho})_\#V$. To prove that Hypothesis (H) holds with $\tilde{V}$, $\BC^{(0)}$, $\BC^c$, $\BC$, in place of $V$, $\BC^{(0)}$, $\BC^c$, and $\BC$ respectively, we now just need to show that (H4) holds. We start by showing that we can compare the coarse excess of $\tilde{V}$ relative to $\BC^c$ to that of $V$ relative to $\BC^c$. Indeed, taking $\epsilon = \epsilon(\BC^{(0)},\rho)\in (0,1)$ and $\gamma = \gamma(\BC^{(0)},\rho)\in (0,1)$ sufficiently small (and, when Hypothesis $(\dagger)$(ii) holds, $\beta = \beta(\BC^{(0)},\rho)\in (0,1)$ sufficiently small, the same argument will hold) so that we may apply Theorem \ref{thm:fine_representation} with $\tau = \rho/64$, we get
	\begin{align*}
	E^2_{\tilde{V},\BC^c} & := \int_{B_1}\dist^2(X,\spt\|\BC^c\|)\ \ext\|\tilde{V}\|\\
	& = \rho^{-n-2}\int_{B_\rho(Z)}\dist^2(X-Z,\spt\|\BC^c\|)\ \ext\|V\|\\
	& \geq \rho^{-n-2}\int_{B_\rho(Z)\cap \{|x|>\rho/16\}}\dist^2(X-Z,\spt\|\BC^c\|)\ \ext\|V\|\\
	& \geq \rho^{-n-2}\sum_i\int_{H_i^c\cap B_\rho(Z)\cap \{|x|>\rho/16\}}|h_i+u_i - \xi^{\perp_i}|^2 + \text{(similar terms over other half-hyperplanes)}\\
	& \geq \rho^{-n-2}\sum_i\int_{H_i^c\cap B_{\rho}(Z)\cap \{|x|>\rho/16\}}\left[\frac{1}{4}|h_i|^2  - |u_i|^2\right] - C\rho^{-2}|\xi|^2 + \text{(terms over other half-hyperplanes)}
	\end{align*}
	where in the last inequality we have used the fact that for any real numbers $a,b,c$ we have $|a+b-c|^2 \geq \frac{1}{4}|a|^2 - |b|^c - |c|^2$. Now, using the bounds from Remark 2 and (\ref{eqn:delta_bound}), we see that for any $\delta>0$, if $\epsilon = \epsilon(\BC^{(0)},\delta)$ and $\gamma = \gamma(\BC^{(0)},\delta)$ are sufficiently small,
	\begin{align*}
	E^2_{\tilde{V},\BC^c} & \geq 2^{-n-4}\bar{C}_1\left(\sum_i|\lambda_i|^2 + \sum_{i,j}|\tilde{\lambda}^j_i|^2 + \sum_i|\mu_i|^2\right) - \rho^{-n-2}E_{V,\BC}^2 - C\rho^{-2}\cdot\delta E^2_{V,\BC^c}\\
	& \geq 2^{-n-4}\bar{C}_1\cdot c_2^2E^2_{V,\BC^c} - \rho^{-n-2}\gamma E^2_{V,\BC^c} - C\rho^{-2}\delta E^2_{V,\BC^c}
	\end{align*}
	where $C = C(n)$ and $\bar{C}_1 = \bar{C}_1(n)\equiv \int_{B^n_{1/2}\backslash\{x^2>1/16\}}|x^2|^2\ \ext\H^n(x^2,y)$ is as defined previously (in the definition of $M_0$); here, the extra factor of $2^{-n-2}$ on the first term in the first inequality arises from a lower bound on how the integral of the linear function scales when comparing its integral over $B_{\rho/2}(Z)$ to $B_{\rho/2}(0,\zeta)$ (note that these are essentially integrals of $|x_2|^2$). Thus we see that, if we choose $\delta = \delta(\rho,n)$ sufficiently small, and then $\epsilon = \epsilon(\BC^{(0)},\delta) = \epsilon(\BC^{(0)},\rho)$ and $\gamma = \gamma(\BC^{(0)},\delta,\rho) = \gamma(\BC^{(0)},\rho)$ sufficiently small, we get
	\begin{equation}\label{E:fine-est-2}
		E_{\tilde{V},\BC^c} \geq CE_{V,\BC^c}
	\end{equation}
	for some $C = C(n)$.
	
	Using (\ref{E:fine-est-2}) we can now prove (H4). Firstly, note from (\ref{eqn:cone_translates}) that
	\begin{align*}
		E_{\tilde{V},\BC}^2 & = \rho^{-n-2}\int_{B_\rho(Z)}\dist^2(X,\spt\|(\tau_Z)_\#\BC\|)\ \ext\|V\|\\
		& \leq 2\rho^{-n-2}\int_{B_\rho(Z)}\dist^2(X,\spt\|\BC\|) + |\xi|^2\ \ext\|V\|\\
		& \leq 2\rho^{-n-2}E_{V,\BC}^2 + C\rho^{-2}|\xi|^2
	\end{align*}
	where $C = C(n)$, and provided that $\epsilon = \epsilon(\BC^{(0)},\rho)$ and $\gamma = \gamma(\BC^{(0)},\rho)$ are sufficiently small to ensure that $|\xi|<\rho/64$,
	\begin{align*}
		\int_{B_{1/2}\backslash\{|x|<1/16\}}\dist^2(X,\spt\|\tilde{V}\|)\ \ext\|\BC\| & = \rho^{-n-2}\int_{B_{\rho/2}(Z)\backslash\{|x-\xi|<\rho/16\}}\dist^2(X,\spt\|V\|)\ \ext\|(\tau_Z)_\#\BC\|\\
		& \leq \rho^{-n-2}\int_{B_{33/64}(0,\zeta)\backslash\{|x|<3\rho/64\}}\dist^3(X,\spt\|V\|)\ \ext\|(\tau_Z)_\#\BC\|\\
		& \leq \rho^{-n-2}\int_{B_{9/16}(0)\backslash\{|x|<\rho/32\}}\dist^2(X,\spt\|V\|)\ \ext\|\BC\| + C\rho^{-2}|\xi|^2\\
		& \leq C\rho^{-n-2}Q_{V,\BC}^2 + C\rho^{-2}|\xi|^2.
	\end{align*}
	where we have used (\ref{eqn:cone_translates}) and the fact that $V$ is graphical in $\{|x|>\rho/32\}$ by Theorem \ref{thm:fine_representation} and choice of $\epsilon,\gamma,\beta$, in the second last inequality. Thus combining the above two inequalities, we get
	\begin{equation}\label{E:fine-est-extra1}
	Q_{\tilde{V},\BC}^2 \leq C\rho^{-n-2}Q_{V,\BC}^2 + C\rho^{-2}|\xi|^2
	\end{equation}
	and thus using the fact that from (\ref{eqn:delta_bound}) that for any $\delta>0$ if $\epsilon,\gamma$ are sufficiently small (depending on $\delta$ also) we have $|\xi|^2 \leq \delta E_{V,\BC^c}^2$, and as (H4) holds for $V$ we get $Q_{V,\BC}^2 < \gamma E_{V,\BC^c}^2$, and also from (\ref{E:fine-est-2}) we have $E_{V,\BC^c} \leq CE_{\tilde{V},\BC^c}$, we therefore have
	$$Q_{\tilde{V},\BC}^2 \leq E^2_{\tilde{V},\BC^c}\left[C\rho^{-n-2}\gamma + C\rho^{-2}\delta\right].$$
	Hence for any $\tilde{\gamma}>0$, if we choose $\delta = \delta(n,\rho,\tilde{\gamma})$, $\epsilon = \epsilon(\BC^{(0)},\rho,\delta) = \epsilon(\BC^{(0)},\rho,\tilde{\gamma})$ and $\gamma = \gamma(\BC^{(0)}, \rho,\delta) = \gamma(\BC^{(0)},\rho,\tilde{\gamma})$ sufficiently small, we get
	$$Q_{\tilde{V},\BC}^2 < \tilde{\gamma}E_{\tilde{V},\BC^c}^2$$
	i.e. (H4) holds for $\tilde{V}$.
	
	Next we need to verify that Hypothesis ($\star$) is satisfied with $\tilde{V}$ in place of $V$ (with the same $\BC^c$), for suitable $M$. Let $\tilde{\BC}\in \FL_I$ be close to $\BC^c$ as varifolds, so that over each half-hyperplane in $\BC^c$ we can represent the corresponding half-hyperplane in $\tilde{\BC}$ as a single-valued (perhaps with multiplicity 2) linear function in the normal direction; let us write $\hat{\lambda}_i$ for the gradient of this linear function over the $i^{\text{th}}$ half-hyperplane in $\BC^c$. Then, reasoning as what led us to (\ref{E:fine-est-2}), we get that for sufficiently small $\epsilon,\gamma$ depending on $\BC^{(0)}$ and $\rho$:
	\begin{align}\allowdisplaybreaks\label{E:fine-est-extra-1a}
		\nonumber E^2_{\tilde{V},\tilde{\BC}} & = \rho^{-n-2}\int_{B_\rho(Z)}\dist^2(X-Z,\spt\|\tilde{\BC}\|)\ \ext\|V\|\\
		\nonumber & \geq 2^{-n-4}\bar{C}_1\sum_i|\lambda_i - \hat{\lambda}_i|^2 - \rho^{-n-2}E_{V,\BC}^2 - C\rho^{-2}|\xi|^2\\
		\nonumber & \geq 2^{-n-4}\bar{C}_1\dist^2_\H(\spt\|\BC\|\cap B_1,\spt\|\tilde{\BC}\|\cap B_1) - \rho^{-n-2}E_{V,\BC}^2 - C\rho^{-2}|\xi|^2\\
		\nonumber & \geq 2^{-n-4}\bar{C}_1\cdot(12\w_n)^{-1}\int_{B_1}\dist^2(X,\spt\|\tilde{\BC}\|)\ \ext\|V\|\\
		\nonumber & \hspace{18em} - \left(\rho^{-n-2}+2^{-n-4}\bar{C}_1(12\w_n)^{-1}\right)E^2_{V,\BC} - C\rho^{-2}|\xi|^2\\
		\nonumber & = 2^{-n-4}\bar{C}_1(12\w_n)^{-1}E_{V,\tilde{\BC}}^2 - \left(\rho^{-n-2}+2^{-n-4}\bar{C}_1(12\w_n)^{-1}\right)E_{V,\BC}^2 - C\rho^{-2}|\xi|^2\\
		& \geq \left[(2^{n+8}\w_n)^{-1}\bar{C}_1 - M\gamma\left(\rho^{-n-2}+(2^{n+7}\w_n)^{-1}\bar{C}_1\right)\right]E_{V,\tilde{\BC}}^2 - C\rho^{-2}\delta E_{V,\BC^c}^2\\
		\nonumber & \geq \tilde{C}_1 E_{V,\tilde{\BC}}^2
	\end{align}
	where here $C = C(n)$ and $\delta>0$ can be made arbitrarily small for suitable $\epsilon,\delta,\beta$ (depending on $\delta$); we remark that the third inequality here follows from the triangle inequality in the form $\dist^2(X,\spt\|\tilde{\BC}\|) \leq 2\dist^2(X,\spt\|\BC\|) + 2\dist^2_\H(\spt\|\BC\|\cap B_1,\spt\|\tilde{\BC}\|\cap B_1)$, and that the constant $\tilde{C}_1$ is given by
	$$\tilde{C}_1 := (2^{n+8}\w_n)^{-1}\bar{C}_1 - M\gamma\left(\rho^{-n-2}+ (2^{n+7}\w_n)^{-1}\bar{C}_1\right) - C\rho^{-2}\delta\cdot M.$$
	
	But we also have, using Remark 1, again using (\ref{eqn:cone_translates}) and the triangle inequality, assuming that $\epsilon$ is sufficiently small to ensure $|\xi|\leq \rho/2$,
	\begin{align*}\allowdisplaybreaks
		E_{\tilde{V},\BC^c}^2 & = \rho^{-n-2}\int_{B_\rho(Z)}\dist^2(X-Z,\spt\|\BC^c\|)\ \ext\|V\|\\
		& \leq 2\rho^{-n-2}\int_{B_\rho(Z)}\dist^2(X,\spt\|\BC^c\|)\ \ext\|V\| + C\rho^{-2}|\xi|^2\\
		& \leq 4\rho^{-n-2}\int_{B_\rho(Z)}\dist^2(X,\spt\|\BC\|)\ \ext\|V\|\\
		& \hspace{5em} + 4\rho^{-n-2}\int_{B_\rho(Z)}\dist^2_\H(\spt\|\BC\|\cap B_{2\rho},\spt\|\BC^c\|\cap B_{2\rho})\ \ext\|V\| + C\rho^{-2}|\xi|^2\\
		& \leq 4\rho^{-n-2}E_{V,\BC}^2 + 4\rho^{-n-2}\|V\|(B_{2\rho}(0,\zeta))\cdot\dist^2_\H(\spt\|\BC\|\cap B_{2\rho},\spt\|\BC^c\|\cap B_{2\rho}) + C\rho^{-2}|\xi|^2\\
		& \leq 4\rho^{-n-2}E_{V,\BC}^2 + 2^{n+4}\cdot (6\w_n)\cdot c_1^2 \cdot E^2_{V,\BC^c} + C\rho^{-2}\delta E_{V,\BC^c}^2\\
		& \leq \tilde{C}_2 E_{V,\BC^c}^2
	\end{align*}
	where $C = C(n)$ and
	$$\tilde{C}_2 := 2^{n+7}\w_n c_1^2 + 4\rho^{-n-2}\gamma + C\rho^{-2}\delta.$$
	The above inequalities held whenever $\tilde{\BC}$ was close enough to $\BC^c$ as varifolds. But note that, for suitable $C = C(M,n) = C(n)$, we have\footnote{This can be seen as follows. If this did not hold for any such $C$, we could find sequences $\BC^{(0)}_k$, $\BC^c_k$, $V_k$, such that Hypothesis ($\star$) hold, yet the infimum $\inf_{\tilde{\BC}\in \FL_I}E^2_{V,\tilde{\BC}}$ is attained as some $\hat{\BC}_k\in \FL_I$ for each $k$ such that the distance between $\BC^c_k$ and $\hat{\BC}_k$ is $>kE_{V_k,\BC^c_k}$; in particular $E^2_{V_k,\hat{\BC}_k} < E^2_{V_k,\BC_k^c}$. But we know from Hypothesis ($\star$) that $E_{V_k,\BC_k^c}^2 \leq M E_{V_k,\hat{\BC}_k}^2$ for all $k$. Hence the coarse blow-up sequences of $V_k$ relative to $\BC_k^c$ and $\hat{\BC}_k$, say $v_k$ and $\hat{v}_k$, obey $\|\hat{C}_k\hat{v} - v_k\|_{C^0(B_{3/4}\cap\{|x|>1/8\})}\to \infty$ as $k\to \infty$, where $\hat{C}_k$ is a constant obeying $\hat{C}_k\in [M^{-1},1]$ for all $k$ sufficiently large. This is clearly a contradiction, as we know that both $\hat{v}_k$ and $v_k$ converge uniformly in $B_{3/4}\cap\{|x|>1/8\}$ to $C^{1,1/2}$ functions.} for any such $V$ and $\BC^c$,
\begin{equation}\label{E:fine-est-3}
\inf_{\tilde{\BC}\in \FL_I}E_{V,\tilde{\BC}}^2 = \inf_{\tilde{\BC}\in \FL_I\cap \FL_{CE_{V,\BC^c}}(\BC^c)}E_{V,\tilde{\BC}}^2
\end{equation}
and thus as Hypothesis $(\star)$ holds for $V$ and $\BC^c$, i.e. $E_{V,\BC^c}^2 \leq M\inf_{\tilde{\BC}\cap \FL_I}E_{V,\tilde{\BC}}^2$, the above inequalities give
$$E_{\tilde{V},\BC^c} \leq \tilde{C}_2 E^2_{V,\BC^c} \leq \tilde{C}_2\cdot M\inf_{\tilde{\BC}\in \FL_I}E^2_{V,\tilde{\BC}} \leq \tilde{C}_2M E_{V,\tilde{\BC}}^2 \leq \tilde{C}_2\cdot \tilde{C}_1^{-1}\cdot M E_{\tilde{V},\tilde{\BC}}^2.$$
As this was true for any $\tilde{\BC}$ sufficiently close to $\BC^c$, again using (\ref{E:fine-est-3}) we see that
\begin{equation}\label{E:fine-est-extra2}
E_{\tilde{V},\BC^c}^2 \leq \tilde{C}_2 \tilde{C}_1^{-1}\cdot M\inf_{\tilde{\BC}\in \FL_I}E^2_{\tilde{V},\tilde{\BC}}.
\end{equation}
Thus, we see from the expressions for $\tilde{C}_1$ and $\tilde{C}_2$ that if we choose $\epsilon,\gamma$ sufficiently small, depending on $\BC^{(0)}, \rho$, that $\tilde{C}_2\tilde{C}_1^{-1}\leq (2^{n+8}\w_n c_1^2)(2^{-n-9}\w_n^{-1}\bar{C}_1)^{-1} = 2^{2n+17}\w_n^2 c_1^2 \bar{C}_1^{-1}$, and hence $\tilde{C}_2\tilde{C}_1^{-1}M \leq \frac{3}{2}M_0^4$ (recall that $M = \frac{3}{2}M_0^3$). Hence we have, for $\epsilon,\gamma,\beta$ sufficiently small depending only on $\BC^{(0)}$, $\rho$, that Hypothesis ($\star$) holds for $\tilde{V}$ and $\BC^c$, with $M = \frac{3}{2}M_0^4$. In particular, when Hypothesis ($\dagger$)(i) holds, this completes the proof that we can apply Theorem \ref{thm:fine_representation} with $\tilde{V}$ in place of $V$ (fixing the cones), and thus completes the proof of (b).

Now let us continue to assume that Hypothesis ($\dagger$)(i) and prove (a) in this setting; the same argument will hold when Hypothesis ($\dagger$)(ii) holds once we have shown that Hypothesis ($\dagger$)(ii) still holds for $\tilde{V}$ for suitably chosen $\epsilon,\gamma,\beta$. Note that from (\ref{eqn:delta_bound}) and Remark 3, we see that for each $\theta>0$, there is an $\epsilon = \epsilon(\BC^{(0)},\theta), \gamma = \gamma(\BC^{(0)},\theta)$ such that when the hypothesis hold with these we have for any $X\in\spt\|V\|\cap B_{3/4}^{n+1}(0)\cap\{|x|>\theta\}\equiv W_\theta$, that, if $H$ is the half-hyperplane in $\BC$ closest to $X$,
\begin{equation}\label{E:fine-est-extra4}
\dist(X,\tau_Z (H)) = |\dist(X,H) - \xi^{\perp_H}|
\end{equation}
simply because the bound in Remark 3 ensures that the closest half-hyperplane to $X-Z$ and $X$ is the same. But now, by a similar argument seen in Lemma \ref{lemma:L2_coarse}, for any $\rho_0>0$, for $\epsilon,\delta$ sufficiently small (allowed to depend on $\rho_0$ now also) we deduce the existence of a constant $c = c(n)$ and a subset $S\subset W_{\rho_0/4}\cap B_{\rho_0}(Z)$ with $\H^n(S)\geq \frac{1}{2}\w_n\rho_0^n$ such that for any $X \in S$ we have $E_{V,\BC^c}|\xi| \leq c|\xi^\perp|$ (here, $\xi^\perp$ is the orthogonal projection onto $\BC$). Integrating this inequality over $S$ we then get
\begin{equation}\label{E:fine-est-extra5}
E_{V,\BC^c}^2|\xi|^2 \leq c\rho_0^{-n}\int_{W_\theta\cap B_{\rho_0}(Z)}|\xi^\perp|^2.
\end{equation}
But also, for each half-hyperplane $H$ in $\spt\|\BC\|$, the angle between $H$ and the corresponding half-hyperplane $H^c$ in $\BC^c$ is bounded above by $cE_{V,\BC^c}$ (see Remark 1), and thus we have $|\xi^{\perp_{H^c}}-\xi^{\perp_H}| \leq cE_{V,\BC^c} |\xi|$, and as such combining this with (\ref{E:fine-est-extra5}) we get
\begin{equation}\label{E:fine-est-extra6}
	|\xi^{\perp_{H^c}}|^2 + E^2_{V,\BC^c}|\xi|^2 \leq c\rho_0^{-n}\int_{W_\theta\cap B_{\rho_0}(Z)}|\xi^\perp|^2.
\end{equation}
But now using (\ref{E:fine-est-extra4}), we get
$$|\xi^{\perp_{H^c}}|^2 + E_{V,\BC^c}^2|\xi|^2 \leq c\rho_0^{-n}\int_{B_{\rho_0}(Z)}\dist^2(X,\spt\|(\tau_Z)_\#\BC\|) + c\rho_0^{-n}\int_{B_1}\dist^2(X,\spt\|\BC\|)\ \ext\|V\|.$$
But now using the fact that we know by the above arguments that Theorem \ref{thm:fine_representation}(d) holds with $(\eta_{Z,1/4})_\#V$ in place of $V$ (provided $\epsilon,\gamma$ are sufficiently small independent of $Z$) we get
\begin{align*}
\rho_0^{-n-3/2}\int_{B_{\rho_0}(Z)}\dist^2(X,\spt\|(\tau_Z)_\#\BC\|) \ \ext\|V\| &  \leq C\int_{B_1}\dist^2(X,\spt\|(\tau_Z)_\#\BC\|)\ \ext\|V\|\\
& \leq C\int_{B_1}\dist^2(X,\spt\|\BC\|)\ \ext\|V\| + C|\xi|^2
\end{align*}
where we have used (\ref{eqn:cone_translates}) in the second inequality, and thus we end up with
\begin{equation}\label{E:fine-est-extra7}
|\xi^{\perp_{H^c}}|^2 + E_{V,\BC^c}^2|\xi|^2 \leq C\rho_0^{-n}E_{V,\BC}^2 + C\rho^{3/2}|\xi|^2.
\end{equation}
This was true for any half-hyperplane $H^c$ in $\spt\|\BC^c$. However, as the rays to these half-hyperplanes span $\R^2$, we can sum the above over the different $H^c$ to see that
$$|\xi|^2 \leq C\rho_0^{-n}E_{V,\BC}^2 + C\rho^{3/2}|\xi|^2.$$
Thus, choosing $\rho_0 = \rho_0(n)$ sufficiently small, we get $|\xi|^2\leq CE_{V,\BC}^2$, as desired. Moreover, we know that once we have this, one may return to (\ref{E:fine-est-extra7}) to obtain for each $H^c$,
\begin{equation}\label{E:fine-est-extra8}
	|\xi^{\perp_{H^c}}|^2 + E_{V,\BC^c}^2 |\xi^{\top_{H^c}}|^2 \leq CE_{V,\BC}^2.
\end{equation}

To prove the corollary in the setting where Hypothesis ($\dagger$)(ii) holds, we will need to use the fact that we now know that the corollary holds when Hypothesis ($\dagger$)(i) holds; in particular, we need to use the fact that $|\xi|\leq CE_{V,\BC}$ when $\epsilon = \epsilon(\BC^{(0)}), \gamma = \gamma(\BC^{(0)})$ are sufficiently small and Hypothesis ($\dagger$)(i) holds,, which is of course much stronger than the bound $|\xi|\leq CE_{V,\BC^c}$ provided by Lemma \ref{lemma:L2_coarse}.

So we first claim that for any $\delta>0$, there exist $\epsilon$, $\gamma$, and $\beta$, only depending on $\BC^{(0)}$ and $\delta$ (and $M$) such that if the hypothesis of the corollary hold with these choices of $\epsilon,\gamma$, and $\beta$, then we have
\begin{equation}\label{E:fine-est-4}
|\xi|^2\leq \delta (Q_V^*)^2.
\end{equation}
Indeed, suppose not. Then one may find sequences $\epsilon_k,\gamma_k$, and $\beta_k\downarrow 0$, and sequences $V_k,\BC^c_k,\BC_k$, obeying the hypothesis of the corollary with $\epsilon_k,\gamma_k\beta_k$ in place of $\epsilon,\gamma,\beta$, respectively, yet there is some $Z_k = (\xi_k,\zeta_k)\in \spt\|V_k\|\cap (\R^2\times B^{n-1}_{3/8})$ with $\Theta_{V_k}(Z_k) \geq\frac{5}{2}$ and $|\xi_k|^2\geq \delta (Q_{V_k}^*)^2$. If there is some $t>0$ such that $(Q^*_{V_k})^2 \geq tE_{V_k,\BC^c_k}^2$ for infinitely many $k$, then one may apply the coarse blow-up argument described at the start of the proof to deduce the contradiction in exactly the same manner. Otherwise, we may pass to a subsequence and find a sequence $t_k\downarrow 0$ such that $(Q^*_{V_k})^2 < t_k E^2_{V_k,\BC^c_k}$ for all $k$. So choose for each $k$ a cone $\tilde{\BC}_k\in \FL_1\cup \FL_2$ such that
$$Q_{V_k,\tilde{\BC}_k}^2 < \frac{3}{2}(Q^*_{V_k})^2.$$
We then know, as $Q_{V_k,\tilde{\BC}_k}^2 < \frac{3}{2}M t_k \inf E_{V_k,\hat{\BC}_k}$, where the infimum is taken over all $\hat{\BC}_k\in \FL_2$; in particular, for all $k$ sufficiently large, we have $\tilde{\BC}_k\in \FL_1$. We are therefore in the setting where one can produce a fine blow-up of $V_k$ relative to the sequences $(\tilde{\BC}_k)_k\subset \FL_1$, $(\BC_k^c)\subset \FL_2$ (along which Hypothesis ($\dagger$)(i) holds). To describe this (see a more detailed description in Section \ref{sec:fine_properties}), take a sequence $\tau_k\downarrow 0$ sufficiently slowly. Then one may pass to a subsequence to apply for each $k$ the results of Theorem \ref{thm:fine_representation} to $V_k$, $\tilde{\BC}_k$, $\BC^c_k$ with $\tau = \tau_k$ to generate a function $u_k$ describing $V_k$ relative to $\tilde{\BC}_k$ in the region $B_{1-\tau_k}^{n+1}(0)\cap \{|x|>\tau_k\}$. The estimates provided by Theorem \ref{thm:fine_representation} and Corollary \ref{cor:fine_estimates} give that the sequence $v_k:= E_{V_k,\tilde{\BC}_k}^{-1}u_k$ (suitably parameterised over the fixed cone $\BC^{(0)}$) converges, strongly in $L^2(B_1)$ and locally in $C^{1,1/2}$ in $B_1\cap\{|x|>0\}$, to a function $v\in C^{1,1/2}(\BC\res B_1\cap \{|x|>0\})$. Moreover, since by Hypothesis ($\dagger$)(ii) holds with $\beta = \beta_k$ for $V_k$, $\BC^c_k$, $\BC_k$, we know that $Q_{V_k,\BC_k}^2 < \beta_k (Q_{V_k}^*)^2$, and thus as here we have, from Theorem \ref{thm:fine_representation} that $Q_{V_k,\tilde{\BC}_k}^2 \leq CE_{V_k,\tilde{\BC}}^2$ for some $C = C(n)$ fixed, we see that $E_{V_k,\tilde{\BC}_k}^{-1}Q_{V_k,\BC_k}\to 0$. In particular, this tells us that $v$ must be supported on formed by linear functions over each of the half-hyperplanes in $\BC^{(0)}$ which vanish along the axis $\{0\}^2\times\R^{n-1}$. But from Theorem \ref{thm:fine_representation} (more so Corollary \ref{cor:fine_estimates}) for the construction of the fine blow-up $v$, we have (in the same way as for the construction of the coarse blow-up class) that
$$\int_{\spt\|\BC^{(0)}\|\cap B_1}\frac{|v-\xi^\perp|^2}{|X-Z|^{n+3/2}} \leq C < \infty$$
where $Z_k \to Z = (0,0,\zeta)$ and $E_{V_k,\tilde{\BC_k}}^{-1}\xi_k\to \xi$, where by assumption we have $|\xi|^2\geq \delta (Q_{V_k}^*)^2 > \frac{2\delta}{3}E_{V_k,\tilde{\BC}_k}^2$, and $|\xi|^2 \leq CE_{V_k,\tilde{\BC}_k}$ by Corollary \ref{cor:fine_estimates}; hence $|\xi|\in [2\delta/3,C]$, i.e. $|\xi|>0$. But the above integral being finite implies that $\xi^{\perp_{H_i}} = 0$ for each half-hyperplane in $\spt\|\BC^{(0)}\|$, which, as the unit vectors in the directions of the rays of the cross-section $\BC^{(0)}_0$ span $\R^2$, implies that $|\xi| = 0$, a contradiction to the fact that $|\xi|>0$. Hence we have established (\ref{E:fine-est-4}).

We know that we can prove that $\tilde{V}$ obeys Hypothesis (H) and Hypothesis ($\star$) with respect to $\BC^c$ and $\BC$ in exactly the same way as above. We need (\ref{E:fine-est-4}) to show that Hypothesis ($\dagger$)(ii) will still hold, for suitably chosen $\epsilon,\gamma,\beta$, depending only on $\BC^{(0)}$ and $\rho$. We know from Remark 3 that for $\epsilon,\gamma,\beta$ sufficiently small depending on $\BC^{(0)}$, if $\tilde{\BC}\in \FL_1\cup \FL_2$ is any level 1 or level 2 cone, we will have
$$\dist^2_\H(\spt\|\tilde{\BC}\|\cap B_1,\spt\|\BC\|\cap B_1)\geq C(Q^*_V)^2$$
where $C = C(n)$. Thus, for sufficiently small $\epsilon,\gamma,\beta$ so that Theorem \ref{thm:fine_representation} holds with $\tau = \rho/64$, we can estimate similarly as before to get
\begin{align*}
	E_{\tilde{V},\tilde{\BC}}^2 & = \rho^{-n-2}\int_{B_\rho(Z)}\dist^2(X-Z,\spt\|\tilde{\BC}\|)\ \ext\|V\|\\
	& \geq \rho^{-n-2}\int_{B_\rho(Z)\cap \{|x|>\rho/16\}}\dist^2(X-Z,\spt\|\tilde{\BC}\|)\ \ext\|V\|\\
	& \geq \rho^{-n-2}\sum_i\int_{H^c_i\cap B_\rho(Z)\cap \{|x|>\rho/16\}}\dist^2(h_i(\tilde{X})+u_i(\tilde{X})-\zeta^{\perp_i},\spt\|\tilde{\BC}\|)\ \ext\H^n(\tilde{X})\\
	& \hspace{23em} + \text{(terms over other half-hyperplanes)}\\
	& \geq \rho^{-n-2}\sum_i\int_{H_i^c\cap B_\rho(Z)\cap \{|x|>\rho/16\}}\dist^2(h_i(\tilde{X}),\spt\|\tilde{\BC}\|)\ \ext\H^n(\tilde{X})- C\rho^{-n-2}E_{V,\BC}^2 - C\rho^{-2}|\xi|^2\\
	& \hspace{23em} + \text{(terms over other half-hyperplanes)}\\
	& \geq C\dist^2(\spt\|\tilde{\BC}\|\cap B_{1},\spt\|\BC\|\cap B_1) - \rho^{-n-2}E_{V,\BC}^2 - C\rho^{-2}|\xi|^2\\
	& \geq C(Q^*_V)^2 - \rho^{-n-2}\beta(Q_V^*)^2 - C\rho^{-2}\delta(Q_V^*)^2
\end{align*}
where in the last inequality we have used Hypothesis ($\dagger$)(ii) (which holds for $V,\BC$) and (\ref{E:fine-est-4}) (for any $\delta>0$, provided we allow $\epsilon,\gamma,\beta$ to depend on $\delta$). Hence, choosing $\delta = \delta(n,\rho)$ sufficiently small, and $\epsilon,\gamma,\beta$ sufficiently small accordingly (depending only on $\BC^{(0)}$, $\rho$) we see that
$$E_{\tilde{V},\tilde{\BC}}^2 \geq C(Q_V^*)^2,\ \ \ \ \text{and thus}\ \ \ \ Q^2_{\tilde{V},\tilde{\BC}}\geq C(Q_V^*)^2$$
where $C = C(n)$. This argument held for any $\tilde{\BC}\in \FL_1\cap \FL_2$ sufficiently close to $\BC$, but arguing as before we know that $\inf_{\hat{\BC}\in \FL_1\cup \FL_2}Q_{V,\hat{\BC}} = \inf_{\hat{\BC}\in (\FL_1\cup \FL_2)\cap \FL_{CE_{V,\BC}}(\BC)}Q_{V,\hat{\BC}}$, and this taking the infimum over all such $\tilde{\BC}$, we see
$$(Q_{\tilde{V}}^*)^2 \geq C(Q^*_V)^2.$$
But then we know from (\ref{E:fine-est-extra1}) that
$$Q_{\tilde{V},\BC}^2 \leq C\rho^{-n-2}Q_{V,\BC}^2 + C\rho^{-2}|\xi|^2$$
and so combining this with the above estimates we have
$$Q_{\tilde{V},\BC}^2 \leq C\beta\rho^{-n-2}(Q_V^*)^2 + C\rho^{-2}\delta (Q_V^*)^2 \leq C(\beta\rho^{-n-2} + \rho^{-2}\delta)(Q^*_{\tilde{V}})^2$$
and so for any $\tilde{\beta}$, choosing $\delta = \delta(n,\rho,\tilde{\beta})$ sufficiently small, and then $\epsilon, \gamma,\beta$ sufficiently small depending on $\BC^{(0)},\rho$ and $\tilde{\beta}$ accordingly, we have
\begin{equation}\label{E:fine-est-5}
Q_{\tilde{V},\BC}^2 < \tilde{\beta}(Q_{\tilde{V}}^*)^2.
\end{equation}
So choosing $\tilde{\beta} = \tilde{\beta}(\BC^{(0)})$ as in Theorem \ref{thm:fine_representation}, we get that $\tilde{V}$, $\BC$, satisfy Hypothesis ($\dagger$)(ii), and thus we can apply Theorem \ref{thm:fine_representation} to complete the proof.
\end{proof}

\textbf{Remark 4:} Notice that whilst Corollary \ref{cor:fine_estimates} establishes $|\xi|\leq CE_{V,\BC}$ for $Z = (\xi,\zeta)$ with $\Theta_V(Z)\geq 5/2$, it also establishes the finer inequality (\ref{E:fine-est-extra8}). Notice that if $H$ is a half-hyperplane in $\BC$ which is generated from the half-hyperplane $H^c$ in $\BC^c$, and moreover $H$ is represented over $H^c$ by a linear function with gradient $\lambda$, then
\begin{equation}\label{E:remark-4}
\xi^{\perp_H} = \xi^{\perp_{H^c}} - \lambda \xi^{\top_{H^c}}.
\end{equation}
Such an equality, combined with Remark 3 and (\ref{E:fine-est-extra8}) will be crucial for showing that any multiplicity two half-hyperplane which splits in $\BC$ generates in the blow-up two \textit{separated} single-valued functions.

Armed now with Corollary \ref{cor:fine_estimates}(b), we can now prove that the fine excess, $E_{V,\BC}$, does not accumulate along the spine, giving the corresponding result to Corollary \ref{corollary:non-concentration} in the coarse blow-up setting; in particular, we will get strong $L^2$ convergence to the fine blow-up.

\begin{lemma}\label{lemma:non-concentration_fine}
	Let $I\in \{1,2\}$, $\delta\in (0,1/10)$ and $\BC^{(0)}\in \FL_S\cap \FL_I$. Then there exists $\epsilon_1 = \epsilon_1(\BC^{(0)},\delta)\in (0,1)$, $\gamma_1 = \gamma_1(\BC^{(0)},\delta)\in (0,1)$, and $\beta_1 = \beta_1(\BC^{(0)})\in (0,1)$ such that the following is true: if $V,\BC^{(0)},\BC^c,\BC$ satisfy Hypothesis (H), Hypothesis $(\star)$, and Hypothesis $(\dagger)$ with $\epsilon_1,\gamma_1,\beta_1$, and $\frac{3}{2}M_0^3$ in place of $\epsilon,\gamma,\beta$, and $M$ respectively, then for each $\sigma\in [\delta,1/4)$,
	$$\int_{B_{3/4}\cap \{|x|<\sigma\}}\dist^2(X,\spt\|\BC\|)\ \ext\|V\| \leq C\sigma^{1/2}E_{V,\BC}^2$$
	where $C = C(n)$ is independent of $\delta$.
\end{lemma}

\begin{proof}
	From Corollary \ref{cor:fine_estimates}(i) we have that, for each $Z = (\xi,\zeta)\in \spt\|V\|\cap B_{3/8}$ with $\Theta_V(Z)\geq \frac{5}{2}$ and any $X\in \R^{n+1}$ that if $\epsilon,\gamma,\beta$ are sufficiently small depending on $\BC^{(0)}$ (recall (\ref{eqn:cone_translates})),
	$$|\dist(X,\spt\|(\tau_Z)_\#\BC\|) - \dist(X,\spt\|\BC\|)| \leq C|\xi| \leq CE_{V,\BC}.$$
	Given this, we can now argue as in \cite[Corollary 3.2]{simoncylindrical}, using Lemma \ref{lemma:gaps}.
\end{proof}

\subsection{Constructing the Fine Blow-Up Class}\label{sec:fine_properties}

Using the results of Section \ref{sec:fine_estimates} we now construct the class of fine blow-ups.

Fix $M_1 = M_1(n)\in (1,\infty)$, $I\in \{1,2\}$, and $\BC^{(0)}\in \FL_S\cap \FL_I$ throughout. Let $(\epsilon_k)_k$, $(\gamma_k)_k$, and $(\beta_k)_k$ be (decreasing) sequences of positive numbers converging to $0$. Consider sequences of varifolds $(V_k)_k\subset\S_2$, $(\BC^c_k)_k\subset\FL_I$, and $(\BC_k)_k\subset\FL$ such that, for each $k=1,2,\dotsc$, $V_k$, $\BC^{(0)}$, $\BC^c_k$, $\BC_k$ obey Hypothesis (H), Hypothesis $(\star)$, and Hypothesis $(\dagger)$ with $\epsilon_k,\gamma_k,\beta_k$, and $M_1$, in place of $\epsilon,\gamma,\beta$, and $M$, respectively. Thus, for each $k=1,2,\dotsc$, we assume:
\begin{enumerate}
	\item [($1_k$)] $V_k\in \CN_{\epsilon_k}(\BC^{(0)})$;
	\item [($2_k$)] $\BC_k\in \FL_{\epsilon_k}(\BC^{(0)})$ and $\BC^c_k\in \FL_{\epsilon_k}(\BC^{(0)})\cap \FL_I$;
	\item [($3_k$)] $E^{-2}_{V_k,\BC_k^c}Q^2_{V_k,\BC_k} < \gamma_k$;
	\item [($4_k$)] $E^2_{V_k,\BC^c_k} < M_1 \inf_{\tilde{\BC}\in \FL_I}E^2_{V_k,\tilde{\BC}}$;
	\item [($5_k$)] One of (i) or (ii) below holds:
	\begin{enumerate}
		\item [(i)] $\BC_k\in \FL_{I-1}$;
		\item [(ii)] $I=2$, $\BC_k\in \FL_0$, and $(Q^*_{V_k})^{-2}Q^2_{V_k,\BC_k}<\beta_k$.
	\end{enumerate}
\end{enumerate}

Write $p_k$, $q_k$ for the number of multiplicity one half-hyperplanes in $\BC_k$ respectively, and set $r_k:=q^{(0)}-q_k$, where $q^{(0)}:= I$, for the change in the level between $\BC^{(0)}$ and $\BC_k$, i.e. the number of splitting multiplicity two half-hyperplanes; we know that $r_k\geq 1$ for all $k$ sufficiently large by Remark 2, and therefore that $p_k= p^{(0)}+2r_k$, where $p^{(0)}:= 5-2I$; moreover, we may pass to a subsequence to ensure that $r_k\in \{1,2\}$ is a constant (and hence $p_k,q_k$, are constant also), i.e. $r_k\equiv r$ for all $k$, and that the hyperplane(s) in $\BC^{(0)}$ which split are the same for all $k$. In particular, we write $\BC^{(0)} = \sum^{p^{(0)}}_{i=1}|H_i^{(0)}| + 2\sum^{q^{(0)}-r}_{i=1}|G^{(0)}_i| + 2\sum^{r}_{i=1}|\tilde{G}_i^{(0)}|$, where the $\tilde{G}_i^{(0)}$ are the multiplicity two half-hyperplanes in $\BC^{(0)}$ which split in $\BC_k$, i.e. are close to two multiplicity one half-hyperplanes in $\BC_k$. Similarly, we write $\BC^c_k = \sum^{p^{(0)}}_{i=1}|H^c_{i,k}| + 2\sum^{q^{(0)}-r}_{i=1}|G^c_{i,k}| + 2\sum^r_{i=1}|\tilde{G}^c_{i,k}|$, and
$$\BC_k = \sum^{p^{(0)}}_{i=1}|H^k_i| + \sum^{r_k}_{i=1}\left(|\tilde{H}_i^{k,1}| + |\tilde{H}^{k,2}_i|\right) + 2\sum^{q^{(0)}-r}_{i=1}|G^k_i|$$
where $\tilde{H}^{k,1}_i$, $\tilde{H}^{k,2}_i$ are the two multiplicity one half-hyperplanes in $\BC_k$ close to the multiplicity two half-hyperplane $\tilde{G}^{(0)}_i$ in $\BC^{(0)}$. Moreover, for some (decreasing) sequence $(\tau_k)_k$ converging to zero sufficiently slowly, on $\{|x|>\tau_k\}$ we write $(\lambda^k_i)_{i=1}^{p^{(0)}}$, $(\tilde{\lambda}_i^{k,j})_{i=1,\dotsc,r;\; j=1,2}$, and $(\mu^k_i)_{i=1}^{q^{(0)}-r}$, for the gradients of the respective half-hyperplanes in $\BC_k$ relative to the corresponding half-hyperplanes in $\BC^c_k$; these constants therefore determine linear functions $(h^k_i)_{i=1}^{p^{(0)}}$, $(\tilde{g}^{k,j}_i)_{i=1,\dotsc,r;\;j=1,2}$, and $(g^k_i)_{i=1}^{q^{(0)}-r}$ whose graphs are the respective half-hyperplanes in the region $\{|x|>\tau_k\}$. Our fine blow-ups will be defined relative to the $\BC^c_k$, however we will use the fixed domain $\spt\|\BC^{(0)}\|$ as a parameter space for our functions, so that they have a fixed domain of definition (just as in Section \ref{sec:coarse-construction}); however, we shall suppress this extra notation for the sake of ease of presentation, and interchange between functions defined on $\BC^k_c$ and $\BC^{(0)}$ freely.

Now let $(\delta_k)_k$ be a decreasing sequence of positive numbers converging to $0$. Changing the definitions of $(\delta_k)_k$, $(\tau_k)_k$ if necessary to ensure that they do not go to zero too quickly, we may then deduce from the results of Section \ref{sec:fine_estimates} that the following assertions hold:

\begin{enumerate}
	\item [(A$_k$)] For every point $Y\in S(\BC^{(0)})\cap B_{1/2}$, we have for all $k$ sufficiently large,
	$$B_{\delta_k}(Y)\cap \{Z:\Theta_{V_k}(Z)\geq 5/2\} \neq\emptyset;$$
	\item [(B$_k$)] For each $\sigma\in [\delta_k,1/4)$ we have
	$$\int_{B_{3/4}\cap\{|x|<\sigma\}}\dist^2(X,\spt\|\BC_k\|)\ \ext\|V_k\| \leq C\sigma^{1/2} E^2_{V_k,\BC_k};$$
	\item [(C$_k$)] There are $p^{(0)}+2r$ single-valued functions, $(u^k_i)_{i=1}^{p^{(0)}}$, $(\tilde{u}^{k,j}_i)_{i=1,\dotsc,r;\;j=1,2}$, and $q^{(0)}-r$ two-valued functions, $(v^k_i)_{i=1}^{q^{(0)}-r}$, where $u^k_i\in C^2(H^c_{i,k}\cap B_{3/4}\cap \{|x|>\tau_k\};(H^c_{i,k})^\perp)$, $\tilde{u}^{k,j}_i\in C^2(\tilde{G}^c_i\cap B_{3/4}\cap \{|x|>\tau_k\};(\tilde{G}^c_{i,k})^\perp)$, and $v^k_i\in C^{1,1/2}(G^c_{i,k}\cap B_{3/4}\cap \{|x|>\tau_k\};\A_2((G^c_{i,k})^\perp))$, each with stationary graph, such that
	$$V_k\res (B_{3/4}\cap \{|x|>\tau_k\}) = \sum^{p^{(0)}}_{i=1}|\graph(h^k_i+u^k_i)| + \sum_{i,j}|\graph(\tilde{g}^{k,j}_i+\tilde{u}^{k,j}_i)| + \sum^{q^{(0)}-r}_{i=1}\mathbf{v}(g_i^k + v^k_i);$$
	\item [(D$_k$)] For each point $Z = (\xi,\zeta)\in \spt\|V_k\|\cap B_{3/8}$ with $\Theta_{V_k}(Z)\geq 5/2$ we have
	$$|\xi| \leq CE_{V_k,\BC_k};$$
	\item [(E$_k$)] We have
	$$c_2 E_{V_k,\BC_k^c} \leq \max_{i,j}\{|\lambda_i^k|,|\tilde{\lambda}_i^{k,j}|,|\mu_i^k|\} \leq c_1E_{V_k,\BC_k^c}$$
	and moreover for some $i\in \{1,\dotsc,r\}$ we have
	$$|\tilde{\lambda}^{k,1}_i - \tilde{\lambda}^{k,2}_i|\geq 2c_3 E_{V_k,\BC_k^c};$$
	\item [(F$_k$)] For each $\rho\in (0,1/4]$, we can find $K = K(\rho)\in \Z_{\geq 1}$ such that for all $k\geq K$ the following holds: for each $Z = (\xi,\zeta)\in \spt\|V_k\|\cap B_{3/8}$ with $\Theta_{V_k}(Z)\geq 5/2$,
	\begin{align*}
		&\sum^{p^{(0)}}_{i=1}\int_{H_i^{(0)}\cap B_{\rho/2}(Z)\cap \{|x|>\tau_k\}}\frac{|u_i^k - \xi^{\perp_{H_i^k}}|^2}{|(h^k_i(r\w_i,y) + u^k_i(r\w_i,y), r\w_i,y)-Z|^{n+3/2}}\\
		& \hspace{1.5em} + \sum_{i,j}\int_{\tilde{G}^{(0)}_i\cap B_{\rho/2}(Z)\cap \{|x|>\tau_k\}}\frac{|\tilde{u}^{k,j}_i - \xi^{\perp_{\tilde{H}^{k,j}_i}}|^2}{|(\tilde{h}^{k,j}_i(r\tilde{\w}_i,y) + \tilde{u}^{k,j}_i(r\tilde{\w}_i,y), r\tilde{\w}_i,y) - Z|^{n+3/2}}\\
		& \hspace{3em} + \sum^{q^{(0)}-r}_{i=1}\int_{G^{(0)}_i\cap B_{\rho/2}(Z)\cap \{|x|>\tau_k\}}\frac{|v^k_i - \xi^{\perp_{G^k_i}}|^2}{|(g^k_i(r\vartheta_i,y) + v^k_i(r\vartheta_i,y),r\vartheta_i,y) - Z|^{n+3/2}}\\
		& \hspace{20em} \leq C\rho^{-n-3/2}\int_{B_\rho(Z)}\dist^2(X,\spt\|(\tau_Z)_\#\BC_k\|)\ \ext\|V_k\|;
	\end{align*}
	\item [(G$_k$)] For each $\rho\in (0,1/4]$, we can find $K = K(\rho)\in \Z_{\geq 1}$ such that for all $k\geq K$ the following holds: for each $Z = (\xi,\zeta)\in \spt\|V_k\|\cap B_{3/8}$ with $\Theta_{V_k}(Z)\geq 5/2$,
	\begin{align*}
		&\sum^{p^{(0)}}_{i=1}\int_{H^{(0)}_i\cap B_{\rho/2}(Z)\cap \{|x|>\tau_k\}} R_Z^{2-n}\left(\frac{\del(u^k_i/R_Z)}{\del R_Z}\right)^2\\
		& \hspace{1.5em} + \sum_{i,j}\int_{\tilde{G}^{(0)}_i\cap B_{\rho/2}(Z)\cap \{|x|>\tau_k\}}R_Z^{2-n}\left(\frac{\del(\tilde{u}^{k,j}_i/R_Z)}{\del R_Z}\right)^2\\
		& \hspace{3em} + \sum^{q^{(0)}-r}_{i=1}\int_{G_i^{(0)}\cap B_{\rho/2}(Z)\cap \{|x|>\tau_k\}}R_Z^{2-n}\left(\frac{\del(v^k_i/R_Z)}{\del R_Z}\right)^2\\
		& \hspace{20em} \leq C\int_{B_\rho(Z)}\dist^2(X,\spt\|(\tau_Z)_\#\BC_k\|)\ \ext\|V_k\|
	\end{align*}
	where $R_Z(X):= |X-Z|$
\end{enumerate}
In all the above, $C = C(n)$ is a fixed dimensional constant. To see why the above inequalities hold, note that (A$_k$) follows from Lemma \ref{lemma:gaps}, (B$_k$) follows from Lemma \ref{lemma:non-concentration_fine}, (C$_k$) follows from Theorem \ref{thm:fine_representation}(a), (D$_k$) follows from Corollary \ref{cor:fine_estimates}, (E$_k$) follows from Remark 1, Remark 2, and Remark 3 (with the modified form of Remark 3 following from (B$_k$)), (F$_k$) follows from Corollary \ref{cor:fine_estimates} (in the same was as the corresponding inequality for the coarse blow-up classes did there, from Corollary \ref{corollary:non-concentration}), and (F$_k$) follows from Theorem \ref{thm:fine_representation}(e), applied to $(\eta_{Z,\rho})_\# V_k$, which is possible for all sufficiently large $k$ by the argument in the proof of Corollary \ref{cor:fine_estimates}. Note that the constant $M_1$ will only change by a factor of $M_0$ in the proofs of these statements, and so for all $k$ sufficiently large we are still able to apply Theorem \ref{thm:fine_representation} to $(\eta_{Z,\rho})_\#V_k$. We extend $u^k_i$, $\tilde{u}_i^{k,j}$, and $v^k_j$ to all of $H_i^{(0)}\cap B_{3/4}$, $\tilde{G}^{(0)}_j\cap B_{3/4}$, and $G_i^{(0)}\cap B_{3/4}$, respectively, by defining them to be zero outside their domain of definition.

By (E$_k$), we can find numbers $(\ell_i)_{i=1}^{p^{(0)}},(\tilde{\ell}_i^j)_{i=1,\dotsc,r;\; j=1,2}$, and $(m_i)_{i=1}^{q^{(0)}-r}$, obeying
$$c_2\leq \max_{i,j}\{|\ell_i|, |\tilde{\ell}^j_i|, |m_i|\} \leq c_1\ \ \ \ \text{and}\ \ \ \ \min_i |\tilde{\ell}^1_i - \tilde{\ell}^2_i| \geq 2c_3$$
such that, after passing to an appropriate subsequence, we have $E^{-1}_{V_k,\BC^c_k}\lambda^k_i \to \ell_i$, $E^{-1}_{V_k,\BC_k^c}\tilde{\lambda}^{k,j}_i\to \tilde{\ell}^j_i$, and $E^{-1}_{V_k,\BC^c_k}\mu_i \to m_i$. By (C$_k$) and elliptic estimates for single-valued and two-valued stationary graphs (those seen in Section \ref{sec:two-valued_stationary_graphs} or Theorem \ref{thm:wick1}), we know that there exist $p^{(0)}+2r$ single-valued $C^2$ harmonic functions, $(\phi_i)_{i=1}^{p^{(0)}}$, $(\tilde{\phi}^j_i)_{i=1,\dotsc,r;\;j=1,2}$, and $q^{(0)}-r$ two-valued $C^{1,1/2}$ harmonic functions, $(\psi_i)_{i=1}^{q^{(0)}-r}$, which patched together form a function on $\BC^{(0)}\res B_{3/4}^{n+1}(0)\cap \{|x|>0\}$, such that, after perhaps passing to another subsequence,
$$E_{V_k,\BC_k}^{-1}u_i^k \to \phi,\ \ \ \ E^{-1}_{V_k,\BC_k} \tilde{u}^{k,j}_i\to \tilde{\phi}^j_i,\ \ \ \ \text{and}\ \ \ \ E^{-1}_{V_k,\BC_k}v_i^k\to \psi_i$$
where the convergence is in $C^{1,1/2}(K)$ for each compact subset $K\subset \spt\|\BC^{(0)}\|\cap B_{3/4}^{n+1}(0)\cap \{|x|>0\}$. From (B$_k$) it follows that, in the same way as in Section \ref{sec:coarse-construction} for the construction of the coarse blow-up class, that for each $\sigma\in (0,1/4)$
$$\int_{B_{3/4}^{n+1}(0)}|\phi|^2 + |\tilde{\phi}|^2 + |\psi|^2 \leq C\sigma^{1/2}$$
and moreover that the convergence to $\phi_i$, $\tilde{\phi}_i^j$, and $\psi_i$, is in fact strongly in $L^2$ on $B_{3/4}^{n+1}(0)$; here we have written $\phi = (\phi_1,\dotsc,\phi_{p^{(0)}})$, $\tilde{\phi} = (\tilde{\phi}^1_1, \tilde{\phi}^2_1, \dotsc, \tilde{\phi}^1_r, \tilde{\phi}^2_r)$, and $\psi = (\psi_1,\dotsc,\psi_{q^{(0)}-r})$.

\begin{defn}
	Fix $\BC^{(0)}\in \FL_S\cap \FL_I$, where $I\in \{1,2\}$, and $M>1$. Fix $q<I$, and set $p:= 5-2I$ and $r:= I-q$. Then any triple of functions $(\phi,\tilde{\phi},\psi) \equiv ((\phi)_{i=1}^p, (\tilde{\phi}^j_i)_{i=1,\dotsc,r;\;j=1,2}, (\psi_i)_{i=1}^q)$ constructed as above with $M_1 = M$ for sequences of varifolds $(V_k)_k$, $(\BC^c_k)_k$, $(\BC_k)_k$ obeying $\BC_k\in \FL_q$ for all $k$ is called a \textit{fine blow-up} \textit{of} $(V_k)_k$ \textit{off} $\BC^{(0)}$ \textit{relative to the sequences} $(\BC^c_k)_k$ and $(\BC_k)_k$. We write $\FB^F_{p,q;M}(\BC^{(0)})$ for the collection of all such fine blow-ups with $\BC_k\in \FL_q$ for all $k$.
\end{defn}

\textbf{Remark:} The crucial point to note here is that, since $q<I$, the number of two-valued functions used to describe functions in the fine blow-up class $\FB^F_{p,q;M}(\BC^{(0)})$ is strictly fewer than the number used to describe functions in the coarse blow-up class $\FB(\BC^(0))$.

\subsection{Initial Properties of the Fine Blow-Up Class}\label{sec:fine_initial_properties}

In this section we shall prove initial properties satisfied by the fine blow-up classes $\FB^F_{p,q;M}(\BC^{(0)})$. We will be able to show that they satisfy properties $(\FB1)$, $(\FB2)$, $(\FB3)$, $(\FB4)$, and $(\FB6)$ from Section (\ref{sec:coarse_regularity}), as well as a modified version of $(\FB5)$. To be more precise regarding this latter point, we will be able to show that the functions described in $(\FB5)$ do not belong to $\FB^F_{p,q;M}(\BC^{(0)})$, but instead to $\FB^F_{p,q;M_0M}(\BC^{(0)})$, where $M_0 = M_0(n)$ is the constant defined at the start of Section \ref{sec:fine_construction}. As explained in the discussion preceding Theorem \ref{thm:fine-reg}, provided all the other properties in Section \ref{sec:coarse_regularity} hold for the classes $(\FB^F_{p,q;\tilde{M}}(\BC^{(0)}))_{\tilde{M}>1}$, this is enough to deduce the boundary regularity of functions in each $\FB^F_{p,q;M}(\BC^{(0)})$. It should be noted that of course when $(p,q) = (5,0)$, the situation is much simpler and the regularity conclusions follow from $(\FB3)$, as the functions in the blow-up class $\FB^F_{5,0;M}(\BC^{(0)})$ consist of single-valued harmonic functions, for which the boundary regularity will follow from standard elliptic theory once the regularity of the boundary values and continuity at the boundary is established, and so the only reason for this additional care is that when $(p,q) = (3,1)$, the functions in the blow-up class $\FB^F_{3,1;M}(\BC^{(0)})$ contain a two-valued function for which the boundary regularity theory is more involved.

Thus, let us now fix $\BC^{(0)}\in \FL_S\cap \FL_I$, where $I\in \{1,2\}$, as well as non-negative integers $p,q$ obeying $q<I$ and $p+2q=5$. For $(\phi,\tilde{\phi},\psi)\in \FB^F_{p,q;M}(\BC^{(0)})$, let us write $(V_k)_k$, $(\BC_k^c)$, $(\BC_k)_k$, $(\epsilon_k)_k$, $(\gamma_k)_k$, $(\beta_k)_k$, $(\tau_k)_k$, and $(\delta_k)_k$ for the sequences generating $(\phi,\tilde{\phi},\psi)$ as described in Section \ref{sec:fine_properties}. 

Note that $(\FB1)$ and $(\FB2)$ hold simply from the discussion in Section \ref{sec:fine_properties}; moreover, note $(\FB6)$ follows by essentially the same diagonal argument used in Section \ref{sec:coarse_properties} to prove that $(\FB6)$ held for the coarse blow-up class, and so we do not repeat it here. Also, again $(\FB4)$ will follow from passing in the limit in (G$_k$) and applying $(\FB5\text{II})$ (once we know its validity), as the inequality from the (G$_k$) for the function in $(\FB5\text{II})$, which will lie in $\FB^F_{p,q;M_0M}(\BC^{(0)})$, is exactly what we want for $v\in \FB^F_{p,q;M_0M}(\BC^{(0)})$; so, once again we are left with establishing $(\FB3)$ and $(\FB5)$.

Let us know look at establishing the variant of $(\FB5\text{I})$. So fixing $v\in \FB^F_{p,q;M}(\BC^{(0)})$, $Z\in S(\BC^{(0)})\cap B_{3/8}$ (one can work in $S(\BC^{(0)})\cap B_1$ simply by scaling), and $\sigma\in (0,1/2)$, it suffices (by the same argument as in Section \ref{sec:coarse_properties} for the coarse blow-up class) to show that we may take a fine blow-up of the sequence $\tilde{V}_k:= (\eta_{Z,\sigma})_\# V_k$ relative to the same generating sequences $(\BC^c_k)_k$ and $(\BC_k)_k$ (changes to the sequences $\epsilon_k,\gamma_k,\beta_k$ are irrelevant). The argument for this is identical to that seen in Corollary \ref{cor:fine_estimates} (in fact simpler, as here every cone $\BC^c_k, \BC_k$, is invariant under translations by $Z$). As such, we see that $\tilde{V}_k\in \CN_{\tilde{\epsilon}_k}(\BC^{(0)})$ for some $\tilde{\epsilon}_k\downarrow 0$, from (\ref{E:fine-est-3}) that $Q_{\tilde{V}_k,\BC_k}^2 \leq \tilde{\gamma}_k E_{\tilde{V}_k,\BC_k}^2$ for some $\tilde{\gamma}_k\downarrow 0$, from (\ref{E:fine-est-extra2}) that $E^2_{\tilde{V}_k,\BC_k^c} \leq M_0 M\cdot\inf_{\tilde{\BC}\in \FL_I}E^2_{\tilde{V}_k,\tilde{\BC}}$, and (\ref{E:fine-est-5}) that (when Hypothesis ($\dagger$)(ii) holds for all $k$) $Q_{\tilde{V}_k,\BC_k}^2 < \tilde{\beta}_k(Q^*_{\tilde{V}_k})^2$, for some $\tilde{\beta}_k\downarrow 0$; as such we may perform a fine blow-up of $(\tilde{V}_k)_k$ relative to $(\BC^c_k)_k$ and $(\BC_k)_k$ to see that $v_{Z,\sigma}\in \FB^F_{p,q;M_0M}(\BC^{(0)})$, as desired.

For $(\FB3)$, note that for each $Y\in S(\BC)\cap B_{3/8}$, one may apply Lemma \ref{lemma:gaps} to deduce the existence of $Z_k = (\xi_k,\zeta_k)\in \spt\|V_k\|\cap B_{3/8}$ obeying $\Theta_{V_k}(Z_k)\geq 5/2$ and $Z_k\to Y$. In particular, (D$_k$) tells us that $|\xi_k|\leq CE_{V_k,\BC_k}$, and hence we deduce, after passing to a subsequence, the existence of a limit $E^{-1}_{V_k,\BC_k}\xi_k\to \kappa(Y)$; we will see momentarily that this limit is independent of the approximation sequence $(Z_k)_k$ and so only depends on $Y$. The only caveat now in this setting when compared to that for the coarse blow-up is that in (D$_k$) the projections of $\xi$ are projections of $\xi$ onto the half-hyperplanes in the $(\BC_k)_k$ sequence, and not onto those in $(\BC^c_k)_k$ (or $\BC^{(0)}$), and thus writing $\kappa^\perp$ in the integral over $\spt\|\BC^{(0)}\|$ is now misleading, as the value of $\kappa^\perp$ can (and indeed will) differ between $\tilde{\phi}^1$ and $\tilde{\phi}^2$ on each half-hyperplane which splits. Indeed, recall from Remark 4 and (\ref{E:fine-est-extra8}) that we have
$$\xi^{\perp_{\tilde{H}^{k,j}_i}} = \xi^{\perp_{\tilde{G}^c_{i,k}}} - \tilde{\lambda}^{k,j}_i\xi^{\top_{\tilde{G}^c_{i,k}}}$$
and 
$$\left|\xi^{\perp_{\tilde{G}^c_{i,k}}}\right|^2 + E_{V_k,\BC_k^c}^2\left|\xi^{\top_{\tilde{G}^c_{i,k}}}\right|^2 \leq CE_{V_k,\BC_k}^2$$
and thus we see that we may pass to a subsequence so that for each $i,j$ we have $E_{V_k,\BC_k}^{-1}\xi^{\perp_{\tilde{G}^c_{i,k}}}\to \kappa_{i}^j(Y)$ and $E_{V_k,\BC_k}^{-1}E_{V_k,\BC_k^c}\xi^{\top_{\tilde{G}^c_{i,k}}}\to \tilde{\kappa}_i^j(Y)$, and thus
$$E_{V_k,\BC_k}^{-1}\xi^{\perp_{\tilde{H}^{k,j}_i}} \to \kappa_{i}^j(Y) - \tilde{\ell}_i^j\tilde{\kappa}_i^j(Y) =:\Lambda^j_i(Y).$$
In particular, note that by Remark 3,
$$|\Lambda^1_i(Y)- \Lambda^2_i(Y)| = |\tilde{\kappa}_i^j(Y)|\cdot|\tilde{\ell}^1_i - \tilde{\ell}^2_i| \geq 2c_3|\tilde{\kappa}^i_j(Y)|$$
and thus these will differ whenever $\tilde{\kappa}^i_j(Y)|\neq 0$. Thus, if we denote the boundary values by a function $\Lambda$, by applying (F$_k$) with $Z = Z_k$, we deduce, in the same way as we did in Section \ref{sec:coarse-construction}, that
$$\int_{B_{\rho/2}(Y)\cap \spt\|\BC^{(0)}\|}\frac{|\Phi-\Lambda(Y)|^2}{|X-Z|^{n+3/2}}\ \ext\H^n \leq C\rho^{-n-3/2}\int_{B_\rho(Y)\cap \spt\|\BC^{(0)}\|}|\Phi - \kappa^\perp(Y)|^2$$
where $\Phi$ is the function on $\spt\|\BC^{(0)}\|$ determined by $(\phi,\tilde{\phi},\psi)$ and the value of $\Lambda$ depends on which function $\Phi$ takes; in particular, finiteness of the above integral and the fact that the unit vectors in the directions of the rays in $\spt\|\BC^{(0)}_0\|$ span all of $\R^2$ is what provides that $\Lambda(Y)$ only depends on $Y$ and not the approximating sequence $(Z_k)_k$. We also get, as before, that $\sup_{B_{5/16}\cap S(\BC^{(0)})}|\Lambda|^2 \leq C\int_{B_{1/2}}|\Phi|^2$, where $C = C(n)$; this is the first half of $(\FB3)$. All that remains is to show smoothness of $\Lambda$. As explained in Section \ref{sec:coarse_properties}, the above integral inequality is enough to deduce a (potentially multi-valued) $C^{0,\alpha}$ Campanato estimate for suitable $\alpha\in (0,1)$, and thus proves that $\Phi$ is in fact $C^{0,\alpha}$ up-to-the-boundary. 

To prove that $\Lambda$ is smooth, we can follow a similar argument as to that seen in Section \ref{sec:coarse_properties}, except now making use of Theorem \ref{thm:fine_representation}(a), (c) instead of the coarse estimates in Lemma \ref{lemma:L2_coarse}. The only slight difference is that one now needs to include the gradient functions $h^k_i,\tilde{g}^{k,j}_i$, and $g^k_i$ in the argument when passing from integrals with respect to $\ext\|V_k\|$ to over the half-hyperplanes in $\BC^c_k$. This is dealt with in an analogous (and in fact simpler, as we do not need to consider different variations as the rays of the cross-section still span $\BC^{(0)}$) way to the calculations seen in \cite[(12.17) -- (12.22)]{wickstable}, and so we do not include the argument here. Thus we have established $(\FB3)$.

The final property left is $(\FB5\text{II})$. This however follows in an identical fashion to that seen in the coarse blow-up setting: one modifies the sequence of cones $\BC_k$ based on the function $\psi$ (notation as in Section \ref{sec:coarse_properties}), now by a factor of $E_{V_k,\BC_k}$ and takes a fine blow-up relative to the sequence of modified cones. Once again, we must check that suitable forms of ($1_k$) -- ($5_k$) hold for these new sequences, but these can be checked in the same manner as we have seen already, so we omit the details.

Hence we see that $\FB^F_{p,q;M}(\BC^{(0)})$ always obeys $(\FB1) - (\FB6)$ (with $(\FB5)$ suitably modified as discussed). In particular, when $(p,q) = (5,0)$, we are able to immediately deduce the boundary regularity of the functions in $\FB_{5,0;M}^F(\BC^{(0)})$, and so we deduce:

\begin{prop}\label{prop:fine}
	The conclusions of Theorem \ref{thm:fine-reg} hold for $\FB^F_{5,0;M}(\BC^{(0)})$ whenever $\BC^{(0)}\in \FL_S\cap (\FL_1\cup \FL_2)$.
\end{prop}

\section{The Fine $\epsilon$-Regularity Theorem}\label{sec:L1fine}

The aim of this section is to prove two $\epsilon$-regularity results, one at the varifold level and the other at the coarse blow-up level. The key result is the one at the varifold level, which will be referred to as a \textit{fine} $\epsilon$-regularity theorem. The fine $\epsilon$-regularity theorem will serve two purposes for us. The first purpose will be to deduce the second $\epsilon$-regularity result of this section, namely to prove that $(\FB7)$ holds for any coarse blow-up class $\FB(\BC^{(0)})$ when $\BC^{(0)}\in \FL_S\cap \FL_1$ is a \textit{level 1} cone; in particular, we can then deduce that Theorem \ref{thm:coarse_reg} holds for $\FB(\BC^{(0)})$, and so we have the desired boundary regularity of the coarse blow-ups relative to level 1 cones. Armed now with this, the second purpose will be to prove (which will be done in Section \ref{sec:L1coarse_decay}) that Theorem \ref{thm:A} holds whenever $\BC^{(0)}\in \FL_S\cap \FL_1$ is a level one cone; exactly how this works will be discussed in Section \ref{sec:L1coarse_decay}.

We start by proving an excess decay result in the setting of the fine blow-up class. We note that this lemma is also true in the setting when $\BC^{(0)}\in \FL_2$ and Hypothesis ($\dagger$)(ii) holds (this will be used later to prove Theorem \ref{thm:A} when $\BC^{(0)}\in \FL_S\cap \FL_2$).

\begin{lemma}[Fine Excess Decay: Level 0]\label{lemma:fine_excess_decay}
	Let $\BC^{(0)}\in \FL_S\cap \FL_I$, where $I\in \{1,2\}$. Fix $\theta\in (0,1/4)$. Then, there exist numbers $\epsilon_2 = \epsilon_2(\BC^{(0)},\theta)\in (0,1/2)$, $\gamma_2 = \gamma_2(\BC^{(0)},\theta)\in (0,1/2)$, and $\beta_2 = \beta_2(\BC^{(0)},\theta)\in (0,1/2)$ such that the following is true: if $V\in \S_2$, $\BC^c\in \FL_I$, and $\BC\in \FL_0$ satisfy Hypothesis (H), Hypothesis ($\star$), and Hypothesis ($\dagger$)\footnote{As we are assuming $\BC\in \FL_0$ here, when $I=2$ we are implicitly assuming that Hypothesis ($\dagger$)(ii) holds.} with $\epsilon_2,\gamma_2,\beta_2$, and $\frac{3}{2}M_0$, in place of $\epsilon,\gamma,\beta$, and $M$, respectively, then there exists an orthogonal rotation $\Gamma$ of $\R^{n+1}$ and a cone $\BC^\prime\in \FL_0$ such that the following hold:
	\begin{enumerate}
		\item [(a)] $|\Gamma-\id|\leq \kappa E_{V,\BC}$;
		\item [(b)] $\dist^2_\H(\spt\|\BC\|\cap B_1,\spt\|\BC^\prime\|\cap B_1)\leq \kappa E_{V,\BC}^2$;
		\item [(c)] $$\theta^{-n-2}\int_{B_\theta}\dist^2(X,\spt\|\Gamma_\#\BC^\prime\|)\ \ext\|V\| + \theta^{-n-2}\int_{\Gamma(B_{\theta/2}\backslash\{|x|<\theta/16\})}\dist^2(X,\spt\|V\|)\ \ext\|\Gamma_\#\BC^\prime\| \leq \kappa\theta^2 E_{V,\BC}^2;$$
		\item [(d)] For any $\tilde{\BC}\in \FL_I$ with $\tilde{\BC}\in \FL_{1/10}(\BC^{c})$, we have:
		$$\left(\theta^{-n-2}\int_{B_\theta}\dist^2(X,\spt\|\tilde{\BC}\|)\ \ext\|\Gamma^{-1}_\#V\|\right)^{1/2} \geq \sqrt{2^{-n-4}\bar{C}_1}\dist_\H(\spt\|\BC\|\cap B_1,\spt\|\tilde{\BC}\|\cap B_1) - \kappa E_{V,\BC};$$
	\end{enumerate}
	here, $\kappa = \kappa(n)$, and $\bar{C}_1 = \bar{C}_1(n):= \int_{B^n_{1/2}\cap\{x^2>1/16\}}|x^2|^2\ \ext\H^n(x^2,y)$ is as before.
\end{lemma}

\textbf{Remark:} The only unfamiliar property here from what we have seen before is (d). This will be used to verify Hypothesis $(\star)$ still holds in the proof of the fine $\epsilon$-regularity theorem later.

\begin{proof}
	The proof will be similar to the excess decay lemma from Section \ref{sec:L0excess}, namely Lemma \ref{lemma:ed-level-0}, except now we need to take more care in verifying that the conditions required to perform a fine blow-up are still satisfied when we take appropriate rotations of our varifolds.
	
	We again argue by contradiction; so suppose that the lemma does not hold (for $\kappa = \kappa(n)$ to be chosen): therefore we may find sequences $\epsilon_k,\gamma_k,\beta_k\downarrow 0$, $V_k$, $\BC_k^c$, and $\BC_k\in \FL_0$ satisfying Hypothesis (H), Hypothesis $(\star)$, and Hypothesis $(\dagger)$ with $\epsilon_k,\gamma_k,\beta_k$, and $\frac{3}{2}M_0$, in place of $\epsilon,\gamma,\beta$, and $M$, respectively (i.e. ($1_k$) -- ($5_k$) from Section \ref{sec:fine_properties}), such that the lemma does not hold for this choice of $\theta$ and $\BC^{(0)}$. We need to show that all conclusions of the lemma are satisfied for infinitely many $k$.
	
	For $i=1,\dotsc,n-1$, let $Y_i:= \frac{1}{2}\theta e_{2+i}\in S(\BC^{(0)})$. Lemma \ref{lemma:gaps} tells us that for each $k\geq 1$ and $i\in \{1,\dotsc,n-1\}$ we may find sequences $Z_{i,k} = (\xi_{i,k},\zeta_{i,k})\in \spt\|V_k\|\cap B_1$ such that $\Theta_{V_k}(Z_{i,k})\geq 5/2$ and $Z_{i,k}\to Y_i$. As in Lemma \ref{lemma:ed-level-0}, we may assume without loss of generality that $\{Z_{1,k},\dotsc,Z_{n-1,k}\}$ span an $(n-1)$-dimensional subspace of $\R^{n+1}$, which we call $\Sigma_k$. We may then choose rotations $\Gamma^\prime_k$ of $\R^{n+1}$ such that $\Gamma^\prime_k(\Sigma_k) = S(\BC^{(0)})\equiv\{0\}^2\times\R^{n-1}$ and $\Gamma^\prime_k(\hat{Z}_{i,k}) \to e_{2+i}$, where $\hat{Z}_{i,k}:= Z_{i,k}/Z_{i,k}$, for each $i=1,\dotsc,n-1$. So far, this puts little restriction on how rotates the cross-section $\R^2\times\{0\}^{n-1}$, and so to ensure $|\Gamma-\id|$ is small we need to reset any significant change in this subspace. Thus, if $\pi_{12}:\R^{n+1}\to \R^2\times\{0\}^{n-1}$ is the orthogonal projection, choose a rotation $\Gamma+k^{\prime\prime}$ of $\R^{n+1}$ such that $\left.\Gamma_k^{\prime\prime}\right|_{\{0\}^2\times\R^{n-1})} = \id$ and $\Gamma_k^{\prime\prime}(\pi_{12}\Gamma_k^\prime(e_1)/|\pi_{12}\Gamma_k^\prime(e_1)|) = e_1$. Now, if we set $\Gamma_k := \Gamma^{\prime\prime}_k\circ\Gamma_k^\prime$, we have
	\begin{equation}\label{E:fine-decay-0a}
	\Gamma_k(\Sigma_k) = \{0\}^2\times\R^{n-1}\ \ \ \ \text{and}\ \ \ \ \Gamma_k(\hat{Z}_{i,k}) \to e_{2+i}\ \ \text{for }i\in\{1,\dotsc,n-1\}.
	\end{equation}
	Moreover, as $|\xi_{i,k}|\leq CE_{V_k,\BC_k}$ by Corollary \ref{cor:fine_estimates}, we have
	\begin{equation}\label{E:fine-decay-0}
	|\Gamma_k-\id| \leq CE_{V_k,\BC_k}.
	\end{equation}
	Now set $\tilde{V}_k:= (\eta_{0,7/8}\circ\Gamma_k)_\#V_k$. It readily follows that for any $\BC\in \FL$ that
	\begin{equation}\label{E:fine-decay-1}
	\dist_\H(\spt\|(\Gamma_k^{-1})_\#\BC\|\cap B_1,\spt\|\BC\|\cap B_1) \leq CE_{V_k,\BC_k}
	\end{equation}
	from which it immediately follows, from the triangle inequality in the form $\dist^2(X,\spt\|(\Gamma_k^{-1})_\#\BC\|) \leq 2\dist^2(X,\spt\|\BC\|) + 2\dist_\H^2(\spt\|(\Gamma_k^{-1})_\#\BC\|\cap B_1,\spt\|\BC\|\cap B_1)$ for suitable choices of $\BC$, that
	\begin{equation}\label{E:fine-decay-2}
	E_{\tilde{V}_k,\BC_k}^2 \leq CE^2_{V_k,\BC_k}\ \ \ \ \text{and}\ \ \ \ E^2_{\tilde{V}_k,\BC_k^c}\leq CE^2_{V_k,\BC^c_k}.
	\end{equation}
	where for the second inequality we have used the fact that $E^2_{V_k,\BC_k}\leq Q^2_{V_k,\BC_k} \leq \gamma_k E^2_{V_k,\BC_k^c}$; here, $C = C(n)$. We claim further that, for some $\tilde{C} = \tilde{C}(n)$,
	\begin{equation}\label{E:fine-decay-3}
	\tilde{C}E_{V_k,\BC_k^c}^2 \leq E^2_{\tilde{V}_k,\BC^c_k}.
	\end{equation}
	To see this, note that the coarse blow-up, $v$, of $(V_k)_k$ relative to $(\BC_k^c)_k$ is homogeneous of degree one (indeed, its graph is in fact a union of half-hyperplanes meeting along $\{0\}^2\times\R^{n-1}$); this is simply because $E_{V_k,\BC_k^c}^{-2}Q_{V_k,\BC_k}^2<\gamma_k\to 0$. Moreover, by Remark 2, there is a definitive constant $c = c(n)$ for which $\int_{B_1}|v|^2\geq c$; thus as $v$ is homogeneous of degree one, we have $\sigma^{-n-2}\int_{B_\sigma}|v|^2 = \int_{B_1}|v|^2 \geq c$ for each $\sigma\in (0,1)$. Thus we have for all $k$ sufficiently large, using Corollary \ref{corollary:non-concentration} (recall $\tau_k\downarrow 0$ is suitably chosen)
	\begin{align*}
		\int_{B_\sigma}\dist^2(X,\spt\|\BC^c_k\|)\ \ext\|V_k\| & \geq \int_{B_\sigma\backslash\{|x|<\tau_k\}}\dist^2(X,\spt\|\BC^c_k\|)\ \ext\|V_k\|\\
		& \geq \int_{B_\sigma\cap\{|x|>\tau_k\}}|u_k|^2\\
		& \geq \int_{B_\sigma}|u_k|^2 - C\tau_k^{1/2}E_{V_k,\BC_k}^2\\
		& \geq \frac{1}{2}\sigma^{n+2}c E_{V_k,\BC_k}^2 - C\tau_k^{1/2}E_{V_k,\BC_k}^2
	\end{align*}
	which tells us, for sufficiently large $k$ (depending on $\sigma$), we have $\int_{B_\sigma}\dist^2(X,\spt\|\BC_k\|)\ \ext\|V_k\| \geq \frac{c}{4}\sigma^{n+2}E^2_{V_k,\BC_k}$. In particular, using (\ref{E:fine-decay-1}), the triangle inequality, and Hypothesis (H4), we deduce for a suitable choice of $\sigma = \sigma(n)\in (0,1)$ that (\ref{E:fine-decay-3}) holds for suitable $\tilde{C} = \tilde{C}(n)$.
	
	We now claim that, after passing to a subsequence, that Hypothesis ($1_k$) -- ($5_k$) of Section \ref{sec:fine_properties} are satisfied with $\tilde{V}_k$ in place of $V_k$ (and keeping the same $\BC^c_k$, $\BC_k$), for suitable sequences $\tilde{\epsilon}_k,\tilde{\gamma}_k$, and $\tilde{\beta}_k\to 0$ in place of $\epsilon_k,\gamma_k$, and $\beta_k$, respectively and with $\frac{3}{2}M_0^2$ in place of $M_1$. Of course ($2_k$) still holds as we have not changed the cone sequences, and ($1_k$) follows from the second inequality in (\ref{E:fine-decay-2}). For ($3_k$), note that by the first inequality in (\ref{E:fine-decay-2}) and (\ref{E:fine-decay-3}) we have $E^2_{\tilde{V}_k,\BC_k} \leq C^\prime\gamma_k E^2_{\tilde{V}_k,\BC_k^c}$, which is one half of ($3_k$); the other half is dealt with similarly to what we have seen, using (\ref{E:fine-decay-0}) and Theorem \ref{thm:fine_representation}.
	
	Let us now look at ($4_k$). First note that, just as argued before in (\ref{E:fine-est-3}), for some $C = C(n)$ we have
	$$\inf_{\tilde{\BC}}E^2_{\tilde{V}_k,\tilde{\BC}} = \inf_{\tilde{\BC}\in\FL_I\cap \FL_{CE_{\tilde{V}_k,\BC_k^c}}(\BC^c_k)}E^2_{\tilde{V}_k,\tilde{\BC}}$$
	and thus for $\tilde{\BC}\in \FL_I$ as in the infimum on the right hand side, for all sufficiently large $k$ we have:
	\begin{align*}
		E_{\tilde{V}_k,\tilde{\BC}}^2 & \geq (7/8)^{-n-2}\int_{B_{1/2}}\dist^2(X,\spt\|(\Gamma_k^{-1})_\#\spt\|\tilde{\BC}\|)\ \ext\|V_k\|\\
		& \geq \frac{1}{2}(7/8)^{-n-2}\int_{B_{1/2}}\dist^2(X,\spt\|\tilde{\BC}\|)\ \ext\|V_k\|\\
		& \hspace{12em} - (7/8)^{-n-2}\cdot 3\w_n\cdot \dist^2_\H(\spt\|(\Gamma_k^{-1})_\#\tilde{\BC}\|\cap B_1,\spt\|\tilde{\BC}\|\cap B_1)\\
		& \geq \frac{1}{2}(7/8)^{-n-2}\cdot 2^{-n-9}\w_n^{-1} \bar{C}_1 E^2_{V_k,\tilde{\BC}} - CE_{V_k,\BC_k}^2\\
		& \geq \frac{1}{2}(7/8)^{-n-2}\cdot 2^{-n-9}\w_n^{-1}\bar{C}_1 M^{-1}E_{V_k,\BC_k^c}^2 - C\gamma_k E_{V_k,\BC_k^c}^2\\
		& \geq 2^{-n-11}\w_n^{-1}\bar{C}_1 M^{-1}E_{V_k,\BC^c_k}^2
	\end{align*}
	where for third inequality, we have used (\ref{E:fine-est-extra-1a}) (with $\rho=1/2$ and $Z= 0$) to lower bound the first term and then used (\ref{E:fine-decay-1}) to lower bound the second term. Then we also have
	\begin{align*}
		E^2_{\tilde{V}_k,\BC^c_k} & = (7/8)^{-n-2}\int_{B_{7/8}}\dist^2(X,\spt\|(\Gamma_k^{-1})_\#\BC_k^c\|)\ \ext\|V_k\|\\
		& \leq 2(7/8)^{-n-2}\int_{B_1}\dist^2(X,\spt\|\BC^c_k\|)\ \ext\|V_k\|\\
		& \hspace{12em} + 2(7/8)^{-n-2}\cdot 3\w_n\cdot \dist^2_\H(\spt\|(\Gamma_k^{-1})_\#\BC_k^c\|\cap B_1,\spt\|\BC^c_k\|\cap B_1)\\
		& \leq 2^{n+3}E_{V_k,\BC_k^c}^2 + CE^2_{V_k,\BC_k}\\
		& \leq 2^{n+4}E_{V_k,\BC^c_k}^2
	\end{align*}
	for all $k$ sufficiently large, again using (\ref{E:fine-decay-0}) and Hypothesis (H4). Thus combining the above two inequalities, we see that, for all $k$ sufficiently large,
	$$E^2_{\tilde{V}_k,\BC^c_k} \leq 2^{2n+14}\w_n\bar{C}_1^{-1}\cdot M\cdot E_{\tilde{V}_k,\tilde{\BC}}^2.$$
	As this constant factor is always at most $M_0$, by definition of $M_0$, this verifies that ($4_k$) holds for $\tilde{V}_k$ (and $\BC^c_k$) for all $k$ sufficiently large, with $M_1 = \frac{3}{2}M_0^2$ (as $M = \frac{3}{2}M_0$ by assumption).
	
	Finally, if $I=1$ then ($5_k$) is automatically satisfied. Otherwise, $I=2$ and $\BC_k\in \FL_0$ and so we need to show for all $k$ sufficiently large that $Q_{\tilde{V}_k,\BC_k}^2 \leq \tilde{\beta}_k (Q^*_{\tilde{V}_k})^2$ for some $\tilde{\beta}_k\downarrow 0$. To see this, it suffices to show that for some $C = C(n)$ we have $Q_{\tilde{V}_k,\BC_k}^2 \leq CQ_{V_k,\BC_k}^2$ and $(Q^*_{V_k})^2\leq C (Q^*_{\tilde{V}_k})^2$, as then we would have $Q_{\tilde{V}_k,\BC_k}^2 \leq C^2\beta_k (Q_{\tilde{V}_k}^*)^2$. Indeed, for the first of these two inequalities, we already have by the first inequality in (\ref{E:fine-decay-2}) that $E^2_{\tilde{V}_k,\BC_k} \leq CQ_{V_k,\BC_k}^2$; the other term in $Q^2_{\tilde{V}_k,\BC_k}$ can be dealt with in a similar way to what we have currently seen by using (\ref{E:fine-decay-1}) and Theorem \ref{thm:fine_representation}. For the second of these inequalities, this follows in a similar manner to (\ref{E:fine-decay-3}), namely by (\ref{E:fine-decay-0}) and showing that there is a $\tilde{c} = \tilde{c}(n)$ such that for any $\tilde{\BC}\in \FL_1\cup \FL_2$,
	$$\int_{B_{7/8}}\dist^2(X,\spt\|\tilde{\BC}\|)\ \ext\|V_k\| + \int_{B_{7/16}\backslash\{|x|<(7/8)/16\}}\dist^2(X,\spt\|V_k\|)\ \ext\|\tilde{\BC}\| \geq \tilde{c}(Q^*_{V_k})^2$$
	which follows by taking cones $\tilde{\BC}_k\in \FL_1\cup \FL_2$ with $Q_{V_k,\tilde{\BC}_k}^2 < \frac{3}{2}(Q_{V_k}^*)^2$ and taking a blow-up (coarse blow-up if $\tilde{\BC}_k\in \FL_2$ for infinitely many $k$, otherwise if $\tilde{\BC}_k\in \FL_1$ for infinitely many $k$, we decide first whether $E^{-1}_{V_k,\BC^c_k}Q_{V_k,\tilde{\BC}_k}< \tilde{\gamma}$ for infinitely many $k$ or not, where $\tilde{\gamma} = \tilde{\gamma}(\BC^{(0)})$ is the constant from Theorem \ref{thm:fine_representation} with $\tau = 7/16^2$ and on the slightly larger ball $B_{15/16}$ rather than $B_{3/4}$: if it does, we can take a fine blow-up relative to $(\BC^c_k)_k$ and $(\tilde{\BC}_k)_k$, in the region $B_{15/16}\cap \{|x|>7/16^2\}$, using ($7_k$) to deduce that the limit $\Phi$ must be homogeneous of degree one and from Remark 3 must obey $\int_{B_1}|\Phi|^2\geq c$ for some $c= c(n)$; otherwise, we have a fixed lower bound $E^{-1}_{V_k,\BC_k^c}Q_{V_k,\tilde{\BC}_k}\geq\tilde{\gamma}$, only depending on $\BC^{(0)}$, and so one may take a coarse blow-up of $V_k$ relative to $\BC^k_c$, and argue in the same way to prove the result.
	
	Thus, we may generate a fine blow-up $\Phi = (\phi,\tilde{\phi})$ from taking the fine blow-up of $(\tilde{V}_k)_k$ relative to $(\BC^c_k)_k$ and $(\BC_k)_k$; since $\Phi\in \FB^F_{5,0;\frac{3}{2}M_0^2}(\BC^{(0)})$, $\Phi$ has no two-valued component $\psi$. It follows by (\ref{E:fine-decay-0a}) and (F$_k$) (applied with $\tilde{V}_k$ in place of $V$ and $\frac{8}{7}\Gamma_k(Z_{i,k})$ in place of $Z$) that $\phi(Y_i) = \tilde{\phi}(Y_i) = 0$ for each $i=1,\dotsc,n-1$ (simply because after rotating by $\Gamma_k$ the points of density $\geq 5/2$ converging to each $Y_i$, i.e. $\Gamma_k(Z_{i,k})$,  have no component in the $\R^2\times\{0\}^{n-1}$ variables); moreover, by translating to assume without loss of generality $\Theta_{V_k}(0)\geq 5/2$ for all $k$ (which can be arranged using Corollary \ref{cor:fine_estimates}), we also have $\phi(0) = \tilde{\phi}(0) = 0$. Since $Y_i = \frac{\theta}{2}e_{2+i}$, from the regularity conclusion along the boundary established in Proposition \ref{prop:fine}, this tells us that there exists points $S_{i,j}$, $\tilde{S}^1_{i,j}$, $\tilde{S}^2_{i,j}\in B_{\theta/2}\cap (\{0\}^2\times\R^{n-1})$ such that for each $i=1,\dotsc,n-1$ and $j$,
	$$\frac{\del\phi_j}{\del y^i}(S_{j,i}) = 0,\ \ \ \ \frac{\del\tilde{\phi}^1_j}{\del y^i}(\tilde{S}^1_{j,i}) = 0,\ \ \ \ \frac{\del\tilde{\phi}^2_j}{\del y^i}(\tilde{S}^2_{j,i}) = 0.$$
	The estimate provided by Proposition \ref{prop:fine} therefore gives that
	$$|D_y\Phi(0)|^2 \leq C\theta^2\int_{B_{1/2}}|\Phi|^2$$
	where $C = C(n)$ (note that we can get a factor of $\theta^2$ here, as our functions in $\Phi$ are single-valued harmonic functions up-to-the-boundary here). Now we define new linear functions just as in the proof of Lemma \ref{lemma:ed-level-0} as follows: if a half-hyperplane in $\spt\|\BC^{(0)}\|$ were $H= \{(x^2,y)\in \R\times\R^{n-1}: x^2>0\}$, and the fine blow-up over this half-hyperplane (or a component of, if there are two) is represented by a function $\varphi$, then we write $L_\varphi(x^2,y):= D\varphi(0)\cdot (x^2,y)$, and $P_\varphi(x^2,y):= \frac{\del\varphi}{\del x^2}(0)x^2$; since $L_\varphi - P_{\varphi} = D_y\varphi$, the above bounds would then tell us
	$$|L_\varphi(x^2,y)-P_\varphi(x^2,y)|^2 \leq C|y|^2\theta^2\int_{B_{1/2}}|\Phi|^2$$
	and consequently, from Proposition \ref{prop:fine}, that
	\begin{equation}\label{E:fine-decay-4}
	\theta^{-n-2}\int_{B_{2\theta}\cap H}|\varphi-P_\varphi|^2\leq C\theta^2.
	\end{equation}
	Thus, if this $\varphi$ component of the cross blow-up was generated by blowing up a function defined over a half-hyperplane $H_\varphi$ in $\BC_k$ which had gradient $\lambda_\varphi$ relative to the corresponding half-hyperplane $H^c_\varphi$ in $\BC^c_k$, we would define a new half-hyperplane, $\tilde{H}^k_\varphi$, which instead has gradient over $H^c_\varphi$ given by
	\begin{equation}\label{E:fine-decay-5}
	\tilde{\lambda}_\varphi:= \lambda_\varphi + E_{\tilde{V}_k,\BC_k}\cdot\frac{\del\varphi}{\del x^2}(0).
	\end{equation}
	Of course, all the above discussion was done under the assumption that $H = \{(x^2,y)\in \R^n:x^2>0\}$ was the half-hyperplane in $\spt\|\BC^{(0)}\|$ in question, but the whole discussion follows through for each half-hyperplane by working in coordinates relative to the given half-hyperplane. Thus, we generate a new sequence of cones, $\BC^\prime_k$, which from the above definition we see will still be level $0$ cones (by Remark 2 and Remark 3 and the bound on the derivatives of the blow-up provided by Proposition \ref{prop:fine}), for which the discussion above readily gives (from (\ref{E:fine-decay-4}) and the strong $L^2$ convergence to the fine blow-up)
	\begin{equation}\label{E:fine-decay-6}
	\theta^{-n-2}\int_{B_\theta}\dist^2(X,\spt\|\BC^\prime_k\|)\ \ext\|(\Gamma_k)_\#V_k\| \leq \left(\frac{7}{8}\right)^{n+2} C\theta^2 E_{\tilde{V}_k,\BC_k}^2.
	\end{equation}
	Applying now the first inequality in (\ref{E:fine-decay-2}) shows that the first term in (c) has the correct bound for infinitely many $k$. Furthermore, by the first inequality in (\ref{E:fine-decay-2}), (\ref{E:fine-decay-5}), as well as the supremum bound on the derivatives at the boundary of the fine blow-up provided by Proposition \ref{prop:fine}, we clearly have
	\begin{equation}\label{E:fine-decay-7}
	\dist_\H(\spt\|\BC_k^\prime\|\cap B_1,\spt\|\BC_k\|\cap B_1) \leq CE_{V_k,\BC_k}
	\end{equation}
	which shows (b) holds for all $k$ sufficiently large. To see the bound on the second term in (c), this follows readily from the graphical representation provided by (C$_k$) in Section \ref{sec:fine_properties}, as it enables us to bound it by the first term in (c), i.e. we have
	$$\hspace{-1em}\theta^{-n-2}\int_{\Gamma^{-1}_k(B_{\theta/2}\backslash\{|x|<\theta/16\}}\dist^2(X,\spt\|V_k\|)\ \ext\|(\Gamma_k^{-1})_\#\BC_k^\prime\| \leq C\theta^{-n-2}\int_{B_\theta}\dist^2(X,\spt\|(\Gamma_k^{-1})_\#\BC^\prime_k\|)\ \ext\|V_k\|$$
	where $C = C(n)$; this shows that (c) holds for infinitely many $k$.
	
	Thus all that is left to show is that (d) holds for infinitely many $k$ to arrive at a contradiction and complete the proof. Indeed, by (C$_k$) again (with $\tilde{V}_k$ in place of $V_k$), as well as the first inequality in (\ref{E:fine-decay-2}) as well as (\ref{E:fine-decay-6}) and (\ref{E:fine-decay-7}), if we set $\tilde{\theta} = 8\theta/7$) and fix any $\tilde{\BC}\in \FL_{1/10}(\BC^c_k)$, writing $(h_i)_{i=1}^{5-2I}$ and $(\tilde{h}^j_i)_{i=1,\dotsc,I;\;j=1,2}$ for the linear functions over $\BC^c_k$ determining $\tilde{\BC}$, we have:
	\begin{align*}
		\tilde{\theta}^{-n-2}&\int_{B_{\tilde{\theta}}}\dist^2(X,\spt\|\tilde{\BC}\|)\ \ext\|\tilde{V}_k\|\\
		& \geq \frac{1}{2}\tilde{\theta}^{-n-2}\sum_{i}\int_{B_{\tilde{\theta}}\cap H^c_{i,k}\cap\{|x|>\tilde{\theta}/16\}}|h^k_i + u^k_i - h_i|^2\\
		& \hspace{13em} + \frac{1}{2}\tilde{\theta}^{-n-2}\sum_{i,j}\int_{B_{\tilde{\theta}/2}\cap \tilde{G}^c_{i,k}\cap \{|x|>\tilde{\theta}/16\}}|\tilde{g}^{k,j}_i + \tilde{u}^{k,j}_i - \tilde{h}^j_i|^2\\
		& \geq \frac{1}{4}\tilde{\theta}^{-n-2}\sum_i\int_{B_{\tilde{\theta}/2}\cap H_{i,k}^c\cap \{|x|>\tilde{\theta}/16\}}|h^k_i - h_i|^2 + \frac{1}{4}\tilde{\theta}^{-n-2}\sum_{i,j}\int_{B_{\tilde{\theta}/2}\cap \tilde{G}^{c}_{i,k}\cap \{|x|>\tilde{\theta}/16\}}|\tilde{g}^{k,j}_i - \tilde{h}^j_i|^2\\
		& \hspace{3em} - \frac{1}{2}\tilde{\theta}^{-n-2}\sum_i\int_{B_{\tilde{\theta}/2}\cap H^c_{i,k}\cap \{|x|>\tilde{\theta}/16\}}|u^k_i|^2 - \frac{1}{2}\tilde{\theta}^{-n-2}\sum_{i,j}\int_{B_{\tilde{\theta}/2}\cap \tilde{G}^{c}_{i,k}\cap \{|x|>\tilde{\theta}/16\}}|\tilde{u}^{k,j}_i|^2\\
		& \geq 2^{-n-4}\bar{C}_1\dist^2_\H(\spt\|\BC_k\|\cap B_1, \spt\|\tilde{\BC}\|\cap B_1) - \frac{1}{2}\tilde{\theta}^{-n-2}\int_{B_{\tilde{\theta}}}\dist^2(X,\spt\|\BC_k\|)\ \ext\|\tilde{V}_k\|\\
		& \geq 2^{-n-4}\bar{C}_1\dist^2_\H(\spt\|\BC_k\|\cap B_1,\spt\|\tilde{\BC}\|\cap B_1) - \tilde{\theta}^{-n-2}\int_{B_{\tilde{\theta}}}\dist^2(X,\spt\|\BC^\prime_k\|)\ \ext\|\tilde{V}_k\| - CE_{\tilde{V}_k,\BC_k}^2\\
		& \geq 2^{-n-4}\bar{C}_1\dist^2_\H(\spt\|\BC_k\|\cap B_1,\spt\|\tilde{\BC}\|\cap B_1) - C\theta^2E_{\tilde{V}_k,\BC_k}^2 -CE_{\tilde{V}_k,\BC_k}^2\\
		& \geq 2^{-n-4}\bar{C}_1\dist^2_\H(\spt\|\BC_k\|\cap B_1,\spt\|\tilde{\BC}\|\cap B_1) - CE_{V_k,\BC_k}^2
	\end{align*}
	where here $C = C(n)$ and $\bar{C}_1 = \int_{B_{1/2}\cap \{x^2>1/16\}}|x^2|^2\ \ext\H^n(x^2,y)$ is the usual constant; of course, we have abused our notation and have written $u^k_i$, $\tilde{u}^{k,j}_i$ for the functions representing $\tilde{V}_k$ in the application of (C$_k$). This readily gives the validity of (d) for all $k$ sufficiently large, and thus the proof is completed.
\end{proof}


Armed now with Lemma \ref{lemma:fine_excess_decay}, we will now be able to prove the first $\epsilon$-regularity result of this section, namely the fine $\epsilon$-regularity theorem for level 1 cones (note that currently we are not able to say anything for level 2 cones, as we currently do not understand the boundary regularity of functions in $\FB_{3,1;M}^F(\BC^{(0)})$).

\begin{theorem}[Varifold Fine $\epsilon$-Regularity Theorem: Level 1 Setting]\label{thm:L1fine_reg}
	Let $\BC^{(0)}\in \FL_S\cap \FL_1$ and $\alpha\in (0,1)$. Then, there exist constants $\epsilon_1 = \epsilon_1(\BC^{(0)},\alpha)\in (0,1)$ and $\gamma_1 = \gamma_1(\BC^{(0)},\alpha)\in (0,1)$ such that the following holds: if $V\in \S_2$, $\BC^c\in \FL_1$, and $\BC\in \FL_0$ are such that $\Theta_V(0)\geq 5/2$, $V\in \CN_{\epsilon_1}(\BC^{(0)})$, $\BC^c,\BC\in \FL_{\epsilon_1}(\BC^{(0)})$, $E^2_{V,\BC^c}< \frac{3}{2}\inf_{\tilde{\BC}\in \FL_1}E^2_{V,\tilde{\BC}}$, and $Q_{V,\BC}^2 < \gamma_1 E_{V,\BC^c}^2$, then there is a cone $\BC^\prime\in \FL_S\cap \FL_0$ with
	$$\dist_\H(\spt\|\BC^\prime\|\cap B_1,\spt\|\BC\|\cap B_1) \leq CQ_{V,\BC}$$
	and an orthogonal rotation $\Gamma:\R^{n+1}\to \R^{n+1}$ with $|\Gamma-\id|\leq CQ_{V,\BC}$, such that $\BC^\prime$ is the unique tangent cone to $\Gamma^{-1}_\#V$ at $0$, and
	$$\sigma^{-n-2}\int_{B_\sigma}\dist^2(X,\spt\|\BC^\prime\|)\ \ext\|\Gamma^{-1}_\#V\|\leq C\sigma^{2\alpha}Q_{V,\BC}^2\ \ \ \ \text{for all }\sigma\in (0,1/2)$$
	and furthermore, $V$ has the structure of a $C^{1,\alpha}$ classical singularity of vertex density $5/2$; more precisely, there is a $C^{1,\alpha}$ function $u$ defined over $\spt\|\BC^c\|$, in the manner described in Theorem \ref{thm:A}, obeying $V\res B_{1/2} = \mathbf{v}(u)$ and over the multiplicity two half-hyperplane in $\BC^c$, $u$ is expressible as two (disjoint) $C^{1,\alpha}$ single-valued functions; thus, $V\res B_{1/2}$ has no (density 2) branch points, and $\sing(V)\cap B_{1/2} = \{\Theta_V = 5/2\}\cap B_{1/2}$ is the set of points determined by the boundary values of $u$. Here, $C = C(n)$.
\end{theorem}

\begin{proof}
	Let $\kappa = \kappa(n)$ be the constant from Lemma \ref{lemma:fine_excess_decay}. Then, choose $\theta = \theta(n)\in (0,1/4)$ such that $\kappa\theta^{2(1-\alpha)} < 1$. Now let $\epsilon_2 = \epsilon_2(\BC^{(0)},\theta) = \epsilon_2(\BC^{(0)},\alpha)$ and $\gamma_2 = \gamma_2(\BC^{(0)},\theta) = \gamma_2(\BC^{(0)},\alpha)$ be the constants from Lemma \ref{lemma:fine_excess_decay} with this choice of $\theta$ (note that we have no constant $\beta_2$ in this situation as $\BC^{(0)})\in \FL_1$ and so Hypothesis ($\dagger$)(i) is satisfied). Now fix $\epsilon_1\in (0,\epsilon_2)$ and $\gamma_1\in (0,\gamma_2)$; these will eventually be chosen depending only on $\BC^{(0)}$.
	
	Suppose that the hypotheses of the theorem hold with $\epsilon_1$ and $\gamma_1$. For the sake of brevity in our notation, let us write
	$$\hspace{-1em}Q_{V,\BC}(\Gamma,\rho)^2:= \rho^{-n-2}\int_{B_\rho}\dist^2(X,\spt\|\Gamma_\#\BC\|)\ \ext\|V\| + \rho^{-n-2}\int_{\Gamma(B_{\rho/2}\backslash\{|x|<\rho/16\})}\dist^2(X,\spt\|V\|)\ \ext\|\Gamma_\#\BC\|$$
	i.e. $Q_{V,\BC}(\Gamma,\rho)^2 \equiv Q_{(\eta_{0,\rho}\circ\Gamma^{-1})_\#V,\BC}^2$.
	We first claim that we can apply Lemma \ref{lemma:fine_excess_decay} iteratively to obtain sequences of orthogonal rotations $\Gamma_k:\R^{n+1}\to \R^{n+1}$ and cones $\BC_k\in \FL_0$ with $\Gamma_0 = \id$, $\BC_0 = \BC$, and $\BC_k\in \FL_0$, such that
	\begin{equation}\label{E:fine-reg-1}
		|\Gamma_k-\Gamma_{k-1}|^2 \leq C\theta^{2k\alpha}Q^2_{V,\BC};
	\end{equation}
	\begin{equation}\label{E:fine-reg-2}
		\dist^2_\H(\spt\|\BC_k\|\cap B_1,\spt\|\BC_{k-1}\|\cap B_1) \leq C\theta^{2k\alpha}Q_{V,\BC}^2;
	\end{equation}
	\begin{equation}\label{E:fine-reg-3}
		Q^2_{V,\BC_k}(\Gamma_k,\theta^k) \leq 4^{-1}\theta^{2\alpha}Q^2_{V,\BC_{k-1}}(\Gamma_{k-1},\theta^{k-1}) \leq \cdots \leq \theta^{2k\alpha}Q_{V,\BC}^2;
	\end{equation}
	and such that for all $\tilde{\BC}\in \FL_1\cap \FL_{1/10}(\BC^c)$ we have
	\begin{equation}\label{E:fine-reg-4}
		\begin{split}
			\left((\theta^k)^{-n-2}\int_{B_{\theta^k}}\dist^2(X,\spt\|\tilde{\BC}\|)\ \ext\|(\Gamma^{-1}_k)_\#V\|\right)^{1/2} &\\
			& \hspace{-17em} \geq \sqrt{2^{-n-4}\bar{C}_1} \dist_\H(\spt\|\BC_{k-1}\|\cap B_1,\spt\|\tilde{\BC}\|\cap B_1) - \tilde{\kappa}Q_{V,\BC_{k-1}}(\Gamma_{k-1},\theta^{k-1});
		\end{split}
	\end{equation}
	here, $\tilde{\kappa} = \tilde{\kappa}(n)\in 0,\infty)$ and $C = C(n)$. The verification of these will be similar to that seen in the proof of Theorem \ref{thm:A} for level 0 cones from the (coarse) excess decay lemma (Lemma \ref{lemma:ed-level-0}) we saw in Section \ref{sec:L0excess}. Note that for this choice of $\epsilon_1$ and $\gamma_1$, we may apply directly Lemma \ref{lemma:fine_excess_decay} to $V$, $\BC^c$, and $\BC$ to see that properties (\ref{E:fine-reg-1}) -- (\ref{E:fine-reg-4}) hold for $k=1$ (also note that, by properties of multiplicity two classes we still have $(\w_n\theta^n)^{-1}\|V\|(B_\theta) < 5/2 + 1/8$). So now let us suppose that $k\geq 2$ and that (\ref{E:fine-reg-1}) -- (\ref{E:fine-reg-4}) hold for $1,2,\dotsc,k-1$. We wish to apply Lemma \ref{lemma:fine_excess_decay} with $V_{k-1}:= (\eta_{0,\theta^{k-1}}\circ\Gamma_{k-1}^{-1})_\#V$ and $\BC_{k-1}$ in place of $V$ and $\BC$, respectively (with the same $\BC^c$), as this would then establish the validity of (\ref{E:fine-reg-1}) -- (\ref{E:fine-reg-4}) for $k$. Let us write $\theta_k:= \theta^k$.
	
	To begin, firstly note that simply by the triangle inequality, and the fact that $(\w_n\theta_k^n)^{-1}\|V\|(B_{\theta_k}) < 5/2 + 1/8$, we have
\begin{equation}\label{E:fine-reg-5}
\begin{split}
E^2_{V_{k-1},\BC^c} & = \theta_{k-1}^{-n-2}\int_{B_{\theta_{k-1}}}\dist^2(X,\spt\|\BC^c\|)\ \ext\|(\Gamma^{-1}_{k-1})_\#V\|\\
& \leq 2\theta_{k-1}^{-n-2}\int_{B_{\theta_{k-1}}}\dist^2(X,\spt\|\BC_{k-1}\|)\ \ext\|(\Gamma^{-1}_{k-1})_\#V\|\\
& \hspace{15em} + 6\w_n\dist^2_\H(\spt\|\BC_{k-1}\|\cap B_1,\spt\|\BC^c\|\cap B_1).
\end{split}
\end{equation}
Now, if one applies the validity of (\ref{E:fine-reg-2}) with $1,2,\dotsc,k-1$ in place of $k$ and the triangle inequality, we get (noting that $\sum^{k-1}_{i=1}(\theta^\alpha)^i = (\theta^\alpha - (\theta^\alpha)^k)/(1-\theta^\alpha) \leq 4^{-\alpha}(1-4^{-\alpha})^{-1}$)
\begin{align*}
	\dist_\H(\spt\|\BC_{k-1}\|&\cap B_1, \spt\|\BC^c\|\cap B_1)\\
	& \leq \dist_\H(\spt\|\BC_0\|\cap B_1,\spt\|\BC^c\|\cap B_1) + \sum^{k-1}_{i=1}\dist_\H(\spt\|\BC_i\|\cap B_1,\spt\|\BC_{i-1}\|\cap B_1)\\
	& \leq \dist_\H(\spt\|\BC_0\|\cap B_1,\spt\|\BC^c\|\cap B_1) + CQ_{V,\BC}\sum^{k-1}_{i=1}(\theta^\alpha)^i\\
	& \leq \dist_\H(\spt\|\BC_0\|\cap B_1,\spt\|\BC^c\|\cap B_1) + C_\alpha Q_{V,\BC},
\end{align*}
where $C_\alpha = C_\alpha(n,\alpha)$ is independent of $k$. Applying this with (\ref{E:fine-reg-3}) and substituting into (\ref{E:fine-reg-5}, we get
$$E^2_{V_{k-1},\BC^c} \leq 6\w_n\dist^2_\H(\spt\|\BC\|\cap B_1,\spt\|\BC^c\|\cap B_1) + CQ_{V,\BC}^2$$
and thus from Remark 1 of Section \ref{sec:fine_construction} and Hypothesis (H4), this gives, for $\gamma_1$ smaller than a constant depending only on $n$ and $\alpha$ (which is crucially independent of $k$)
\begin{equation}\label{E:fine-reg-6}
	E^2_{V_{k-1},\BC^c} \leq 12\w_n c_1^2E_{V,\BC^c}^2.
\end{equation}
But also, again from (\ref{E:fine-reg-2}) and the triangle inequality one has
\begin{align*}
	\dist_\H(\spt\|\BC_{k-2}\|&\cap B_1,\spt\|\BC^c\|\cap B_1)\\
	& \geq \dist_\H(\spt\|\BC_0\|\cap B_1,\spt\|\BC^c\|\cap B_1) - \sum^{k-2}_{i=1}\dist_\H(\spt\|\BC_{i-1}\|\cap B_1,\spt\|\BC_i\|\cap B_1)\\
	& \geq \dist_\H(\spt\|\BC_0\|\cap B_1,\spt\|\BC^c\|\cap B_1) - CQ_{V,\BC}
\end{align*}
where the constant $C = C(n,\alpha)$ is once again essentially unchanged from that in (\ref{E:fine-reg-3}), and thus using this with (\ref{E:fine-reg-3}) and (\ref{E:fine-reg-4}), with $k-1$ in place of $k$, we get
$$E_{V_{k-1},\BC^c} \geq \sqrt{2^{-n-4}\bar{C}_1}\dist_\H(\spt\|\BC\|\cap B_1,\spt\|\BC^c\|\cap B_1) - CQ_{V,\BC}$$
and thus by Remark 2 of Section \ref{sec:fine_construction} and Hypothesis (H4), we have
\begin{equation}\label{E:fine-reg-7}
	E_{V_{k-1},\BC^c} \geq (C_1- \tilde{C}\gamma_1)E_{V,\BC^c}
\end{equation}
where $C_1,\tilde{C}$ are fixed constants, independent of $k$, such that $C_1$ is dependent only on $n$ and $\tilde{C}$ depends on $n$ and $\alpha$. Thus, if $\gamma_1$ is such that $2\tilde{C}\gamma_1 < C_1$, then combining the assumed Hypothesis (H4) for $V,\BC^c,\BC$ (with $\gamma_1$) with (\ref{E:fine-reg-3}) (with $k-1$ in place of $k$) and (\ref{E:fine-reg-7}), we get
	\begin{equation}\label{E:fine-reg-8}
	Q_{V_{k-1},\BC_{k-1}}^2 \leq 4^{-k+1}Q_{V,\BC}^2 \leq 4^{-k+1}\gamma_1 E_{V,\BC^c}^2 \leq 4^{-k+1}\gamma_1\cdot (C_1/2)^{-2}\cdot E_{V_{k-1},\BC^c}^2.
	\end{equation}
	Of course, we already have from (\ref{E:fine-reg-6}) we already know that $E^2_{V_{k-1},\BC^c} \leq C\epsilon_1$, and thus if $C\epsilon_1< \epsilon_2$ and $\gamma_1(C_1/2)^{-2} < \gamma_2$, we have that $V_{k-1},\BC^c,\BC_{k-1}$ satisfy the Hypothesis (H) assumption of Lemma \ref{lemma:fine_excess_decay} with the correct parameters. So now let us turn to establishing Hypothesis ($\star$) holds with $M = \frac{3}{2}M_0$. Again, using (\ref{E:fine-decay-4}) with $k-1$ in place of $k$ as well as (\ref{E:fine-decay-2}) with $1,2,\dotsc,k-1$ in place of $k$, we have:
	$$\left(\int_{B_1}\dist^2(X,\spt\|\tilde{\BC}\|)\ \ext\|V_{k-1}\|\right)^{1/2} \geq \sqrt{2^{-n-4}\bar{C}_1}\dist_\H(\spt\|\BC\|\cap B_1,\spt\|\tilde{\BC}\|) - CQ_{V,\BC}$$
	and so we have
	\begin{align*}
		\int_{B_1}\dist^2(X,\spt\|\tilde{\BC}\|)\ \ext\|V_{k-1}\| & \geq \frac{1}{2}\cdot 2^{-n-4}\bar{C}_1\dist^2_\H(\spt\|\BC\|\cap B_1,\spt\|\tilde{\BC}\|\cap B_1) - CQ_{V,\BC}^2\\
		& \geq 2^{-n-5}\bar{C}_1\cdot (6\w_n)^{-1}\int_{B_1}\dist^2(X,\spt\|\tilde{\BC}\|)\ \ext\|V\| - \hat{C}Q^2_{V,\BC}\\
		& \geq 2^{-n-8}\bar{C}_1\w_n^{-1}(3/2)^{-1}E_{V,\BC^c}^2 - \hat{C}Q_{V,\BC}^2\\
		& \geq 2^{-n-9}\bar{C}_1\w_n^{-1}\cdot (12\w_n c_1^2)^{-1}E^2_{V_{k-1},\BC^c} - \hat{C}Q_{V,\BC}^2\\
		& \geq (2^{-n-13}\bar{C}_1\w_n^{-2}c_1^{-2} -C^\prime\gamma_1)E^2_{V_{k-1},\BC^c}
	\end{align*}
	where here in the third inequality we used our assumption that Hypothesis $(\star)$ holds for $V,\BC^c$ with $M= \frac{3}{2}$, in the fourth inequality we have used (\ref{E:fine-reg-6}), and in the last inequality we have used our assumption of Hypothesis (H4) on $V,\BC^c,\BC$ with $\gamma = \gamma_1$, followed by (\ref{E:fine-reg-7}) (thus $C^\prime = C^\prime(n,\alpha)$ is independent of $k$). Hence, if we choose $\gamma_1 = \gamma_1(\BC^{(0)})$ sufficiently small, we will ensure that (as the above was true for any such $\tilde{\BC}\in \FL_1\cap \FL_{1/10}(\BC^c)$
	$$E^2_{V_{k-1},\BC^c} \leq (2^{n+14}\bar{C}_1^{-1}\w_n^2c_1^2)\cdot \inf_{\tilde{\BC}\in \FL_1}E_{V_{k-1},\tilde{\BC}}$$
	and as this constant is $<\frac{3}{2}M_0$, we see that Hypothesis ($\star$) holds with $V_{k-1}$ and $\BC^c$ with $M= \frac{3}{2}M_0$. Hence, we see that as long as $\epsilon_1 = \epsilon_1(\BC^{(0)},\alpha)$ and $\gamma_1 = \gamma_1(\BC^{(0)},\alpha)$ are sufficiently small (independent of $k$, as the constants are reset in each application of Lemma \ref{lemma:fine_excess_decay}) we can apply Lemma \ref{lemma:fine_excess_decay} to $V_{k-1},\BC^c$, and $\BC_{k-1}$ to obtain a orthogonal rotation $\Gamma:\R^{n+1}\to \R^{n+1}$ and cone $\BC_k\in \FL_0$ such that, with $\Gamma_k := \Gamma_{k-1}\circ\Gamma$, (\ref{E:fine-decay-1}) -- (\ref{E:fine-decay-4}) hold; this completes the inductive proof that (\ref{E:fine-decay-1}) -- (\ref{E:fine-decay-4}) hold for all $k\in\{1,2,\dotsc\}$.
	
	Now write $\lambda^k_1,\dotsc,\lambda^k_3$, $\tilde{\lambda}^{k,1}_1, \tilde{\lambda}^{k,2}_1$ for the gradients of the half-hyperplanes in $\BC_k$ relative to the corresponding half-hyperplanes in $\BC_k^c$ in the usual way. Then from Remark 2 from Section \ref{sec:fine_construction} (applied with $V_k$ and $\BC_k$ in place of $V$ and $\BC$, respectively) we get, using (\ref{E:fine-reg-8}) and (\ref{E:fine-reg-7}), that
\begin{equation}\label{E:fine-reg-9}
	|\tilde{\lambda}^{k,1}_1-\tilde{\lambda}^{k,2}_1|\geq CE_{V,\BC^c}
\end{equation}
	where $C = C(n,\alpha)$; we stress that this is a fixed lower bound independent of $k$.
	
	Now, (\ref{E:fine-reg-2}) tells us that $(\spt\|\BC_k\|\cap B_1)_k$ is a Cauchy sequence (with respect to Hausdorff distance), and moreover as each $\BC_k\in \FL_0$ is level 0, so formed of multiplicity one half-hyperplanes, and moreover since we have fixed lower bounds on the Hausdorff distance between any pair of half-hyperplanes in $\BC_k$ (see \ref{E:fine-reg-9} and Hypothesis (H3)) we can find $\BC_*\in \FL_0$ such that $\BC_k\weakly \BC_*$; moreover, by the triangle inequality and (\ref{E:fine-reg-2}), we have for each $k\in \{1,2,\dotsc\}$,
	\begin{equation}\label{E:fine-reg-10}
	\dist^2_\H(\spt\|\BC_*\|\cap B_1,\spt\|\BC_k\|\cap B_1) \leq C\theta_k^{2\alpha}Q_{V,\BC}^2
	\end{equation}
	where $C = C(n,\alpha)$. But then from (\ref{E:fine-reg-3}) and (\ref{E:fine-reg-10}), and our mass upper bounds on $V_k$ in $B_1$ (from the multiplicity two class) we get for each $k\in \{1,2,\dotsc\}$,
	\begin{equation}\label{E:fine-reg-11}
	\int_{B_1}\dist^2(X,\spt\|\BC_*\|)\ \ext\|V_k\| \leq C\theta_k^{2\alpha}Q_{V,\BC}^2
	\end{equation}
	and
	\begin{equation}\label{E:fine-reg-12}
		\int_{B_{1/2}\backslash\{|x|<1/16\}}\dist^2(X,\spt\|V_k\|)\ \ext\|\BC_k\| \leq C\theta_k^{2\alpha}Q_{V,\BC}^2.
	\end{equation}
	Since all the $(V_k)_k$ belong to a multiplicity two class (Theorem \ref{thm:M2C}) we then have that $V_k\weakly \BC_*$; indeed, every subsequence of $(V_k)_k$ has a further subsequence (by the compactness property of multiplicity two classes) which converges in $B^{n+1}_1(0)$ to some varifold $V_*\in \S_2$; (\ref{E:fine-reg-11}) then tells us that $\spt\|V_*\|\cap B_1\subset \spt\|\BC_*\|\cap B_1$, and (\ref{E:fine-reg-11}) along with the weak convergence $\BC_k\weakly \BC_*$ gives that $\spt\|\BC_*\|\cap (B_{1/2}\backslash\{|x|<1/16\}) \subset \spt\|V\|\cap (B_{1/2}\backslash \{|x|<1/16\})$, which along with the mass upper bound $\|V_k\|(B_1) \leq (5/2+1/4)\w_n$ gives that $V_* = \BC_*$; as this limit was independent of the subsequences taken, this tells us that $V_k\weakly \BC_*$ without needing to pass to any subsequence.
	
	Let us now pass this information back to $V$. We know from (\ref{E:fine-reg-1}) that $(\Gamma_k)_k$ form a Cauchy sequence of rotations, and thus $\Gamma_k\to \Gamma$, where again for every $k\geq 0$, $|\Gamma-\Gamma_k| \leq C\theta^{k\alpha}Q_{V,\BC}$ for some $C = C(n,\alpha)$. Using this, the triangle inequality, and (\ref{E:fine-reg-10}) we therefore have for each $k\geq 1$,
	$$\theta_k^{-n-2}\int_{B_{\theta_k}}\dist^2(X,\spt\|\BC_*\|)\ \ext\|(\Gamma^{-1})_\#V\| \leq C\theta_k^{2\alpha}Q_{V,\BC}^2$$
	from which a standard scale-interpolation argument gives that, for each $\sigma\in (0,1/2)$,
	\begin{equation}\label{E:fine-reg-13}
		\sigma^{-n-2}\int_{B_\sigma}\dist^2(X,\spt\|\BC_*\|)\ \ext\|(\Gamma^{-1})_\#V\| \leq \tilde{C}\sigma^{2\alpha}Q_{V,\BC}^2
	\end{equation}
	where $\tilde{C} = C\theta^{-n-2-2\alpha}$ ($C$ the constant from (\ref{E:fine-reg-10})); in particular, $\tilde{C} = \tilde{C}(n,\alpha)$. Note that (\ref{E:fine-reg-13}) tells us two pieces of information, namely (i) $\BC_*$ is the unique tangent cone to $(\Gamma^{-1})_\#V$ at $0$ (so in particular $\BC_*\in \FL_S\cap \FL_0$), and (ii) there is a $\sigma\in (0,1/2)$ such that $V\res B^{n+1}_\sigma(0)$ has the structure of classical singularity of vertex density $5/2$; this latter fact follows from Theorem \ref{thm:A} in the case $I = 0$, as $\BC_*\in \FL_S\cap \FL_0$. The issue however is that this $\sigma$ will depend on $\BC_*$, $V$, and the point considered (which in this case is $0$) and so is not uniform in any manner from which one could deduce Theorem \ref{thm:L1fine_reg} at this moment. To get around this, we need to apply the above argument but with different base points $Z$ obeying $\Theta_V(Z)\geq 5/2$. Let us summarise everything we have proved as a consequence of our arguments so far:
	
	\textbf{Summary 1:} Given any $M_1\in [1,\infty)$, we have seen that (by simple modifications to our arguments) there are constants $\tilde{\epsilon}_1 = \tilde{\epsilon}_1(\BC^{(0)},M_1,\alpha)\in (0,1)$ and $\tilde{\gamma}_1 = \tilde{\gamma}_1(\BC^{(0)},M_1,\alpha)\in (0,1)$ such that if $V\in \S_2$, $\BC^c\in \FL_1$, and $\BC\in \FL_0$ are such that $\Theta_V(0)\geq 5/2$, $V\in \CN_{\tilde{\epsilon}_1}(\BC^{(0)})$, $\BC^c\in \FL_{\tilde{\epsilon}_1}(\BC^{(0)})$, $E_{V,\BC^c}<\frac{3}{2}M_1\inf_{\tilde{\BC}\in \FL_1}E_{V,\tilde{\BC}}^2$, and $Q^2_{V,\BC}<\tilde{\gamma}_1E^2_{V,\BC^c}$, then we may find $\theta = \theta(n,M_1,\alpha)\in (0,1/4)$ and orthogonal rotations $\Gamma,\Gamma_k:\R^{n+1}\to \R^{n+1}$ with $\Gamma_0 = \id$ obeying (from (\ref{E:fine-reg-1}))
	\begin{equation}\label{E:fine-reg-14}
		|\Gamma-\Gamma_k| \leq C\sigma_k^\alpha Q_{V,\BC};
	\end{equation}
	and a cone $\BC_0\in \FL_0$ such that (from (\ref{E:fine-reg-10}) and (\ref{E:fine-reg-13})) 
	\begin{equation}\label{E:fine-reg-15}
		\dist_\H(\spt\|\BC_0\|\cap B_1,\spt\|\BC\|\cap B_1) \leq CQ_{V,\BC};
	\end{equation}
	\begin{equation}\label{E:fine-reg-16}
		\sigma^{-n-2}\int_{B_\sigma}\dist^2(X,\spt\|\BC_0\|)\ \ext\|(\Gamma^{-1})_\#V\| \leq C\sigma^{2\alpha}Q_{V,\BC}^2\ \ \ \ \text{for all }\sigma\in (0,1/2);
	\end{equation}
	and for $k=1,2,\dotsc$ (from (\ref{E:fine-reg-11}) and (\ref{E:fine-reg-12}))
	\begin{equation}\label{E:fine-reg-17}
		\theta_k^{-n-2}\int_{B_{\theta_k}}\dist^2(X,\spt\|\BC_0\|)\ \ext\|(\Gamma^{-1}_k)_\#V\| \leq C\theta_k^{2\alpha}Q^2_{V,\BC};
	\end{equation}
	\begin{equation}\label{E:fine-reg-18}
		\theta_k^{-n-2}\int_{B_{\theta_k/2}\backslash\{|x|<\theta_k/16\}}\dist^2(X,\spt\|(\Gamma_k^{-1})_\#V\|)\ \ext\|\BC_0\| \leq C\theta_k^{2\alpha}Q_{V,\BC}^2;
	\end{equation}
	and moreover from (\ref{E:fine-reg-6}) and (\ref{E:fine-reg-7})
	\begin{equation}\label{E:fine-reg-19a}
		C^{-1}E_{V,\BC^c} \leq E_{(\eta_{0,\theta^k}\circ\Gamma_k)_\#V,\;\BC^c} \leq CE_{V,\BC^c};
	\end{equation}
	here, we have $C = C(n,M_1,\alpha)$. 
	
	\textbf{Summary 2:} From the proof of Corollary \ref{cor:fine_estimates}, we know that for any fixed $\tilde{\epsilon},\tilde{\gamma}\in (0,1/2)$, then there exist $\tilde{\epsilon}_2 = \tilde{\epsilon}_2(\tilde{\epsilon},\tilde{\gamma},\BC^{(0)},M_1,\alpha)\in (0,1/2)$ and $\tilde{\gamma}_2 = \tilde{\gamma}_2(\tilde{\epsilon},\tilde{\gamma},\BC^{(0)},M_1,\alpha)$ such that if the above hypotheses hold with $\tilde{\epsilon}_2,\tilde{\gamma_2}$ in place of $\tilde{\epsilon}_1$ and $\tilde{\gamma}_1$ (with the same $M_1$), then for any $Z\in \spt\|V\|\cap B_{9/16}$ with $\Theta_V(Z)\geq 5/2$, if we set $V_Z:= (\eta_{Z,1/8})_\#V$, then the above hypotheses hold for $V_Z$ (with the same $\BC^c$ and $\BC$) with $\tilde{\epsilon}, \tilde{\gamma}$, and $M_1M_0$ in place of $M_1$; moreover, we have (see (\ref{E:fine-est-2}))
	\begin{equation}\label{E:fine-reg-19b}
		E_{V_Z,\BC^c} \geq CE_{V,\BC^c}
	\end{equation}
	and combining (\ref{E:fine-est-3}) with Corollary \ref{cor:fine_estimates}(a),
	\begin{equation}\label{E:fine-reg-20}
		Q_{V_Z,\BC}\leq CQ_{V,\BC};
	\end{equation}
	again, here $C = C(n,M_1,\alpha)$ can depend on the value of $M_1$; we emphasis that in this situation there is no need to introduce Hypothesis ($\dagger$) (as Hypothesis ($\dagger$)(i) is always satisfied).
	
	Now let us fix $M_1\in [1,\infty)$. Let $\tilde{\epsilon}_1 = \tilde{\epsilon}_1(\BC^{(0)},M_0M_1,\alpha)$ and $\tilde{\gamma}_1 = \tilde{\gamma}_1(\BC^{(0)},M_0M_1,\alpha)$ be as in Summary 1 above. Now set $\tilde{\epsilon}_2 = \tilde{\epsilon}_2(\tilde{\epsilon}_1,\tilde{\gamma}_1,\BC^{(0)},M_1,\alpha)$ and $\tilde{\gamma}_2 = \tilde{\gamma}_2(\tilde{\epsilon}_1,\tilde{\gamma}_1,\BC^{(0)},M_1,\alpha)$ be as in Summary 2. Now fix $\epsilon_3\in (0,\tilde{\epsilon}_2]$ and $\gamma_3 \in (0,\tilde{\gamma}_2]$, and suppose the hypothesis above (as in Summary 1) hold with $\epsilon_3,\gamma_3$, and $M_1$. Hence, in view of Summary 1 and Summary 2, we see that the conclusions of Summary 1 hold for each base point $Z\in \spt\|V\|\cap B_{9/16}$ with $\Theta_V(Z)\geq 5/2$, i.e. there is a $\theta = \theta(n,M_1,\alpha)$ such that for each such $Z$ we can find orthogonal rotations $\Gamma_Z,\Gamma^k_Z:\R^{n+1}\to \R^{n+1}$ with $\Gamma_0^Z = \id$ such that
	\begin{equation}\label{E:fine-reg-21}
		|\Gamma_Z-\Gamma^k_Z|\leq C\theta_k^\alpha Q_{V_Z,\BC};
	\end{equation}
	a cone $\BC_Z\in \FL_0$ such that
	\begin{equation}\label{E:fine-reg-22}
		\dist_\H(\spt\|\BC_Z\|\cap B_1,\spt\|\BC\|\cap B_1) \leq CQ_{V_Z,\BC};
	\end{equation}
	\begin{equation}\label{E:fine-reg-23}
		\sigma^{-n-2}\int_{B_\sigma}\dist^2(X,\spt\|\BC_Z\|)\ \ext\|(\Gamma_Z^{-1})_\#V\| \leq C\sigma^{2\alpha}Q_{V_Z,\BC}^2;
	\end{equation}
	and for each $k=1,2,\dotsc$,
	\begin{equation}\label{E:fine-reg-24}
		\theta_k^{-n-2}\int_{B_{\theta_k}}\dist^2(X,\spt\|\BC_Z\|)\ \ext\|(\Gamma_Z^k)^{-1}_\#V\| \leq C\theta_k^{2\alpha}Q_{V_Z,\BC}^2;
	\end{equation}
	\begin{equation}\label{E:fine-reg-25}
		\theta_k^{-n-2}\int_{B_{\theta_k}}\dist^2(X,\spt\|(\Gamma_Z^k)^{-1}_\#V\|)\ \ext\|\BC_Z\| \leq C\theta_k^{2\alpha}Q_{V_Z,\BC}^2;
	\end{equation}
	which also obey
	\begin{equation}\label{E:fine-reg-26}
		C^{-1}E_{V_Z,\BC^c} \leq E_{(\eta_{0,\theta^k}\circ\left(\Gamma^k_Z\right)^{-1})_\#V_Z,\;\BC^c}\leq CE_{V_Z,\BC^c};
	\end{equation}
	and also
	\begin{equation}\label{E:fine-reg-27}
		E_{V_Z,\BC^c}\geq CE_{V,\BC^c};
	\end{equation}
	\begin{equation}\label{E:fine-reg-28}
		Q_{V_Z,\BC}\leq CQ_{V,\BC};
	\end{equation}
	here, $C = C(n,M_1,\alpha)$. Thus, we can conclude that every $Z\in \spt\|V\|\cap B_{9/16}$ with $\Theta_V(Z)\geq 5/2$ has a unique tangent cone $(\Gamma_Z)_\#\BC_Z\in \FL_0$; in particular, by Theorem \ref{thm:A} in the level 0 case, each such $Z$ is a classical singularity of $V$ and moreover $V\res B_{9/16}$ has no points of density $>5/2$.
	
	Now let us take $\tilde{\epsilon}_3\in (0,\tilde{\epsilon}_2(\epsilon_3,\gamma_3,\BC^{(0)},M_1M_0,\alpha)]$ and $\tilde{\gamma}_3\in (0,\tilde{\gamma}_2(\epsilon_3,\gamma_3,\BC^{(0)},M_1M_0,\alpha)]$. We first want to follow the proof of (\ref{E:thm-A-level-0-1}) in Theorem \ref{thm:mainL0}, now based on (\ref{E:fine-reg-21}) and (\ref{E:fine-reg-23}), to show that every slice $\R^2\times\{y\}$, for $y\in \{0\}^2\times B^{n-1}_{9/16}(0)$, has exactly one point of density $5/2$, and moreover that the points of density $5/2$ form a $C^{1,\alpha}$ submanifold. So suppose $Z_1,Z_2\in \spt\|V\|\cap B_{9/16}$ obey $\Theta_{V}(Z_1) = 5/2$ and $\Theta_V(Z_2)=5/2$; set $\sigma:= |Z_1-Z_2|$, and choose $k$ such that $\theta^{k+1}< 16\sigma \leq \theta^k$ (of course, by Lemma \ref{lemma:coarse_graphical_rep} we can without loss of generality assume $|Z_1|,|Z_2|<\theta/32$, and so $|Z_1-Z_2|<\theta/16$). Then if we set $\tilde{V}:= (\eta_{0,\theta^k}\circ(\Gamma_{Z_2}^k)^{-1})_\#V_{Z_2}$ and $\tilde{Z}:= (\eta_{0,\theta^k/8}\circ(\Gamma_{Z_2}^k)^{-1})(Z_1-Z_2)$, then clearly $\tilde{V}_{\tilde{Z}}:= (\eta_{\tilde{Z},1/8})_\#\tilde{V} = (\eta_{0,\theta^k}\circ(\Gamma_{Z_2}^k)^{-1})_\#V_{Z_1}$, and $\Theta_{\tilde{V}}(\tilde{Z}) = \Theta_{V}(Z_1) = 5/2$. We wish to verifying our hypotheses hold for $\tilde{V}$ for suitable parameters. Indeed, by (\ref{E:fine-reg-26}) and Summary 2 (with $M_1 = 1$) we have
	$$E^2_{\tilde{V},\BC^c} \leq CE^2_{V_{Z_2},\BC^c} \leq \frac{3}{2}CM_0\inf_{\tilde{\BC}\in \FL_1}E^2_{V_{Z_2},\tilde{\BC}}$$
	where $C = C(n,\alpha)$. Also, by (\ref{E:fine-reg-24}), (\ref{E:fine-reg-25}), (\ref{E:fine-reg-28}), (\ref{E:fine-reg-27}), (\ref{E:fine-reg-26}),
	\begin{equation}\label{E:fine-reg-29}
	Q_{\tilde{V},\BC_{Z_2}} \leq C\theta_k^{\alpha}Q_{V_{Z_2},\BC} \leq C\theta_k^{\alpha}Q_{V,\BC} \leq C\theta_k^\alpha \tilde{\gamma}_3 E_{V,\BC^c} \leq C\theta_k^\alpha \tilde{\gamma}_3 E_{V_{Z_2},\BC^c} \leq C^\prime \tilde{\gamma}_3 E_{\tilde{V},\BC^c}.
	\end{equation}
	Moreover, we clearly have from (\ref{E:fine-reg-26}) that $E_{\tilde{V},\BC^c} \leq CE_{V_{Z_2},\BC^c} \leq C\tilde{\epsilon}_3$. Thus, we have our hypotheses are satisfied with $\tilde{V}$ in place of $V$, $C\tilde{\epsilon}_3$ in place of $\epsilon$, $C^\prime\tilde{\gamma}_3$ in place of $\gamma$, and $M_1 = CM_0$ (with $\BC_{Z_2}$ in place of $\BC$, with the same $\BC^c$); hence we may apply our deductions proceeding Summary 1 and Summary 2 to $\tilde{V}$ to find a cone $\tilde{\BC}_{\tilde{Z}}\in \FL_0$ and a rotation $\tilde{\Gamma}_{\tilde{Z}}$ such that
	\begin{equation}\label{E:fine-reg-30}
		\dist_\H(\spt\|\tilde{\BC}_{\tilde{Z}}\|\cap B_1,\spt\|\BC_{Z_2}\|\cap B_1)\leq CQ_{\tilde{V}_{\tilde{Z}},\BC_{Z_2}};
	\end{equation}
	\begin{equation}\label{E:fine-reg-31}
		|\tilde{\Gamma}_{\tilde{Z}}-\id|\leq CQ_{\tilde{V}_{\tilde{Z}},\BC_{Z_2}};
	\end{equation}
	moreover, $(\tilde{\Gamma}_{\tilde{Z}})_\#\tilde{\BC}_{\tilde{Z}}$ is the unique tangent cone to $\tilde{V}$ at $\tilde{Z}$, but unravelling the transformations reveals
	\begin{equation}\label{E:fine-reg-32}
	(\tilde{\Gamma}_{\tilde{Z}})_\#\tilde{\BC}_{\tilde{Z}} = \left[(\Gamma_{Z_2}^k)^{-1}\circ\Gamma_{Z_1}\right]_\#\BC_{Z_1}.
	\end{equation}
	But then (\ref{E:fine-reg-30}), (\ref{E:fine-reg-31}), and (\ref{E:fine-reg-21}) gives:
	$$\dist_\H(\spt\|(\Gamma_{Z_1})_\#\BC_{Z_1}\|\cap B_1, \spt\|(\Gamma_{Z_2})_\#\BC_{Z_2}\|\cap B_1) \leq CQ_{\tilde{V},\BC_{Z_2}} + C\theta_k^\alpha Q_{V_{Z_2},\BC}$$
	which by (\ref{E:fine-reg-29}) and (\ref{E:fine-reg-20}) gives
	\begin{equation}\label{E:fine-reg-33}
	\dist_\H(\spt\|(\Gamma_{Z_1})_\#\BC_{Z_1}\|\cap B_1, \spt\|(\Gamma_{Z_2})_\#\BC_{Z_2}\|\cap B_1) \leq C|Z_1-Z_2|^\alpha Q_{V,\BC}.
	\end{equation}
	Moreover, we can show that there is at most one point of density $5/2$ in each slice $\R^2\times\{y\}$ with $y\in \{0\}^2\times B^{n-1}_{1/2}$; moreover, we can show that if we take distinct points $Z_1,Z_2\in \spt\|V\|\cap (\R^2\times\{y\})$ and $\Theta_V(Z_1)=5/2$, then in fact $\Theta_{V}(Z_2) = 1$; in particular, $Z_2\in \reg(V)$. Indeed, to see this choose $k$ such that $\theta^{k+1}<|Z_1-Z_2|\leq \theta^k$ and use (\ref{E:fine-reg-24}) -- (\ref{E:fine-reg-28}) with $Z = Z_1$ to see that the assumptions of the theorem hold for $V_*:= (\eta_{Z_1,\theta^k}\circ(\Gamma_{Z_1}^k)^{-1})_\#V$ for suitable $\epsilon,\gamma$, and with $M_1 = CM_0$, and thus in particular by Lemma \ref{thm:fine_representation} taken with $\tau = \theta^2/2$, we see that if $Z_* := (\Gamma^k_{Z_1})^{-1}(Z_2-Z_1)/\theta$ that $|Z_*|\in \{|x|>\theta^2\}$, and so we must have $\Theta_{V_*}(Z_*) = 1$; as $\Theta_{V_*}(Z_*) = \Theta_V(Z_2)$, this produces the desired conclusion. But from Lemma \ref{lemma:gaps} (or in fact we could now use Theorem \ref{thm:A} in the level 0 case) we see that in fact every such slice $\R^2\times\{y\}$ must contain a point of density $5/2$. Thus, if we define $\phi:\{0\}^2\times B^{n-1}_{1/2}(0)\to \R^2$ to be $\phi(y) = Z_y$, for $Z_y$ the unique point in $\spt\|V\|\cap (\R^2\times\{y\})$ with $\Theta_V(Z_y) =5/2$, we see from Theorem \ref{thm:A} that $\graph(\phi) = \{\Theta_V = 5/2\}\cap B_{1/2}$ is a $C^{1,\alpha}$ submanifold, and moreover that the unique tangent plane at a $y$ mapped to $Z_y$ is the spine of $(\Gamma_{Z_y}) _\#\BC_{Z_y}$, which is $\Gamma_{Z_y}(\{0\}^2\times \R^{n-1})$; in particular, we get from (\ref{E:fine-reg-33}) that $[D\phi]_{0,\alpha} \leq CQ_{V,\BC}$. Moreover, the corresponding bounds on $\sup|D\phi|$ follow from (\ref{E:fine-reg-21}) with $k=0$, and the bounds on $\sup|\phi|$ follow from Corollary \ref{cor:fine_estimates}(a).
	
	The last thing we must justify is the graph structure away from the points of density $5/2$. However, given our estimates this follows in essentially the same way as in the last stages of \cite[Proof of Theorem 16.1]{wickstable}, and so we do not repeat the details here. This therefore completes the proof of the theorem, for suitable choices of $\epsilon,\gamma$.
\end{proof}

As a consequence of Theorem \ref{thm:L1fine_reg}, we are now able to prove that property $(\FB7)$ holds for the coarse blow-up class $\FB(\BC^{(0)})$ when $\BC^{(0)}\in \FL_S\cap \FL_1$ is a level 1 cone.

\begin{corollary}\label{cor:L1B7}
	Let $\BC^{(0)}\in \FL_S\cap \FL_1$. Then, the coarse blow-up class $\FB(\BC^{(0)})$ as defined in Section \ref{sec:coarse_construction} obeys property $(\FB7)$ of Section \ref{sec:coarse_regularity}; in particular, $\FB(\BC^{(0)})$ is a proper blow-up class in the sense of Definition \ref{defn:proper-blow-up-class}, and so satisfies the conclusions of Theorem \ref{thm:coarse_reg}.
\end{corollary}

\begin{proof}
	Suppose the for contradiction that $(\FB7)$ does not hold. Thus, for each $k\in \{1,2,\dotsc\}$, we could find $v^k\in \FB(\BC^{(0)})$ obeying $v^k_a(0) = 0$, $Dv^k_a(0) = 0$, $\|v^k\|_{L^2} = 1$, and $v_*^k$ such that $v_*^k$ is comprised of linear functions with common boundary and zero average (over each half-hyperplane), which moreover satisfy that over the (unique) multiplicity two half-hyperplane in $\BC^{(0)}$ is represented by two linear functions $\ell^k_1,\ell^k_2$ obeying $\ell^k_1 = -\ell^k_2\not\equiv 0$ (so $\graph(v^*_k)$ is a level 0 cone), and moreover
	$$\int_{B_1}\G(v^k,v^k_*)^2 < \frac{1}{k}.$$
	It suffices to show that in fact infinitely many of the $v^k$ are $C^{1,\alpha}$ up-to-the-boundary in $B_{1/2}$, for some $\alpha = \alpha(n)$ independent of the choice of sequences $v^k, v^k_*$.
	
	We may pass to a subsequence to assume that $v_*^k\to v_*$ (e.g. in $C^1$); by hypothesis, $v_*$ will be zero over each multiplicity one half-hyperplane in $\spt\|\BC^{(0)}\|$, and over the multiplicity two half-hyperplane will be given by two distinct linear functions with zero average: moreover, as $\|v^k_*\|_{L^2} >1-1/k$, we have $\|v_*\|_{L^2} = 1$ (which is what tells us that over the multiplicity two half-hyperplane, the linear functions cannot agree); in particular, $\graph(v_*)$ is a level 0 cone still). Moreover, we have
	$$\int_{B_{1}}\G(v^k,v_*)^2 \to 0.$$
	Now for each $k\in \{1,2,\dotsc\}$ let $(V^k_{j})_j\subset\S_2$ and $(\BC^k_j)_j\subset\FL_1$ be sequences such that the coarse blow-up sequence $v^k_j:= E^{-1}_{V_j^k,\BC^k_j} u^k_{j}$ of $V^k_j$ relative to $\BC^k_j$ gives rise to $v^k$ (as $j\to\infty$); without loss of generality, we may translate to assume that $\Theta_{V^k_j}(0) = 5/2$ (using Lemma \ref{lemma:gaps}). We therefore know that, for any $\sigma\in (0,1)$, given any $\delta>0$, for each $k$ we can find $j_k$ such that, for all $j\geq j_k$,
	$$\int_{B_\sigma(0)}\G(v^k_{j_k},v^k)^2 < \delta^2.$$
	Thus, setting $V_k:= V^k_{j_k}$, $\BC_k:= \BC^k_{j_k}$, and $\tilde{v}_k:= v^k_{j_k} \equiv E^{-1}_{V_k,\BC_k}\tilde{u}_k$, we have
	$$\int_{B_\sigma}\G(\tilde{v}_k,v_*)^2\to 0.$$
	In particular, we may assume that $v_*\in \FB(\BC^{(0)})$ is the coarse blow-up of $V_k$ relative to $\BC_k$. We now claim that we must have, for all $k$ sufficiently large,
	\begin{equation}\label{E:B7-0}
		\int_{B_1}\dist^2(X,\spt\|\BC_k\|)\ \ext\|V_k\| < \frac{3}{2}\inf_{\tilde{\BC}\in \FL_1}\int_{B_1}\dist^2(X,\spt\|\tilde{\BC}\|)\ \ext\|V_k\|.
	\end{equation}
	We argue this again by construction. If this were false, then we could find a subsequence (which we pass to) such that the reverse inequality holds; thus, choosing $\tilde{\BC}_k\in \FL_1$ such that
	$$\inf_{\tilde{\BC}\in \FL_1}\int_{B_1}\dist^2(X,\spt\|\tilde{\BC}\|)\ \ext\|V_k\|> \frac{4}{5}\int_{B_1}\dist^2(X,\spt\|\tilde{\BC}_k\|\ \ext\|V_k\|$$
	we have
	\begin{equation}\label{E:B7-1}
		\int_{B_1}\dist^2(X,\spt\|\tilde{\BC}_k\|)\ \ext\|V_k\| < \frac{5}{6}\int_{B_1}\dist^2(X,\spt\|\BC_k\|)\ \ext\|V_k\|.
	\end{equation}
	Let us now denote by $(\lambda^i_k)_i$ the gradients of the half-hyperplanes in $\tilde{\BC}_k$ relative to the corresponding half-hyperplanes in $\BC_k$. In particular, for each $\sigma\in (1/2,1)$ and sufficiently large $k$,
	\begin{equation}\label{E:B7-2}
		\sum_{i}(1+(\lambda^i_k)^2)^{-1}\int_{H_i^{(0)}\cap B_\sigma\backslash \{|x|<1-\sigma\}}|\tilde{u}^i_k - \lambda^i_k x^{\perp_i}|^2 < \frac{5}{6}E^2_{V_k,\BC_k}
	\end{equation}
	where $H^{(0)}_i$ are the half-hyperplanes of $\spt\|\BC^{(0)}\|$ and $\perp_i$ the orthogonal projection onto $(H^{(0)}_i)^\perp$; we stress here that for $k$ sufficiently large and $\sigma\in (0,1/2)$ fixed, the two-valued piece of $\tilde{u}_k$ will necessarily be two single-valued functions, due to the fact that $\tilde{v}_k\to v_*$, and thus by the above sum we are including the half-hyperplanes in $\spt\|\BC^{(0)}\|$ with their respective multiplicities. In particular, (\ref{E:B7-2}) gives that, $\sum_i(1+(\lambda^i_k)^2)^{-1}(\lambda_k^i)^2\int_{B_{1/2}\backslash\{|x|<1/4\}}|x^2|^2 \leq \frac{5}{3}E^2_{V_k,\BC_k} + 2E_{V_k,\BC_k}^2 = \frac{11}{3}E_{V_k,\BC_k}^2$, and hence $|\lambda^i_k| \leq CE^2_{V_k,\BC_k}$ for all $k$ sufficiently large, where $C = C(n)$. Thus, we may assume that $E^{-1}_{V_k,\BC_k}\to \ell^i$ for some $\ell^i\in \R$. We then get from (\ref{E:B7-2}), dividing by $E^2_{V_k,\BC_k}$, taking $k\to\infty$, and then $\sigma\uparrow 1$, that
	$$\sum_{i,j}\int_{B_1}|(v_*)^j_i - \ell^i|^2 \leq \frac{5}{6};$$
	here, the sum over $i$ is over each distinct half-hyperplane in $\spt\|\BC^{(0)}\|$ and the sum over $j$ is over the number of values of $(v_*)_i$ over a given half-hyperplane. Expanding this, noting that $\sum_j (v_*)_i^j = 0$ (as each component of $v_*$ is average-free) we get
	$$\int_{B_1}|v_*|^2 + \sum_i\int_{B_1}|\ell^i|^2 \leq \frac{5}{6}$$
	which obviously contradicts $\|v_*\|_{L^2} = 1$. Thus, (\ref{E:B7-0}) holds for all $k$ sufficiently large.
	
	Now define a new sequence of cones, $\hat{\BC}_k$, via $v_*$ in the usual way: by modifying the gradients of the half-hyperplanes in $\BC_k$ relative to $\BC^{(0)})$ by $E_{V_k,\BC_k}\cdot (v_*)^j_i$ (depending on the number of values of $v_*$ over the respective half-hyperplane); thus, $\hat{\BC}_k$ is a level 0 cone. Then, for any $\sigma\in (0,1)$, the estimates from Corollary \ref{corollary:non-concentration} give:
	\begin{equation}\label{E:B7-3}
		\int_{B_\sigma}\dist^2(X,\spt\|\hat{\BC}_k\|)\ \ext\|V_k\| \leq 2\int_{B_\sigma}\G(u_k,E_{V_k,\BC_k} v_*)^2 + C\sigma^{1/2}E^2_{V_k,\BC_k}.
	\end{equation}
	Moreover, as $v_*$ is homogeneous of degree one and obeys $\|v_*\|_{L^2(B_1)}=1$, we know that for all $\tau,\sigma\in (0,1)$ with $\tau<\sigma$ we have $\int_{B_\sigma\backslash\{|x|<\tau\}}|v_*|^2 = \sigma^{n+2} - \tau^3/\sigma^3$, and thus for any $\theta\in (0,1/8)$, for all $k$ sufficiently large we have $\int_{B_\sigma\backslash\{|x|<\tau\}}|u_k|^2 \geq (1-\theta)(\sigma^{n+2} - \tau/\sigma)E_{V_k,\BC_k}^2$. Moreover, again by Corollary \ref{corollary:non-concentration}, we know
	$$\int_{B_\sigma}\dist^2(X,\spt\|\BC_k\|)\ \ext\|V_k\| \geq \int_{B_\sigma\backslash\{|x|<\tau\}}|u_k|^2 - C\tau^{1/2}E_{V_k,\BC_k}^2$$
	and thus from this we see, for all sufficiently large $k$,
	\begin{equation}\label{E:B7-4}
		\int_{B_1\backslash B_\sigma}\dist^2(X,\spt\|\BC_k\|)\ \ext\|V_k\| \leq \left(1-(1-\theta)(\sigma^{n+2}-\tau/\sigma) + C\tau^{1/2}\right)E_{V_k,\BC_k}^2.
	\end{equation}
	But then, from the triangle inequality and the definition of $\hat{\BC}_k$ we have
	$$\int_{B_1\backslash B_\sigma} \dist^2(X,\spt\|\hat{\BC}_k\|)\ \ext\|V_k\| \leq 2\int_{B_1\backslash B_\sigma}\dist^2(X,\spt\|\BC_k\|)\ \ext\|V_k\| + C\H^n(B_1\backslash B_\sigma)E^2_{V_k,\BC_k};$$
	combining this with (\ref{E:B7-4}) and (\ref{E:B7-3}), and using the fact that $E_{V_k,\BC_k}^{-1}u_k\to v_*$ in (\ref{E:B7-3}), we see that for any $\delta\in (0,1)$, we may choose $\sigma = \sigma(n,\delta)$ sufficiently close to $1$, $\theta = \theta(n,\delta)$ and $\tau = \tau(n,\delta,\theta)$ sufficiently close to $0$ to get that, for all $k$ sufficiently large,
	$$\int_{B_1}\dist^2(X,\spt\|\hat{\BC}_k\|)\ \ext\|V_k\| \leq \delta E^2_{V_k,\BC_k}.$$
	For any $\alpha\in (0,1)$, if we choose $\delta = \gamma_1/4$, where $\gamma_1 = \gamma_1(\BC^{(0)},\alpha)$ is the constant from Theorem \ref{thm:L1fine_reg}, then we see that for all $k$ sufficiently large,
	\begin{equation}\label{E:B7-5}
		\int_{B_1}\dist^2(X,\spt\|\hat{\BC}_k\|)\ \ext\|V_k\| \leq \frac{\gamma_1}{4}E^2_{V_k,\BC_k}.
	\end{equation}
	This bounds on half of $Q_{V_k,\hat{\BC}_k}$. Similarly to how we have seen before, we can bound the other half of $Q_{V_k,\hat{\BC}_k}$ using the graphical representation provided by Lemma \ref{lemma:L2_coarse}, achieving
	\begin{equation}\label{E:B7-6}
		\int_{B_{1/2}\backslash \{|x|<1/16\}}\dist^2(X,\spt\|V_k\|)\ \ext\|\hat{\BC}_k\| \leq \eta_k E_{V_k,\BC_k}^2
	\end{equation}
	where $\eta_k\downarrow 0$. Thus, combining (\ref{E:B7-5}) and (\ref{E:B7-6}) we have for all $k$ sufficiently large,
	\begin{equation}\label{E:B7-7}
	Q_{V_k,\hat{\BC}_k}^2 < \frac{\gamma_1}{2}E_{V_k,\BC_k}^2.
	\end{equation}
	Hence we can now apply Theorem \ref{thm:L1fine_reg} to see that each $V_k$ is represented by functions which are $C^{1,\alpha}$ up-to-the-boundary, with estimates. In particular, as $V_k := V^k_{j_k}$, we can take $u^k_{j_k}$ to simply be the function which is $C^{1,\alpha}$ up-to-the-boundary whose function agrees with $V^k_{j_k}$. But we could also re-run this argument for $V^k_j$, where $j\geq j_k$ is arbitrary, to see that $v^k_j$ is $C^{1,\alpha}$ up-to-the-boundary, with estimates, for all $j\geq j_k$; hence $v^k$ is $C^{1,\alpha}$ up-to-the-boundary with estimates, for all $k$ sufficiently large; but this is a contradiction to our original assumption, and hence the proof is completed.
\end{proof}

\section{Level 1: Proof of Main Theorem}\label{sec:L1coarse_decay}

We have now proved, in Corollary \ref{cor:L1B7}, that for $\BC^{(0)}\in \FL_S\cap \FL_1$ a level 1 cone, the coarse blow-up class $\FB(\BC^{(0)})$ obeys the regularity conclusions of Theorem \ref{thm:coarse_reg}. We have also seen in Theorem \ref{thm:L1fine_reg} the fine $\epsilon$-regularity theorem for level 1 cones. In this section, we will combine these two results to prove Theorem \ref{thm:A} in the setting where $\BC^{(0)}\in \FL_S\cap \FL_1$ is level 1.

\begin{theorem}\label{thm:mainL1}
	Theorem \ref{thm:A} is true whenever $\BC^{(0)}\in \FL_S\cap \FL_1$.
\end{theorem}

\begin{proof}
	Fix $\BC^{(0)}\in \FL_S\cap \FL_1$. We first claim the following: there exists $\epsilon = \epsilon(\BC^{(0)})\in (0,1)$ and $\theta = \theta(n)\in (0,1)$ such that the following dichotomy holds: if $\BC^c\in \FL_1$ obeys $\dist_\H(\spt\|\BC^c\|\cap B_1,\spt\|\BC^{(0)}\|\cap B_1)<\epsilon$, and if $V\in \S_2$ is such that $\Theta_V(0)\geq 5/2$, $(2^n\w_n)^{-1}\|V\|(B_2(0))\in (2+1/16, 3-1/16)$, and $E_{V,\BC^c}^2<\epsilon$, then either:
	\begin{enumerate}
		\item [(i)] there is a cone $\BC^\prime\in \FL_1$ with $\dist_\H(\spt\|\BC^\prime\|\cap B_1,\spt\|\BC^c\|\cap B_1) \leq CE_{V,\BC^c}$ and moreover $\theta^{-n-2}\int_{B_{\theta}}\dist^2(X,\spt\|\BC^\prime\|)\ \ext\|V\|\leq \frac{1}{2}E_{V,\BC^c}^2$; or,
		\item [(ii)] there is a cone $\BC\in \FL_0$ and a rotation $\Gamma$ with $|\Gamma-\id|\leq CE_{V,\BC^c}$, $\dist_\H(\spt\|\BC\|\cap B_1,\spt\|\BC^c\|\cap B_1)\leq CE_{V,\BC^c}$, and $\rho^{-n-2}\int_{B_\rho}\dist^2(X,\spt\|\BC\|)\ \ext\|\Gamma_\#V\| \leq C\rho^{2\mu} E_{V,\BC^c}^2$ for all $\rho\in (0,\theta/8]$;
	\end{enumerate}
	here $C = C(n)$ and $\mu = \mu(n)$. To prove this, we argue by contradiction. Suppose we have a sequence of varifolds $(V_k)_k\subset \S_2$ and a sequence of level 1 cones $(\BC_k^c)_k\subset \FL_1$ with $\dist_\H(\spt\|\BC^c_k\|\cap B_1,\spt\|\BC^{(0)}\|\cap B_1) <\epsilon_k$ such that $\Theta_{V_k}(0)\geq 5/2$, $(2^n\w_n)^{-1}\|V_k\|(B_2(0))\in (2+1/16,3-1/16)$, and $E_{V_k,\BC^c_k} < \epsilon_k$, where $\epsilon_k\downarrow 0$. If necessary, we can replace $\BC^c_k$ with a sequence of level 1 cones $\tilde{\BC}^c_k$ obeying $\tilde{\BC}^c_k\in \FL_{\epsilon_k}(\BC^{(0)})$; so let us assume this without loss of generality for our cones $\BC^c_k$. Let $v\in \FB(\BC^{(0)})$ be the coarse blow-up of the sequence $(V_k)_k$ relative to $(\BC^c_k)_k$; since $\Theta_{V_k}(0) \geq 5/2$ for all $k$, this implies that $v_a(0) = 0$. So, by Theorem \ref{thm:coarse_reg}, we know that there is some $\phi\in C^1(\BC^{(0)})$ with $\mathbf{v}(\phi)\in \FL_0\cup \FL_1$ such that for every $\sigma\in (0,1/8]$
	\begin{equation}\label{E:thmA-L1-1}
	\sigma^{-n-2}\int_{B_\sigma}\G(v,\phi)^2 \leq C_1\sigma^{2\mu}\int_{B_{1/2}}|v|^2
	\end{equation}
	where here $C_1 = C_1(n)$ and $\mu = \mu(n)$. Since $\int_{B_1}|v|^2 \leq 1$, this with the homogeneity of $\phi$ implies that 
	\begin{equation}\label{E:thmA-L1-2}
		\int_{B_1}|\phi|^2 \leq C_2
	\end{equation}
	where $C_2 = C_2(n)$; we know that $\phi_a(0) = Dv_a(0)\cdot x$ (understood as equality on each respective half-hyperplane). Let us now choose $\theta = \theta(\BC^{(0)})$ such that
	$$\max\{C_1(2\theta)^\mu, (2\theta)^\mu\}<\min\{1/8,\epsilon_1\}$$
	where here $\epsilon_1 = \epsilon_1(\BC^{(0)})$ is the constant from $(\FB7)$ for the class $\FB(\BC^{(0)})$, and $\mu = \mu(n)$ is from (\ref{E:thmA-L1-1}).
	
	We then have two cases. Firstly, if $(2\theta)^{-n-2}\int_{B_{2\theta}}\G(v,\phi_a)^2 <(2\theta)^\mu$. In this case, define $\tilde{\BC}_k\in \FL_1$ in the usual fashion, but modifying the gradients of the half-hyperplanes in $\BC^c_k$ relative to the corresponding half-hyperplanes $H_i$ in $\BC^{(0)}$ by $E_{V_k,\BC^c_k}\cdot D_{H_i}\phi_a$ (for the corresponding value of $D_{H_i}\phi_a$, where by $D_{H_i}$ we mean the derivative in the direction of the ray in the cross-section of $H_i$ giving rise to $H_i$; see also the proof of Theorem \ref{thm:mainL0}). It is then standard to check that (i) holds for infinitely many case in this situation, with $C = C(n)$.
	
	The second case is when the first case fails, i.e. when we have
	\begin{equation}\label{E:L1-3}
		(2\theta)^{-n-2}\int_{B_{2\theta}}\G(v,\phi_a)^2\geq (2\theta)^\mu;
	\end{equation}
	in this case, we must have $\mathbf{v}(\phi)\in \FL_0$, since otherwise $\mathbf{v}(\phi)\in \FL_1$, implying that over the two linear functions in $\phi$ over the multiplicity two half-hyperplane in $\BC^{(0)})$ agree, from which (\ref{E:thmA-L1-1})  (with $\sigma = 2\theta$) would imply, using (\ref{E:thmA-L1-2}), that (\ref{E:L1-3}) does not hold. One may then find rotations $\tilde{\Gamma}_k$ of $\R^{n+1}$ which rotate the spine of $\mathbf{v}(E_{V_k,\BC_k^c}\cdot\phi)$ to $\{0\}^2\times\R^{n-1}$ and obey $\|\tilde{\Gamma}_k-\id\|\to 0$, and such that
	$$\tilde{v}(x):= \|v(2\theta(\cdot))-\phi_a(2\theta(\cdot))\|^{-1}_{L^2(B_1)}(v(2\theta(\cdot))-\phi_a(2\theta(\cdot)))$$
	is the coarse blow-up of (a subsequence of) $W_k:=(\eta_{0,2\theta}\circ\tilde{\Gamma}_k)_\#V_k$ (relative to $\BC^c_k$); but then if $\tilde{\phi}(x):= \|v(2\theta(\cdot))-\phi_a(2\theta(\cdot))\|^{-1}_{L^2(B_1)}(\phi(2\theta(\cdot))-\phi_a(2\theta(\cdot)))$, from (\ref{E:thmA-L1-1}) we have
	$$\int_{B_1}\G(\tilde{v},\tilde{\phi})^2 \leq C_1(2\theta)^\mu < \epsilon_1.$$
	Since $\|\tilde{v}\|_{L^2(B_1)} = 1$ and $\tilde{\phi}_a\equiv 0$ yet $\tilde{\phi}\not\equiv 0$, the assumptions of $(\FB7)$ are satisfied; thus, by the proof of Corollary \ref{cor:L1B7}, we know that for all $k$ sufficiently large the hypotheses of Theorem \ref{thm:L1fine_reg} are satisfied with $W_k$ in place of $V$, and thus we see that in fact (ii) must hold in this case.
	
	We can now apply the established dichotomy iteratively, to deduce that (taking $\BC^c_0 = \BC^{(0)}$) one of the following must hold (set $V_k:= (\eta_{0,\theta^k})_\#V$):
	\begin{enumerate}
		\item [(i)'] there is a sequence of level 1 cones $(\BC^c_k)_k$ with $\dist_\H(\spt\|\BC^c_{k+1}\|\cap B_1,\spt\|\BC^c_{k}\|\cap B_1) \leq CE_{V_k,\BC_k}$ and $E^2_{V_{k+1},\BC_{k+1}}\leq \frac{1}{2}E^2_{V_{k},\BC_{k}}$ for all $k\geq 0$; or,
		\item [(ii)'] there is an integer $I\geq 0$ and a finite sequence of level 1 cones $\BC^c_0 = \BC^{(0)},\BC^c_1,\dotsc,\BC^c_I$, such that (i)' holds for $k=0,1,\dotsc,I-1$ (if $I\geq 1$), and there is a level 0 cone $\BC\in \FL_S\cap \FL_0$ with $\dist_\H(\spt\|\BC\|\cap B_1,\spt\|\BC^c_I\|\cap B_1)\leq CE_{V_{I},\BC_{I}}$ and a rotation $\Gamma$ with $|\Gamma-\id| \leq CE_{V_I,\BC^c_I}$ such that $(\rho\theta^I)^{-n-2}\int_{B_{\rho\theta^I}}\dist^2(X,\spt\|\BC\|)\ \ext\|\Gamma_\#V\| \leq C\rho^{2\mu}E^2_{V_{I},\BC^c_I}$ for all $\rho\in (0,\theta/8]$.
	\end{enumerate}
	From these, we readily deduce that there are constants $C = C(n)\in (0,\infty)$ and $\beta = \beta(\BC^{(0)})\in (0,1)$ such that we have either:
	\begin{enumerate}
		\item [(A)] there is a (unique) level 1 cone $\BC_1\in \FL_S\cap \FL_1$ with $\dist_\H(\spt\|\BC_1\|\cap B_1,\spt\|\BC^{(0)}\|\cap B_1) \leq CE_{V,\BC^{(0)}}$ and $E_{(\eta_{0,\rho})_\#V,\BC_1}^2 \leq C\rho^{2\beta}E_{V,\BC^{(0)}}^2$ for all $\rho\in (0,\theta/8]$; or,
		\item [(B)] there is a (unique) level 0 cone $\BC_0\in \FL_S\cap \FL_0$ and rotation $\Gamma:\R^{n+1}\to \R^{n+1}$ with $\dist_\H(\spt\|\BC_0\|\cap B_1,\spt\|\BC^{(0)}\|\cap B_1) \leq CE_{V,\BC^{(0)}}$, $|\Gamma-\id|\leq CE_{V,\BC^{(0)}}$, and $E_{(\eta_{0,\rho}\circ\Gamma)_\#V,\BC_0}^2\leq C\rho^{2\beta}E_{V,\BC^{(0)}}^2$ for all $\rho\in (0,\theta/8]$.
	\end{enumerate}
	Indeed, (A) holds when (i)' holds and (B) holds when (ii)' holds. In particular, $\Theta_V(0) = 5/2$ and $V$ has a unique tangent cone at $0$ which is either a level 0 or level 1 stationary cone.
	
	Finally, to complete the proof note that the hypotheses of Theorem \ref{thm:A} will still holds if one replaces $V$ with $(\eta_{Z,1/4})_\#V$ for any $Z\in \spt\|V\|\cap B_{3/8}$ obeying $\Theta_V(Z)\geq 5/2$, provided $\epsilon = \epsilon(\BC^{(0)})$ is sufficiently small (this follows from Lemma \ref{lemma:L2_coarse}(i)); thus (A) or (B) above hold at each $Z\in \spt\|V\|\cap B_{3/8}$ obeying $\Theta_V(Z)\geq 5/2$. At this point, the proof can be completed in a similar manner to that seen in Theorem \ref{thm:mainL0}; thus we have completed the proof.
\end{proof}

\textbf{Remark:} Currently, the final power $\beta$ in the above proof could depend on the base cone $\BC^{(0)}$. However, once we have established Theorem \ref{thm:A} in the level 2 case with a power which is independent of the level 2 base cone (the independence of which is immediate in the level 2 case as there is only one level 2 base cone up to rotations) we will be able to deduce that the power can be chosen in the level 1 case to only depend on the dimension $n$ and not on the specific choice of level 1 cone.

\section{The Ultra Fine Blow-Up Class}\label{sec:L2}

In this section we will begin the proof that for each level 2 base cone $\BC^{(0)}\in \FL_S\cap\FL_2$, the fine blow-up class $\FB^F_{3,1;M}(\BC^{(0)})$ obeys property $(\FB7)$ of Section \ref{sec:coarse_regularity}, as thus enjoys the $C^{1,\alpha}$ boundary regularity as seen in Section \ref{sec:coarse_regularity}. Currently, if we were to try and replicate the proof of Corollary \ref{cor:L1B7} in this setting, we would arrive at a situation where the excess relative to a sequence of level 0 cones is significantly smaller than the excess relative to a sequence of level 1 cones which we were taking a fine blow-up sequence of. However, we are unable to transfer this to a situation where Hypothesis $(\dagger)$(ii) holds relative to this sequence of level 0 cones (as thus in a different fine blow-up situation where a suitable variant of the fine $\epsilon$-regularity theorem, Theorem \ref{thm:L1fine_reg}, might apply) as we do not know if the sequence of level 1 cones the fine blow-up was taken relative to were close to $Q^*_V$, i.e. close to the infimum of the excess over all level 1 (and level 2) cones.\footnote{It is then natural to wonder that if we were to modify our definition of the fine blow-up class $\FB^{F}_{3,1}(\BC^{(0)}$, namely by only considering fine blow-ups relative to sequences of level 1 cones $\BC_k\in \FL_1$ which obey $Q_{V_k,\BC_k}^2 < \frac{3}{2}M^\prime(Q^*_{V_k})^2$ for some $M^\prime\geq 1$, if this would suffice to prove the desired $(\FB7)$ property. Whilst such a class may obey the desired $(\FB7)$, the resulting class would not obey a suitable form of $(\FB5)$ to apply the results of \cite{minter2021} to establish the boundary regularity, even given $(\FB7)$. Indeed, without any additional assumptions, when performing rescalings by $\rho$ as in $(\FB5\text{I})$, the constant $M^\prime$ would necessary change to some $M^\prime(\rho)$, which would obey $M^\prime(\rho)\to \infty$ as $\rho\downarrow 0$. Hence we would see that the $\epsilon$ in $(\FB7)$ would depend on $M^\prime$, and hence on $\rho$ for the functions in $(\FB5\text{I})$, leading to no uniform choice of $\epsilon$ for the arguments in \cite{minter2021}. In order to achieve a uniform $\tilde{M}>M^\prime$ which all the functions described in $(\FB5)$ obey, we would need additional hypotheses, which ultimately leads one to the ultra fine blow-up constructed in this section.} This leads us naturally to consider this as a separate situation, which we refer to as an \textit{ultra fine blow-up}.

\subsection{Construction of Ultra Fine Blow-Ups}\label{sec:UF-construction}

Fix $\BC^{(0)}\in \FL_S\cap \FL_2$ a level 2 cone throughout; as usual, this will be our base cone. Let us first outline the hypotheses under which an ultra fine blow-up is constructed; this will be similar to those seen in Section \ref{sec:fine_construction}, and thus at numerous places we will refer back to the arguments there. Let $V\in \S_2$ and $\BC^c,\BC_1,\BC_0\in \FL$.

\textbf{Remark:} Up to rotation, there is only one level 2 cone $\BC^{(0)}\in \FL_S\cap \FL_2$ in $\R^{n+1}$; hence all of our constants in this section will in fact only be dependent on the dimension $n$.
\hspace{0.5em}
\begin{leftbar}
	\textbf{Hypothesis (G):} For appropriately small $\epsilon,\gamma_0,\gamma_1\in (0,1)$, to be determined depending only on $n$, we have:
	\begin{enumerate}
		\item [(G1)] $\BC^c\in \FL_{\epsilon}(\BC^{(0)})\cap \FL_2$;
		\item [(G2)] $V\in \CN_{\epsilon}(\BC^{(0)})$ and $\Theta_V(0) \geq 5/2$;
		\item [(G3)] $\BC_1\in \FL_{\epsilon}(\BC^{(0)})\cap \FL_1$ with $Q^2_{V,\BC_1}<\gamma_1 E_{V,\BC^c}^2$;
		\item [(G4)] $\BC_0\in \FL_{\epsilon}(\BC^{(0)})$ with $Q_{V,\BC_0}^2 < \gamma_0 E^2_{V,\BC_1}$.
	\end{enumerate}
\end{leftbar}

Moreover, for $M = M(n)>1$ a dimensional constant, we will also assume:

\begin{leftbar}
	\textbf{Hypothesis $(\diamond)$:} We have both:
	\begin{enumerate}
		\item [(a)] $E^2_{V,\BC^c} < M\inf_{\tilde{\BC}\in \FL_2}E^2_{V,\tilde{\BC}}$;
		\item [(b)] $E^2_{V,\BC_1} < M\inf_{\tilde{\BC}\in \FL_1}E^2_{V,\tilde{\BC}}$.
	\end{enumerate}
\end{leftbar}

Hypothesis (G) and Hypothesis ($\diamond$) should of course be compared to Hypothesis (H) and Hypothesis ($\dagger$) seen in Section \ref{sec:fine_construction}. These properties imply, for $\epsilon = \epsilon(n)$ and $\gamma_0 = \gamma_0(n)$ sufficiently small, that $\BC_0\in \FL_0$ is a level 0 cone. It should also be noted that Hypothesis (G4) and Hypothesis ($\diamond$) give $E_{V,\BC_0}^2 < \gamma_0M^2 \inf_{\tilde{\BC}\in \FL_1}E_{V,\tilde{\BC}}^2 \leq \gamma_0 M (Q_V^*)^2$, and thus these two properties give that $\BC_0$ obeys a form of Hypothesis ($\dagger$)(ii) from Section \ref{sec:fine_construction}. This already provides some intuition for our methods here: both two-valued functions represented $V$ over the multiplicity two half-hyperplanes in $\BC^c$ should split into two pairs of single-valued functions. We want to represent $V$ as single-valued functions over the half-hyperplanes in $\BC_0$. Note that all the results in Section \ref{sec:fine_construction} hold for $V,\BC^{(0)},\BC^c,\BC_1$ under Hypothesis (G) and Hypothesis ($\diamond$).

Now, for $\epsilon,\gamma_1,\gamma_0$ sufficiently small depending only on $n$, we know that we can write the half-hyperplanes in $\BC_0$ as linear functions over the half-hyperplanes in $\BC_1$, where over the (unique) multiplicity two half-hyperplanes in $\BC_1$ we have two distinct linear functions representing the two nearby half-hyperplanes in $\BC_0$. Let us write $\lambda_1,\lambda_2,\lambda_3$ for the gradients of these linear functions over the multiplicity one half-hyperplanes in $\BC_1$, and $\mu_1,\mu_2$ for the gradients of the two linear functions over the multiplicity two half-hyperplane in $\BC_1$. By essentially the same arguments as in Remark 1 and Remark 2 of Section \ref{sec:fine_construction}, we then have that there exist dimensional constants $c_1 = c_1(n)$ and $c_2 = c_2(n)$ such that
\begin{equation}\label{E:UF-1}
	c_2E_{V,\BC_1} \leq \max_{i,j}\{|\lambda_i|,|\mu_j|\} \leq c_1 E_{V,\BC_1}
\end{equation}
and there is a dimensional constant $c_3 = c_3(n)$ such that
\begin{equation}\label{E:UF-2}
	c_3E_{V,\BC_1}\leq |\mu_1-\mu_2|;
\end{equation}
these are also essentially analogous estimates to those in Remark 3 of Section \ref{sec:fine_construction}, i.e. when Hypothesis ($\dagger$)(ii) held. We note that here we do not need any equivalent hypothesis to Hypothesis ($\dagger$) of Section \ref{sec:fine_construction}, as we are in a ``smallest'' possible setting where no more splitting can occur and all degeneration is removed.

Let us first prove the analogue of Theorem \ref{thm:fine_representation} in this setting.

\begin{theorem}[Ultra Fine Representation]\label{thm:UF-construction}
	Let $\tau\in (0,1/40)$ and $\BC^{(0)}\in \FL_S\cap \FL_2$. Then, there exist constants $\epsilon^* = \epsilon^*(n,\tau)\in (0,1)$, $\gamma_1^* = \gamma^1_*(n,\tau)\in (0,1)$, and $\gamma_0^* = \gamma_0^*(n,\tau)\in (0,1)$ such that the following is true: let $V,\BC^{(0)},\BC^c,\BC_1,\BC_0$ satisfy Hypothesis (G) and Hypothesis ($\diamond$) with $\epsilon^*,\gamma_1^*,\gamma_2^*$, and $\frac{3}{2}M_0^4$ in place of $\epsilon,\gamma_1,\gamma_0$, and $M$, respectively. Then we have:
	\begin{align*}
		\ \ \ \ \textnormal{(a)}\ \ & V\res B_{3/4}\cap \{|x|>\tau\} = \mathbf{v}(u)\res\{|x|>\tau\}\text{, where }u\in C^{2}(\BC_0\res B_{3/4}(0)\cap \{|x|>\tau\});\\
		&\text{equivalently, we can express $V\res B_{3/4}\cap \{|x|>\tau\}$ as a sum of 5 single-valued functions}\\
		&\text{over $\BC_1$ or $\BC^c$, in the same way as in Theorem \ref{thm:fine_representation}(a), namely using $u$ and the corres--}\\
		&\text{ponding linear function defining the half-hyperplane in $\BC_0$ from the given half-hyperplane}\\
		&\text{in the cone $\BC_1$ or $\BC^c$;}\\
		\textnormal{(b)}\ \ & \int_{B_{5/8}(0)}\frac{|X^\perp|^2}{|X|^{n+2}}\ \ext\|V\|\leq CE_{V,\BC_0}^2;\\
		\textnormal{(c)}\ \ & \int_{B_{5/8}(0)}\sum^{n+1}_{j=3}|e_j^\perp|^2\ \ext\|V\|\leq CE_{V,\BC_0}^2;\\
		\textnormal{(d)}\ \ & \int_{B_{5/8}(0)}\frac{\dist^2(X,\spt\|\BC_0\|)}{|X|^{n+3/2}}\ \ext\|V\| \leq CE_{V,\BC_0}^2;
	\end{align*}
	here, $C = C(n)$ is a constant which is in particular independent of $\tau$.
\end{theorem}

\textbf{Remark:} We will not need the Hardt--Simon inequality in our proof this time, as our ultra fine blow-ups will be comprised of single-valued harmonic functions and thus we can use classical results from elliptic PDE theory to establish their boundary regularity.

\begin{proof}
Let us first prove (a); this essentially follows by an appropriate modification of Lemma \ref{lemma:splitting}, phrased in terms of cones in $\FL$ as opposed to hyperplanes and under hypotheses similar to those seen in Hypothesis (G). Indeed, if (a) we note true, then we could find sequences $\epsilon_k,\gamma_1^k,\gamma_0^k\downarrow 0$ and sequences of varifolds $V_k,\BC^c_k,\BC_1^k,\BC_0^k$ such that under the hypotheses of the lemma with $\epsilon_k,\gamma^1_k,\gamma_0^k$, in place of $\epsilon,\gamma_1,\gamma_0$, (a) does not hold. In particular, we have $E^{-1}_{V_k,\BC^k_1}Q_{V_k,\BC^k_0}\to 0$. Thus, if $\Phi$ denotes the fine blow-up of $(V_k)_k$ relative to $(\BC^c_k)_k$ and $(\BC^k_1)_k$, we see that $\Phi$ consists of $5$ linear functions, which have disjoint graphs in the region $\{|x|>0\}$. But this would imply that, from the local uniform convergence of the fine blow-up sequence to $\Phi$ on $\{|x|>0\}$, then $V$ has no multiplicity two singular points on the region $B_{3/4}\cap\{|x|>\tau\}$, and so the any two-valued function is in fact simply two single-valued functions; this then gives (a).

To prove (b) -- (d), just as in the proof of Theorem \ref{thm:fine_representation}, we need to extend the definition of $u$ from $B_{3/4}\cap\{|x|>\tau\}$ to a domain $U\subset \spt\|\BC^{(0)}\|$ such that, if $G := \graph(u|_U)$, then
\begin{equation}\label{E:UF-rep-1}
\int_{B_{3/4}\backslash G}r^2\ \ext\|V\| + \int_{U\cap B_{3/4}}r^2|Du|^2 \leq CE_{V,\BC_0}^2.
\end{equation}
This can be done in much the same way as seen in Theorem \ref{thm:fine_representation}, except now we have 8 different possibilities: for $C_1 = C_1(n)$, $C_2 = C_2(n)$, and $C_3 = C_3(n)$ sufficiently small to be chosen, first ask whether (i) $E^2_{V,\BC_0}(\tilde{T}_\rho(\zeta)) <C_1 E^2_{V,\BC_1}(T_\rho(\zeta))$, then ask whether (ii) $E^2_{V,\BC_1}(\tilde{T}_\rho(\zeta))< C_2 E^2_{V,\BC^c}(T_\rho(\zeta))$, and finally ask (iii) if the answer to (i) was ``yes'', then ask if $E^2_{V,\BC_0}(T_\rho(\zeta))<C_3$, otherwise if the answer to (i) was ``no'' and the answer to (ii) was ``yes'', ask whether $E^2_{V,\BC_1}(T_\rho(\zeta))<C_3$, otherwise if the answer to (i) was ``no'' and the answer to (ii) was ``no`` ask whether $E^2_{V,\BC^c}<C_3$; note that here we have written $E^2_{V,\BC}(T_\rho(\zeta)):= \rho^{-n-2}\int_{T_\rho(\zeta)}\dist^2(X,\spt\|\BC\|)\ \ext\|V\|$ for the excess over the region $T_\rho(\zeta)$, and similarly defined for $E^2_{V,\BC}(\tilde{T}_\rho(\zeta))$. In any of the 8 possibilities, if the answer to (iii) is ``yes'', then we include it in the definition of $U$; otherwise, we do not. Then one can check that this definition of $U$ gives rise to (\ref{E:UF-rep-1}), and thus the proof can be completed in the same manner as in Lemma \ref{lemma:L2_coarse}.
\end{proof}

\begin{proof}
	We only discuss (a), as the rest follows similarly to previous arguments. To prove (a), note that if it were not true, we could find sequences $\epsilon^k,\gamma_1^k,\gamma_0^k\to 0$ and sequences $V_k,\BC^c_k,\BC_1^k,\BC_0^k$ such that it were not true. In particular we have
	$$E_{V_k,\BC_1^k}^{-1}E_{V_k,\BC_0^k}\to 0$$
	and so if we take the fine blow-up $\Phi$ of $(V_k)_k$ relative to $(\BC^c_k)_k$, $(\BC^k_1)_k$, we would see that necessarily we have that $\Phi$ consists of 5 linear pieces (this is again due to (\ref{E:UF-1}) and (\ref{E:UF-2})). But this would imply that, from the local uniform convergence of the fine blow-ups away from the spine, that on $\{|x|>\tau\}$ $V$ has no multiplicity two pieces, and so the two-valued graph over the multiplicity two half-hyperplane in $\BC^k_1$ splits into two single-valued stationary graphs. This then gives the result.
\end{proof}

We then have the corresponding corollary to Corollary \ref{cor:fine_estimates} in this setting:

\begin{corollary}\label{cor:UF-estimates}
	Let $\BC^{(0)}\in \FL_S\cap \FL_2$. Then there exist constants $\epsilon^* = \epsilon^*(n)$, $\gamma^*_1 = \gamma^*_1(n)$, and $\gamma^*_0 = \gamma^*_0(n)\in (0,1)$ such that the following holds: if $V,\BC^{(0)},\BC^c,\BC_1,\BC_0$ satisfy Hypothesis (G) and Hypothesis ($\diamond$) with $\epsilon^*,\gamma_1^*,\gamma_0^*$, and $M = \frac{3}{2}M_0^3$ in place of $\epsilon,\gamma_1,\gamma_0$, and $M$, respectively, then for each $Z = (\xi,\zeta)\in \spt\|V\|\cap (\R^2\times B_{3/8}^{n-1}(0))$ with $\Theta_V(Z)\geq 5/2$, we have the following:
	\begin{enumerate}
		\item [(a)] $|\xi|\leq CE_{V,\BC_0}$;
		\item [(b)] For any $\rho\in (0,1)$, if we allow $\epsilon^*,\gamma_1^*,\gamma_0^*$ to depend on $\rho$ also, we have
		$$\int_{B_{5\rho/8}(Z)}\frac{\dist^2(X,\spt\|(\tau_Z)_\#\BC_0\|)}{|X-Z|^{n+3/2}}\ \ext\|V\|(X) \leq C\rho^{-n-3/2}\int_{B_\rho(Z)}\dist^2(X,\spt\|(\tau_Z)_\#\BC_0\|)\ \ext\|V\|(X);$$
	\end{enumerate}
	here, $C = C(n)$ is independent of $\rho$.
\end{corollary}

\begin{proof}
	The proof follows in much the same way as in Corollary \ref{cor:fine_estimates}. Indeed, we first argue that for any $\delta\in(0,1)$, there exists $\epsilon^* = \epsilon^*(n,\delta)$, $\gamma_1 = \gamma_1(n,\delta)$, and $\gamma_0 = \gamma_0(n,\delta)$ sufficiently small such that if Hypothesis (G) and Hypothesis ($\diamond$) hold for $V,\BC^{(0)},\BC^c,\BC_1,\BC_0$ with $\epsilon_*,\gamma_1^*,\gamma_0^*$ in place of $\epsilon,\gamma_1,\gamma_0$, respectively (and $M = \frac{3}{2}M_0^3$), then
	\begin{equation}\label{E:UF-1a}
		|\xi|\leq \delta E_{V,\BC_1}.
	\end{equation}
	Indeed, from Corollary \ref{cor:fine_estimates}(a) we already know under the present assumptions that there is some $C = C(n)$ for which $|\xi|\leq CE_{V,\BC_1}$, so this is an improved estimate. Indeed, to show this one may argue by contradiction in the same way as in proving (\ref{eqn:delta_bound}), except now taking a fine blow-up as described in Section \ref{sec:fine_construction} and using (F$_k$) from Section \ref{sec:fine_properties}. Given Theorem \ref{thm:UF-construction}, the proof now essentially follows in an identical fashion to that seen in the arguments from (\ref{eqn:delta_bound}) -- (\ref{E:fine-est-extra8}).
\end{proof}

Given Corollary \ref{cor:UF-estimates}, we know can estimate the following non-concentration of excess result in the usual fashion:

\begin{lemma}\label{lemma:UF-non-concentration}
	Let $\delta\in (0,1/10)$ and $\BC^{(0)}\in \FL_S\cap \FL_2$. Then, there exist constants $\epsilon^* = \epsilon^*(n,\delta)\in (0,1)$, $\gamma_1^* = \gamma_1^*(n,\delta)\in (0,1)$, and $\gamma^*_0 = \gamma^*_0(n,\delta)$ such that the following is true: if $V,\BC^{(0)},\BC^c,\BC_1,\BC_0$ satisfy Hypothesis (G) and Hypothesis ($\diamond$) with $\epsilon^*,\gamma_1^*,\gamma_0^*$, and $\frac{3}{2}M_0^3$ in place of $\epsilon,\gamma_1,\gamma_0$, and $M$, respectively, then:
	$$\int_{B_{3/4}\cap \{|x|<\sigma\}}\dist^2(X,\spt\|\BC_0\|)\ \ext\|V\| \leq C\sigma^{1/2}E_{V,\BC_0}^2$$
	for each $\sigma\in [\delta,1/4)$, where $C = C(n)$ is independent of $\sigma$.
\end{lemma}

\begin{proof}
	Given Corollary \ref{cor:UF-estimates} this is now identical to the proof in Lemma \ref{lemma:non-concentration_fine}.
\end{proof}

\subsection{Constructing the Ultra Fine Blow-Up Class}\label{sec:UF-class}

Using the results of Section \ref{sec:UF-construction} we now construct the class of ultra fine blow-ups, in a similar fashion to that seen in Section \ref{sec:fine_properties} for the fine blow-up classes.

Fix $M_1 = M_1(n)\in (1,\infty)$ and $\BC^{(0)}\in \FL_S\cap \FL_2$. Let $(\epsilon_k)_k$, $(\gamma_1^k)_k$, $(\gamma_0^k)_k$ be (decreasing) sequences of positive numbers converging to $0$. Consider sequences of varifolds $(V_k)_k\subset\S_2$, $(\BC^c_k)_k\subset\FL_2$, $(\BC^k_1)_k\subset \FL_1$, and $(\BC^k_0)_k\subset \FL_0$ such that for each $k\geq 1$, $V_k,\BC^{(0)},\BC^c,\BC^k_1,\BC^k_0$ obey Hypothesis (G) and Hypothesis $(\diamond)$ with $\epsilon_k,\gamma^1_k,\gamma^k_0$, and $M_1$ in place of $\epsilon,\gamma_1,\gamma_0$, and $M$, respectively. Thus, for each $k=1,2,\dotsc$, we suppose the following:
\begin{enumerate}
	\item [(1$_k$)] $V_k\in \CN_{\epsilon_k}(\BC^{(0)})$;
	\item [(2$_k$)] $\BC_k\in \FL_{\epsilon_k}(\BC^{(0)})\cap \FL_2$, $\BC^k_1\in \FL_{\epsilon_k}(\BC^{(0)})\cap \FL_1$, and $\BC^k_0\in \FL_{\epsilon_k}(\BC^{(0)})$;
	\item [(3$_k$)] $E^{-2}_{V_k,\BC^c_k}Q_{V_k,\BC_k^1}^2<\gamma^k_1$;
	\item [(4$_k$)] $E^{-2}_{V_k,\BC^k_1}E_{V_k,\BC^k_0}^2<\gamma^k_0$;
	\item [(5$_k$)] $E_{V_k,\BC^c_k}^2 < M_1\inf_{\tilde{\BC}\in \FL_2} E_{V_k,\tilde{\BC}}^2$;
	\item [(6$_k$)] $E_{V_k,\BC^k_1}^2 < M_1\inf_{\tilde{\BC}\in \FL_1}E_{V_k,\tilde{\BC}}^2$;
\end{enumerate}

Now let $(\delta_k)_k$ and $(\tau_k)_k$ be decreasing sequences of positive numbers converging to $0$. Let us write $H^k_1,\dotsc,H^k_4$ for the distinct half-hyperplanes in $\spt\|\BC^k_1\|$, so that
$$\BC^k_1 = 2|H^k_4| + \sum^3_{i=1}|H_i^k|.$$
We then write $h^k_1,\dotsc,h^k_5$ for the linear functions over the half-hyperplanes in $\BC^k_1$ whose graphs coincide, in the region $\{|x|>\tau_k\}$, with the half-hyperplanes in $\BC^k_0$; here, $h^k_4$ and $h^k_5$ are defined on $H_k^4$; for the sake of notational simplicity, we introduce $H^k_5:= H^k_4$. Note that we may also pass to a subsequence to ensure that it is the same multiplicity two half-hyperplane in $\BC^c_k$ which splits in $\BC^1_k$. For $i=1,2,\dotsc,5$, write $\lambda^k_i$ for the gradient of the linear function $h^k_i$. Write also $\w_i$ for the unit vector in $\R^2$ which determines the ray in the cross-section of $H^{k}_i$. We also write $\BC^k_0 = \sum^5_{i=1}|H_{0,i}^k|$, where $H_{0,i}^k$ are the half-hyperplanes determining $\BC^k_0$. As before in Section \ref{sec:fine_properties} and Section \ref{sec:coarse-construction}, we will be using the fixed domain $\BC^{(0)}$ as a parameter space for our functions, and we do not make a distinction in our notation between functions defined on half-hyperplanes in $\BC^k_1$ and those defined over half-hyperplanes in $\BC^{(0)}$.

By passing to an appropriate subsequence (and modifying the sequences $(\delta_k)_k$ and $(\tau_k)_k$ is needed), we may then deduce from the results in Section \ref{sec:UF-construction} that the following assertions hold:

\begin{enumerate}
	\item [(A$_k$)] For every point $Y\in S(\BC^{(0)})\cap B_{1/2}$, we have for all $k$ sufficiently large,
	$$B_{\delta_k}(Y)\cap \{Z:\Theta_V(Z)\geq 5/2\} \neq\emptyset;$$
	\item [(B$_k$)] For each $\sigma\in [\delta_k,1/4)$ we have:
	$$\int_{B_{3/4}\cap \{|x|<\sigma\}}\dist^2(X,\spt\|\BC^0_k\|)\ \ext\|V_k\| \leq C\sigma^{1/2}E^2_{V_k,\BC_0^k};$$
	\item [(C$_k$)] There are $5$ single-valued $C^2$ functions $u^k_1,\dotsc,u^k_5$, where $u^k_i\in C^2(H^k_{i}\cap B_{3/4}\cap \{|x|>\tau_k\}; (H^k_i)^\perp)$, each with stationary graph, such that
	$$V_k\res (B_{3/4}\cap \{|x|>\tau_k\}) = \sum^5_{i=1}|\graph(h^k_i + u^k_i)|;$$
	\item [(D$_k$)] For each point $Z = (\xi,\zeta)\in \spt\|V_k\|\cap B_{3/8}$ with $\Theta_{V_k}(Z)\geq 5/2$, we have
	$$|\xi|\leq CE_{V_k,\BC^k_0};$$
	\item [(E$_k$)] We have:
	$$c_2E_{V_k,\BC^k_1}\leq \max_i|\lambda^k_i| \leq c_1E_{V_k,\BC^k_1}$$
	and
	$$|\lambda^k_4 - \lambda^k_5|\geq c_3 E_{V_k,\BC^k_1};$$
	\item [(F$_k$)] For each $\rho\in (0,1/4]$, we can find $K = K(\rho)\in \Z_{\geq 1}$ such that for all $k\geq K$ the following holds: for each $Z = (\xi,\zeta)\in \spt\|V_k\|\cap B_{3/8}$ with $\Theta_{V_k}(Z)\geq 5/2$,
	\begin{align*}
	\sum^5_{i=1}\int_{H^k_{i}\cap B_{\rho/2}(Z)\cap\{|x|>\tau_k\}}&\frac{|u^k_{i}-\xi^{\perp_{H^k_{0,i}}}|^2}{|(h^k_{i}(r\w_{i},y) + u^k_{i}(r\w_{i},y),r\w_{i},y)-Z|^{n+3/2}}\\
	& \hspace{10em} \leq C\rho^{-n-3/2}\int_{B_\rho(Z)}\dist^2(X,\spt\|(\tau_Z)_\#\BC^k_0\|)\ \ext\|V_k\|;
	\end{align*}
	moreover, we have
	$$\xi^{\perp_{H_{0,i}^k}} = \xi^{\perp_{H_i^k}} - \lambda^{k}_i\xi^{\top_{H^k_i}}.$$
\end{enumerate}
here $C = C(n)$. Once again, (A$_k$) holds from Lemma \ref{lemma:gaps}, (B$_k$) holds by Lemma \ref{lemma:UF-non-concentration}, (C$_k$) holds by Theorem \ref{thm:UF-construction}, (D$_k$) holds by Corollary \ref{cor:UF-estimates}, (E$_k$) holds from (\ref{E:UF-1}) and (\ref{E:UF-2}), and (F$_k$) holds from Corollary \ref{cor:UF-estimates}. We may extend $u^k_I$ to all of $H^{(0)}_i\cap B_{3/4}$ by extending them by $0$ outside their domains of definition.

From (E$_k$) it follows that we can numbers $(\ell_i)_{i=1}^5$ obeying
$$c_2\leq \max_i|\ell_i|\leq c_1\ \ \ \ \text{and}\ \ \ \ |\ell_4-\ell_5|\geq 2c_3$$
such that, after passing to an appropriate subsequence, we have $E_{V_k,\BC_1^k}^{-1}\lambda^k_i \to \ell_i$. Moreover, by (C$_k$) and elliptic estimates for single-valued stationary graphs, we know that there exist $5$ single-valued harmonic functions, $f = (f_1,\dotsc,f_5)$ which patched together give a form on $\BC^{(0)}\res B^{n+1}_{3/4}\cap \{|x|>0\}$ (with the number of functions defined over a given half-hyperplane in $\BC^{(0)}$ equal to the multiplicity of the respective half-hyperplane in $\BC^{(0)}$) such that, after passing to another subsequence,
$$E_{V_k,\BC^k_0}^{-1}u_i^k\to f_i$$
where the convergence is in $C^2(K)$ for each compact subset $K\subset\spt\|\BC^{(0)}\|\cap B_{3/4}\cap \{|x|>0\}$. From (B$_k$) it follows that, in the same way as in Section \ref{sec:coarse-construction} that, for each $\sigma\in (0,1/4)$,
$$\int_{B_{3/4}}|f|^2 \leq C\sigma^{1/2}$$
and moreover that the convergence $E_{V_k,\BC^0_k}^{-1}u^k_i\to f_i$ is strongly in $L^2(B_{3/4})$.

\begin{defn}
	Fix $\BC^{(0)}\in\FL_S\cap \FL_2$ and $M>1$. Then we say that any quintuple of functions $f = (f_1,\dotsc,f_5)$ constructed as above with $M_1 = M$ for sequences of varifolds $(V_k)_k$, $(\BC^c_k)_k$, $(\BC^k_1)_k$, $(\BC^k_0)_k$, is called an \textit{ultra fine blow-up of} $(V_k)_k$ off $\BC^{(0)}$ relative to the sequences of cones $(\BC^c_k)_k$, $(\BC^k_1)_k$, and $(\BC^k_0)_k$, We write $\FB^{\FF}_M(\BC^{(0)})$ for the collection of all possible ultra fine blow-ups when we take $M_1 = M$ in (5$_k$) and (6$_k$).
\end{defn}

\subsection{Boundary Regularity of Ultra Fine Blow-Ups}\label{sec:UF-reg}

In Section \ref{sec:UF-class}, we constructed the ultra fine blow-up class $\FB^\FF_M(\BC^{(0)})$. We now need to understand the boundary regularity theory of functions in this class, so that we may in turn prove a suitable $\epsilon$-regularity theorem at the varifold level which in turn can be used to establish property $(\FB7)$ of Section \ref{sec:coarse_regularity} holds for the fine blow-up class $\FB^F_{3,1;M}(\BC^{(0)})$.

The boundary regularity theory of the ultra fine class is the simplest situation, as the $f \in \FB^\FF_M(\BC^{(0)})$ is comprised of $5$ single-valued harmonic functions defined, each defined on a half-hyperplane; thus, if one can prove that each harmonic function is continuous up-to-the-boundary of the half-hyperplane and that its boundary values are $C^{2,\alpha}$, one may invoke standard $C^{2,\alpha}$ boundary regularity theory of harmonic functions to deduce that $f$ is $C^{2,\alpha}$ up-to-the-boundary.

Note that properties (A$_k$), (D$_k$), (E$_k$) and (F$_k$) of Section \ref{sec:UF-class} give the following: for each $Y \in B_{3/8}\cap S(\BC^{(0)})$, we can find a sequence $Z_k\to Y$, where $Z_k = (\xi_k,\zeta_k)\in \spt\|V_k\|\cap B_{1/2}$ obeys $\Theta_{V_k}(Z_k) \geq 5/2$, and moreover that $\xi_k^{\perp_{H^k_{0,i}}} = \xi_k^{\perp_{H_i^k}} - \lambda^k_i\xi_k^{\top_{H^k_i}}$ and
$$\left|\xi_k^{\perp_{H_i^k}}\right|^2 + E_{V_k,\BC^k_1}^2\left|\xi_k^{\perp_{H_i^k}}\right|^2 \leq CE_{V_k,\BC^0_k}^2$$
and thus, up to passing to a subsequence, we have $E_{V_k,\BC_0^k}^{-1}\xi_k^{\perp_{H^k_i}}\to \kappa_i^\perp(Y)$ and $E_{V_k,\BC_0^k}^{-1}E_{V_k,\BC^k_1}\xi_k^{\top_{H_i^k}}\to \kappa_i^\top(Y)$ (where in the usual fashion we shall see momentarily that $\kappa^\perp_i(Y), \kappa^\top_i(Y)$ are only dependent on $Y$ and not on the approximating sequences $(Z_k)_k$), and so
$$E_{V_k,\BC_0^k}^{-1}\xi^{\perp_{0,i}^k}\to \kappa^\perp_i(Y) - \ell_i \kappa^\top_i(Y)$$
and thus, for each $\rho\in (0,1/4]$, we have
$$\sum^5_{i=1}\int_{H^{(0)}_i\cap B_{\rho/2}(Y)}\frac{|f_i-(\kappa^\perp_i(Y) - \ell_i\kappa^\top_i(Y))|^2}{|X-Y|^{n+3/2}} \leq C\rho^{-n-3/2}\int_{B_\rho(Z)}|f_i - (\kappa_i^\perp(Y) - \ell_i\kappa^\top_i(Y))|^2;$$
here, $H^{(0)}_i$ are the half-hyperplanes in $\BC^{(0)}$, counted with multiplicity. Such an inequality gives us, by Campanato style arguments for single-valued functions as discussed before, that $f_i\in C^{0,\alpha}(\overline{H_i^{(0)}}\cap B_{1/8};(H_i^{(0)})^\perp)$, with boundary values given by $\kappa_i:= \kappa_i^\perp - \ell_i\kappa_i^\top$ (we stress here that $\perp$ and $\top$ in $\kappa$ are purely notational, and do not represent projections of some fixed $\kappa$). Thus, all that remains to show is that $\kappa_i$ is a $C^{2,\alpha}$ function along $S(\BC^{(0)})\cap B_{1/8}$ (with estimates on its $C^{2,\alpha}$ norm in terms of $\int_{B_{1/2}}|\kappa|^2$). This can be done in much the same way as seen in the corresponding results for coarse and fine blow-ups seen in Section \ref{sec:coarse_properties} and Section \ref{sec:fine_initial_properties}, using now results from Section \ref{sec:UF-construction}; as such, we shall not duplicate the calculations here.

Thus we have now seen:

\begin{prop}\label{prop:UF-reg}
	For each $\alpha\in (0,1)$, elements of the ultra fine blow-up $\FB^\FF_M(\BC^{(0)})$ are harmonic functions which are $C^{1,\alpha}$ up-to-the-boundary.
\end{prop}

\subsection{The Ultra Fine $\epsilon$-Regularity Theorem}

Equipped now with the boundary regularity of the fine blow-up class, our next step is to prove an $\epsilon$-regularity theorem for varifolds under the assumptions seen in Section \ref{sec:UF-construction}. We will prove such a result in much the same way as seen in Section \ref{sec:L1fine}, by first proving a suitable excess improvement lemma.

\begin{lemma}[Ultra Fine Excess Improvement]\label{lemma:UF-excess-decay}
	Let $\BC^{(0)}\in \FL_S\cap \FL_2$ and $\theta\in (0,1/4)$. Then, there exist numbers $\bar{\epsilon} = \bar{\epsilon}(n,\theta)\in (0,1/2)$, $\bar{\gamma}_1 = \bar{\gamma}_1(n,\theta)\in (0,1/2)$, and $\bar{\gamma}_0 = \bar{\gamma}_0(n,\theta)\in (0,1/2)$ such that the following is true: if $V\in \S_2$, $\BC^c\in \FL_2$, $\BC_1\in \FL_1$, and $\BC_0\in \FL_0$ satisfy Hypothesis (G) and Hypothesis $(\diamond)$ with $M = \frac{3}{2}M_0$, then there exists an orthogonal rotation $\Gamma$ of $\R^{n+1}$ and a cone $\BC^\prime\in \FL_0$ such that the following hold:
	\begin{enumerate}
		\item [(a)] $|\Gamma-\id|\leq \kappa E_{V,\BC_0}$;
		\item [(b)] $\dist^2_\H(\spt\|\BC_0\|\cap B_1,\spt\|\BC^\prime\|\cap B_1) \leq \kappa E_{V,\BC_0}^2$;
		\item [(c)] $$\hspace{-0.5em}\theta^{-n-2}\int_{B_\theta}\dist^2(X,\spt\|\Gamma_\#\BC^\prime\|)\ \ext\|V\| + \theta^{-n-2}\int_{\Gamma(B_{\theta/2}\backslash\{|x|<\theta/16\})}\dist^2(X,\spt\|V\|)\ \ext\|\Gamma_\#\BC^\prime\|\leq \kappa \theta^2 E_{V,\BC_0}^2;$$
		\item [(d)] For any $\tilde{\BC}\in \FL_1$ with $\tilde{\BC}\in \FL_{1/10}(\BC_1)$ we have
		$$\left(\theta^{-n-2}\int_{B_\theta}\dist^2(X,\spt\|\tilde{\BC}\|)\ \ext\|\Gamma^{-1}_\#V\|\right)^{1/2} \geq \sqrt{2^{-n-4}\bar{C}_1}\dist_\H(\spt\|\BC_0\|\cap B_1,\spt\|\tilde{\BC}\|\cap B_1) - \kappa E_{V,\BC_0};$$
	\end{enumerate}
	here, $\kappa = \kappa(n)$ and $\bar{C}_1 = \bar{C}_1(n)\equiv \int_{B^n_{1/2}\cap\{x^2>1/16\}}|x^2|^2\ \ext\H^n(x^2,y)$ is as before.
\end{lemma}

\begin{proof}
	The proof follows by the same arguments as seen in Lemma \ref{lemma:fine_excess_decay}; indeed, the verification of Hypothesis (G1) -- (G3) and Hypothesis $(\diamond)$(a) is identical to as before (and indeed we can simply take the $\BC^c_k$ and $\BC_1^k$ sequences to be fixed), and the verification of Hypothesis (G4) and Hypothesis $(\diamond)$(b) is also the same, except now whenever in the corresponding argument of Lemma \ref{lemma:fine_excess_decay} a coarse blow-up was used (e.g. to prove $\tilde{C}E^2_{V_k,\BC^k_1}\leq E_{\tilde{V}_k,\BC^k_1}^2$, for some $\tilde{C} = \tilde{C}(n)$, which is the corresponding inequality to (\ref{E:fine-decay-3})), we instead use a fine blow-up relative to $(\BC^c_k)_k$ and $(\BC^k_1)_k$ (which obey Hypothesis $(\dagger)$ of Section \ref{sec:fine_construction}) and results from Theorem \ref{thm:UF-construction}, Corollary \ref{cor:UF-estimates}, and Lemma \ref{lemma:UF-non-concentration}. Thus, in the end we take an ultra fine blow-up and use the boundary regularity from Proposition \ref{prop:UF-reg} to generate the new sequence of (level 0) cones along which (a) -- (d) above hold.
\end{proof}

Now we are able to prove the ultra fine $\epsilon$-regularity theorem for varifolds, which will then be used to verify property $(\FB7)$ holds for the fine blow-up class $\FB^F_{3,1}(\BC^{(0)})$.

\begin{theorem}[Varifold Ultra Fine $\epsilon$-Regularity Theorem]\label{thm:UF-reg}
	Let $\BC^{(0)}\in \FL_S\cap \FL_2$ and $\alpha\in (0,1)$. Then there exist constants $\epsilon^* = \epsilon^*(n,\alpha)\in (0,1)$, $\gamma^*_1 = \gamma^*_1(n,\alpha)\in (0,1)$, and $\gamma^*_0 = \gamma^*_0(n,\alpha)\in (0,1)$ such that the following holds: if $V\in \S_2$, $\BC^c\in \FL_2$, $\BC_1\in \FL_1$, and $\BC\in \FL_0$ are such that $\Theta_V(0)\geq 5/2$, $V\in \CN_{\epsilon^*}(\BC^{(0)})$, $\BC^c,\BC_1,\BC_0\in \FL_{\epsilon^*}(\BC^{(0)})$, $E^2_{V,\BC^c}< M\inf_{\tilde{\BC}\in \FL_2}E_{V,\tilde{\BC}}^2$, $E^2_{V,\BC_1}<\frac{3}{2}\inf_{\tilde{\BC}\in \FL_1}E_{V,\tilde{\BC}}^2$, $Q_{V,\BC_1}^2<\gamma^*_1 E^2_{V,\BC^c}$, and $E^2_{V,\BC_0} < \gamma^*_0 E_{V,\BC_1}^2$, then there is a cone $\BC^\prime\in \FL_S\cap \FL_0$ with
	$$\dist_\H(\spt\|\BC^\prime\|\cap B_1,\spt\|\BC\|\cap B_1)\leq CQ_{V,\BC_0}$$
	and an orthogonal rotation $\Gamma:\R^{n+1}\to \R^{n+1}$ with $|\Gamma-\id|\leq CQ_{V,\BC}$ such that $\BC^\prime$ is the unique tangent cone to $\Gamma_\#V^{-1}$ at $0$, and
	$$\sigma^{-n-2}\int_{B_\sigma}\dist^2(X,\spt\|\BC^\prime\|)\ \ext\|\Gamma^{-1}_\#V\| \leq C\sigma^{2\alpha}Q^2_{V,\BC}\ \ \ \ \text{for all }\sigma\in (0,1/2);$$
	furthermore, $V$ has the structure of a $C^{1,\alpha}$ classical singularity of vertex density $5/2$; more precisely, there is a $C^{1,\alpha}$ function $u$ defined over $\spt\|\BC^c\|$, in the manner described in Theorem \ref{thm:A}, obeying $V\res B_{1/2} = |\graph(u)|$, and over any multiplicity two half-hyperplane in $\BC^c$, $u$ is given by two (disjoint) $C^{1,\alpha}$ single-valued functions, which meet only at the boundary; thus $V\res B_{1/2}$ has no (density 2) branch points, and $\sing(V)\cap B_{1/2} = \{\Theta_V = 5/2\}\cap B_{1/2}$ is the set of points determined by the boundary values of $u$.
\end{theorem}

\begin{proof}
	Given Lemma \ref{lemma:UF-excess-decay}, the proof is now similar to the proof of Theorem \ref{thm:L1fine_reg} and thus we do not repeat the arguments; indeed, one may take $\BC^c$, $\BC_1$ fixed, and just show that Hypothesis (G4) and Hypothesis ($\diamond$)(b) will hold inductively along applications of Lemma \ref{lemma:UF-excess-decay} (which give rise to sequences of cones $\BC_0^k$ and rotations $\Gamma_k$ as in the proof of Theorem \ref{thm:L1fine_reg}). We note that here we have changed our assumptions slightly, namely in Hypothesis $(\diamond)$ we have taken different constants in (a) and (b) for our assumption here, but this does not impact the previous arguments and so this assumption is still valid for the validity of Lemma \ref{lemma:UF-excess-decay}.
\end{proof}

\subsection{Property $(\FB7)$ for the Fine Blow-Up class $\FB^F_{1,3;M}(\BC^{(0)})$}

We can now use the ultra fine $\epsilon$-regularity theorem for varifolds, Theorem \ref{thm:UF-reg}, to prove that the fine blow-up class $\FB^F_{3,1;M}(\BC^{(0)})$ obeys property $(\FB7)$ from Section \ref{sec:coarse_regularity}, for any $M>1$.

\begin{corollary}\label{cor:UF-B7}
	Let $\BC^{(0)}\in \FL_S\cap \FL_2$. Then for each $M>1$, the fine blow-up class $\FB^F_{3,1;M}(\BC^{(0)})$, as defined in Section \ref{sec:fine_construction}, satisfied property $(\FB7)$ of Section \ref{sec:coarse_regularity} (with $\epsilon$ depending on $M$); in particular, it obeys the conclusions of Theorem \ref{thm:fine-reg}.
\end{corollary}

\begin{proof}
	The proof follows the same strategy as seen in the proof of Corollary \ref{cor:L1B7}, although let us sketch the proof to note some differences. Suppose for contradiction that $(\FB7)$ does not hold for some fixed $M>1$. Then we can find $v^k\in \FB^F(\BC^{(0)})$ obeying $v^k_a(0) = 0$, $Dv_a^k(0) = 0$, $\|v^k\|_{L^2} = 1$, and $v^k_*$ such that $v^k_*$ is comprised of linear functions with common boundary and zero average which obey 
	$$\int_{B_1}\G(v^k,v^k_*)^2 < \frac{1}{k};$$
	here, we stress that the average of any single-valued function is simply the function itself, whilst the average of any two-valued function is the usual average; thus, ``average'' here does not refer to ``average'' over functions defined on a given half-hyperplane, but instead of the individual functions themselves. Thus, over the multiplicity two half-hyperplane in $\BC^{(0)}$ for which $v^k$ is given by two single-valued functions, the two (single-valued) linear functions in $v^k_*$ over the same (multiplicity two) half-hyperplane are both zero; of course, over the multiplicity one half-hyperplane, $v^k_*$ is zero also. On the final remaining half-hyperplane, we are assuming that $v^k_*$ is given by two linear functions, $\ell_1,\ell_2$, which obey $\ell_1 \equiv -\ell_2 \not\equiv 0$. In particular, it is \textit{not} the case here that $\graph(v^k_*)$ is a level 0 cone; but when we pass to the cone level, $v_*^k$ will modify the sequence of level 1 cones giving rise to the fine blow-up, so will still give rise to a level 0 cone as it will split the multiplicity two piece into two.
	
	As in the proof of Corollary \ref{cor:L1B7}, we may pass to a subsequence to ensure that $v^k_*\to v_*$ (e.g. in $C^1$); $v^*$ then obeys $\|v^*\|_{L^2} = 1$ and $\int_{B_1}\G(v^k,v_*)^2\to 0$. Now let $(V_{k,j})_j\subset\S_2$, $(\BC^c_{k,j})_j\subset\FL_2$, and $(\BC_{k,j})_j\subset\FL_1$ be such that the fine blow-up sequence $v_{k,j}:= E_{V_{k,j},\BC_{k,j}}^{-1}u_{k,j}$ of $V_{k,j}$ relative to $\BC^c_{k,j}$ and $\BC_{k,j}$ gives rise to $v^k$. Again, for any $\sigma\in (0,1)$, we may find for each $k$ an index $j_k$ such that if $V_k:= V_{k,j_k}$, $\BC^c_k:= \BC^c_{k,j_k}$, $\BC_{k}:= \BC_{k,j_k}$, and $v_k:= v_{k,j_k}$, then
	$$\int_{B_\sigma}\G(v_k,v_*)^2 \to 0.$$
	Note that we know, by definition of the fine blow-up, that for all $k$ sufficiently large (for $j_k$ chosen appropriately, depending on $k$, $n$, and $\alpha$) we have $Q_{V_k,\BC_k}^2<\gamma^*_1 E_{V_k,\BC^c_k}$ and $E^2_{V_k,\BC^c_k}< M\inf_{\tilde{\BC}\in \FL_2}E^2_{V_k,\tilde{\BC}}$, where $\gamma^*_1 = \gamma^*_1(n,\alpha)$ is the constant from Theorem \ref{thm:UF-reg} (for $\alpha\in (0,1)$ fixed). Thus, we only need to verify the assumptions of Theorem \ref{thm:UF-reg} which correspond to Hypothesis (G4) and Hypothesis ($\diamond$)(b).
	
	To begin with, we claim that for all sufficiently large $k$,
	$$E^2_{V_k,\BC^k_1}<\frac{3}{2}\inf_{\tilde{\BC}\in \FL_1}E_{V_k,\tilde{\BC}}^2$$
	i.e. Hypothesis ($\diamond$)(b) holds with $M=3/2$. This follows in essentially the same manner as (\ref{E:B7-0}), and so we do not repeat the argument here. 
	
	We now generate the sequence of level 0 cones, $\BC^k_0\in \FL_0$, in the usual fashion: we modify the gradient of the half-hyperplanes in $\BC^k_1$ relative to $\BC^c_k$ by $E_{V_k,\BC_1^k}\cdot v_*$, for the corresponding value of $v_*$; in particular, note that it is only the multiplicity two half-hyperplane in $\BC_1^k$ which is modified, splitting into two distinct (multiplicity one) half-hyperplanes, by construction; thus $\BC_0^k$ is level 0. We can then follow (\ref{E:B7-3}) -- (\ref{E:B7-7}) to show that for all $k$ sufficiently large we have $Q_{V_k,\BC_0^k}^2<\frac{\gamma^*_0}{2}E^2_{V_k,\BC^k_1}$, where $\gamma^*_0 = \gamma^*(n,\alpha)$ is the constant from Theorem \ref{thm:UF-reg}. Thus, the assumptions of Theorem \ref{thm:UF-reg} are satisfied for all $k$ sufficiently large for $V_k,\BC^c_k,\BC^k_1,\BC^k_0$ (and in a uniform manner,  by which we mean they are also satisfied for all $V_{k,j}$, $\BC^c_{k,j}$, $\BC^{k,j}_1$, and correspondingly created $\BC^{k,j}_0$, for all $j\geq j_k$, as we just need the parameters $\epsilon_k$, $\gamma_k$, in the construction of the fine blow-up to be sufficiently small for this), and thus we can conclude using Theorem \ref{thm:UF-reg} in the same manner as in Corollary \ref{cor:L1B7}.
\end{proof}

\section{Boundary Regularity of Level 2 Coarse Blow-Ups and Completion of Main Theorem}

The aim of this section is to complete the proof of Theorem \ref{thm:A} in the case where the base cone $\BC^{(0)}$ is level 2. To do this, we first need to prove that the coarse blow-up class $\FB(\BC^{(0)})$ satisfies property $(\FB7)$ of Section \ref{sec:coarse_regularity}, and prove the version of the fine $\epsilon$-regularity theorem for varifolds, Theorem \ref{thm:L1fine_reg}, in the level 2 setting. The first step towards both of these results is the boundary regularity for the fine blow-up class $\FB^F_{3,1;M}(\BC^{(0)})$ established in Corollary \ref{cor:UF-B7}, which, coupled with Lemma \ref{lemma:fine_excess_decay}, will be used to prove a fine excess decay lemma when the base cone is level 2.

\subsection{Fine Excess Decay for Level 2 Cones}

Using Corollary \ref{cor:UF-B7} and Lemma \ref{lemma:fine_excess_decay}, we can now prove a fine excess decay lemma when the base cone is level 2.

\begin{lemma}\label{lemma:L2-fine-excess-decay-1}
	Let $\BC^{(0)}\in \FL_S\cap \FL_2$ and fix $\theta\in (0,1/4)$. Then, Lemma \ref{lemma:fine_excess_decay} holds, for some decay rate $\alpha = \alpha(n)$, without the assumption $\BC\in \FL_0$, i.e., there exist numbers $\bar{\epsilon} = \bar{\epsilon}(n,\theta)\in (0,1/2)$, $\bar{\gamma} = \bar{\gamma}(n,\theta)\in (0,1/2)$, and $\bar{\beta} =\bar{\beta}(n,\theta)\in (0,1/2)$ such that the following holds: if $V\in \S_2$, $\BC^c\in \FL_2$, and $\BC\in \FL_0\cup \FL_1$ satisfy Hypothesis (H), Hypothesis ($\star$), and Hypothesis ($\dagger$) of Section \ref{sec:fine_construction} with $\bar{\epsilon},\bar{\gamma},\bar{\beta}$, and $\frac{3}{2}M_0$ in place of $\epsilon,\gamma,\beta,$ and $M$, respectively, then there exists an orthogonal rotation $\Gamma$ of $\R^{n+1}$ and a cone $\BC^\prime\in \FL_0\cup \FL_1$ such that the following hold:
	\begin{enumerate}
		\item [(a)] $|\Gamma-\id|\leq \kappa E_{V,\BC}$;
		\item [(b)] $\dist^2_\H(\spt\|\BC\|\cap B_1,\spt\|\BC^\prime\|\cap B_1)\leq \kappa E_{V,\BC}^2$;
		\item [(c)] $$\hspace{-0.5em}\theta^{-n-2}\int_{B_\theta}\dist^2(X,\spt\|\Gamma_\#\BC^\prime\|)\ \ext\|V\| + \theta^{-n-2}\int_{\Gamma(B_{\theta/2}\backslash\{|x|<\theta/16\})}\dist^2(X,\spt\|V\|)\ \ext\|\Gamma_\#\BC^\prime\| \leq \kappa\theta^{2\alpha}E^2_{V,\BC};$$
		\item [(d)] For any $\tilde{\BC}\in \FL_2$ with $\tilde{\BC}\in \FL_{1/10}(\BC^c)$, we have:
		$$\left(\theta^{-n-2}\int_{B_\theta}\dist^2(X,\spt\|\tilde{\BC}\|)\ \ext\|\Gamma^{-1}_\#V\|\right)^{1/2} \geq \sqrt{2^{-n-4}\bar{C}_1}\dist_\H(\spt\|\BC\|\cap B_1,\spt\|\tilde{\BC}\|\cap B_1) - \kappa E_{V,\BC};$$
		here, $\kappa = \kappa(n)$, $\alpha = \alpha(n)$, and $\bar{C}_1 = \bar{C}_1(n) \equiv \int_{B_{1/2}^n\cap\{x^2>1/16\}}|x^2|^2\ \ext\H^n(x^2,y)$ is as before.
	\end{enumerate}
\end{lemma}

\begin{proof}
	We argue by contradiction in the same manner as in the proof of Lemma \ref{lemma:fine_excess_decay}: if the lemma does not hold for $\kappa = \kappa(n)\in (0,\infty), \alpha= \alpha(n)\in (0,1)$ to be chosen, then we may find sequences $\epsilon_k,\gamma_k,\beta_k\downarrow 0$, $V_k,\BC^c_k$, and $\BC_k$ satisfying Hypothesis (H), Hypothesis ($\star$), and Hypothesis $(\dagger)$ with $\epsilon_k,\gamma_k,\beta_k$, and $\frac{3}{2}M_0$ in place of $\epsilon,\gamma,\beta$, and $M$, respectively, such that the lemma does not hold for this choice of $\theta$ (and $\BC^{(0)}$). We already know from Lemma \ref{lemma:fine_excess_decay} that if Hypothesis ($\dagger$)(ii) holds for infinitely many $k$, then the lemma holds; so we may assume without loss of generality that Hypothesis ($\dagger$)(i) holds for all (but finitely many) $k$, i.e. that $\BC_k\in \FL_1$.
	
	But then if we follow the proof of Lemma \ref{lemma:fine_excess_decay} in this situation (which is entirely analogous to the situation where $\BC^{(0)}\in \FL_S\cap \FL_1$ is level 1), we may take a fine blow-up of (a rotation of) $V_k$ relative to $\BC^c_k$ and $\BC_k$; call this fine blow-up $v$. But we know from Corollary \ref{cor:UF-B7} that $v$ is $C^{1,\alpha}$ up-to-the-boundary on each half-hyperplane of $\spt\|\BC^{(0)}\|$, for some $\alpha = \alpha(n)$, with decay estimates, and moreover that the boundary values of the two-valued piece in $v$ are in fact given by a (multiplicity two) single-valued function. Thus, the new cone which the fine blow-up determines (in the same manner as Lemma \ref{lemma:fine_excess_decay}) will again by level 1, and the same proof as in Lemma \ref{lemma:fine_excess_decay} shows that the result holds for infinitely many $k$, providing the contradiction and thus proving the result.
\end{proof}

We now prove a stronger fine excess decay statement which removes the assumption of Hypothesis $(\dagger)$ from Lemma \ref{lemma:L2-fine-excess-decay-1}. In doing so, our excess decay lemma will change slightly -- we will no longer have one decay scale, but two possible decay scales; we remark that such a change will not greatly impact our previously arguments when using excess decay statements.

\begin{lemma}[Level 2: Fine Excess Decay]\label{lemma:L2-fine-excess-decay-2}
	Let $\BC^{(0)}\in \FL_S\cap \FL_2$ and fix $\theta_1,\theta_2\in (0,1/4)$ such that $\theta_2 < \theta_1/8$. Then, there exist numbers $\tilde{\epsilon} = \tilde{\epsilon}(n,\theta_1,\theta_2)\in (0,1/2)$, $\tilde{\gamma} = \tilde{\gamma}(n,\theta_1,\theta_2)\in (0,1/2)$ such that the following holds: if $V\in \S_2$, $\BC^c\in \FL_2$, and $\BC\in \FL_0\cup \FL_1$ satisfy Hypothesis (H) and Hypothesis $(\star)$ of Section \ref{sec:fine_construction} with $\tilde{\epsilon},\tilde{\gamma}$, and $\frac{3}{2}M_0$ in place of $\epsilon,\gamma$, and $M$, respectively, then there exists an orthogonal rotation $\Gamma$ of $\R^{n+1}$ and a cone $\BC^\prime\in \FL_0\cup\FL_1$ such that we have:
	\begin{enumerate}
		\item [(a)] $|\Gamma-\id|\leq \kappa Q_{V,\BC}$;
		\item [(b)] $\dist^2_\H(\spt\|\BC\|\cap B_1,\spt\|\BC^\prime\|\cap B_1)\leq \kappa Q_{V,\BC}^2$;
	\end{enumerate}
	and for some $j\in \{1,2\}$,
	\begin{enumerate}
		\item [(c)] $$\hspace{-0.5em}\theta_j^{-n-2}\int_{B_{\theta_j}}\dist^2(X,\spt\|\Gamma_\#\BC^\prime\|)\ \ext\|V\| + \theta_j^{-n-2}\int_{\Gamma(B_{\theta_j}/2\backslash\{|x|<\theta_j/16\})}\dist^2(X,\spt\|V\|)\ \ext\|\Gamma_\#\BC^\prime\|\leq \nu_j\theta^{2\alpha}_j Q^2_{V,\BC};$$
		\item [(d)] For any $\tilde{\BC}\in \FL_2$ with $\tilde{\BC}\in \FL_{1/10}(\BC^c)$, we have
		$$\left(\theta_j^{-n-2}\int_{B_{\theta_j}}\dist^2(X,\spt\|\tilde{\BC}\|)\ \ext\|\Gamma_\#^{-1}V\|\right)^{1/2} \geq \sqrt{2^{-n-4}\bar{C}_1}\dist_\H(\spt\|\BC\|\cap B_1,\spt\|\tilde{\BC}\|\cap B_1) - \kappa Q_{V,\BC};$$
	\end{enumerate}
	here, $\kappa = \kappa(n,\theta_1)$, $\alpha = \alpha(n)$, $\bar{C}_1 = \bar{C}_1(n)$ (is the usual constant), $\nu_1 = \nu_1(n)$, and $\nu_2 = \nu_2(n,\theta_1)$.
\end{lemma}

\textbf{Note:} We will see from the proof that our bounds must be in terms of $Q_{V,\BC}$ and not $E_{V,\BC}$.

\begin{proof}
	The proof of this, given Lemma \ref{lemma:L2-fine-excess-decay-1}, follows in the same manner as \cite[Lemma 13.2 and Lemma 13.3]{wickstable} do from \cite[Lemma 13.1]{wickstable}; we outline this argument here for the sake of completeness. 
	
	It $\BC\in \FL_1$, then there is nothing to prove: Hypothesis $(\dagger)$ trivially holds in this instance, and so the result (with $j=2$ in (c) and (d)) follows from Lemma \ref{lemma:L2-fine-excess-decay-1} taken with $\theta = \theta_2$; fix the constants from Lemma \ref{lemma:L2-fine-excess-decay-1} gives in this instance, namely $\epsilon_2 = \epsilon_2(n,\theta_2)$, $\gamma_2 = \gamma_2(n,\theta_2)$ (there is no $\beta$ in this case).
	
	So now let us suppose $\BC\in \FL_0$. Firstly, choose a cone $\tilde{\BC}\in \FL_1$ for which
	$$Q_{V,\tilde{\BC}}^2 \leq \frac{3}{2}(Q_V^*)^2.$$
	Now let $\beta_1 = \beta_1(n,\theta_1)$ be as in Lemma \ref{lemma:L2-fine-excess-decay-1} for $\theta = \theta_1$. Then if we have $Q^2_{V,\BC}< \beta_1 (Q_V^*)^2$, then Hypothesis $(\dagger)$ holds for $V,\BC^c,\BC$ (provided $\epsilon_1 = \epsilon_1(n,\theta_1)$ and $\gamma = \gamma_1(n,\theta_1)$ are sufficiently small as in Lemma \ref{lemma:L2-fine-excess-decay-1}), and hence the result follows from Lemma \ref{lemma:L2-fine-excess-decay-1}. Otherwise, we must have
	$$Q_{V,\BC}^2\geq \beta_1(Q^*_V)^2$$
	and thus we would have
	$$Q^2_{V,\tilde{\BC}} \leq \frac{3}{2\beta_1}Q^2_{V,\BC} < \frac{3\gamma}{2\beta_1}E^2_{V,\BC^c}$$
	where we have used Hypothesis (H). Thus, we see that if $\gamma = \gamma(n,\theta_1)$ is sufficiently small, then Hypothesis (H) will hold for $V,\BC^c,\tilde{\BC}$, and so as Hypothesis $(\star)$ still holds (as $\BC^c$ has no changed) and Hypothesis $(\dagger)$ is trivially satisfied in this instance (as $\tilde{\BC}\in \FL_1$), we would be able to apply Lemma \ref{lemma:L2-fine-excess-decay-1} (with $\theta=\theta_1$) to $V,\BC^c,\tilde{\BC}$, provided $\epsilon = \epsilon(n,\theta_1)$ was sufficiently small, to deduce the result, but with $\tilde{\BC}$ in place of $\BC$. But as $E^2_{V,\tilde{\BC}}\leq Q^2_{V,\tilde{\BC}} \leq \frac{3}{2\beta_1}Q_{V,\BC}^2$, the inequalities (a) -- (d) in terms of $E_{V,\tilde{\BC}}$ can readily be written in terms of $Q_{V,\BC}$, up to the constants changing by terms involving factors of $\frac{3}{2\beta_1}$, which depends on $\theta_1$. Moreover, any distance terms involving $\tilde{\BC}$ can be replaced by distance terms involving just $\BC$ by using the fact that here we have $\dist_\H^2(\spt\|\BC\|\cap B_1,\spt\|\tilde{\BC}\|\cap B_1) \leq \tilde{C}(Q_{V,\BC}^2 + Q_{V,\tilde{\BC}}^2)$, where $\tilde{C} = \tilde{C}(n)$. Thus the result follows, by taking $\epsilon = \epsilon(n,\theta_1,\theta_2) \leq \min\{\epsilon_1,\epsilon_2\}$ and $\gamma = \gamma(n,\theta_1,\theta_2)\leq\min\{\gamma_1,\gamma_2\}$ suitably small.
\end{proof}

\subsection{The Fine $\epsilon$-Regularity Theorem: Level 2 Setting}

Now that we have the full fine excess decay lemma, namely Lemma \ref{lemma:L2-fine-excess-decay-2}, when the base cone is level 2, we may now prove the variant of the fine $\epsilon$-regularity theorem for varifolds, i.e. Theorem \ref{thm:L1fine_reg}, in the level 2 setting.

\begin{theorem}[Varifold Fine $\epsilon$-Regularity Theorem: Level 2 Setting]\label{thm:L2_fine_reg}
Let $\BC^{(0)}\in \FL_S\cap \FL_2$. Then, there exist constants $\epsilon_1 = \epsilon_1(n)\in (0,1)$, $\gamma_1 = \gamma_1(n)\in (0,1)$ such that the following is true: if $V\in \S_2$, $\BC^c\in \FL_1$, and $\BC\in \FL$ are such that $\Theta_V(0)\geq 5/2$, $V\in \CN_{\epsilon_1}(\BC^{(0)})$, $\BC^c,\BC\in \FL_{\epsilon_1}(\BC^{(0)})$, $E^2_{V,\BC^c}< \frac{3}{2}\inf_{\tilde{\BC}\in \FL_2}E_{V,\tilde{\BC}}^2$, and $Q_{V,\BC}^2<\gamma_1 E_{V,\BC^c}^2$, then there is a cone $\BC^\prime\in \FL_S\cap (\FL_0\cup \FL_1)$ with
$$\dist_\H(\spt\|\BC^\prime\|\cap B_1,\spt\|\BC\|\cap B_1)\leq CQ_{V,\BC}$$
and an orthogonal rotation $\Gamma:\R^{n+1}\to \R^{n+1}$ with $|\Gamma-\id|\leq CQ_{V,\BC}$ such that $\BC^\prime$ is the unique tangent cone to $\Gamma^{-1}_\#V$ at $0$, and
$$\sigma^{-n-2}\int_{B_\sigma}\dist^2(X,\spt\|\BC^\prime\|)\ \ext\|\Gamma^{-1}_\#V\| \leq C\sigma^{2\alpha}Q_{V,\BC}^2\ \ \ \ \text{for all }\sigma\in (0,1/2).$$
Furthermore, there is a $C^{1,\alpha}$ function $u$ defined over $\spt\|\BC^c\|$, in the manner described in Theorem \ref{thm:A}, obeying $V\res B_{1/2} = \mathbf{v}(u)$, and over one multiplicity two half-hyperplane in $\BC^{(0)}$, $u$ is in fact given by two (disjoint) $C^{2}$ single-valued functions. Here, $C = C(n)\in (0,\infty)$ and $\alpha = \alpha(n)\in (0,1)$.
\end{theorem}

\begin{proof}
	The proof follows the same lines as that seen in Theorem \ref{thm:L1fine_reg}, however some modifications are needed. 
	
	The first is that in our fine excess decay lemma, Lemma \ref{lemma:L2-fine-excess-decay-2}, we have two possible decay scales as opposed to the single scale (and moreover that we no longer have a decay factor of $\theta^2$ but of $\theta^{2\alpha}$; here we will be able to get any power $<\alpha$, where $\alpha = \alpha(n)$ is as in Lemma \ref{lemma:L2-fine-excess-decay-2}). The modifications to deal with this difference are simple: fix any $\alpha^\prime \in (0,\alpha)$, and first choose $\theta_1 = \theta_1(n,\alpha^\prime)$ such that $\nu_1\theta^{2(\alpha-\alpha^\prime)}_1 < 1$, where $\nu_1 = \nu_1(n)$ is as in Lemma \ref{lemma:L2-fine-excess-decay-2}. Then choose $\theta_2 = \theta_2(n,\alpha^\prime)$ obeying $\theta_2<\theta_1/8$ and $\nu_2\theta_2^{2(\alpha-\alpha^\prime)} < 1$, where $\nu_2 = \nu_2(n,\theta_1) = \nu_2(n)$ is as in Lemma \ref{lemma:L2-fine-excess-decay-2}. Then one may follow (\ref{E:fine-reg-1}) --(\ref{E:fine-reg-13}) in an identical fashion, up to changing $\alpha$ to $\alpha^\prime$ and instead of our sequence of scales being $\theta,\theta^2,\theta^3,\dotsc$, we have a sequence of scales of the form $\sigma_k = \theta_1^{n_k}\theta_2^{k-n_k}$, for some $n_k\in\{0,1,2,\dotsc,k\}$, i.e. change $\theta_k\equiv \theta^k$ in the proof of Theorem \ref{thm:L1fine_reg} to this $\sigma_k$. Thus, (\ref{E:fine-reg-1}) -- (\ref{E:fine-reg-13}) follow in the same fashion; thus Summary 1 from the proof of Theorem \ref{thm:L1fine_reg} holds.
	
	However, Summary 2 from the proof of Theorem \ref{thm:L1fine_reg} does not currently hold, as in the proof of Corollary \ref{cor:fine_estimates} we needed to assume Hypothesis $(\dagger)$ holds, which we are currently not assuming (this was necessary to control the various excess quantities when shifting the base point). If Hypothesis $(\dagger)$ does hold for $\BC$, with $\beta = \beta_0$, where $\beta_0 = \beta_0(n)$ is from Corollary \ref{cor:fine_estimates}\footnote{We remark that in the proof of Corollary \ref{cor:fine_estimates}, it was shown that for any $\tilde{\epsilon}$, $\tilde{\gamma}\in (0,1/2)$ and $M_1>1$, there was $\epsilon_0 = \epsilon_0(n,\tilde{\epsilon},\tilde{\gamma},M_1), \gamma_0 = \gamma_0(n,\tilde{\epsilon},\tilde{\gamma},M_1)$, and $\beta_0 = \beta_0(n,M_1)$ such that if Hypothesis (H), Hypothesis $(\star)$, and Hypothesis $(\dagger)$(ii) held with $\epsilon_0,\gamma_0,\beta_0$, and $M_1$ in place of $\epsilon,\gamma,\beta$, and $M$, then for any $Z\in \spt\|V\|\cap B_{3/8}$ with $\Theta_V(Z)\geq 5/2$, the varifold $\tilde{V} := (\eta_{Z,1/4})_\#V$ would satisfy Hypothesis (H) and Hypothesis $(\star)$ with $\tilde{\epsilon}$, $\tilde{\gamma}$, and $M_1M_0$ in place of $\epsilon,\gamma$, and $M$, and moreover the inequalities (\ref{E:fine-reg-19b}) and (\ref{E:fine-reg-20}) will hold, i.e. Summary 2 from Theorem \ref{thm:L1fine_reg} still holds. We stress that, whilst Hypothesis $(\dagger)$ is not guaranteed for $\tilde{V}$ and any chosen constant $\tilde{\beta}$ (which would need additional smallness assumptions on $\beta_0$ in terms of $\tilde{\beta}$), for the conclusions listed $\beta_0$ can be chosen independent of $\tilde{\epsilon}$ and $\tilde{\gamma}$.}, then we know that the hypothesis which lead to Summary 1 will be true for $(\eta_{Z,1/4})_\#V$, where $Z\in \spt\|V\|\cap B_{3/8}$ is such that $\Theta_V(Z)\geq 5/2$, i.e. Summary 2 will still hold; hence the proof can be completed in this case in the same manner as in Theorem \ref{thm:L1fine_reg}. Otherwise, if Hypothesis ($\dagger$) does not hold for $\BC$ with this choice of $\beta$, then choosing $\tilde{\BC}\in \FL_1$ with $Q_{V,\tilde{\BC}}^2 < \frac{3}{2}(Q^*_V)^2$, we note that $Q_{V,\tilde{\BC}}^2 \leq \frac{3}{2\beta}Q_{V,\BC}^2$, and that Hypothesis ($\dagger$) does hold for $\tilde{\BC}$; moreover, for suitably small $\epsilon = \epsilon(n)$, $\gamma = \gamma(n)$, Hypothesis (H) and Hypothesis $(\star)$ will hold for $V,\BC^c$, and $\tilde{\BC}$. Thus we can run the proof of Theorem \ref{thm:L1fine_reg} with $\tilde{\BC}$ in place of $\BC$, and replace the final inequalities, which will be in terms of $\tilde{\BC}$ and $Q_{V,\tilde{\BC}}$, by those in terms of $\BC$ using $Q_{V,\tilde{\BC}}^2 \leq \frac{3}{2\beta}Q_{V,\BC}^2$ and $\dist_\H^2(\spt\|\BC\|\cap B_1,\spt\|\tilde{\BC}\|\cap B_1) \leq \tilde{C}(Q_{V,\BC}^2 + Q_{V,\tilde{\BC}}^2)$, where $\tilde{C} = \tilde{C}(n)$. Of course, in the discussion in the proof of Theorem \ref{thm:L1fine_reg} after (\ref{E:fine-reg-33}) we can no longer show that the any point of density not equal to $5/2$ is regular (as our cone can have a multiplicity two half-hyperplane), but the same argument will now show that any other singular point must either be a density 2 branch point or density 2 classical singularity, from Theorem \ref{thm:wick1} and Theorem \ref{thm:wick2}. Moreover, of course the fact that over one multiplicity two half-hyperplane the two-valued function splits into two single-valued functions follows immediately from the excess decay result, e.g. (\ref{E:fine-reg-16}), and (\ref{E:fine-reg-33}). Thus the proof is complete.
\end{proof}

Using Theorem \ref{thm:L2_fine_reg} we are now able to prove that the coarse blow-up class, $\FB(\BC^{(0)})$, where $\BC^{(0)}\in \FL_S\cap \FL_2$, obeys property $(\FB7)$ from Section \ref{sec:coarse_regularity}.

\begin{corollary}\label{cor:L2-B7}
	Let $\BC^{(0)}\in \FL_S\cap\FL_2$. Then, the coarse blow-up class $(\FB7)$ obeys property $(\FB7)$; in particular, Theorem \ref{thm:coarse_reg} holds for $\FB(\BC^{(0)})$.
\end{corollary}

\begin{proof}
	Given Theorem \ref{thm:L2_fine_reg}, the proof now follows in an identical fashion to that seen in Corollary \ref{cor:L1B7}.
\end{proof}

\subsection{Level 2: Proof of Main Theorem}

\begin{proof}[Proof of Theorem \ref{thm:A} when $\BC^{(0)}\in \FL_S\cap \FL_2$ is level 2]
	Fix $\BC^{(0)}\in \FL_S\cap \FL_2$. Given Theorem \ref{thm:L2_fine_reg} and Corollary \ref{cor:L2-B7}, Theorem \ref{thm:A} now follows in an identical fashion to that seen in Theorem \ref{thm:mainL1}; in the initial dichotomy, alternative (ii) will allow for level 0 and level 1 cones.
\end{proof}

Hence we have now shown that Theorem \ref{thm:A} holds in all cases, and thus have completed the proof of Theorem \ref{thm:A}.

\section{Concluding Remarks and Future Questions}

Let us now outline some possible future research directions arising from the current work. Firstly, we have seen in Theorem \ref{thm:A} the existence of an $\alpha = \alpha(n)\in (0,1)$ such that under the assumptions of Theorem \ref{thm:A}, the nearby varifold $V$ is expressible as a $C^{1,\alpha}$ graph over $\spt\|\BC_0\|$ in the sense described in Theorem \ref{thm:A}. In the case where $\BC_0\in \FL_0$ is level 0, we know from the work of \cite{krummel2014regularity} that in fact $V$ must be smooth (in fact, real analytic) up-to-the-boundary, with the points of density $5/2$ in $V$ forming a real analytic $(n-1)$-dimensional submanifold, i.e. $V$ is locally a ($C^\infty$) classical singularity, in the sense of \cite{wickstable}. Thus, in this case we get an improved regularity conclusion. This naturally raises the question of whether one could establish an optimal regularity conclusion in the general situation (which can at most be $\alpha=1/2$ when two-valued functions are used, see \cite{simon2016frequency}):

\begin{enumerate}
	\item [(Q1)] What is the optimal value of $\alpha$ in Theorem \ref{thm:A}? Can we take $\alpha=1/2$?
\end{enumerate}

It should be noted however that currently we do not have examples of a varifold $V\in \S_2$ which has a point $X$ with a tangent cone $\BC\in \FL_S$ such that $X$ is a limit point of density $2$ branch points in $V$. If such a situation is in fact impossible, then the proof of Theorem \ref{thm:A} would be significantly simplified, as the two-valued functions used for any graphical representation would actually simply be two single-valued functions with stationary graphs; thus the boundary regularity for the blow-up class is significantly simplified. Hence it is natural to ask:

\begin{enumerate}
	\item [(Q2)] Given $V\in \S_2$ and $X\in \spt\|V\|$ with $\Theta_V(X) = 5/2$, is it possible for $X$ to be a limit point of density 2 branch points in $V$, whilst at the same time having a tangent cone $\BC\in \FL_S$? If not, can one construct examples of this behaviour?
\end{enumerate}

Another hurdle which needed to be overcome in the current work was the absence of any general $C^{1,\alpha}$ boundary regularity statements for two-valued $C^{1,\alpha}$ harmonic functions. One could therefore ask whether such a boundary regularity statement might be true in general under weaker assumptions than those seen here, perhaps more in line with classical boundary regularity statements from the theory of elliptic PDEs.

\begin{enumerate}
	\item [(Q3)] Let $H = \{x\in \R^n: x^1>0\}$ and $\alpha\in (0,1)$. Suppose $u\in C^{1,1/2}(H;\A_2(\R))\cap C^{0,\alpha}(\overline{H};\A_2(\R))$ is a symmetric two-valued $C^{1,1.2}$ function in $H$. Suppose also that $\left. u\right|_{\del H} = \{0,0\}$. Then, is $u\in C^{1,\beta}(\overline{H};\A_2(\R))$, for some $\beta\in (0,1/2]$? (With estimates on $\|u\|_{C^{1,\beta}}$ in terms of $\|u\|_{L^2}$.)
\end{enumerate}

Another point of note is that we saw in Theorem \ref{thm:coarse_reg} that it is possible to prove that the boundary branch set for each coarse blow-up $v\in \FB(\BC)$ is well-behaved, in the sense that in fact it is possible to reflect the symmetric part of $v$ across the boundary of the half-hyperplane and still have a $C^{1,1/2}$ harmonic function on all of $\R^n$; thus boundary branch points are just interior branch points of the reflected function, and thus we may apply the interior regularity results of \cite{simon2016frequency} and \cite{krummelwick1} to say more. Is it possible to prove similar results hold at the varifold level, as in Theorem \ref{thm:A}? In particular:

\begin{enumerate}
	\item [(Q4)] Let $V$ be as in Theorem \ref{thm:A} and let $\B$ denote its (multiplicity two) branch set. Let us write $\sing_{5/2}(V):= \{X\in \spt\|V\|: \Theta_V(X) = 5/2\}$. Then must we have $\dim_\H(\sing_{5/2}(V)\cap \overline{\B}) \leq n-2$\;?
\end{enumerate}

Finally, we remark that given Theorem \ref{thm:A}, Theorem \ref{thm:wick1}, and Theorem \ref{thm:wick2}, it now seems reasonable to extend the results of \cite{polyhedral} to more general polyhedral cones for the class $\S_2$, namely those polyhedral cones with $4$-way and $5$-way junctions.

\bibliographystyle{alpha} 
\bibliography{references}

\end{document}